\documentclass[oneside,a4paper]{amsart}

\usepackage[T1]{fontenc}
\usepackage[utf8]{inputenc}
\usepackage[british]{babel}
\usepackage{amsmath, amsthm, amsfonts, amssymb}

\usepackage{hyperref}
\usepackage{url}
\usepackage[all]{xy}
\usepackage{placeins} 

\usepackage[usenames,dvipsnames]{color}
\usepackage{tikz}
\usetikzlibrary{arrows,backgrounds,
    decorations.pathreplacing,
    decorations.pathmorphing
}
\usetikzlibrary{matrix,arrows}

\tikzset{snake it/.style={decorate, decoration=snake, segment length=4pt}}

\tikzset{twist/.style={decorate, decoration={snake, amplitude=1.5pt, segment length=5pt}}}
\tikzset{small_twist/.style={decorate, decoration={snake, amplitude=0.5pt, segment length=3.5pt}}}
\definecolor{red}{rgb}{0.9, 0.4, 0.31}
\definecolor{serre_blue}{rgb}{0, 0.3, 0.7}
\definecolor{serre_green}{rgb}{0, 0.7, 0}
\definecolor{serre_orange}{rgb}{0.99, 0.6, 0}
\definecolor{chambering_blue}{rgb}{0.04, 0.28, 0.34}

\definecolor{halo_green}{rgb}{0.57, 0.8, 0}

\theoremstyle{plain}

\newtheorem{thm}{Theorem}[section]
\newtheorem{lem}[thm]{Lemma}
\newtheorem{prop}[thm]{Proposition}
\newtheorem{cor}[thm]{Corollary}
\newtheorem*{thm*}{Theorem}
\newtheorem*{thm_a*}{Theorem A}
\newtheorem*{thm_b*}{Theorem B}
\newtheorem*{thm_c*}{Theorem C}
\newtheorem*{thm_d*}{Theorem D}
\newtheorem*{thm_e*}{Theorem E}
\newtheorem*{cobordism_hypothesis*}{Cobordism Hypothesis}

\theoremstyle{definition}

\newtheorem{defin}[thm]{Definition}
\newtheorem{rem}[thm]{Remark}
\newtheorem*{rem*}{Remark}
\newtheorem*{defin*}{Definition}
\newtheorem*{conjecture*}{Conjecture}
\newtheorem*{notation*}{Notation}
\newtheorem*{convention*}{Convention}
\newtheorem*{reading_presentation*}{Reading the presentation}
\newtheorem*{reading_graphic*}{Reading the graphic}
\newtheorem*{interpreting_pasting_diagrams*}{Interpreting pasting diagrams}
\newtheorem*{coherence_for_monoidal_bicategories*}{Coherence for monoidal bicategories}
\newtheorem*{coherence_for_symmetric_monoidal_bicategories*}{Coherence for symmetric monoidal bicategories}
\newtheorem*{so2_action_on_fully_dualizable_objects*}{SO(2)-action on the space of fully dualizable objects}
\newtheorem*{cobordism_hypothesis_for_oriented_bordism*}{The Cobordism Hypothesis for two-dimensional oriented bordism}
\newtheorem*{thm_italics*}{Theorem}

\theoremstyle{remark}

\makeatletter
  \def\subsection{\@startsection{subsection}{1}%
  \z@{.7\linespacing\@plus\linespacing}{.5\linespacing}%
  {\normalfont\bfseries\centering}}
\makeatother

\let\oldtocsection=\tocsection
\let\oldtocsubsection=\tocsubsection
\let\oldtocsubsubsection=\tocsubsubsection
\renewcommand{\tocsection}[2]{\hspace{0em}\oldtocsection{#1}{#2}}
\renewcommand{\tocsubsection}[2]{\hspace{1em}\oldtocsubsection{#1}{#2}}
\renewcommand{\tocsubsubsection}[2]{\hspace{2em}\oldtocsubsubsection{#1}{#2}}

\begin{document}
\title[On dualizable objects in monoidal bicategories]{On dualizable objects in monoidal bicategories, framed surfaces and the Cobordism Hypothesis}
\author[Piotr Pstr\k{a}gowski]{Piotr Pstr\k{a}gowski}
\address{University of Bonn}
\email{pstragowski.piotr@gmail.com}

\begin{abstract}
We prove coherence theorems for dualizable objects in monoidal bicategories and for fully dualizable objects in symmetric monoidal bicategories, describing \emph{coherent dual pairs} and \emph{coherent fully dual pairs}. These are property-like structures one can attach to an object that are equivalent to the properties of dualizability and full dualizability.

We extend diagrammatic calculus of  surfaces of Christopher Schommer-Pries to the case of surfaces equipped with a framing. We present two equivalence relations on so obtained \emph{framed planar diagrams}, one which can be used to model isotopy classes of framings on a fixed surface and one modelling diffeomorphism-isotopy classes of surfaces.

We use the language of framed planar diagrams to derive a presentation of the framed bordism bicategory, completely classifying all two-dimensional framed topological field theories with arbitrary target. We then use it to show that the framed bordism bicategory is equivalent to the free symmetric monoidal bicategory on a coherent fully dual pair. In lieu of our coherence theorems, this gives a new proof of the Cobordism Hypothesis in dimension two.
\end{abstract}

\maketitle 

\tableofcontents

\addtocontents{toc}{\protect\setcounter{tocdepth}{-1}}

\section*{Introduction}

The theory of dual objects in symmetric monoidal categories has recently come into new light after the discovery of its intrinsic connection to topology of higher-dimensional manifolds. The Cobordism Hypothesis of Baez-Dolan, now a theorem of Jacob Lurie, can be roughly stated as follows.
\begin{cobordism_hypothesis*}[Baez-Dolan]
The framed bordism $n$-category is the free symmetric monoidal $n$-category on a fully dualizable object.
\end{cobordism_hypothesis*}
The framed bordism $n$-category can be informally described as having objects framed $0$-manifolds, morphisms framed $1$-bordisms between these and so on up to $n$-morphisms, which are given diffeomorphism-isotopy classes of framed $n$-bordisms. Composition is induced from gluing manifolds along their boundaries, the symmetric monoidal structure from disjoint union. 

The progress on the hypothesis was for a long time rather slow, possibly due to several reasons. On one hand, the theory of weak $n$-categories was not sufficiently well developed, on the other, no precise statement was known. This has now changed, several good definitions of a more general notion of an $(\infty,n)$-category are available and the groundbreaking work \cite{lurietqfts} of Jacob Lurie made it clear that one should focus on the \emph{property} of being fully dualizable. 

A common feature of most of these new definitions of higher categories is their \emph{non-algebraicity}, that is, in such a category one might for example know that for any two morphisms there is a composite, but this composite is not necessarily unique. Instead, one has a contractible space of choices, none of which is distinguished.

The situation is different in low dimensions, where a fully weak algebraic definition of a 2-category, that of a \emph{bicategory}, has been available since the 1960s. Similarly, an algebraic notion of a weak 3-category, that of a \emph{tricategory}, and related notions of \emph{monoidal} and \emph{symmetric monoidal bicategories}, have been already studied in late 1990s. Perhaps more importantly, in all of these cases there is an array of coherence theorems available, making it possible to perform explicit calculations.

It seems natural to construct bordism categories in this framework. This was done by Christopher Schommer-Pries in \cite{chrisphd} using symmetric monoidal bicategories, who also devised methods of obtaining presentations of such categories. One can then ask whether a direct proof of the Cobordism Hypothesis in dimension two is possible in this context, the main aim of current work is to show that the answer is affirmative. 

Our proof of the Cobordism Hypothesis can be essentially split into two parts, an algebraic and a geometric one. In the former, we describe explicitly monoidal bicategories satisfying universal properties like the one predicted by the Cobordism Hypothesis, these results can be understood as coherence theorems connected to the earlier work \cite{gurski_biequivalences_in_tricategories} of Nick Gurski. 

In the latter, geometric part, we extend the work of Christopher Schommer-Pries to the case of framed surfaces and obtain a complete classification of two-dimensional framed topological field theories. This classification together with one of the coherence theorems reduces the Cobordism Hypothesis to a comparison between two finitely generated symmetric monoidal bicategories, which we perform directly .

Before delving into the content of specific chapters, let us describe our results in a little bit more detail. We start by defining a collection of data, a \emph{dual pair}, that one can attach to a dualizable object in a monoidal bicategory such that this data alone is already enough to witness that the object is indeed dualizable. We term such collections of data satisfying some additional equations \emph{coherent} and prove the following theorem. 

\begin{thm_a*}[Coherence for dualizable objects in monoidal bicategories, Theorem \ref{thm:coherence_for_dualizable_objects}]
Let $\mathbb{M}$ be a monoidal bicategory. The forgetful homomorphism
\[ \pi: \mathcal{C}ohDualPair(\mathbb{M}) \rightarrow \mathbb{K}(\mathbb{M}^{d}) \]
between, respectively, the bicategory of coherent dual pairs in $\mathbb{M}$ and the groupoid of dualizable objects in $\mathbb{M}$, is a surjective on objects equivalence.
\end{thm_a*} 
In particular, an object can be completed to a coherent dual pair if and only if it is dualizable and any two coherent dual pairs living over a given object are equivalent. In other words, the theorem implies that a coherent dual pair is a property-like structure equivalent to dualizability. 

We then study the case of fully dualizable objects in symmetric monoidal bicategories. Again, we describe a collection of data, a \emph{fully dual pair}, that one can attach to a fully dualizable object and we prove that this data alone is enough to witness that the object is fully dualizable. We introduce additional equations describing \emph{coherent fully dual pairs} and prove the following.

\begin{thm_b*}[Coherence for fully dualizable objects in symmetric monoidal bicategories, Theorem \ref{thm:coherence_for_fully_dualizable_objects}]
Let $\mathbb{M}$ be a symmetric monoidal bicategory. The forgetful homomorphism 
\[ \pi: \mathcal{C}ohFullyDualPair(\mathbb{M}) \rightarrow \mathbb{K}(\mathbb{M}^{fd}) \]
between, respectively, the bicategory of coherent fully dual pairs in $\mathbb{M}$ and the groupoid of fully dualizable objects in $\mathbb{M}$, is a surjective on objects equivalence.
\end{thm_b*}
In both cases the additional equations we impose are inspired by the generating relations of the oriented bordism bicategory, this is important in applications to topological field theories.

In the geometric part of the current work we extend the diagrammatic calculus of oriented planar diagrams of Christopher Schommer-Pries to the case of framed surfaces, describing \emph{framed planar diagrams}. These are combinatorial objects that one can attach to an oriented surface equipped with a generic map into $\mathbb{R}^{2}$ and a sufficiently well-behaved framing, conversely, given its framed planar diagram, one can reconstruct the surface together with its map to the plane and the isotopy class of its framing. 

We introduce two equivalence relations on the set of framed planar diagrams, that of \emph{strong equivalence} and \emph{equivalence}, and prove the following. 

\begin{thm_c*}[Calculus of framed surfaces, Theorems \ref{thm:description_of_framings_on_surface} and \ref{thm:framed_planar_decomposition}]
Two framed planar diagrams $\aleph, \aleph ^\prime$ with the same underlying oriented planar diagram are strongly equivalent if and only if they describe isotopic framings on the oriented surface that can be reconstructed from them. 

Two framed planar diagrams $\aleph, \aleph ^\prime$ with possibly different underlying oriented planar diagrams are equivalent if and only if they describe diffeomorphic-isotopic framed surfaces.
\end{thm_c*}
We then use framed planar diagrams together with the relation of equivalence to obtain a presentation of the framed bordism bicategory, which we also carefully construct. Our proof is formally analogous to the one used in \cite{chrisphd} to establish presentations of unoriented and oriented bordism bicategories.

\begin{thm_d*}[Presentation of the framed bordism bicategory, Theorem \ref{thm:presentation_of_the_framed_bordism_bicategory}]
There exists a presentation $(G_{bord}^{fr}, \mathcal{R}_{bord}^{fr})$ of a symmetric monoidal bicategory consisting of a finite list of generating objects, morphisms, $2$-cells and a finite list of relations between them such that there is an equivalence
\[ \mathbb{F}_{bord}^{fr} \simeq \mathbb{B}ord_{2}^{fr} \]
between the free symmetric monoidal bicategory on the presentation and the framed bordism bicategory.
\end{thm_d*}
The presentation is a little bit larger than the oriented one, but we describe it explicitly. Note that due to the universal property of a freely generated symmetric monoidal bicategory, the theorem implies a complete classification of two-dimensional framed topological field theories with arbitrary target.

We then use the presentation to show directly that the framed bordism bicategory is equivalent to the free symmetric monoidal bicategory on a coherent fully dual pair and explain how this result together with coherence for fully dualizable objects implies the following formulation of the Cobordism Hypothesis in dimension two. 

\begin{thm_e*}[The Cobordism Hypothesis of Baez-Dolan, Theorem \ref{thm:cobordism_hypothesis_in_dimension_two}]
Let $\mathbb{B}ord_{2}^{fr}$ be the framed bordism bicategory and let $\mathbb{M}$ be any other symmetric monoidal bicategory. Then, the evaluation at the positively framed point homomorphism induces an equivalence
\[ \textbf{SymMonBicat}(\mathbb{B}ord_{2}^{fr}, \mathbb{M}) \rightarrow \mathbb{K}(\mathbb{M}^{fd}) \]
between the bicategory of framed topological field theories with values in $\mathbb{M}$ and the groupoid of fully dualizable objects in $\mathbb{M}$.
\end{thm_e*}
Lastly, we describe some possible extensions of the current work related to an explicit form of the $SO(2)$-action on the space of fully dualizable objects and to the Cobordism Hypothesis for two-dimensional oriented topological field theories.

\subsection*{Description of specific chapters}

Chapter 1 is introductory and illustrates the main ideas of the current work. We present the well-known proof of the Cobordism Hypothesis in dimension one as a consequence of two separate statements, coherence for dualizable objects in monoidal categories and classification of one-dimensional framed topological field theories.

Chapter 2 concerns the theory of dualizable objects in monoidal bicategories. We define a collection of data one can attach to a dualizable object, a \emph{dual pair}, whose existence is trivially equivalent to dualizability. We then describe additional equations on the components of a dual pair and term those pairs that satisfy them \emph{coherent}. We organize dual pairs in a fixed monoidal bicategory into a bicategory of their own and show that it is a $2$-groupoid. We then prove a strictification result that any dualizable object can be completed to a coherent dual pair and a coherence theorem which says that the forgetful homomorphisms from the bicategory of coherent dual pairs into the 2-groupoid of dualizable objects is an equivalence.

Chapter 3 develops the theory of fully dualizable objects in symmetric monoidal bicategories. We start by recalling the notion of a \emph{Serre autoequivalence}, which is a canonical automorphism one can attach to any fully dualizable object in a symmetric monoidal bicategory, and verify its basic properties. We then introduce the notions of fully dual pairs and coherent fully dual pairs, organize them into bicategories and show that these bicategories are $2$-groupoids. We prove strictification and coherence theorems analogous to the dualizable case. 

Chapter 4 concerns the general theory of bordism bicategories, as developed by Christopher Schommer-Pries. It does not contain any new results. We recall the construction of bordism bicategories and the related idea of haloed manifolds. We give a brief account of the calculus of planar diagrams and describe how it leads to the presentation of the unoriented bordism bicategory. 

Chapter 5 sees us constructing the framed bordism bicategory $\mathbb{B}ord_{2}^{fr}$, which intuitively can be described as having as objects framed $0$-manifolds,  framed $1$-bordisms as morphisms, and $2$-cells given by diffeomorphism-isotopy classes of framed $2$-bordisms. We prove some basic results concerning framings on haloed manifolds and introduce the generators of the framed bordism bicategory. 

Chapter 6 is devoted to diagrammatic calculus of framed surfaces. We introduce an extension of the notion of an oriented planar diagram into what we call a \emph{framed planar diagram}. We describe how a surface equipped with a generic map into $\mathbb{R}^{2}$ and a sufficiently generic framing leads to a framed planar diagram and consequently, how given such a diagram one can reconstruct a framed surface mapping into $\mathbb{R}^{2}$. We then introduce two different relations on framed planar diagrams, that of strong equivalence and equivalence. The former is a relation on framed planar diagrams with the same underlying oriented diagram, we show that two such diagrams are strongly equivalent if and only if they describe isotopic framings. We then build up on this result by showing that two arbitrary framed planar diagrams are equivalent if and only if the surfaces they describe are diffeomorphic-isotopic.

Chapter 7 concerns the classification of two-dimensional framed topological field theories. We first introduce \emph{relative framed planar diagrams}, which model framed $2$-bordisms and \emph{framed linear diagrams}, which model framed $1$-bordisms. We then apply the results of the previous chapter to derive a presentation of the framed bordism bicategory.

Chapter 8 is devoted to the proof of the Cobordism Hypothesis in dimension two. More precisely, we first show that the statement is equivalent to proving that the free symmetric monoidal bicategory on a coherent fully dual pair is equivalent to the framed bordism bicategory. We then prove that the strict homomorphism from the former to the latter classifying the fully dual pair of positively and negatively framed points is an equivalence.

In the Appendix \ref{appendix_free_monoidal_bicategories} we develop a variant on the theory of free monoidal bicategories, modeled on computadic symmetric monoidal bicategories of Christopher Schommer-Pries. The variant described has the advantage of allowing generating $1$-cells and generating $2$-cells whose sources and targets are only consequences of other generating data, and are not necessarily generating themselves. We also prove a group of technical results concering the truncation of generating data, these are rather simple-minded in nature, but are used in our proofs of coherence theorems to simplify bookkeeping.

\subsection*{Notation and terminology} 

A \emph{monoidal bicategory} is by definition a tricategory with one unnamed object. We use definition of the latter from \cite{gurski_algebraic_theory_of_tricategories}, which differs from the one given in \cite{coherence_for_tricategories_gordon_power_street} by the fact that it is fully algebraic and so all functors that are postulated to be equivalences come in form of adjoint equivalences. The related notions of \emph{symmetric monoidal bicategories}, \emph{monoidal homomorphisms}, \emph{symmetric monoidal homomorphisms}, \emph{transformations} and \emph{modifications} can be found in \cite{chrisphd}. 

We will always refer to functors between bicategories as \emph{homomorphisms}, they are assumed to be strong, that is, their constraint 2-cells are isomorphisms, but not necessarily identities. If they are identities, we will talk about \emph{strict homomorphisms}. The word \emph{functor} itself will be reserved for ordinary functors between categories.

If $\mathbb{M}$ is a monoidal bicategory, we will denote its monoidal product by $\otimes$ and by $I$ its monoidal unit. If it doesn't lead to confusion, we will also denote the monoidal product by juxtaposition. The associator will be denoted by $a$ and if $\mathbb{M}$ is symmetric, we will denote the symmetry by $b$. Unless working with \textbf{Gray}-monoids, we will not explicitly invoke constraint $2$-cells by name, instead relying on various coherence theorems to pin them down.

We will alternatively use the name of the object in question or the symbol $1$ to denote identity one-cells, constraint 2-cells witnessing naturality of some homomorphisms of bicategories will be denoted by the name of the homomorphism. 

If $\mathbb{B}$ is bicategory, by $\mathbb{K}(\mathbb{B})$ we denote the \emph{underlying 2-groupoid} of $\mathbb{B}$, that is, the bicategory obtained by disregarding all the morphisms that are not equivalences and all the $2$-cells that are not isomorphisms. By $\text{Ho}(\mathbb{B})$ we will denote its \emph{homotopy category}, which is obtained from $\mathbb{B}$ by identifying isomorphic $1$-cells. The categories of morphisms will be denoted by $\mathbb{B}(-,-) := \textnormal{Hom}_{\mathbb{B}}(-,-)$. 

When working with \textbf{Gray}-monoids, we will use the "first the maps on the left" convention. More specifically, when writing $f_{1} \otimes \cdots \otimes f_{n}$ we will always mean $(1 \otimes \ldots 1 \cdots f_{n}) \circ (1 \otimes \cdots f_{n-1} \otimes 1) \circ \cdots \circ (f_{1} \otimes 1 \otimes \cdots \otimes 1)$. We will denote the interchange isomorphism via $\Sigma _{f, g}: (1 \otimes g) (f \otimes 1) \Rightarrow (f \otimes 1)(1 \otimes g)$. We reserve the right to suppress it in the presence of a different 2-cell if it is clear from the context that it should be inserted.


Unless explicitly stated otherwise, by the word \emph{manifold} we mean a compact manifold, possibly with corners.

\subsection*{Coherence issues}

\begin{interpreting_pasting_diagrams*}
We will use the theory of coherence for bicategories to interpret our pasting diagrams and we will not name or draw any constraint 2-cells coming from the structure of a bicategory.  
\end{interpreting_pasting_diagrams*}

\begin{coherence_for_monoidal_bicategories*}
There are two related coherence results in the theory of monoidal bicategories which we will use. The first one is a classical theorem of Robert Gordon, Anthony Power and Ross Street from \cite{coherence_for_tricategories_gordon_power_street}, which in our language reads as follows.

\begin{thm_italics*}[Coherence for tricategories]
Any monoidal bicategory is equivalent to a \textbf{Gray}-monoid. 
\end{thm_italics*}
Additionally, we will sometimes implicitly use its following more direct consequence derived by Nick Gurski in \cite{gurski_algebraic_theory_of_tricategories}. 

\begin{thm_italics*}[Simple-minded coherence]
Let $\mathbb{M}$ be a monoidal bicategory and let $E$ be a locally discrete category-enriched graph. Any diagram of 2-cells in the image of a strict monoidal homomorphism $\mathbb{F}(E) \rightarrow \mathbb{M}$, where $\mathbb{F}(E)$ is the free monoidal bicategory on $E$, commutes. 
\end{thm_italics*}
To distinguish it from the classical coherence theorem, we have chosen to refer to the result above by \emph{simple-minded coherence}. This is motivated by its immediate applicability, as it can be used to assert that certain diagrams commute. We will think of any $2$-cells which are obviously in the image of a homomorphism as above as canonical, in which case we reserve a right not to label them. 
\end{coherence_for_monoidal_bicategories*}

\begin{coherence_for_symmetric_monoidal_bicategories*}
We will use two coherence results from the theory of symmetric monoidal bicategories. The first one, taken from \cite{gurski_osorno_coherence_for_symmetric_monoidal_bicategories}, is due to Gurski and Ang\'{e}lica Osorno, and reads as follows. 

\begin{thm_italics*}[Coherence of Gurski-Osorno]
In a free symmetric monoidal bicategory $\mathbb{F}(X)$ on a set of objects $X$, any diagram of $2$-cells commutes.
\end{thm_italics*}
This gives a new class of canonical coherence $2$-cells between composites of associators and symmetries. Due to classical statement on symmetric monoidal categories, two such composites are isomorphic if and only if they represent the same underlying permutation, the coherence theorem above implies that the isomorphism between them is unique. 

We will use extensively the theory of \emph{unbiased semistrict symmetric monoidal bicategories}, which are a semistrict version of symmetric monoidal bicategories introduced in \cite{chrisphd}. These are precisely strict enough to admit a variant of string diagram calculus, similar to the one for usual bicategories, except regions are labeled by tensor products of objects and one has special edges reprenting symmetries. We will use these diagrams to perform calculations. 

We will also make repeated use of the following result of Christopher Schommer-Pries, which allows us to replace any symmetric monoidal bicategory by an unbiased semistrict one in a very controlled manner. 

\begin{thm_italics*}[Coherence for unbiased semistrict symmetric monoidal bicategories]
Let $(G, \mathcal{R})$ be a generating datum for a symmetric monoidal bicategory, $\mathbb{F}(G, \mathcal{R})$ the free symmetric monoidal bicategory on the given datum and $\mathbb{U}(G, \mathcal{R})$ the free semistrict unbiased symmetric monoidal bicategory. Then, the canonical induced homomorphism
\[ \mathbb{F}(G, \mathcal{R}) \rightarrow \mathbb{U}(G, \mathcal{R}) \]
is an equivalence of symmetric monoidal bicategories.
\end{thm_italics*}

As an alternative to unbiased semistrict string diagrams, we will sometimes informally use spaced-out ordinary string diagrams for symmetric monoidal categories to draw composites of $2$-cells in a symmetric monoidal bicategory. This will not bother us, but an analogous approach to diagrammatic calculus of $2$-cells can be completely formalized in the case of so called \emph{quasistrict symmetric monoidal bicategories}, see \cite{bartlett_wire_diagrams}.
\end{coherence_for_symmetric_monoidal_bicategories*}

\subsection*{Acknowledgements}

This paper is a version of my Master's thesis, which I have written as a master student at Bonn University. 

First and foremost, I would like to thank my supervisor Chrisopher Schommer-Pries for his guidance, patience and many helpful comments. While working with him I have learned a lot of mathematics and, more importantly, a lot on what mathematics is about. I enjoyed our frequent meetings which made me never feel alone in my pursuit. 

I would also like to extend gratitude to Peter Teichner for his support, for first telling me about this topic and also for teaching me the Thom-Pontryagin construction. I am also very grateful to Stefan Schwede, whose mentoring helped me go through many difficult times during my master studies. 

I want to thank my family and friends, without whom this undertaking would not be possible at all, and the staff of Caf\'{e} Blau\footnote{Franziskanerstrasse 9, 53113 Bonn}, for creating a wonderful environment in which I could do my work and drink tea. 

\addtocontents{toc}{\protect\setcounter{tocdepth}{3}}

\section{The Cobordism Hypothesis in dimension one}

In this chapter we will illustrate the main ideas of the current work by presenting the well-known proof of the Cobordism Hypothesis in dimension one. It is not our goal to review all of the material necessary to understand subsequent chapters, but only to give a feel of the problems involved. None of the results presented here are new.

\subsection{Dualizable objects in monoidal categories}

In this section we will review the basic theory of dualizable objects in monoidal categories. We will introduce dual pairs, organize them into a category and prove a coherence theorem for dualizable objects. The result can be understood as a baby-version of our coherence theorems for dualizable objects in monoidal bicategories, as it is formally analogous. 

\begin{defin}
Let $(\mathbb{M}, \otimes, I)$ be a monoidal category. A \textbf{dual pair} $(X, Y, e, c)$ in $\mathbb{M}$ consists of an object $X$, which we call the \textbf{left dual}, an object $Y$, which we call the \textbf{right dual}, together with \textbf{evaluation} $e: X \otimes Y \rightarrow I$ and \textbf{coevaluation} $c: I \rightarrow Y \otimes X$ maps that satisfy \textbf{triangle equations} pictured below.

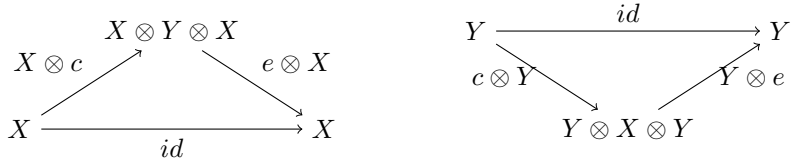
\begin{center}
	\begin{tikzpicture}
		\node (X1) at (0, 0) {$ X $};
		\node (XYX1) at (2,1.3) {$ X \otimes Y \otimes X $};
		\node (X2) at (4, 0) {$ X $};
		
		\draw [->] (X1) to node[auto] {$ X \otimes c $} (XYX1);
		\draw [->] (XYX1) to node[auto] {$ e \otimes X $} (X2);
		\draw [->] (X1) to node[below] {$ id $} (X2);
		
		\node (Y1) at (6,1.3) {$ Y $};
		\node (YXY1) at (8, 0) {$ Y \otimes X \otimes Y $};
		\node (Y2) at (10, 1.3) {$ Y $};
		
		\draw [->] (Y1) to node[left] {$ c \otimes Y $} (YXY1);
		\draw [->] (YXY1) to node[right] {$ Y \otimes e $} (Y2);
		\draw [->] (Y1) to node[auto] {$ id $} (Y2);
	\end{tikzpicture}
\end{center}
An object $X \in \mathbb{M}$ is \textbf{left dualizable} if it can be completed into a dual pair as a left dual.
\end{defin}

\begin{rem}
Similarly, one defines an object to be right dualizable if it can be completed to a dual pair as a right dual. Their theory is formally analogous and so we will focus exclusively on left dualizability, referring to it simply as \emph{dualizability}. 
\end{rem}

One could say that an object $X$ is \emph{monoidally invertible} if there exists some $Y$ such that $X \otimes Y$ and $Y \otimes X$ are both isomorphic to the monoidal unit $I$. It is possible to show that such an object is in particular dualizable, but the latter property is much weaker. 

In practice, it is often the objects that satisfy some finiteness conditions that are dualizable. For example, in the category $\mathcal{V}ect_{k}$ of vector spaces with the monoidal structure coming from the tensor product, it is exactly the finite-dimensional vector spaces that are dualizable. On the other hand, only one-dimensional vector spaces are monoidally invertible. 

The triangle equations can be drawn using string diagrams, which reveal the connection with topology of $1$-manifolds. If we denote evaluation with the right elbow and coevaluation with the left elbow, we can picture the equations as in \textbf{Figure \ref{fig:triangle_equations_in_string_diagram_notation}}.

\begin{figure}[htbp!]
	\begin{tikzpicture}[thick, scale=0.8, every node/.style={scale=0.8}]
	\node (X1) at (2, -0) {$ X $};
	\node (Y1) at (2, -0.8) {$ Y $};
	\node (X2) at (2, -1.6) {$ X $};
	
	\node (X3) at (4, -0.8) {$ X $};
	\node (X4) at (5, -0.8) {$ X $};
	
	\draw (X1) to [out = 0, in = 0] (Y1);
	\draw (Y1) to [out = 180, in = 180] (X2);
	\draw (X3) to (X4);
	
	\node at (3.2, -0.83) {$ = $};

	\node (rY1) at (8, -0) {$ Y $};
	\node (rX1) at (8, -0.8) {$ X $};
	\node (rY2) at (8, -1.6) {$ Y $};
	
	\node (rY3) at (10, -0.8) {$ Y $};
	\node (rY4) at (11, -0.8) {$ Y $};
	
	\draw (rY1) to [out=180, in=180] (rX1);
	\draw (rX1) to [out=0, in=0] (rY2);
	\draw (rY3) to (rY4);
	
	\node at (9.2, -0.83) {$ = $};
	\end{tikzpicture}
\caption{Triangle equations in string diagram form}
\label{fig:triangle_equations_in_string_diagram_notation}
\end{figure}
Observe that being dualizable is a \emph{property}. More specifically, it entails existence of some amount of data, but does not specify it. A given dualizable object can be usually completed to many different dual pairs and thus it is a valid question to ask whether all such pairs are in fact in some sense "the same". 

If that was the case, we would say that the structure of a dual pair is \emph{property-like}, as the "space" of dual pairs over a given object would always be either empty or contractible. The terminology is motivated by the fact that in such a case the association 
\[ X \leadsto (\textrm{Collection of dual pairs over} \ X) \]
takes in essence values in $\{ 0, 1 \}$, where $0$ denotes the empty space and $1$ denotes a space consisting of one-point. Thus, it \emph{does} behave like a property, as the latter could be in principle defined as this sort of function.

Note that although our use of this terminology is somewhat informal, the category-theoretic practice of calling structures satisfying some form of uniqueness \emph{property-like} dates back many years, see for example \cite{kelly_lack_on_property_like_structures}.

To be able to compare different dual pairs, we have to first organize them into a category, which is not difficult.

\begin{defin}
If $\mathbb{M}$ is a monoidal category and $(X_{1}, Y_{1}, e_{1}, c_{1})$ and $(X_{2}, Y_{2}, e_{2}, c_{2})$ are dual pairs in $\mathbb{M}$, then a \textbf{morphisms of dual pairs} consists of arrows $x: X_{1} \rightarrow X_{2}$, $y: Y_{1} \rightarrow Y_{2}$ that are natural with respect to (co)evaluation maps. In other words, we have $e_{2} \circ (x \otimes y) = e_{1}$ and $c_{2}= (y \otimes x) \circ c_{1}$. 
\end{defin}
We will denote the category of dual pairs in $\mathbb{M}$ by $\mathcal{D}ualPair(\mathbb{M})$. Another possible approach to defining it would be to consider the \emph{free monoidal category} on a dual pair, this is the one which is the easiest to generalize to bicategories and so we present it in some detail. Observe that the collection of dual pairs is functorial along strict homomorphisms of monoidal categories, hence the following characterization is possible.

\begin{prop}
There exists a monoidal category $\mathbb{F}_{d}$, the \textbf{free monoidal category on a dual pair}, such that for any other monoidal category there is a bijection 
\[ \textbf{MonCat}_{strict}(\mathbb{F}_{d}, \mathbb{M}) \simeq \mathcal{D}ualPair(\mathbb{M}) \]
between strict monoidal homomorphisms $\mathbb{F}_{d} \rightarrow \mathbb{M}$ and the set of dual pairs in $\mathbb{M}$, natural with respect to strict homomorphisms.
\end{prop}
By Yoneda lemma, such a monoidal category is unique up to a canonical invertible strict homomorphism. Observe that this endows $\mathcal{D}ualPair(\mathbb{M})$, which is here a priori only a set, with a natural structure of a category, as the left hand side of the bijection forms a category under monoidal transformations. It is not difficult to check that this structure coincides with the one we already described.

The existence of the category $\mathbb{F}_{d}$ can be verified on abstract grounds, however it is perhaps more instructive to construct it explicitly. Since a strict homomorphism $\mathbb{F}_{d} \rightarrow \mathbb{M}$ is supposed to in particular correspond to two objects $X, Y \in \mathbb{M}$, we define $Ob(\mathbb{F}_{d})$ to be the set of all \emph{binary words} in $X, Y$. This set can be described inductively by saying that $I$ is a binary word, $X, Y$ are binary words and that if $w_{1}, w_{2}$ are binary words, so is $w_{1} \otimes w_{2}$. 

To define morphisms, one first defines \emph{generating morphisms} to be the maps $e: X \otimes Y \rightarrow I, c: I \rightarrow Y \otimes X$ and all the constraint isomorphisms required by the structure of a monoidal category. That is, for each triple $A, B, C \in \mathbb{F}_{d}$ of objects we have generating associators $a: (A \otimes B) \otimes C \rightarrow A \otimes (B \otimes C)$ and their inverses, similarly we have left and right unitors. 

One then constructs the set of \emph{expressions in generating morphisms} by freely closing the latter under tensor product and composition, the same way the set of binary words freely closed $X, Y, I$ under tensor product of objects. The arrows in $\mathbb{F}_{d}$ are equivalence classes of such expressions, the equivalence relation is the smallest one ensuring that the constraint cells are natural, the axioms of a monoidal category hold and that $e, c$ satisfy triangle equations.

\begin{rem}
The approach to construction of $\mathbb{F}_{d}$ we presented above readily generalizes, in a similar way one could construct the free braided monoidal category on a dual pair $\mathbb{F}_{d}^{b}$ and the free symmetric monoidal category on a dual pair $\mathbb{F}_{d}^{s}$, both having an analogous universal property with respect to strict braided homomorphisms.
\end{rem}

The question of whether the structure of a dual pair is property-like can now be answered in affirmative in the form of the following theorem.

\begin{thm}[Coherence for dualizable objects in monoidal categories]
\label{thm:coherence_for_dualizable_objects_in_monoidal_categories}
Suppose 
\[ (x, y): (X_{1}, Y_{1}, e_{1}, c_{1}) \rightarrow (X_{2}, Y_{2}, e_{2}, c_{2})\] 
is a map of dual pairs. Then $x: X_{1} \rightarrow X_{2}$ is an isomorphism and, conversely, for any such isomorphism there is a unique $y: Y_{1} \rightarrow Y_{2}$ that completes it to a map of dual pairs. More concisely, the forgetful functor
\begin{gather*}
\mathcal{D}ualPair(\mathbb{M}) \rightarrow \mathbb{K}(\mathbb{M}^{d}) \\
(X, Y, e, c) \mapsto X
\end{gather*}
from the category of dual pairs into the groupoid of dualizable objects in $\mathbb{M}$ is a surjective on objects equivalence of categories.
\end{thm}
Before going into the proof, let us observe that one of the main results of the current work will be an analogous result for dualizable objects in \emph{monoidal bicategories}. In this context, the most naive definition of a \emph{dual pair} will not satisfy the conditions of the above theorem and we will have to restrict to a class of dual pairs satisfying some additional coherence equations.

\begin{proof}
If $a: X_{2} \rightarrow X_{1}$ and $b: Y_{1} \rightarrow Y_{2}$ are morphisms, then we define their duals $\widetilde{a}, \widetilde{b}$ using the string diagrams
\begin{center}
.
\end{center}
The triangle equations imply that this establishes inverse bijections $\mathbb{M}(Y_{1}, Y_{2}) \simeq \mathbb{M}(X_{2}, X_{1})$. We claim that $\widetilde{y}$ is inverse to $x$. Indeed, one composite is equal to identity by the easy computation

\begin{center}
	%
,
\end{center}
where we have used the fact that $y \otimes x$ commutes with coevaluation. For the other we compute
\begin{center}
	%
,
\end{center}
using that $x \otimes y$ commutes with evaluation. This also shows that there is at most one such $y$, as it must then be necessarily equal to $\widetilde{x^{-1}}$, as taking dual twice recovers the morphism.

Conversely, one verifies easily that if $x:X_{1} \rightarrow X_{2}$ is invertible, then setting $y = \widetilde{x^{-1}}$ does indeed complete it to a morphism of dual pairs. 
\end{proof}

\subsection{Framed bordism category}
In this section we will finish the proof of the Cobordism Hypothesis in dimensione one. As we have already established the needed categorical preliminaries, what is missing is a classification result for one-dimensional framed topological field theories. Let us first recall the relevant notions.

\begin{defin}
If $M$ is an $n$-manifold and $k \geq n$, then a \textbf{$k$-framing} $v$ of $M$ is a trivialization of the vector bundle $T_{M}^{k} = T_{M} \oplus \mathbb{R}^{k-n}$, that is, a collection $v_{1}, \ldots, v_{k}$ of sections of $T_{M}^{k}$ such that at all points $x \in M$ the vectors $v_{1}(x), \ldots, v_{k}(x)$ form a basis of $T_{M, x}^{k}$. An \textbf{isotopy} of $k$-framings is a family of such sections parametrized by the interval. 
\end{defin}
Observe that $k$-framings can be always restricted to the boundary, as if $M$ is a manifold with boundary, then $T_{\partial M} \oplus \mathbb{R} \simeq T_{M} |_{\partial M}$ \emph{almost} canonically, with minor ambiguity related to the fact that to obtain this splitting one has to first choose a collar. 

In dimension one, the framed bordism category can be easily constructed as an ordinary symmetric monoidal category. 

\begin{defin}
The \textbf{framed bordism category} $\mathbb{B}ord_{1}^{fr}$ is a symmetric monoidal category with objects $1$-framed $0$-manifolds and morphisms given by diffeomorphism-isotopy classes of $1$-framed bordisms. Composition is given by gluing of bordisms along the boundary, the symmetric monoidal structure comes from the disjoint union of manifolds.
\end{defin}
Let us describe maps in the framed bordism category in a little bit more detail. A morphism $M \rightarrow N$ is given by an equivalence class of framed $1$-manifolds $w$ together with diffeomorphisms $w \simeq M \sqcup N$ such that the framing of $W$ restricts to the given framing on the boundary. Two such $1$-manifolds are considered to be equivalent if there is a diffeomorphism between them preserving the decomposition of the boundary and preserving framings up to isotopy relative to the boundary. 

Observe that to restrict the framing of $w$ to its boundary, we need to choose isomorphisms $T_{M} \oplus \mathbb{R} \simeq T_{w} |_{M}$ and $T_{N} \oplus \mathbb{R} \simeq T_{w} |_{N}$, we agree to use the isomorphism obtained from the inward vector for the domain $M$ and from the outward vector for the codomain $N$. Even then, these isomorphisms are still defined only up to a positive scalar, but with a little care one can show that these difficulties can be ignored and one still gets a well-defined category. 

Our task will be now to classify all one-dimensional framed topological field theories or, in other words, symmetric monoidal homomorphisms out of the framed bordism category. By a \emph{presentation} of $\mathbb{B}ord_{1}^{fr}$ we will mean an equivalence with some freely generated symmetric monoidal category $\mathbb{F}$. Observe that due to the universal property of the latter, such an equivalence would lead to a compact description of the category of symmetric monoidal homomorphisms out of $\mathbb{B}ord_{1}^{fr}$ with arbitrary target, as such a description is available for homomorphisms with source $\mathbb{F}$. Thus, a presentation implies a \emph{complete classification} we are after, at least up to equivalence.

\begin{figure}[htbp!]
	\begin{tikzpicture}
	
	\begin{scope}
		\draw[very thick] (0, 0) [out=0, in=90] to (1, -0.5) [out=-90, in=0] to (0, -1);
		
		\draw [
	    		decoration={
	    	    	brace,mirror,
	      		raise=.55cm
		},
		decorate
		] (-0.3, -1) -- (1.3,-1);
		\node at (0.5, -2)  {Right elbow};
	\end{scope}

	\begin{scope}[xshift=3cm]
		\draw [very thick] (1, 0) [out=180, in=90] to (0, -0.5) [out=-90, in=180] to (1, -1);
		
		\draw [
	    		decoration={
	    	    	brace,mirror,
	      		raise=.55cm
		},
		decorate
		] (-0.3, -1) -- (1.3,-1);
		\node at (0.5, -2)  {Left elbow};
	\end{scope}
	\end{tikzpicture}
\caption{Left and right elbow $1$-bordisms}
\label{fig:left_and_right_elbow_1_framed_bordisms}
\end{figure}
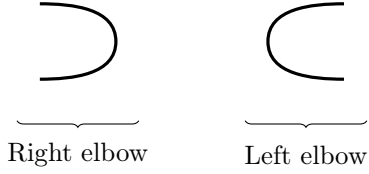

Let us first fix some notation. By the \textbf{positively framed point} $pt_{+}$ we mean a framed $0$-manifold consisting of a single point $pt$ together with the trivialization $1 \in \mathbb{R} \simeq \mathbb{R} \oplus T_{pt}$ of its one-tangent space, by the \textbf{negatively framed point} we mean a single point together with the opposite trivialization $-1 \in \mathbb{R}$. The \textbf{left and right elbow} are the $1$-bordisms presented in \textbf{Figure \ref{fig:left_and_right_elbow_1_framed_bordisms}}.

Topologically, they are both intervals, but they decomposition of the boundary is chosen so that the right elbow is a $1$-bordism $pt_{+} \sqcup pt_{-} \rightarrow \emptyset$ and the left elbow is a $1$-bordism $\emptyset \rightarrow pt_{-} \sqcup pt_{+}$. There is a unique class of isotopies of framings on them compatible with the given framings on the boundary and so this characterization specifies well-defined arrows in $\mathbb{B}ord_{1}^{fr}$. 

Note that the left and right elbow are precisely the elementary bordisms in dimension $1$, corresponding to respectively a critical point of index $0$ or $1$. Observing that their domains and codomains coincide with those of (co)evaluation maps leads to the the following classical result. 

\begin{thm}[Presentation of the framed bordism category]
\label{thm:presentation_of_the_framed_bordism_category}
The framed bordism category $\mathbb{B}ord_{1}^{fr}$ is freely generated, as a symmetric monoidal category, by the dual pair consisting of the positively and negatively framed points and left and right elbows. More precisely, the induced strict homomorphism
\[ \mathbb{F}^{s}_{d} \rightarrow \mathbb{B}ord_{1}^{fr} \]
from the free symmetric monoidal category on a dual pair is an equivalence.
\end{thm}

\begin{proof}
We have to verify that the induced homomorphism is essentially surjective on objects and fully faithful.

Essential surjectivity is equivalent to saying that any framed $0$-manifold $A$ is isomorphic to some disjoint union of positively and negatively framed points. Since up to isotopy$1$-framings are classified by their orientation, $A$ is certainly diffeomorphic to such a disjoint union with a diffeomorphism that preserves framings up to isotopy. Once such an isotopy is chosen, it can be spread out along a framed $1$-bordism with underlying manifold $A \times I$, the resulting map in $\mathbb{B}ord_{1}^{fr}$ will be an isomorphism with an explicit inverse given by the framed $1$-manifold $A \times I$ constructed from the inverse isotopy.

To establish fullness, we have to verify that any $1$-bordism $w$ between disjoint unions of positively and negatively framed points is in the image of the homomorphism. This can be done by choosing a Morse function $w \rightarrow I$ with disjoint critical values, the preimages of sufficiently fine covering of the interval will then decompose $w$ into elementary bordisms, which are precisely the left and right elbow.

Faithfulness is equivalent to proving that any two different decompositions of $w$ into left and right elbows can be related by a sequence of applications of triangle equations. However, the latter correspond to Morse birth-death singularities and the result follows from classical Cerf theory. 
\end{proof}

\begin{cor}[The Cobordism Hypothesis in dimension one] 
Let $\mathbb{B}ord_{1}^{fr}$ be the framed bordism category, let $\mathbb{M}$ be arbitrary symmetric monoidal category. Then, the evaluation at the positive point induces an equivalence 
\[ \textbf{SymMonCat}(\mathbb{B}ord_{1}^{fr}, \mathbb{M}) \rightarrow \mathbb{K}(\mathbb{M}^{d}) \]
between framed one-dimensional topological field theories with values in $\mathbb{M}$ and the groupoid of dualizable objects in $\mathbb{M}$.
\end{cor}

\begin{proof}
Observe that the statement is true for $\mathbb{F}^{s}_{d}$ by the coherence theorem for dualizable objects. It then follows also for $\mathbb{B}ord_{1}^{fr}$, as the two symmetric monoidal categories are equivalent.
\end{proof}

The important lesson here is that to establish the Cobordism Hypothesis in dimension one, we needed \emph{both} results, that is, coherence for dualizable objects and the presentation theorem. The first part can be understood as establishing the needed knowledge about the property of dualizability, the second classifies framed topological field theories.

Our approach to a proof in in dimension two is essentially the same. We will spend the next two chapters on the study of the properties of dualizability and full dualizability in the setting of symmetric monoidal bicategories, our goal will be a coherence theorem which will identify the symmetric monoidal bicategory satisfying the conditions of the Cobordism Hypothesis. The rest of the current work will be devoted to classification of two-dimensional framed topological field theories, mirroring results from \cite{chrisphd} on oriented and unoriented variants.

\section{Dualizable objects in monoidal bicategories}

This chapter is entirely devoted to the theory of dualizable objects in monoidal bicategories. The latter is a property of objects, generalizing the usual notion in ordinary monoidal categories. Our main result is description of a property-like structure that one can attach to an object in a monoidal bicategory,  namely that of a \emph{coherent dual pair}, equivalent to the property of being dualizable. 

\subsection{Dual pairs}

In this section we will describe the structure of a \emph{dual pair}, whose existence is trivially equivalent to dualizability. We then introduce a class of dual pairs satisfying additional equations, which we call \emph{coherent}, and prove a strictification result which says that any dualizable object can be completed to a coherent dual pair. We do the latter under the technical assumption that the ambient monoidal bicategory is a \textbf{Gray}-monoid, we will show later how one can deduce the general case from this partial result. 

We start by generalizing the notion of dualizability to monoidal bicategories. To do so, we use the homotopy category $\text{Ho}(\mathbb{M})$, which inherits an ordinary monoidal structure.

\begin{defin}
Let  $\mathbb{M}$ be a monoidal bicategory. An object $L \in \mathbb{M}$ is \textbf{left dualizable} if it is left dualizable as an object of $\text{Ho}(\mathbb{M})$. 
\end{defin}

\begin{rem}
Like in the case of monoidal categories, there is an analogous notion of right dualizability. We will exclusively focus on left dualizability and refer to it simply as \emph{dualizability}. 
\end{rem}

In the case of monoidal bicategories, like in the ordinary case, the property of being dualizable can be rephrased as ensuring existence of some auxiliary data. It is rather easy to specify the latter, as it arises as a direct categorification of the notion of a dual pair, where the triangle equations are naively replaced by invertible $2$-cells. 
\begin{defin}
A \textbf{dual pair} in a monoidal bicategory $\mathbb{M}$ is a tuple $(L, R, e, c, \alpha, \beta)$, where $L, R \in \mathbb{M}$ are objects, $e: L \otimes R \rightarrow I$, $c: I \rightarrow R \otimes L$ are 1-cells and $\alpha, \beta$ are isomorphisms 

\begin{center}
	\begin{tikzpicture}[thick, scale=0.8, every node/.style={scale=0.8}]		
		\node (L1) at (-1,0.3) {$ L $};
		\node (LI) at (-0.3,1.3) {$ L \otimes I $};
		\node (LRL1) at (2,2) {$ L \otimes (R \otimes L) $};
		\node (LRL2) at (5,2) {$ (L \otimes R) \otimes L $};
		\node (IL) at (7.3,1.3) {$ I \otimes L $};
		\node (L2) at (8,0.3) {$ L $};
		\draw [->] (L1) to node [auto] {$ r^{\bullet} $} (LI);
		\draw [->] (LI) to node [auto] {$ L \otimes c $} (LRL1);
		\draw [->] (LRL1) to node [auto] {$ a^{\bullet} $} (LRL2);
		\draw [->] (LRL2) to node [auto] {$ e \otimes L $} (IL);
		\draw [->] (IL) to node [auto] {$ l $} (L2);
		\draw [->] (L1) to [out=-10, in=-170] node [auto] {$ L $} (L2);
		\node at (3.5,1) {$ \Downarrow \alpha $};
	\end{tikzpicture}
\end{center}

\begin{center}
	\begin{tikzpicture}[thick, scale=0.8, every node/.style={scale=0.8}]		
		\node (R1) at (-1,0.3) {$ R $};
		\node (IR) at (-0.3,1.3) {$ I \otimes R $};
		\node (RLR1) at (2,2) {$ (R \otimes L) \otimes R $};
		\node (RLR2) at (5,2) {$ R \otimes (L \otimes R) $};
		\node (RI) at (7.3,1.3) {$ R \otimes I $};
		\node (R2) at (8,0.3) {$ R $};
		\draw [->] (R1) to node [auto] {$ l^{\bullet} $} (IR);
		\draw [->] (IR) to node [auto] {$ c \otimes R $} (RLR1);
		\draw [->] (RLR1) to node [auto] {$ a $} (RLR2);
		\draw [->] (RLR2) to node [auto] {$ R \otimes e $} (RI);
		\draw [->] (RI) to node [auto] {$ r $} (R2);
		\draw [->] (R1) to [out=-10, in=-170] node [auto] {$ R $} (R2);
		\node at (3.5,1) {$ \Downarrow \beta $};
	\end{tikzpicture}.
\end{center}
We will refer to $e$ as \textbf{evaluation}, to $c$ as \textbf{coevaluation} and to $\alpha, \beta$ as \textbf{cusp isomorphisms}.
\end{defin}

\begin{notation*}
By abuse of language we will sometimes refer to the whole dual pair $(L, R, e, c, \alpha, \beta)$ just by referring to the underlying objects, we will then denote it by $\langle L, R \rangle _{d}$.
\end{notation*}
The terminology of \emph{cusp isomorphisms} comes from the oriented bordism bicategory, which we will describe later in the text. In that monoidal bicategory one has a duality between positive and negative points. The (co)evaluation maps are the usual left/right elbows and cusp isomorphisms are given by surfaces that exhibit a cusp when drawn as embedded in space. 

\begin{prop}
An object $L \in \mathbb{M}$ is dualizable if and only if it can be made part of a dual pair $\langle L, R \rangle _{d}$ for some choice of $R$, (co)evaluation maps and cusp isomorphisms.
\end{prop}
\begin{proof}
Clearly if an object is a part of a dual pair then it is dualizable. Conversely, assume that $L$ is dualizable. This means that one can complete it to a dual pair in $\text{Ho}(\mathbb{M})$, giving us the right dual and (co)evaluation maps. The fact that the triangle equations hold in the homotopy category immediately implies existence of some isomorphisms that witness it, which we can take to be cusp isomorphisms.
\end{proof}

A basic property of dualizable objects is an obvious categorification of the usual result on bijections between certain Hom-sets for dual pairs in ordinary monoidal categories, obtained by taking duals of morphisms. In monoidal bicategories, taking duals leads to an \emph{equivalence} of categories. This is one of main technical tools we will use in our study of dualizable objects. 

\begin{lem}[On equivalences of Hom-categories]
\label{lem:hom_equivalences}
Let $A, B \in \mathbb{M}$ be objects and let $\langle L, R \rangle _{d}$ be a dual pair. Then, the functor $\mathbb{M}(A \otimes R, B) \rightarrow \mathbb{M}(A, B \otimes L)$ given by $f \mapsto (f \otimes L) \circ  (A \otimes c)$ is an equivalence of categories with an explicit pseudoinverse $g \mapsto (B \otimes e) \circ (g \otimes R)$.

Similarly, we have an equivalence $\mathbb{M}(L \otimes A, B) \rightarrow \mathbb{M}(A, R \otimes B)$ given by $f \mapsto (R \otimes f) \circ (c \otimes A)$ with an explicit pseudoinverse $g \mapsto (e \otimes B) \circ (L \otimes g)$.
\end{lem}

\begin{proof}
The needed natural isomorphisms between the relevant composites and identites are induced by cusp isomorphisms.
\end{proof}
Observe that in the statement of the lemma above we have supressed the obvious associators and unitors that should be included in the definition of the functors so that they are well-formed. By simple-minded coherence for monoidal bicategories, any choices one would have to make do not matter, as all the different composites of associators and unitors are uniquely isomorphic.

\begin{defin}
\label{defin_coherent_dual_pair}
A dual pair $\langle L, R \rangle _{d} = (L, R, e, c, \alpha, \beta)$ is said to be \textbf{coherent} if the \textbf{Swallowtail composites} depicted in \textbf{Figures \ref {fig:axiom_1_coherent_dual_pair}, \ref{fig:axiom_2_coherent_dual_pair}} are both identities.

\begin{figure}[htbp!]

	
\caption{Swallowtail composite (E)}
\label{fig:axiom_1_coherent_dual_pair}
\end{figure}

\begin{figure}[htbp!]
	%
	
\caption{Swallowtail composite (C)}
\label{fig:axiom_2_coherent_dual_pair}
\end{figure}
\end{defin}
The Swallowtail equations were first introduced by Dominic Verity in \cite{verity_enriched_categories_internal_categories_and_change_of_base} in the context of \textbf{Gray}-categories and later generalized to tricategories by Nick Gurski, see \cite{gurski_biequivalences_in_tricategories}.

We will now prove the strictification theorem, which says that given a fixed choice of $L, R$ and (co)evaluation maps $e, c$, if the cusp isomorphisms $\alpha, \beta$ can be chosen at all, then one can also choose them in a "particularly nice" way. By "nice" we mean a choice that would make the resulting dual pair coherent, it then follows that any dualizable object can be completed to a coherent dual pair.

\begin{thm}[Strictification for dual pairs]
\label{thm:strictification_for_dual_pairs}
If a dual pair $\langle L, R \rangle _{d}$ satisfies either of the Swallowtail equations, then it satisfies both of them. 

Moreover, if we keep the objects and (co)evaluation maps fixed, then for any choice of a cusp isomorphism $\alpha$ there is a unique cusp isomorphism $\beta$ such that together they satisfy both Swallowtail equations. In particular any dual pair can be made coherent by only a change of $\beta$. 
\end{thm}

\begin{cor}
Any dualizable object can be made part of a coherent dual pair.
\end{cor}
The strictification theorem was first proven by Nick Gurski in \cite{gurski_biequivalences_in_tricategories} under the additional assumption that the (co)evaluation maps $e, c$ are equivalences. The proof given there is very different, using the Yoneda embedding to reduce to the case of tricategory of bicategories, where a direct construction is possible.

Notice that our statement is very precise in saying how much freedom one has in choosing a pair of cusp isomorphisms satisfying Swallowtail equations, saying that one can first choose $\alpha$ in an arbitrary way and this will fix a unique $\beta$ with the required properties. This additional determination of data will be important later in our proof of the coherence theorem. 

In this section we only prove the theorem under an additional assumption that the ambient monoidal bicategory $\mathbb{M}$ is a \textbf{Gray}-monoid,  we will later show how to deduce the general result from this special case by organizing coherent dual pairs into a bicategory. 

\begin{proof}[Proof (of strictification for dual pairs in \textbf{Gray}-monoids).]
We first observe that the Swallowtail identities take a very simple form if the background monoidal bicategory $\mathbb{M}$ is a \textbf{Gray}-monoid. Namely, Swallowtail identity (C) reads

\begin{center}

\end{center}
and similarly identity $(E)$ is 

\begin{center}
	%
.
\end{center}
We will show that equation $(E)$ implies equation $(C)$, the other case is analogous. By \textbf{Lemma \ref{lem:hom_equivalences}}, it is enough to show that 

\begin{center}
	%
,
\end{center}
as the functor $(R \otimes e) \circ (- \otimes R)$ is an equivalence on the respective Hom-categories. Using the naturality of the interchanger we can rewrite the left hand side of the equation as 

\begin{center}
	%
.
\end{center}
Similarly, the right hand side can be rewritten as 

\begin{center}
	%
\end{center}
by first using the Swallowtail identity $(E)$ to replace $R \alpha R$ by $LR \beta$ and some coherence and then using the naturality of the interchanger exactly like we did during the manipulation of the left hand side. Thus, we have shown that under the assumption of $(E)$, equation $(C)$ is equivalent to one of the form 

\begin{center}
	%
.
\end{center}
The equality of 

\begin{center}
	%
.
\end{center}
follows from the simple fact that in any ordinary monoidal category, if  $A \rightarrow I$ is an isomorphism with the monoidal unit, then $A \simeq A \otimes I \rightarrow A \otimes A = A \simeq I \otimes A \rightarrow A \otimes A$. Here we apply it to the the isomorphism $\beta$ in the monoidal category $\mathbb{M}(R, R)$. 

Now we only have to verify that the coherence cells on both sides coincide. Unfortunately, we cannot apply coherence for monoidal bicategories directly, as the relevant composites involve non-trivially the unit. However, one rewrites 

\begin{center}
	%
\end{center}
using the compatibility of the interchanger with composition. Coherence now applies to any of the inside regions of the right hand sides and we are done.
\end{proof}

\subsection{Organizing dual pairs into a bicategory}
Our main goal in this section will be to organize dual pairs in some fixed ambient monoidal bicategory $\mathbb{M}$ into a bicategory in itself. We establish a basic property of these bicategories, which is that they are weak $2$-groupoids. As an application we show how to deduce the full statement of the strictification theorem from the special case of \textbf{Gray}-monoids, which we have already proven. 

To organize dual pairs into a bicategory, we will use the language of freely generated monoidal bicategories. There are many ways one can formalize this notion, we will use one modeled on the language of \emph{computadic symmetric monoidal bicategories} from \cite{chrisphd}. In essence, one allows generating objects, $1$-cells and $2$-cells, where the sources and targets of $1$-cells and $2$-cells are allowed not to be precisely generating lower cells, but only their consequences. One also allows relations, but only on the level of $2$-cells. 

The relevant definitions of \emph{generating datum}, \emph{freely generated monoidal bicategories} and of \emph{bicategories of shapes} can be found in \textbf{Appendix \ref{appendix_free_monoidal_bicategories}}. Familiarity with these notions is not strictly necessary for understanding of the results, as we will also give explicit descriptions of the bicategories of dual pairs. However, we do warn the reader that some simple-minded technical results from the appendix do occur in proofs. 

\begin{defin}
Let $G_{d}$ by a free generating datum for a monoidal bicategory consisting of

\begin{itemize}
\item two generating objects $L, R$,
\item two generating maps $e: L \otimes R \rightarrow I$ and $c: I \rightarrow R \otimes L$,
\item four generating $2$-cells $\alpha, \beta, \alpha^{-1}$ and $\beta^{-1}$, where $\alpha, \beta$ have sources and targets like exactly like the cusp isomorphisms and $\alpha^{-1}, \beta^{-1}$ the opposite.
\end{itemize}

Let $\mathcal{R}_{d}$ be a class of relations on the free monoidal bicategory $\mathbb{F}(G_{d})$ consisting of the relations that $\alpha, \alpha^{-1}$ and $\beta, \beta^{-1}$ are inverse to each other, let $\mathcal{R}_{cd}$ contain these relations and additionally the Swallowtail identities. 

We call the freely generated monoidal bicategory $\mathbb{F}_{d} := \mathbb{F}(G_{d}, \mathcal{R}_{d})$ the \textbf{free monoidal bicategory on a dual pair} and $\mathbb{F}_{cd} := \mathbb{F}(G_{d}, \mathcal{R}_{cd})$ the \textbf{free monoidal bicategory on a coherent dual pair}.
\end{defin}
One can now proceed to define bicategories of dual pairs in any monoidal bicategory $\mathbb{M}$ in the standard way, by looking at the bicategory of homomorphisms out of the corresponding freely generated monoidal bicategory. In this case, it is advantageous to consider only strict homomorphisms, as then the relevant bicategories admit very compact descriptions using the language of $(G, \mathcal{R})$-shapes, where $(G, \mathcal{R})$ is a generating datum.

\begin{defin}
 If $\mathbb{M}$ is a monoidal category, we define the \textbf{bicategory of dual pairs in $\mathbb{M}$} and the \textbf{bicategory of coherent dual pairs in $\mathbb{M}$} as the bicategories of shapes
\begin{gather*}
\mathcal{D}ualPair(\mathbb{M}) := \mathbb{M}(G_{d}, \mathcal{R}_{d}), \\
\mathcal{C}ohDualPair(\mathbb{M}) := \mathbb{M}(G_{d}, \mathcal{R}_{cd}).
\end{gather*}

Equivalently, the bicategory of dual pairs is precisely the bicategory of strict monoidal homomorphisms $\mathbb{F}_{d} \rightarrow \mathbb{M}$ and similarly the bicategory of coherent dual pairs is the bicategory of strict homomorphisms $\mathbb{F}_{cd} \rightarrow \mathbb{M}$. 
\end{defin}
The equivalence of the two definitions given above follows from \textbf{Appendix \ref{appendix_free_monoidal_bicategories}}, as the bicategories of $(G, \mathcal{R})$-shapes are precisely constructed to model the corresponding bicategories of strict homomorphisms. One can also give the following explicit description of the relevant bicategories, which we present here for conveniance of the reader and also to fix the notation.

Let $\mathbb{M}$ be the ambient monoidal bicategory, assumed to be a \textbf{Gray}-monoid for the purposes of drawing diagrams, although other than that what we say holds in the general case. The objects of the bicategory $\mathcal{D}ualPair(\mathbb{M})$ are precisely the dual pairs in $\mathbb{M}$. 

If $\langle L, R \rangle _{d}, \langle L^\prime, R^\prime \rangle _{d} \in \mathcal{D}ualPair(\mathbb{M})$ are objects, where $\langle L, R \rangle _{d} = (L, R, e, c, \alpha, \beta)$, $\langle L^\prime, R^\prime \rangle _{d} = (L^\prime, R^\prime, e^\prime, c^\prime, \alpha^\prime, \beta^\prime)$, then a morphism $(s, t) _{d}: \langle L, R \rangle _{d} \rightarrow \langle L^\prime, R^\prime \rangle _{d}$ of dual pairs consists of data of $1$-cells $s: L \rightarrow L^\prime$, $t: R \rightarrow R^\prime$ and constraint isomorphisms $\gamma, \delta$ of type

\begin{center}
,
\end{center}
whose purpose is to witness naturality of $s, t$ with respect to (co)evaluation maps. Additionally, these constraint isomorphisms are required to satisfy further naturality condition with respect to $2$-cells. In our case, we require compatibility with cusp isomorphisms $\alpha$ and $\beta$, pictured below.

\begin{center}
	%
.
\end{center}
We will now describe $2$-cells of $\mathcal{D}ualPair(\mathbb{M})$. If $(s_{1}, t_{1}) _{d} = (s_1, t_1, \gamma_1, \delta_1)$ and $(s_{2}, t_{2}) _{d} = (s_2, t_2, \gamma_2, \delta_2)$ are morphisms of type $\langle L, R \rangle _{d} \rightarrow \langle L^\prime, R^\prime \rangle _{d}$, then a 2-cell $\Gamma: (s_1, t_1) _{d} \rightarrow (s_{2}, t_{2}) _{d}$ consists of data of $2$-cells 

\begin{center}
	%
\end{center}
in $\mathbb{M}$. Additionally, we require these $2$-cells to satisfy naturality properties with respect to isomorphisms witnessing naturality of maps $(s_{1}, t_{1})_{d}, (s_{2}, t_{2}) _{d}$. More precisely, we assume that equations expressing compatibility with $\delta$ and $\gamma$, pictured below, hold.

\begin{center}
	%
.
\end{center}
The bicategory $\mathcal{C}ohDualPair(\mathbb{M})$ is the full subbicategory of $\mathcal{D}ualPair(\mathbb{M})$ spanned by those dual pairs that are coherent. 

The following technical result is an immediate consequence of the cofibrancy theorem for freely generated monoidal bicategories, which we state in the appendix as \textbf{Theorem \ref{thm:cofibrancy_theorem}}. The latter can be understood as saying that the difference between strict and non-strict homomorphisms out of a freely generated monoidal bicategory is neglible, as it implies that the relevant bicategories of homomorphisms are equivalent.
\begin{prop}
The strict injective homomorphisms
\begin{gather*}
\mathcal{D}ualPair(\mathbb{M}) \hookrightarrow \textbf{MonBicat}(\mathbb{F}_{d}, \mathbb{M}) \\
\mathcal{C}ohDualPair(\mathbb{M}) \hookrightarrow \textbf{MonBicat}(\mathbb{F}_{cd}, \mathbb{M})
\end{gather*}
identifying the bicategories on the left as the full subbicategories on the strict homomorphisms and monoidally-strict natural transformations are equivalences of bicategories.
\end{prop}
The invariance of equivalence type of bicategories of homomorphisms then immediately yields the following corollary.
\begin{cor}
The equivalence type of bicategories of dual and coherent dual pairs in $\mathbb{M}$ depends only on the equivalence type of $\mathbb{M}$. 
\end{cor}
We now proceed to establish a basic property of bicategories of dual pairs, namely that they are groupoids. This is analogous to the similar result for dual pairs in ordinary monoidal categories, which we have discussed in the introduction. 

\begin{prop}
\label{prop:bicategory_of_dual_pairs_is_a_groupoid}
The bicategory $\mathcal{D}ualPair(\mathbb{M})$ is a 2-groupoid.
\end{prop}

\begin{proof}
Since the equivalence type of $\mathcal{D}ualPair(\mathbb{M})$ depends only on the equivalence type of $\mathbb{M}$, for the purposes of the proof we may assume that $\mathbb{M}$ is a \textbf{Gray}-monoid. 

We first show that all morphisms are equivalences. Observe that a 1-cell $(s, t) _{d} = (s, t, \gamma, \delta)$ is an equivalence if and only if $s, t$ are, this follows from the corresponding statement in the homomorphism bicategory. Any map of dual pairs in a monoidal bicategory is in particular a map of dual pairs in $\text{Ho}(\mathbb{M})$, so by coherence for dualizable objects in monoidal categories, which we presented in the previous chapter as \textbf{Theorem \ref{thm:coherence_for_dualizable_objects_in_monoidal_categories}}, both its components must be isomorphisms in the homotopy category. It follows that $s, t$ are equivalences.

We are now left with showing that all 2-cells are isomorphisms, by the part above it is enough to do so for all $2$-cells between autoequivalences. We will do so by exhibiting for any endomorphism $(s, t)_{d}: \langle L, R \rangle _{d} \rightarrow \langle L, R \rangle _{d}$ a "canonical" pseudoinverse to the component $s: L \rightarrow L$, together with "canonical" witnessing isomorphisms between the respective composites and the identity. By "canonical" we mean that this structure will be natural with respect to maps of such endomorphisms.

We construct it as follows. The pseudoinverse of $s$ is given by $s^{-1} := (e \otimes L) \circ (L \otimes t \otimes L) \circ (L \otimes c)$, that is, the the dual of $t$. It is easy to observe that the witnessing isomorphisms $\theta: s s^{-1} \simeq id_{L}$ and $\vartheta: s^{-1} s \simeq id_{L}$ pictured below

\begin{center}

\end{center}
are natural with respect to maps between such endomorphisms. This follows from the naturality properties of the interchanger and the naturality of $\gamma, \delta$ constraint 2-cells with respect to such maps. 

Once we know that the component $s: L \rightarrow L$ of an endomorphism of a dual pair is canonically exhibited as a witnessed equivalence and that 2-cells between such endomorphisms commute with these witnessing isomorphisms $\theta, \vartheta$, the fact that such 2-cells must be invertible follows by a variation on the argument from \textbf{Theorem \ref{thm:coherence_for_dualizable_objects_in_monoidal_categories}}.

In detail, suppose we are now given two endomorphisms $(s, t) _{d} = (s, t, \gamma, \delta)$ and $(s^\prime, t^\prime) _{d} = (s^\prime, t^\prime, \gamma^\prime, \delta^\prime)$ of $\langle L, R \rangle _{d}$. Consider the composite 

\begin{center}
	\begin{tikzpicture}[thick, scale=0.75, every node/.style={scale=0.75}]
		\node (L0) at (-3, 0) {$ L $};
		\node (L1) at (-1.5,0) {$ L $};
		\node (LRL1) at (1.5, 0) {$ L \otimes R \otimes L $};
		\node (LRL2) at (6.5, 0) {$ L \otimes R \otimes L $};
		\node (L2) at (9.5, 0) {$ L $};
		\node (L3) at (11, 0) {$ L $};
		
		\draw [->] (L1) to node[auto] {$ L \otimes c $} (LRL1);
		\draw [->] (LRL2) to node[auto] {$ e \otimes L $} (L2);
		\draw [->] (LRL1) to [out=25, in=155] node[above] {$ L \otimes t_{1} \otimes L $} (LRL2);
		\draw [->] (LRL1) to [out=-25, in=-155] node[below] {$ L \otimes t_{2} \otimes L $} (LRL2);
		\draw [->] (L0) to node[auto] {$ s_{2} $} (L1);
		\draw [->] (L2) to node[auto] {$ s_{1} $} (L3);
		
		\draw [->] (L1) to [out=30, in=150] node[above] {$ L $} (L3);
		\draw [->] (L0) to [out=-30, in=-150] node[below] {$ L $} (L2);
		
		\node at (4, 0) {$ \Downarrow L \otimes \Gamma_{R} \otimes L $};
		\node at (1.3, -1) {$ \simeq \vartheta _{2} $};
		\node at (6.7, 1) {$ \simeq \theta _{1} $};
	\end{tikzpicture}.
\end{center}
If we paste $\Gamma _{L}$ from below, which amounts to postcomposition, by the naturality of $\theta$ the resulting pasting diagram will be an isomorphism involving $\theta_{2},\vartheta _{2}$. Similarly, if we paste $\Gamma _{L}$ from above, which amounts to precomposition, the naturality of $\vartheta$ implies that the resulting diagram will be an isomorphism involving $\theta_{1}, \vartheta_{1}$. It follows that $\Gamma_{L}$ is an isomorphism itself. 

The argument for $\Gamma_{R}$ is analogous, where we would now exhibit the $t$-component of a map of dual pairs as a witnessed equivalence with pseudoinverse given by the dual of $s$. 
\end{proof}
As promised, we will now use the language of bicategories of dual pairs to prove the general case of the strictification theorem. We first recall it for conveniance of the reader. 

\begin{thm*}[Strictification for dual pairs]
If a dual pair $\langle L, R \rangle _{d}$ satisfies one of the Swallowtail equations, then it satisfies both of them. 

Moreover, if we keep the objects and (co)evaluation maps fixed, then for any choice of a cusp isomorphism $\alpha$ there is a unique cusp isomorphism $\beta$ such that together they satisfy both Swallowtail equations. In particular any dual pair can be made coherent by only a change of $\beta$. . 
\end{thm*}
Let us observe that the proof given below makes an essential use of the cofibrancy theorem to reduce everything to the case of a strict homomorphisms. This is most likely not necessary and a more direct argument involving an arbitrary equivalence of $\mathbb{M}$ with a \textbf{Gray}-monoid is perhaps possible, however, it would probably be slightly more difficult. 

\begin{proof}
We've already proved the result under the additional assumption that the ambient monoidal bicategory is a $\textbf{Gray}$-monoid when we first stated it as \textbf{Theorem \ref{thm:strictification_for_dual_pairs}}, so now we only have to show that the general case follows from this one.

Let $\mathbb{M}$ be a monoidal bicategory and let $\widetilde{\mathbb{M}}$ be its cofibrant replacement. That is, let $\widetilde{\mathbb{M}}$ be the monoidal bicategory freely generated by all the objects, 1-cells and 2-cells on $\mathbb{M}$, subject to the relation that two 2-cells are equal if and only if their images under the obvious strict homomorphism $\widetilde{\mathbb{M}} \rightarrow \mathbb{M}$ are equal. The strict homomorphism in question is an equivalence, as it is surjective on objects, morphisms and $2$-cells by construction, which gives essential surjecitivty on objects, morphisms and local fullness. It is also locally faithful, as we have identified any $2$-cells whose images in $\mathbb{M}$ are equal.

By coherence for tricategories, we can choose an equivalence $\mathbb{M} \rightarrow \mathbb{G}$, where $\mathbb{G}$ is a \textbf{Gray}-monoid. Consider the composite $\widetilde{\mathbb{M}} \rightarrow \mathbb{G}$, by cofibrancy theorem for freely generated monoidal bicategories, which is given in the appendix as \textbf{Theorem  \ref{thm:cofibrancy_theorem}}, we can replace it up to natural equivalence by a strict homomorphism. This new homomorphism will be again an equivalence of monoidal bicategories. We will first show that $\widetilde{\mathbb{M}}$ satisfies the conditions of the strictification theorem and later conclude that $\mathbb{M}$ also does. 

We have a commutative diagram of monoidal bicategories and strict homomorphisms, where all the vertical maps are induced by the chosen strict equivalence $\widetilde{\mathbb{M}} \rightarrow \mathbb{G}$. It follows that they are equivalences, too. The horizontal homomorphisms from $\mathcal{C}ohDualPair(-)$ to $\mathcal{D}ualPair(-)$ are inclusions of full subbicategories.

\begin{center}
	\begin{tikzpicture}
	\node at (3, 0) (DualMt) {$ \mathcal{D}ualPair(\widetilde{\mathbb{M}}) $};
	\node at (3, -1.5) (DualG) {$ \mathcal{D}ualPair(\mathbb{G}) $};
	\node at (0, 0) (CohMt) {$ \mathcal{C}ohDualPair(\widetilde{\mathbb{M}}) $};
	\node at (0, -1.5) (CohG) {$ \mathcal{C}ohDualPair(\mathbb{G}) $};
	\node at (-2.4, 0) (Mtd) {$ \widetilde{\mathbb{M}}^{d} $};
	\node at (-2.4, -1.5) (Gd) {$ \mathbb{G}^{d} $};
	
	\draw [->] (CohMt) to (DualMt);
	\draw [->] (CohG) to (DualG);
	\draw [->] (CohMt) to node[auto] {$ \pi_{\widetilde{\mathbb{M}}} $} (Mtd);
	\draw [->] (CohG) to node[auto] {$ \pi_{\mathbb{G}} $} (Gd);
	\draw [->] (DualMt) to (DualG);
	\draw [->] (CohMt) to (CohG);
	\draw [->] (Mtd) to (Gd);
	\end{tikzpicture}
\end{center}

Let $\langle L, R \rangle _{d} \in \mathcal{D}ualPair(\widetilde{\mathbb{M}})$ and suppose that it satisfies one of the Swallowtail identities. Then the same is true for its image in $\mathcal{D}ualPair(\mathbb{G})$ and by strictification theorem for \textbf{Gray}-monoids, which we've already proven, the image is coherent. As $\widetilde{\mathbb{M}} \rightarrow \mathbb{G}$ is an equivalence, it is in particular locally injective on 2-cells and we can use reflection of equations, which is given in the appendix as \textbf{Lemma \ref{lem:reflection_of_equations}}, to deduce that any equations satisfied by the image were already satisfied by $\langle L, R \rangle _{d}$. Hence, $\langle L, R \rangle _{d}$ is also coherent. This shows that any of the Swallowtail identities implies the other for dual pairs in $\widetilde{\mathbb{M}}$. 

Moreover, if $\langle L, R \rangle _{d}$ is a dual pair not necessarily satisfying the Swallowtail equation, then its image in $\mathcal{D}ualPair(\mathbb{G})$ has the property that we can make it coherent by a unique change of cusp isomorphism $\beta$. As $\widetilde{\mathbb{M}} \rightarrow \mathbb{G}$ is locally bijective on 2-cells, this implies there is a unique choice of this cusp isomorphism in $\widetilde{\mathbb{M}}$ such that the modified dual pair will have coherent image in $\mathcal{D}ualPair(\mathbb{G})$. By reflection of equations, this implies that the modified dual pair in $\widetilde{\mathbb{M}}$ is coherent itself. This shows the second part of the strictification result for $\widetilde{\mathbb{M}}$.

Extension of this result to $\mathbb{M}$ is similar, except now we go "the other way around". Consider the analogus diagram induced by the strict equivalence $\widetilde{\mathbb{M}} \rightarrow \mathbb{M}$. 

\begin{center}
	\begin{tikzpicture}
	\node at (3, 0) (DualMt) {$ \mathcal{D}ualPair(\widetilde{\mathbb{M}}) $};
	\node at (3, -1.5) (DualM) {$ \mathcal{D}ualPair(\mathbb{M}) $};
	\node at (0, 0) (CohMt) {$ \mathcal{C}ohDualPair(\widetilde{\mathbb{M}}) $};
	\node at (0, -1.5) (CohM) {$ \mathcal{C}ohDualPair(\mathbb{M}) $};
	\node at (-2.4, 0) (Mtd) {$ \widetilde{\mathbb{M}}^{d} $};
	\node at (-2.4, -1.5) (Md) {$ \mathbb{M}^{d} $};
	
	\draw [->] (CohMt) to (DualMt);
	\draw [->] (CohM) to (DualM);
	\draw [->] (DualMt) to (DualM);
	\draw [->] (CohMt) to node[auto] {$ \pi_{\widetilde{\mathbb{M}}} $} (Mtd);
	\draw [->] (CohM) to node[auto] {$ \pi_{\mathbb{M}} $} (Md);
	\draw [->] (CohMt) to (CohM);
	\draw [->] (Mtd) to (Md);
	\end{tikzpicture}.
\end{center}
The outer vertical arrows are both surjective on objects as any dual pair in $\mathbb{M}$ can be lifted to $\widetilde{\mathbb{M}}$ in a canonical way, since the cofibrant replacement is precisely generated by objects, 1-cells and 2-cells of $\mathbb{M}$. 

Let $\langle L, R \rangle _{d} \in \mathcal{D}ualPair(\mathbb{M})$ and denote its canonical lift by $\langle \widetilde{L}, \widetilde{R} \rangle _{d}$. If the pair $\langle L, R \rangle _{d}$ satisfies one of the Swallowtail identities, by reflection of equations so does the lift.  Since we have proven the stricitification result for $\widetilde{\mathbb{M}}$, it follows that $\langle \widetilde{L}, \widetilde{R} \rangle _{d}$ is coherent, hence so is its image $\langle L, R \rangle _{d}$.

Similarly, one observes that for the lifted dual pair $\langle \widetilde{L}, \widetilde{R} \rangle _{d}$ there is a unique choice of a possibly different cusp isomorphism $\beta$ that would make it coherent. As $\widetilde{\mathbb{M}} \rightarrow \mathbb{G}$ is locally bijective on 2-cells, the same must be true for $\langle L, R \rangle _{d}$.
\end{proof}

\subsection{Coherence for dualizable objects in monoidal bicategories}

We are now ready to prove one of the main results of the current work, which is a coherence theorem for dualizable objects. 

\begin{thm}[Coherence for dualizable objects]
\label{thm:coherence_for_dualizable_objects}
Let $\mathbb{M}$ be a monoidal bicategory.The forgetful homomorphism

\begin{center}
$\pi: \mathcal{C}ohDualPair(\mathbb{M}) \rightarrow \mathbb{K}(\mathbb{M}^{d})$, \\
$\langle L, R \rangle _{d} \mapsto L$
\end{center}
between, respectively, the bicategory of coherent dual pairs in $\mathbb{M}$ and the groupoid of dualizable objects in $\mathbb{M}$, is a surjective on objects equivalence.
\end{thm}

\begin{rem}
We know that the image of $\pi$ must land in the maximal subgroupoid of $\mathbb{M}^{d}$, since $\mathcal{C}ohDualPair(\mathbb{M})$ is a groupoid itself by \textbf{Proposition \ref{prop:bicategory_of_dual_pairs_is_a_groupoid}}.
\end{rem}
The result should be understood as saying that the notion of a coherent dual pairs is a property-like structure equivalent to the property of being dualizable. This means two things, one of which is that an object is dualizable if and only if it can be completed to a coherent dual pair, which we already know. 

The other part, which is hidden in the adjective \emph{property-like}, is that any two coherent dual pairs living over the same dualizable object are equivalent. In fact, even more is true -  the space of coherent dual pairs over any given object is always contractible, this follows easily from the above theorem by taking suitable fibers. 

We will now proceed with the proof, which will take the remainder of this chapter. The general outline is as follows. We already know that the functor $\pi$ is surjective on objects, this is exactly the strictification theorem for dual pairs, that is, \textbf{Theorem \ref{thm:strictification_for_dual_pairs}}. Hence, we are left with showing that $\pi$ is essentially surjective on morphisms and locally fully faithful, as these are exactly the conditions for a map of bicategories to be an equivalence. 

We will first reduce to the case of $\mathbb{M}$ being a \textbf{Gray}-monoid and later prove essential surjectivity and local full faithfulness directly as \textbf{Lemma \ref{lem:essential_surjectivity_on_equivalences}} and \textbf{Lemma \ref{lem:bijectivity_on_isomorphisms_between_equivalences}}. This will end the proof.

\begin{lem}
\label{lem:on_stability_of_properties_of_forgetful_functor_from_dual_pairs}
Let $P$ be a property of homomorphisms of bicategories that is stable under composition with equivalences. Then $\pi_{\mathbb{M}}: \mathcal{C}ohDualPair(\mathbb{M}) \rightarrow \mathbb{K}(\mathbb{M}^{d})$ satisfies $P$ if and only if $\pi_{\mathbb{G}}: \mathcal{C}ohDualPair(\mathbb{G}) \rightarrow \mathbb{K}(\mathbb{G}^{d})$ does, where $\mathbb{G}$ is a \textbf{Gray}-monoid equivalent to $\mathbb{M}$. 
\end{lem}

\begin{rem}
The same is true if we replace $\mathcal{C}ohDualPair(-)$ by $\mathcal{D}ualPair(-)$, or more generally, the same is true for forgetful homomorphism out of general bicategories of $(G, \mathcal{R})$-shapes, where $(G, \mathcal{R})$ is a generating datum for a monoidal bicategory. The proof given below works verbatim.
\end{rem}

\begin{proof}
Let $\mathbb{M}$ be a monoidal category and let $\widetilde{\mathbb{M}}$ be its cofibrant replacement. Choose an equivalence $\mathbb{M} \rightarrow \mathbb{G}$ and deform the composite equivalence $\widetilde{\mathbb{M}} \rightarrow \mathbb{G}$ to a strict one. In the commutative diagram 

\begin{center}
	\begin{tikzpicture}
	\node at (7, 0) (CohG) {$ \mathcal{C}ohDualpair(\mathbb{G}) $};
	\node at (7, -1.5) (Gd) {$ \mathbb{K} (\mathbb{G}^{d} ) $};
	\node at (3.5, 0) (CohMt) {$ \mathcal{C}ohDualPair(\widetilde{\mathbb{M}}) $};
	\node at (0, 0) (CohM) {$ \mathcal{C}ohDualPair(\mathbb{M}) $};
	\node at (3.5, -1.5) (Mtd) {$ \mathbb{K}( \widetilde{\mathbb{M}}^{d} )$};
	\node at (0, -1.5) (Md) {$ \mathbb{K}(\mathbb{M}^{d}) $};
	
	\draw [->] (CohMt) to node[auto] {$ \pi_{\widetilde{\mathbb{M}}} $} (Mtd);
	\draw [->] (CohM) to node[auto] {$ \pi_{\mathbb{M}} $} (Md);
	\draw [->] (CohG) to node[auto] {$ \pi_{\mathbb{G}} $} (Gd);
	\draw [->] (CohMt) to (CohG);
	\draw [->] (Mtd) to (Gd);
	\draw [->] (CohMt) to (CohM);
	\draw [->] (Mtd) to (Md);
	\end{tikzpicture}
\end{center}
of strict homomorphisms the horizontal maps are all equivalences, from which we immediately conclude that $\pi_{\mathbb{G}}$ satisfies $P$ if and only if $\pi_{\mathbb{M}}$ does.
\end{proof}

\begin{lem}
\label{lem:essential_surjectivity_on_equivalences}
The homomorphism $\pi: \mathcal{C}ohDualPair(\mathbb{M}) \rightarrow \mathbb{K}(\mathbb{M}^{d})$ is essentially surjective on morphisms. In other words, any equivalence in $\mathbb{M}^{d}$ is isomorphic to one that is an image under $\pi$ of a morphism in $\mathcal{C}ohDualPair(\mathbb{M})$. 
\end{lem}

\begin{proof}
By \textbf{Lemma \ref{lem:on_stability_of_properties_of_forgetful_functor_from_dual_pairs}} we may assume that $\mathbb{M}$ is a \textbf{Gray}-monoid. We will prove that in this case \emph{every} equivalence in $\mathbb{M}$ may be lifted to one in $\mathcal{C}ohDualPair(\mathbb{M})$. 

Let $\langle L, R \rangle _{d} = (L, R, e, c, \alpha, \beta)$, $\langle L^\prime, R^\prime \rangle _{d} = (L^\prime, R^\prime, e^\prime, c^\prime, \alpha^\prime, \beta^\prime)$ be coherent dual pairs and suppose we are given an equivalence $s: L \rightarrow L^\prime$. We have to complete it to a 1-cell $(s, t) _{d}: \langle L, R \rangle _{d} \rightarrow \langle L^\prime, R^\prime \rangle _{d}$ in $\mathcal{C}ohDualPair(\mathbb{M})$ with components $(s, t) _{d} = (s, t, \gamma, \delta)$ 

We first pass to the homotopy category of $\mathbb{M}$. Observe that by coherence for dual pairs in monoidal categories, which we stated in the previous chapter as \textbf{Theorem \ref{thm:coherence_for_dualizable_objects_in_monoidal_categories}}, there is a unique homotopy class $t$ that will complete $s$ to a map of dual pairs in $\text{Ho}(\mathbb{M})$.  In fact, we see from the proof that it is given by the dual of any pseudoinverse of $s$. We take any representative of this map to be the component $t$. We  are now left with choosing the needed constraint isomorphisms.

The fact that $s, t$ commute with (co)evaluation maps in the homotopy category implies that we can choose the constraint isomorphisms $\gamma, \delta$ \emph{in some way}. However, we need to be careful and choose them in a way that makes them natural with respect to the cusp isomorphisms $\alpha, \beta$. 

The first step in showing that this can be done is to prove that for any choice of $\gamma, \delta$, naturality with respect to $\alpha$ implies naturality with respect to $\beta$. This can be done in several ways, one of which is a straightforward manipulation with diagrams, which is rather lengthy. As throughout the paper we will need several statements of this form, it seems worthwhile to discuss a "cute trick" which applies here.

Suppose we have chosen some $\gamma, \delta$ that are natural with respect to the cusp isomorphism $\alpha$. The data of $(s, t, \gamma, \delta)$ is then already an equivalence $(L, R, e, c) \rightarrow (L^\prime, R^\prime, e^\prime, c^\prime)$ in the bicategory $\mathbb{M}(G_{d} ^\prime)$ of 1-truncated $G_{d}$-shapes.

Using  \textbf{Lemma \ref{lem:promoting_equivalences_of_g_shapes}}, we can transport the structure of a $G_{d}$-shape from $\langle L, R \rangle _{d}$ to $(L^\prime, R^\prime, e^\prime, c^\prime)$. The new $G_{d}$-shape will be a coherent dual pair by \textbf{Lemma \ref{lem:promotion_respects_classes_of_relations}}, as it is equivalent to one, let us denote it by $\widetilde{ \langle L^\prime, R^\prime \rangle _{d}}$. By the precise statement of the lemma, the map $(s, t, \gamma, \delta)$ we began with was in fact an equivalence of dual pairs if and only if $\widetilde{\langle L^\prime, R^\prime \rangle _{d}} = \langle L, R \rangle _{d}$. The latter amounts to agreeing on cusp isomorphisms, as the two dual pairs agree on objects and (co)evaluation maps by construction.

As $(s, t, \gamma, \delta)$ was already natural with respect to cusp isomorphism $\alpha$, the dual pairs $\widetilde{ \langle L^\prime, R^\prime \rangle _{d}}$, $\langle L, R \rangle _{d}$ must agree on it. However, they are both coherent and strictification for dual pairs, that is, \textbf{Theorem \ref{thm:strictification_for_dual_pairs}}, implies that they are equal, as $\alpha$ uniquely determines $\beta$. This proves that the data $(s, t, \gamma, \delta)$ was natural with respect to $\beta$ to begin with and ends our "cute trick".

Observe that the argument given above is completely formal, but at least notationally it simplifies things greatly, as it reduces a question about determination of data in a \emph{morphism} of dual pairs to a question about determination of data in a single dual pair. 

We are left with showing that one can choose constraint isomorphisms $\gamma, \delta$ such that they satisfy the equation of $\alpha$-naturality, which can be rewritten as

\begin{center}
.
\end{center}
By adding the inverse of $\gamma \otimes s$ to the bottom right corner we now see the equation becomes

\begin{center}
	%
,
\end{center}
where we do not redraw the right hand side, as for the purpose of what we're showing it is not needed to know its exact value, as long as it's an isomorphism completely definable in terms of only the two dual pairs in question and $\gamma$. We now apply interchangers to the equation above, to the maps in the bottom right to move the left $s$ from the cell $s \otimes t \otimes s = (L^\prime \otimes t \otimes s) (s \otimes R \otimes L)$ past the $L \otimes c$ before it, so that we obtain 

\begin{center}
	%
.
\end{center}
As $s$ is an equivalence, this equation is equivalent to one of the form 

\begin{center}
	%
,
\end{center}
and \textbf{Lemma \ref{lem:hom_equivalences}} implies that it has a unique solution, as $(e^\prime \otimes L^\prime) \circ (L^\prime \otimes -)$ is an equivalence on Hom-categories.
\end{proof}

We will now finish the proof of the coherence theorem by showing that $\pi$ is locally bijective on 2-cells. Here, it is not necessary to work with \emph{coherent} dual pairs, hence we state the result in a slightly greater generality. 

\begin{lem}
\label{lem:bijectivity_on_isomorphisms_between_equivalences}
The homomorphism $\pi: \mathcal{D}ualPair(\mathbb{M}) \rightarrow \mathbb{K}(\mathbb{M}^{d})$ is locally bijective on $2$-cells.
\end{lem}

\begin{proof}
By \textbf{Lemma \ref{lem:on_stability_of_properties_of_forgetful_functor_from_dual_pairs}}, or more precisely the remark below it, we may assume $\mathbb{M}$ is a \textbf{Gray}-monoid, as the property of being locally bijective on $2$-cells is stable under composition with equivalences.

Suppose we have maps $(s_{1}, t_{1}) _{d}, (s_{2}, t_{2}) _{d}$ in $\mathcal{D}ualPair(\mathbb{M})$  with $(s_{i}, t_{i}) _{d}: \langle L, R \rangle _{d} \rightarrow \langle L^\prime, R^\prime \rangle _{d}$  and components $(s_{i}, t_{i}) _{d} = (s_{i}, t_{i}, \gamma _{i}, \delta _{i})$. Let $\Gamma_{L}: \pi((s_{1}, t_{1}) _{d}) \rightarrow \pi((s_{2}, t_{2})_{d})$ be an isomorphism. We have to show that this lifts uniquely to a 2-cell $\Gamma$ in $\mathcal{D}ualPair(\mathbb{M})$. 

For the purposes of the proof, we may well assume that $\langle L, R \rangle _{d} = \langle L^\prime, R^\prime \rangle _{d} = (L, R, e, c, \alpha, \beta)$. Indeed, the bicategory of dual pairs is a groupoid, so the relevant homomorphism will be locally bijective on $2$-cells if and only if it is bijective on $2$-cells in endomorphism categories.

Since $\pi((s_{i}, t_{i}) _{d}) = s_{i}$, the lift amounts to completing the isomorphism $\Gamma_{L}: s_{1} \rightarrow s_{2}$ by finding an invertible $\Gamma_{R}: t_{1} \rightarrow t_{2}$ such that together they satisfy $\gamma$ and $\delta$-naturality equations. As has by now become our standard approach, we will first show that one of these equations implies the other and later that this equation alone makes $\Gamma_{L}$ uniquely determine $\Gamma_{R}$. This will end the proof of the lemma.

In the first part another "cute trick" applies. Suppose we have chosen some isomorphisms $\Gamma_{L}, \Gamma_{R}$ that are natural with respect to the constraint isomorphism $\delta$. This data alone defines an isomorphism $\Gamma: (s_{1}, t_{1}) \rightarrow (s_{2}, t_{2})$ of maps of $G_{d}^{\prime \prime}$-shapes, where $G_{d}^{\prime \prime}$ is the $0$-truncation of $G_{d}$. 

By \textbf{Lemma \ref{lem:promotion_of_isomorphisms_between_maps_of_gpp_shapes}}, we can transport the structure of a map of $G_{d}^{\prime}$-shapes from $(s_{1}, t_{1}) _{d}$ to $(s_{2}, t_{2})$ along the isomorphism $\Gamma$. This way, we obtain some equivalence $\widetilde{(s_{2}, t_{2})}$ which agrees with $(s_{2}, t_{2}) _{d}$ on its components on objects, moreover, being isomorphic to a map of dual pairs, $\widetilde{(s_{2}, t_{2})}$ is a map of dual pairs too.

Observe now that $(s_{2}, t_{2}) _{d}$ and $\widetilde{(s_{2}, t_{2})}$ are both equivalences of dual pairs and moreover they agree on the constraint isomorphism $\delta$, as $\Gamma$ was natural with respect to it by assumption. However, from the proof of \textbf{Lemma \ref{lem:essential_surjectivity_on_equivalences}} above we know that for any given equivalence of dual pairs, the constraint isomorphism $\delta$ determines the constraint isomorphism $\gamma$.

This implies that we have $(s_{2}, t_{2}) _{d} = \widetilde{(s_{2}, t_{2})}$, so that our isomorphisms $\Gamma_{R}, \Gamma_{L}$ were natural with respect to both $\gamma, \delta$ to begin with. This ends our "cute trick" and shows that naturality with respect to $\delta$ implies naturality with respect to $\gamma$.

We will now prove that $\delta$-naturality makes $\Gamma_{L}$ uniquely determine $\Gamma _{R}$, this will end the proof of the lemma. We can redraw the relevant equation as

\begin{center}

\end{center}
or perhaps more simply as 

\begin{center}
	%
\end{center}
since, again, we do not need to know precisely what is the right hand side. One should be perhaps a little bit cautious here because to obtain an equation of this form from the one just above we needed to apply interchangers to the top and the bottom, as according to our convention $(t \otimes s) = (t \otimes L) (r \otimes R)$ and here we want it the other way around. We can now append $R \otimes \Gamma_{L}^{-1}$ to the left hand side to get 

\begin{center}
	%
.
\end{center}
This equation clearly has a unique solution $\Gamma_{R}$. Indeed, postcomposition with $R \otimes s_{2}$ is an equivalence on Hom-categories, since the map itself is an equivalence, as is tensoring via $- \otimes L$ followed by precomposition with $c$ by our standard lemma.
\end{proof}

\section{Fully dualizable objects in symmetric monoidal bicategories}

In this chapter we focus on the theory of fully dualizable objects in symmetric monoidal bicategories. Our main goal is a coherence result analogous to the one we proved for dualizable objects, namely description of a property-like structure equivalent to full dualizability. 

Note that full dualizability can be in general defined as a property of objects in non-symmetric monoidal bicategories, strictly stronger than "ordinary" dualizability, by introducing so called monoidal bicategories \emph{with duals}. We will not work in this generality and restrict to the symmetric case, as this is, as predicted by the Cobordism Hypothesis, where the theory is the most geometric. 

\subsection{Serre autoequivalence}

In this section we define full dualizability and show that every fully dualizable object in a symmetric monoidal bicategory admits a canonical up to isomorphism autoequivalence, which one calls the Serre autoequivalence. We use it to identify a set of minimal conditions implying full dualizability. Nothing here is really new, with the terminology and results appearing in similar form in \cite{lurietqfts}. 

The idea of full dualizability is to strengten the notion of dualizability by putting some conditions on the (co)evaluation maps. In an ordinary monoidal category, we could require them to be isomorphisms and in this way we obtain the notion of a "monoidally inverse pair".

However, in the bicategorical case there is a property of morphisms which is weaker then being an equivalence, but has some of the same consequences - namely the property of having an adjoint. In the strongest version we might require the (co)evaluation maps to \emph{have all adjoints}, that is, to have both left and right adjoints and also that these left/right adjoints should have both adjoints on their own, and these adjoints should have both adjoints and so on. This is what we call \emph{full dualizability}. 

To make this precise, it is conveniant to define a maximal subbicategory admitting all adjoints, which we do now.

\begin{defin}
If $\mathbb{B}$ is a bicategory, let $\mathbb{B}^{adj}$ be the maximal subbicategory with the property that all its $1$-cells admit both adjoints.
\end{defin}

The bicategory $\mathbb{B}^{adj}$ can be obtained by an iterated process of removing morphisms that do not have a left or a right adjoint. Namely, if $\mathbb{B}$ is a bicategory, let $\mathbb{B}^{(1)} \subseteq \mathbb{B}$ be the full subbicategory on the same objects and those morphisms that admit both a left and right adjoint in $\mathbb{B}$. Then, define inductively $\mathbb{B}^{(i+1)} = (\mathbb{B}^{(i)}) ^{(1)}$. One sees easily that $\mathbb{B}^{adj} = \bigcap _{i \geq 1} \mathbb{B}^{(i)}$. 

\begin{prop}
If $\mathbb{M}$ is monoidal, then the monoidal structure restricts to give one on $\mathbb{M}^{adj}$. The same is true for braidings and syllepses, and hence symmetries.
\end{prop}

\begin{proof}
It's enough to show that the tensor product preserves the property of having adjoints, as then one can naively restrict the monoidal structure to $\mathbb{B}^{adj}$. However, if we have $f_{1} \dashv g_{1}$, $f_{2} \dashv g_{2}$, then $f_{1} \otimes f_{2} \dashv g_{1} \otimes g_{2}$.
\end{proof}

\begin{defin}
An object $L \in \mathbb{M}$ in a symmetric monoidal bicategory is \textbf{fully dualizable} if it is dualizable as an object of $\mathbb{M} ^{adj}$.
\end{defin}

In other words, an object $L \in \mathbb{M}$ is fully dualizable if and only if it can be completed to a dual pair $\langle L, R \rangle _{d}$ with (co)evaluation maps admitting all adjoints. We start by deriving some basic properties of fully dualizable objects, our short-term goal will be to describe a finite set of conditions that guarantee that an object is fully dualizable. 

\begin{defin}
\label{defin:serre_autoequivalence}
Let $\langle L, R \rangle _{d}$ be a dual pair in a symmetric monoidal bicategory $\mathbb{M}$ and suppose we have adjunctions $e^{L} \dashv e, c^{L} \dashv c$. We define the \textbf{Serre autoequivalence} $q$ of $L$ and the \textbf{pseudoinverse to Serre autoequivalence} $q^{-1}$ to be the $1$-cells 

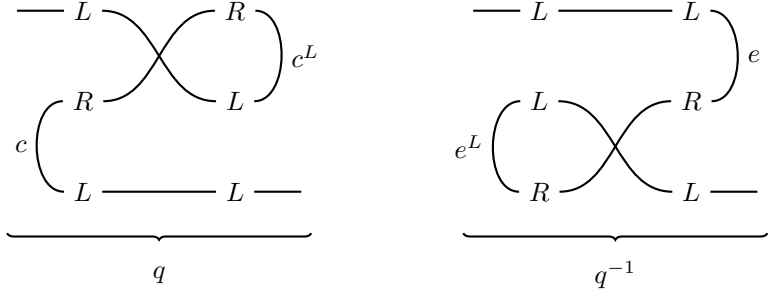
\begin{center}
	\begin{tikzpicture}[thick]
	\node (L) at (0,0) {$ L $};
	\node (R) at  (0, -1.2) {$ R $};
	\node (L2) at (0, -2.4) {$ L $};
	\node (R2) at (2, -0) {$ R $};
	\node (L3) at (2, -1.2) {$ L $};
	\node (L4) at (2, -2.4) {$ L $};
	\node (Domain) at (-1, 0) {$ $};
	\node (Codomain) at (3, -2.4) {$ $};
	
	\draw (L) to [out = 0, in = 180] (L3);
	\draw (R) to [out = 0, in = 180] (R2);
	\draw (L2) to [out = 180, in = 180] node[left] {$ c $} (R);
	\draw (R2) to [out = 0, in = 0] node[right] {$ c^{L} $} (L3);
	\draw (L2) to (L4);
	\draw (Domain) to (L);
	\draw (L4) to (Codomain);
	
	\draw [
    		decoration={
        		brace,mirror,
        		raise=.55cm
    	},
  	decorate
	] (-1,-2.4) -- (3,-2.4);
	\node at (1, -3.5) {$ q $};
	
	\node (rL) at (6, 0) {$ L $};
	\node (rL2) at (6, -1.2) {$ L $};
	\node (rR) at (6, -2.4) {$ R $};
	\node (rL3) at (8, 0) {$ L $};
	\node (rR2) at (8, -1.2) {$ R $};
	\node (rL4) at (8, -2.4) {$ L $};
	\node (rDomain) at (5, 0) {$ $};
	\node (rCodomain) at (9, -2.4) {$ $};
	
	\draw (rL2) [out=180, in = 180] to node[left] {$ e^{L} $} (rR);
	\draw (rL2) [out=0, in=180] to (rL4);
	\draw (rR) [out=0, in=180] to (rR2);
	\draw (rL) to (rL3);
	\draw (rL3) [out=0, in=0] to node[right] {$ e $} (rR2);
	\draw (rDomain) to (rL);
	\draw (rL4) to (rCodomain);
	
	\draw [
    		decoration={
        		brace,mirror,
        		raise=.55cm
    	},
  	decorate
	] (5,-2.4) -- (9,-2.4);
	\node at (7, -3.5) {$ q^{-1} $};
	\end{tikzpicture}
\end{center}
\end{defin}
This terminology is due to Jacob Lurie, see \cite{lurietqfts}. It is motivated by a particular case of a suitably defined $(\infty, 2)$-category of cocomplete dg-categories, where the endofunctor described by Serre autoequivalence has traditionally been called the \emph{Serre functor}. 

The general definition of the Serre autoequivalence given in \cite{lurietqfts} rests on the Cobordism Hypothesis. Namely, if $\mathbb{M}$ is a symmetric monoidal $(\infty, n)$-category, the groupoid of fully dualizable objects in $\mathbb{M}$ is acted upon by $SO(n)$, as it can be identified with the category of $n$-dimensional framed topological field theories with values in $\mathbb{M}$, on which $SO(n)$ acts via change of framing. 

If $n \geq 1$, then the group $\pi_{1}(SO(n))$ is cyclic, generated by a "full twist". Under the action on the groupoid of fully dualizable objects, the generator of $\pi_{1}(SO(n))$ corresponds to a natural endotransformation of the identity of this groupoid. The component of this transformation at any given fully dualizable object is precisely its Serre autoequivalence.

We will proceed by establishing some basic properties of the Serre autoequivalence, using the definition in terms of string diagrams given above. Note that most of these properties would follow trivially if we defined the Serre autoequivalence as the component of some invertible endotransformation of the identity homomorphism, but doing so would require us to assume the Cobordism Hypothesis. As one of the main aims of the current work is to prove the latter, we will not do so.

\begin{prop}
\label{prop:serre_autoequivalence_pseudoinverse_to_its_pseudoinverse}
The Serre autoequivalence and its pseudoinverse are, in fact, pseudoinverse to each other.
\end{prop}

\begin{proof}
We only claim that $qq^{-1} \simeq id_{L}, q^{-1} q \simeq id_{L}$ is \emph{some} way, so it is enough to perform a 1-dimensional computation with string diagrams in the homotopy category $\text{Ho}(\mathbb{M})$, which is an ordinary symmetric monoidal category.  The class of $q^{-1} q$ is given by

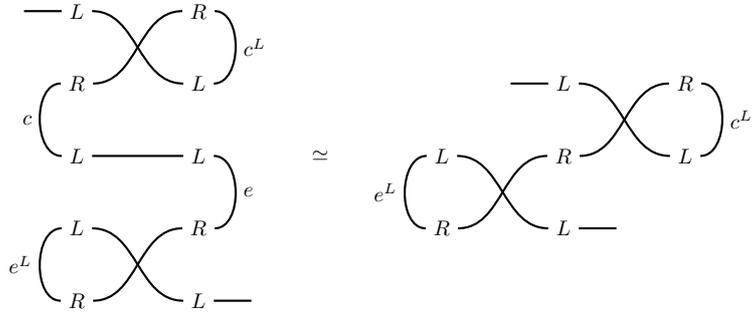
\begin{center}
,
\end{center}
where we obtain the right hand side from the left hand side by using the triangle equations, which are witnessed by cusp isomorphisms. Since in a symmetric monoidal category the composition of braidings depends only on the underlying permutation, we can rewrite the right hand side as

\begin{center}
	%
.
\end{center}
The two sides above are isomorphic by naturality of braidings. The right hand side is the identity, as its representative in $\mathbb{M}$ is left adjoint to the cusp composite, which is itself isomorphic to the identity. As being left adjoint is a relation on isomorphism classes compatible with composition, this ends the proof that $q^{-1} q \simeq id_{L}$. The case of $q q^{-1} \simeq id_{L}$ is completely analogous. 
\end{proof}

\begin{prop}
\label{prop:serre_autoequivalence_is_natural}
The Serre autoequivalence is, up to isomorphisms of 1-cells, natural with respect to equivalences of fully dualizable objects.
\end{prop}

\begin{cor}
Up to isomorphism, the Serre autoequivalence of $L$ does not depend on the choice of its dual, (co)evaluation maps or their left adjoints.
\end{cor}

\begin{proof}
Let $s: L \rightarrow L^\prime$ be an equivalence of fully dualizable objects and suppose that they are part of dual pairs $\langle L, R \rangle _{d}, \langle L^\prime, R^\prime \rangle _{d}$ and that left adjoints $e^{L}, c^{L}, e^{\prime L}, c^{\prime L}$ have been chosen, so that Serre autoequivalences of $L, L^\prime$ are well-defined. 

Choose a map $t: R \rightarrow R^\prime$ that completes it into a map of dual pairs in the homotopy category, this can be done by coherence for dualizable objects in monoidal categories, that is, \textbf{Theorem \ref{thm:coherence_for_dualizable_objects_in_monoidal_categories}}. One computes that, up to isomorphism, $s \circ q$ is 

\begin{center}
.
\end{center}
The right hand side of the above is

\begin{center}
	%
,
\end{center}
as one sees by observing that $c^{\prime L} \circ (t \otimes s)$ is isomorphic to $c^{L}$, as they are both left adjoint to $c$. This ends the proof, as the last composite is precisely $q^\prime \circ s$. 

To deduce the corollary, apply the proposition to the identity morphism of $L$ with different choices of (co)evaluation maps and their left adjoints on both sides. 
\end{proof}

\begin{prop}
\label{prop:left_adjoints_of_coevaluation_maps_in_terms_of_serre}
Let $\langle L, R \rangle _{d}$ be a dual pair and suppose we have adjunctions $c^{L} \dashv c, e^{L} \dashv e$. Then, up to isomorphism, we have

\begin{center}
,
\end{center}
so that the left adjoints to (co)evaluation maps are definable in terms of Serre autoequivalence and its pseudoinverse.
\end{prop}

\begin{proof} We only prove the theorem for $c^{L}$, as the proof for $e^{L}$ is basically the same. In this case, the relevant composite is

\begin{center}
	%
,
\end{center}
where the isomorphism pictured is induced by the cusp. The right hand side is isomorphic to $c^{L}$, as that composite of braidings corresponds to the identity permutation.
\end{proof}

\begin{thm}[Minimal conditions on full dualizability]
\label{thm:minimal_conditions_on_full_dualizability}
Let $(L, R, e, c, \alpha, \beta)$ be a dual pair and suppose that we have adjunctions $c^{L} \dashv c$, $e^{L} \dashv e$, let $\widetilde{e}, \widetilde{c}$ be (co)evaluation maps composed with the symmetry $L \otimes R \simeq R \otimes L$. Then, for all $n \in \mathbb{Z}$ we have adjunctions

\begin{enumerate}
	\item $(q^{-n-1} \otimes R) \circ \widetilde{c} \dashv e \circ (q^{n} \otimes R)$ 
	\item $(R \otimes q^{-n-1}) \circ c \dashv \widetilde{e} \circ (R \otimes q^{n})$ 
	\item $\widetilde{e} \circ (R \otimes q^{n+1}) \dashv (R \otimes q^{-n}) \circ c$ 
	\item $e \circ (q^{n+1} \otimes R) \dashv (q^{-n} \otimes R) \circ \widetilde{c}$. 
\end{enumerate}
In particular, $L$ is already fully dualizable.
\end{thm}

The theorem is very conveniant, as the property of full dualizability seems a priori difficult to verify, as it assures the existence of an infinite strings of adjoints. However, it follows from the above that it is enough to verify that only \emph{two} of the needed adjoints exist.

\begin{proof}
The families $(1)$ and $(2)$ correspond to each other through the symmetry, as do $(3)$ and $(4)$, so it's enough to prove only one from each pair.

For $n = 0$, the adjunctions $(1)$ and $(3)$ are exactly \textbf{Proposition \ref{prop:left_adjoints_of_coevaluation_maps_in_terms_of_serre}}. All the other adjunctions can be obtained by composing with $q \dashv q^{-1}$ or $q^{-1} \dashv q$ and using the fact that the composite of adjoints is adjoint to the composite.
\end{proof}

\subsection{Fully dual pairs}

Having briefly strayed to establish properties of the Serre autoequivalence, in this section we will roughly follow the same path as we did in the case of dualizable objects. 

We first define a structure whose existence is equivalent to full dualizability, that of a \emph{fully dual pair}, and describe a class of them satisfying additional equations, which we call \emph{coherent}. We then organize fully dual pairs into a bicategory and prove that it is a groupoid. Lastly, we prove a strictification result which implies that any fully dualizable object can be completed to a coherent fully dual pair.

Observe that in the case of dualizability, the collection of data needed to witness this property was rather self-evident, and it is this collection that led us to the notion of a dual pair. In the case of full dualizability, things are not so clear. 

One obvious route would be to define a "fully dual pair" to consists of two objects, (co)evaluation maps, cusp isomorphisms, right and left adjoints to (co)evaluation maps together with the relevant (co)units, the adjoints of adjoints together with (co)units, their adjoints and so on. This is a very large collection of data and probably not something we could work with directly.
 
However, as we only work in the symmetric monoidal case, we have \textbf{Theorem \ref{thm:minimal_conditions_on_full_dualizability}} which implies that an object is fully dualizable if and only if it is dualizable and the (co)evaluation maps admit left adjoints. This motivates the following definition.
 
 \begin{defin}
 \label{def:fully_dual_pair}
 A \textbf{fully dual pair} in a symmetric monoidal bicategory $\mathbb{M}$ consists of
 
 \begin{itemize}
 	\item a dual pair $\langle L, R \rangle _{d}$,
	\item morphisms $q, q^{-1}: L \rightarrow L$ together with isomorphisms $\psi: q q^{-1} \simeq id_{L}$, $\phi: q ^{-1} q \simeq id_{L}$,
	\item 2-cells $\mu_{e}: id_{I} \rightarrow e \circ e^{L}$ and $\epsilon_{e}: e^{L} \circ e \rightarrow id_{L \otimes R}$ that are the unit and counit of the adjunction $e^{L} \dashv e$, where $e^{L} = b \circ (R \otimes q^{-1}) \circ c$,
	\item 2-cells $\mu_{c}: id_{R \otimes L} \rightarrow c \circ c^{L}$ and $\epsilon_{c}: c^{L} \circ c \rightarrow id_{I}$ that are the unit and counit of the adjunction $c^{L} \dashv c$, where $c^{L} = e \circ (q \otimes R) \circ b$.
 \end{itemize}
 \end{defin}
 
\begin{notation*}
 By abuse of language we will sometimes refer to the whole fully dual pair $(L, R, e, c, q, q^{-1}, \alpha, \beta, \mu_{e}, \epsilon_{e}, \mu_{c}, \epsilon_{c}, \psi, \phi)$ just by referring to the underlying objects. We will then denote it by $\langle L, R \rangle _{fd}$, do distinguish it from the similar notation for dual pairs.
 \end{notation*}
 The intuition about this definition is as follows. By the theorem on minimal conditions, to witness full dualizability it is enough to give the left adjoints of (co)evaluation maps. We could include these adjoints directly in the definition of a fully dual pair and then in that terms define the Serre autoequivalence and its pseudoinverse, like we did before. 
 
However, for technical reasons - one of which is to make the notion easier to compare with the presentation of the framed bordism bicategory - we take a slightly different route. Instead of adjoints, we postulate the Serre autoequivalence as part of the data of a fully pair and only define the left adjoints in terms of it using the formula given by \textbf{Proposition \ref{prop:left_adjoints_of_coevaluation_maps_in_terms_of_serre}}. Observe that this forces $q, q^{-1}$ that are given as part of a fully dual pair to be, up to isomorphism, Serre autoequivalence and its pseudoinverse of $L$ in the sense of previous definitions.
 
\begin{prop}
\label{prop:completion_to_a_fully_dual_pair}
An object $L \in \mathbb{M}$ in a symmetric monoidal bicategory is fully dualizable if and only if it can be completed to a fully dual pair $\langle L, R \rangle _{fd}$.
\end{prop}
 
\begin{proof}
Clearly any object $L$ that is a part of a fully dual pair is fully dualizable by \textbf{Theorem \ref{thm:minimal_conditions_on_full_dualizability}}, as its (co)evaluation maps have left adjoints. 
 
Conversely, if $L$ is a fully dualizable object, first complete it to a dual pair $\langle L, R \rangle _{d}$. By full dualizability, we can choose some maps $e^{L}, c^{L}$ that are left adjoint to the chosen (co)evaluation maps. In terms of these maps we can define the Serre autoequivalence and its pseudoinverse $q, q^{-1}$ as in \textbf{Definition \ref{defin:serre_autoequivalence}}, we now include them as part of the structure of a fully dual pair.
 
 We are now only left with giving the missing 2-cells. We can find isomorphisms $\psi: q q^{-1} \simeq id_{L}$, $\phi: q ^{-1} q \simeq id_{L}$ as the relevant maps are pseudoinverse by \textbf{Proposition \ref{prop:serre_autoequivalence_pseudoinverse_to_its_pseudoinverse}}. Moreover, by \textbf{Proposition \ref{prop:left_adjoints_of_coevaluation_maps_in_terms_of_serre}}, the composites $(q^{-1} \otimes R) \circ b \circ c$ and $e \circ b \circ (R \otimes q)$ are left adjoint to (co)evaluation maps and this allows us to define the necessary units and counits.
 \end{proof}
 
We will now define coherent fully dual pairs by enforcing some additional equations at the level of $2$-cells. Our main result o this chapter, which we will prove later, is that the structure of a coherent fully dual pair in a symmetric monoidal bicategory is property-like and equivalent to full dualizability.
 \begin{defin}
 \label{defin:coherent_fully_dual_pair}
We say a fully dual pair $\langle L, R \rangle _{fd}$ is \textbf{coherent} if and only if the following conditions hold.
 
 \begin{enumerate}
 \item It is coherent as a dual pair.
 \item Witnessing isomorphisms $\phi, \psi$ make $q, q^{-1}$ into an adjoint equivalence
 \item The cusp-counits equation $(CC1) = (CC2)$, where the composites $(CC1)$, $(CC2)$ are given by pasting diagrams pictured below in \textbf{Figures \ref{fig:cusp_counits_composite_1}, \ref{fig:cusp_counits_composite_2}}, holds.
 \end{enumerate}
 
 \begin{figure}[htbp!]

	
\caption{Cusp-counits composite (CC1)}
\label{fig:cusp_counits_composite_1}
\end{figure}

 \begin{figure}[htbp!]
	%
	
\caption{Cusp-counits composite (CC2)}
\label{fig:cusp_counits_composite_2}
\end{figure}
\end{defin}
As we will see later, the cusp-counits equation is essentially the framed analogue of cusp flip relations from the presentation of the oriented bordism bicategory due to Christopher Schommer-Pries. For details on the latter type of relation, see \cite{chrisphd}. 

Like in the case of dual pairs, we can organize the collection of fully dual and coherent fully dual pairs in a given symmetric monoidal bicategory into bicategories on their own. To do so, we will use the theory of freely generated symmetric monoidal bicategories, introduced in \cite{chrisphd} under the name of \emph{computadic symmetric monoidal bicategories}. 

This theory is formally analogous to freely generated monoidal bicategories we present in \textbf{Appendix \ref{appendix_free_monoidal_bicategories}} and so at the level of what we're doing the passage from the monoidal to symmetric monoidal case is almost completely invisible. This is one of advantages of our approach.

\begin{defin}
Let $G_{fd}$ be the free generating datum for a symmetric monoidal bicategory consisting of

\begin{itemize}
\item two generating objects $L, R$,
\item four generating morphisms $e: L \otimes R \rightarrow I$, $c: I \rightarrow R \otimes L$, $q: L \rightarrow L$ and  $q^{-1}: L \rightarrow L$,
\item twelve generating $2$-cells $\alpha, \alpha^{-1}, \beta \beta ^{-1}, \phi, \phi^{-1}, \psi, \psi^{-1}, \epsilon_{e}, \mu_{e}, \epsilon _{c}$ and $\mu _{c}$, whose sources and targets are exactly like in the definition of a fully dual pair.
\end{itemize}
Let $\mathcal{R}_{fd}$ to be the class of relations on the free symmetric monoidal bicategory $\mathbb{F}(G_{fd})$ consisting of the relations that $\alpha, \alpha^{-1}$, $\beta, \beta ^{-1}$, $\phi, \phi^{-1}$ and  $\psi, \psi^{-1}$ are inverse to each other and that $\epsilon_{e}, \mu_{e}, \epsilon_{c},  \mu_{c}$ satisfy triangle equations.  Let $\mathcal{R} _{cfd}$ contain these relations and additionally the Swallowtail identities, triangle identities for $\phi, \psi$ and the cusp-counits equation.

We call the freely generated symmetric monoidal bicategory $\mathbb{F} _{fd} := \mathbb{F}(G_{fd}, \mathcal{R}_{fd})$ the \textbf{free symmetric monoidal bicategory on a fully dual pair} and $\mathbb{F}_{cfd} := \mathbb{F}(G_{fd}, \mathcal{R}_{cfd})$ the \textbf{free symmetric monoidal bicategory on a coherent fully dual pair}.
\end{defin}

\begin{defin}
 If $\mathbb{M}$ is a symmetric monoidal category, we define the \textbf{bicategory of fully dual pairs in $\mathbb{M}$} and the \textbf{bicategory of coherent fully dual pairs in $\mathbb{M}$} as the bicategories of shapes
\begin{gather*}
\mathcal{F}ullyDualPair(\mathbb{M}) := \mathbb{M}(G_{fd}, \mathcal{R}_{fd}), \\
\mathcal{C}ohFullyDualPair(\mathbb{M}) := \mathbb{M}(G_{fd}, \mathcal{R}_{cfd}).
\end{gather*}

Equivalently, the bicategory of dual pairs is precisely the bicategory of strict symmetric monoidal homomorphisms $\mathbb{F}_{fd} \rightarrow \mathbb{M}$ and similarly the bicategory of coherent dual pairs is the bicategory of strict homomorphisms $\mathbb{F}_{cfd} \rightarrow \mathbb{M}$. 
\end{defin}

Again, one can unpack the above concise definitions to yield a concise description, which we will do mainly to fix notation. Throughout, let $\mathbb{M}$ be the ambient symmetric monoidal bicategory. The objects of the bicategory $\mathcal{F}ullyDualPair(\mathbb{M})$ are precisely the fully dual pairs in $\mathbb{M}$. 

If $\langle L, R \rangle _{fd}, \langle L^\prime, R^\prime \rangle _{fd}$ are fully dual pairs, then a morphism $(s, t) _{fd}: \langle L, R \rangle _{fd} \rightarrow \langle L^\prime, R^\prime \rangle _{fd}$ consists of data of 1-cells $s: L \rightarrow L^\prime$, $t: R \rightarrow R^\prime$ and constraint isomorphisms $\gamma, \delta, \kappa, \tau$ as below.

\begin{center}

\end{center}
The constraint isomorphisms are assumed to satisfy naturality equations with respect to cusp isomorphisms $\alpha, \beta$, (co)units $\epsilon_{e}, \mu_{e}, \epsilon_{c}, \mu_{c}$ and witnessing isomorphisms $\psi$ and $\phi$.

If $(s_{1}, t_{1}) _{fd}, (s_{2}, t_{2}) _{fd}: \langle L, R \rangle _{fd} \rightarrow \langle L^\prime, R^\prime \rangle _{fd}$ are morphisms of fully dual pairs, then a $2$-cell $\Gamma: (s_{1}, t_{1}) _{fd} \rightarrow (s_{2}, t_{2}) _{fd}$ consists of data $2$-cells $\Gamma_{L}, \Gamma_{R}$ 

\begin{center}
	%
\end{center}
in $\mathbb{M}$. These are required to satisfy naturality with respect to constraint isomorphisms $\gamma, \delta, \kappa$ and $\tau$ of maps $(s_{1}, t_{1}) _{fd}, (s_{2}, t_{2}) _{fd}$. 

The bicategory $\mathcal{C}ohFullyDualPair(\mathbb{M}) \subseteq \mathcal{F}ullyDualPair(\mathbb{M})$ is the full subbicategory spanned by those fully dual pairs that are coherent.

An important property of the bicategory dual pairs, namely that it is a groupoid, implies the same property for bicategories of fully dual pairs.

\begin{prop}
The bicategory $\mathcal{F}ullyDualPair(\mathbb{M})$ is a 2-groupoid.
\end{prop}
Observe that this also shows that the bicategory of coherent fully dual pairs is a groupoid, too, as it is defined as a full subbicategory. 

\begin{proof}
A map of fully dual pairs induces a map of underlying dual pairs and similarly a $2$-cell between maps of fully dual pairs induces a $2$-cell between maps of dual pairs. Moreover, using the direct description one verifies immediately that these two notions differ only in constraint cells. Since non-constraint cells of a map or a $2$-cell in $\mathcal{D}ualPair(\mathbb{M})$ are invertible by \textbf{Proposition \ref{prop:bicategory_of_dual_pairs_is_a_groupoid}}, the same is true in the case of fully dual pairs and the result follows.
\end{proof}
The following theorem is a first step to our coherence theorem for fully dualizable objects, which we will prove in the next section.

\begin{thm}[Strictification of fully dual pairs]
\label{thm:strictification_of_fully_dual_pairs}
Any fully dualizable object can be completed to a coherent fully dual pair. 
\end{thm}

\begin{proof}
One can reduce to the case of symmetric monoidal bicategories with underlying \textbf{Gray}-monoid in the following way. If $\mathbb{M}$ is symmetric monoidal, let $\widetilde{\mathbb{M}}$ be its cofibrant replacement as a symmetric monoidal bicategory. That is, it is a symmetric monoidal bicategory freely generated by all the objects, morphisms and $2$-cells of $\mathbb{M}$ subject to relations that two $2$-cells are equal if and only if they coincide under the obvious strict homomorphism $\widetilde{\mathbb{M}} \rightarrow \mathbb{M}$. Moreover, let $\widetilde{\mathbb{M}} _{us}$ be the free unbiased semistrict symmetric monoidal bicategory on the same datum. 

We have a span $\widetilde{\mathbb{M}} _{us} \leftarrow \widetilde{\mathbb{M}} \rightarrow \mathbb{M}$ of symmetric monoidal bicategories and strict homomorphisms, the left arrow is an equivalence due to coherence theorem for unbiased semistrict symmetric monoidal bicategories and the right arrow is an equivalence since it is essentially surjective and fully faithful by construction. One can now deduce the statement for $\mathbb{M}$ from the corresponding statement for $\widetilde{\mathbb{M}} _{us}$ in the same way as we did in strictification theorem for dualizable objects. Since $\widetilde{\mathbb{M}} _{us}$ has an underlying \textbf{Gray}-monoid, we are done. 

Observe that we have in fact proven that one can reduce to the case of unbiased semistrict symmetric monoidal bicategories, we will however not need this stronger statement. From now on, assume that the ambient symmetric monoidal bicategory $\mathbb{M}$ has an underlying \textbf{Gray}-monoid.

Let $L$ be fully dualizable. By \textbf{Proposition \ref{prop:completion_to_a_fully_dual_pair}} we can complete it to a not-necessarily coherent fully dual pair and by strictification for dual pairs, that is, \textbf{Theorem \ref{thm:strictification_for_dual_pairs}}, we can make the chosen data satisfy Swallowtail identities by only a change of the cusp isomorphism $\beta$. 
Since $q, q^{-1}$ are pseudoinverse, there is a unique such $\phi$ such that together with $\psi$ they form an adjoint equivalence. 

We are now only left with enforcing the cusps-counits equation. Under our assumption of $\mathbb{M}$ having an underlying \textbf{Gray}-monoid, we can draw it as 
\begin{center}
	\begin{tikzpicture}[thick, scale=0.8, every node/.style={scale=0.8}]
	\node (L1) at (0,0) {$ L $};
	\node (LRL1) at (1.5, 1.5) {$ L \otimes R \otimes L $};
	\node (L2) at (3,3) {$ L $};
	\node (LRL2) at (4.5, 1.5) {$ L \otimes R \otimes L $};
	\node (L3) at (6, 0) {$ L $};
	
	\draw [->] (L1) to node[auto] {$ L \otimes c $} (LRL1);
	\draw [->] (LRL1) to node[auto] {$ e \otimes L $} (L2);
	\draw [->] (L2) to node[auto] {$ e^{L} \otimes L $} (LRL2);
	\draw [->] (LRL2) to node[auto] {$ L \otimes c^{L} $} (L3);
	\draw [->] (L1) to node[auto] {$ 1 $} (L3);
	\draw [->] (LRL1) to node[auto] {$ 1 $} (LRL2);
	
	\node at (3, 2.15) {$ \Downarrow \epsilon _{e} \otimes L $};
	\node at (3, 0.8) {$ \Downarrow L \otimes \epsilon_{c} $};
	
	\node at (9, 1.6) {$ (CC2) $};
	
	\node at (7.3, 1.6) {$ = $};
	\end{tikzpicture}.
\end{center}
It is important here to observe that the right hand side is an isomorphism and that it does not depend in any way on the (co)units of our fully dual pair. The left hand side is also an isomorphism. Indeed, it is the counit of the induced adjunction $(L \otimes c^{L}) \circ (e^{L} \otimes L) \dashv (e \otimes L) \circ (L \otimes c)$ and since the right adjoint, being isomorphic to the identity, is an equivalence, the adjunction must be an adjoint equivalence and so both (co)units are isomorphisms. It follows that the difference between the right and left hand sides in the equation above is at most some automorphism $\zeta: id_{L} \simeq id_{L}$, so that 

\begin{center}

\end{center}
already holds. The idea is now to "absorb" this $\zeta$ into one of the counits, this will enforce the cusps-counits equation. We can first commute $\zeta$ with $\epsilon_{c}$ to obtain that

\begin{center}
	%
\end{center}
We claim that the upper part of the diagram on the left hand side is already of the form $\epsilon \otimes L$ for some different counit of $e^{L} \dashv e$. This finishes the proof, as replacing the original $\epsilon_{e}, \mu_{e}$ by these new (co)units we obtain a fully dual pair that is coherent. We have 

\begin{center}
	%
,
\end{center}
as $\zeta$ is an automorphism of the identity, and the right hand side is precisely the counit of the adjunction $e^{L} \dashv e$ transferred along the isomorphisms $e \circ (\zeta \otimes R): e \simeq e$ and $id_{e^{L}}: e^{L} \simeq e^{L}$, tensored via $- \otimes L$.
\end{proof}

\subsection{Coherence for fully dualizable objects}

We will devote this section to the proof of the coherence theorem for fully dualizable objects, which identifies the notion of a coherent fully dual pair as a property-like structure equivalent to full dualizability. 

\begin{thm}
\label{thm:coherence_for_fully_dualizable_objects}
Let $\mathbb{M}$ be symmetric monoidal bicategory. The forgetful homomorphism 
\begin{gather*}
\pi: \mathcal{C}ohFullyDualPair(\mathbb{M}) \rightarrow \mathbb{K}(\mathbb{M}^{fd}), \\
\langle L, R \rangle _{fd} \mapsto L
\end{gather*}
between, respectively, the bicategory of coherent fully dual pairs in $\mathbb{M}$ and the groupoid of fully dualizable objects in $\mathbb{M}$, is a surjective on objects weak equivalence.
\end{thm}
Observe that by the cofibrancy theorem for freely generated symmetric monoidal bicategories, the bicategory of coherent fully dual pairs in $\mathbb{M}$ is equivalent to the homomorphism bicategory $\textbf{SymMonBicat}(\mathbb{F}_{cfd}, \mathbb{M})$, where $\mathbb{F}_{cfd}$ is the free symmetric monoidal bicategory on a coherent fully dual pair. 

Thus, the coherence theorem above identifies the latter as the "free symmetric monoidal bicategory on a fully dualizable object". This reduces the proof of the Cobordism Hypothesis to showing that $\mathbb{F}_{cfd}$ is equivalent to the framed bordism bicategory, a fact we will exploit in later chapters. 

Throughout the rest of this section we will be working with fixed ambient symmetric monoidal bicategory $\mathbb{M}$, let us make it a standing assumption that this $\mathbb{M}$ has an underlying \textbf{Gray}-monoid  It is enough to prove the statements that follow in this case, one can then deduce them for general $\mathbb{M}$ as in the proof of \textbf{Theorem \ref {thm:strictification_of_fully_dual_pairs}} and \textbf{Lemma \ref{lem:on_stability_of_properties_of_forgetful_functor_from_dual_pairs}}. 

We will now proceed with the proof. We have already established surjectivity on objects in the form of the strictification theorem and so we are left with essential surjectivity on morphisms, which we prove as \textbf{Lemma \ref{lem:forgetful_from_coherent_fully_dual_pairs_essentially_surjective}}, and local bijectivity on $2$-cells, which is \textbf{Lemma \ref{lem:forgetful_homomorphism_from_fully_dual_pairs_is_fully_faithful}}. We start with a little technical proposition, analogue of which in the dualizable case was phrased directly into the strictification theorem. 

\begin{prop}
\label{prop:counits_of_coevaluation_determine_counits_of_evaluation}
Suppose we have fully dual pairs $\langle L, R \rangle _{fd}$ and $\langle L^\prime, R^\prime \rangle _{fd}$ that share the same objects, the same 1-cells and the same 2-cells except maybe for the (co)units of the adjunction $e^{L} \dashv e$. Then, if they both satisfy the cusps-counits equation, they are necessarily equal, so that we also have $\epsilon _{e} = \epsilon^{\prime} _{e}$ and $\mu _{e} = \mu ^{\prime} _{e}$. 
\end{prop}

\begin{proof}
As both fully dual pairs in question satisfy cusp-counits equation, we also have

\begin{center}
,
\end{center}
as these are exactly the $(CC1)$ composites of the dual pairs in question, and so must be equal to $(CC2)$, which doesn't depend at all on the counits. Both sides are an isomorphism.

Since both $\epsilon_{e}$ and $\epsilon^{\prime}_{e}$ are counits of an adjunction $e^{L} \dashv e$, they can at most differ by an automorphism of any of the components. Let $\zeta: e \simeq e$ be exactly the difference between them, so that we have 

\begin{center}
	%
.
\end{center}
by taking the equation above and replacing $\epsilon^{\prime} _{e}$ by $\epsilon _{e}$ together with the difference $\zeta$. We want to show that this equation implies that $\zeta = id_{e}$, as then we will know that $\epsilon _{e} = \epsilon^{\prime} _{e}$. 

Using the fact that the big triangle composite is an isomorphism, the above equation is equivalent to 

\begin{center}
	%
\end{center}
being the identity, as one sees from pasting the inverse to the triangle from below. The composite of the last two morphisms is left adjoint to one of the cusp composites, so it is an equivalence, and thus the equation above is equivalent to 

\begin{center}
	%
.
\end{center}
By our standard lemma on equivalences of Hom-categories, that is, \textbf{Lemma \ref{lem:hom_equivalences}}, this equation implies that $\zeta$ is the identity, ending the argument.
\end{proof}

\begin{lem}
\label{lem:forgetful_from_coherent_fully_dual_pairs_essentially_surjective}
The homomorphism $\pi: \mathcal{C}ohFullyDualPair(\mathbb{M}) \rightarrow \mathbb{K}(\mathbb{M}^{fd})$ is essentially surjective on morphisms.
\end{lem}

\begin{proof}
We have to show that given coherent fully dual pairs $\langle L, R \rangle _{fd}, \langle L^\prime, R^\prime \rangle _{fd}$ and an equivalence $s: L \rightarrow L^\prime$ in $\mathbb{M}$, one can lift it to an equivalence of fully dual pairs, at least up to isomorphism.

We first choose the component $t: R \rightarrow R^\prime$ so that it commutes with (co)evaluation maps up to isomorphism, this can be done by coherence for dualizable objects in monoidal categories. We are left with choosing the constraint isomorphisms $\gamma, \delta, \kappa, \tau$ in a way that they will be natural with respect to all the 2-cells that are part of the structure of a fully dual pair.

Since both pairs are coherent, among these 2-cells we have three (co)unit pairs, $\epsilon_{e}, \mu_{e}$ of $e^{L} \dashv e$, $\epsilon_{c}, \mu_{c}$ of $c^{L} \dashv c$ and $\psi, \phi$ of $q \dashv q^{-1}$. As counits determine units and vice versa, it is enough to choose our constraint isomorphisms in a way that they are natural with just one of maps from each (co)unit pair, this can be proven by a variation on the "cute trick" we used in the proof of the corresponding result in the dualizable case. 

In detail, suppose we have chosen constraint isomorphisms such that they are natural with respect to either the unit or the counit from one of the three pairs above. This determines an equivalence $\langle L, R \rangle _{fd} \rightarrow \langle L^\prime R^\prime \rangle _{fd}$ of $G_{fd}^{\prime}$-shapes, where $G_{fd}$ is the 1-truncation of $G_{fd}$. We can use this equivalence to transport the structure of a $G_{fd}$ shape from $\langle L, R \rangle _{fd}$ to $\langle L^\prime, R^\prime \rangle _{fd}$ by \textbf{Lemma \ref{lem:promoting_equivalences_of_g_shapes}} and since this equivalence was in fact natural with respect to the chosen (co)unit, the original (co)unit of $\langle L, R \rangle _{fd}$ and the one obtained by transport of structure will necessarily agree. Since (co)units determine each other, the original and transported structure must also agree on the other (co)unit and so our equivalence was natural with respect to both of them to begin with.

By an analogous argument and the technical \textbf{Proposition \ref{prop:counits_of_coevaluation_determine_counits_of_evaluation}} above we see that if we choose constraint isomorphisms that are natural with respect to all structural 2-cells except maybe $\epsilon_{e}, \mu_{e}$, then they will already be natural with respect to all of them.

We have already proven as part of the of the coherence theorem for dual pairs that naturality with respect to $\alpha$ implies naturality with respect to $\beta$ and moreover for any chosen constraint isomorphism $\gamma$ there is a unique constraint $\delta$ such that together they satisfy both equations. 

These few remarks give enough of "determination of data" for our purposes and we can now prove the existence of the necessary isomorphisms directly. The rough outline is as follows. We will first choose $\delta, \gamma$, then choose "some" $\kappa$, correct it to the "right" $\kappa$ and at the very end choose the constraint $\tau$. 

Fix $\delta, \gamma$ such that they are natural with respect to cusp isomorphisms and also choose constraint isomorphism $\kappa$, not necessarily natural with respect to anything. The latter is possible because the Serre autoequivalence is natural up to isomorphism, which is \textbf{Proposition \ref{prop:serre_autoequivalence_is_natural}}. Naturality with respect to $\epsilon_{c}$ is
\begin{center}
.
\end{center}
However, the left hand side is the counit $\epsilon_{c}$ of the pair $\langle L, R \rangle _{fd}$ transported to the pair $\langle L^\prime, R^\prime \rangle _{fd}$, so it differs from $\epsilon ^{\prime} _{c}$ at most by an automorphism of $c^{L \prime}$, as they are both counits of the same adjunction. In other words, we can find some $\zeta: c^{L \prime} \simeq c^{L \prime}$ such that 

\begin{center}
	%
.
\end{center}
already holds. The idea, very similar to the one we used in the proof of \textbf{Proposition \ref{prop:counits_of_coevaluation_determine_counits_of_evaluation}}, is to "absorb" this $\zeta$ into the isomorphism $\kappa$. This new "corrected" $\widetilde{\kappa}$, together with the $\gamma, \delta$ we already have, will satisfy $\epsilon_{c}$ naturality since the above equation holds. 

Observe that we could in principle also try to leave $\kappa$ unchanged and try to absorb $\zeta$ into, for example, $\delta$, which is a little easier. However, we've already chosen $\gamma, \delta$ in a way that they satisfy naturality with respect to cusp isomorphisms and we cannot now change one of them without changing the other, or we may lose this property.

Our goal is now to prove that there exists $\widetilde{\kappa}$ such that we have 

\begin{center}
,
\end{center}
if we then replace $\kappa$ by $\widetilde{\kappa}$, $\epsilon_{c}$-naturality will hold. In fact we will see that such $\widetilde{\kappa}$ must be unique, this will be important later on, in the proof of \textbf{Lemma \ref{lem:forgetful_homomorphism_from_fully_dual_pairs_is_fully_faithful}}.

As $b$ is an equivalence, there exists a unique automorphism $\zeta ^\prime$ of $e^\prime \circ (q^\prime \otimes R^\prime)$ which has the property that it reduces to $\zeta$ under composition with $b$. Thus, an equation equivalent to the one given above is 

\begin{center}
.
\end{center}
As the naturality constraint for $b$ and $\gamma$ are isomorphisms, we can omit them to obtain an equivalent equation. One can then remove $b$ is well, as it is an equivalence. In the end, we are left with showing that the equation 

\begin{center}
	%
.
\end{center}
has a unique solution $\widetilde{\kappa}$. By adding the necessary interchangers to the left bottom square we can move the $t \otimes (-)$ component to the front and whisker $t$ out, obtaining an equivalent equation. We can then leave out $t$ altogether, as it is an equivalence, so that the are left with 

\begin{center}
	%
.
\end{center}
This has clearly a unique invertible solution $\widetilde{\kappa}$ by our standard \textbf{Lemma \ref{lem:hom_equivalences}} on Hom-equivalences, as the left hand side is an isomorphism not depending on $\widetilde{\kappa}$. The unique solution is the "corrected" constraint $\kappa$. 

The argument above shows that we can choose constraint isomorphisms $\gamma, \delta, \kappa$ such that they satisfy naturality with respect to cusp isomorphisms and $\epsilon_{c}$. By our reasoning with transport of structure this also immediately implies that they are natural with respect to $\mu _{c}$.

We will now choose $\tau$, it has to satisfy the equation of $\psi$-naturality, which is

\begin{center}
	\begin{tikzpicture}[thick, scale=0.8, every node/.style={scale=0.8}]
		\node (L1) at (0,0) {$ L $};
		\node (Lp1) at (0,1.5) {$ L^\prime $};
		\node (L2) at (2, 0.5) {$ L $};
		\node (Lp2) at (2, 2) {$ L^\prime $};
		\node (L3) at (4, 0) {$ L $};
		\node (Lp3) at (4, 1.5) {$ L ^\prime $};
		
		\draw [->] (L1) to node[auto] {$ s $} (Lp1);
		\draw [->] (L2) to node[auto] {$ s $} (Lp2);
		\draw [->] (L3) to node[auto] {$ s $} (Lp3);
		\draw [->] (L1) to node[auto] {$ q^{-1} $} (L2);
		\draw [->] (L2) to node[auto] {$ q $} (L3);
		\draw [->] (Lp1) to node[auto] {$ q^{\prime -1} $} (Lp2);
		\draw [->] (Lp2) to node[auto] {$ q^{\prime} $} (Lp3);
		\draw [->] (L1) to [out=-30, in=-150] node[auto] {$ 1 $} (L3);
		
		\node at (1, 1.1) {$ \simeq \tau $};
		\node at (3, 1.1) {$ \simeq \kappa $};
		\node at (2, 0) {$ \simeq \psi $};
		
		\node (rL1) at (6,0) {$ L $};
		\node (rLp1) at (6,1.5) {$ L^\prime $};
		\node (rLp2) at (8, 2) {$ L^\prime $};
		\node (rL3) at (10, 0) {$ L $};
		\node (rLp3) at (10, 1.5) {$ L ^\prime $};
		
		\draw [->] (rL1) to node[auto] {$ s $} (rLp1);
		\draw [->] (rL3) to node[auto] {$ s $} (rLp3);
		\draw [->] (rLp1) to node[auto] {$ q^{\prime -1} $} (rLp2);
		\draw [->] (rLp2) to node[auto] {$ q^{\prime} $} (rLp3);
		\draw [->] (rL1) to [out=-30, in=-150] node[auto] {$ 1 $} (rL3);
		\draw [->] (rLp1) to [out=-30, in=-150] node[auto] {$ 1 $} (rLp3);
		
		\node at (8, 1.5) {$ \simeq \psi ^\prime $};
		\node at (8, 0.2) {$ = $};
		
		\node at (5, 0.7) {$ = $};
	\end{tikzpicture}.
\end{center}
This clearly has a unique invertible solution $\tau$, as everything in sight is either invertible or an equivalence. This $\tau$ will also then satisfy $\phi$-naturality, as $\psi, \phi$ are a (co)unit pair. 

Hence, we have chosen the constraint isomorphisms $\gamma, \delta, \kappa, \tau$ such that they are natural with respect to $\alpha, \beta, \epsilon_{c}, \mu_{c}, \psi, \phi$. We've already given argument above how in this case they are in fact natural with respect to all structure 2-cells and so we are done. 
\end{proof}

\begin{lem}
\label{lem:forgetful_homomorphism_from_fully_dual_pairs_is_fully_faithful}
The homomorphism $\pi: \mathcal{F}ullyDualPair(\mathbb{M}) \rightarrow \mathbb{K}(\mathbb{M}^{fd})$ is locally bijective on $2$-cells.
\end{lem}

\begin{proof}
One reduces to the case of the ambient symmetric monoidal bicategory $\mathbb{M}$ having an underlying \textbf{Gray}-monoid in the usual way, using the cofibrant replacement.

Observe that the datum of a transformation between maps of fully dual pairs is the same as for transformation between maps of dual pairs, since it consists only of components indexed on objects, only the axioms differ. This will allow us to reduce most of the statement to the already proven analogue for dual pairs, which was \textbf{Lemma \ref{lem:bijectivity_on_isomorphisms_between_equivalences}}.

We will first prove a "little claim" on uniqueness of data in any given equivalence of fully dual pairs and later use the promotion of isomorphisms to deduce what we need.

Suppose $(s_{1}, t_{1}) _{fd}, (s_{2}, t_{2}) _{fd}: \langle L, R \rangle _{fd} \rightarrow \langle L^\prime, R^\prime \rangle _{fd}$ are equivalences of fully dual pairs that agree on components on objects and also agree on the constraint isomorphisms $\gamma, \delta$. We claim that they must be in fact equal.

Since both $(s_{1}, t_{1}) _{fd}, (s_{2}, t_{2}) _{fd}$ satisfy naturality with respect to $\epsilon_{c}$, we must have

\begin{center}
,
\end{center}
as then both sides are simply equal to $\epsilon^{\prime} _{c}$. If we paste 

\begin{center}
	%
\end{center}
to the right hand side along their common boundary $c^\prime = c^\prime \otimes 1$, then after an application of the naturality of $(s_{2}, t_{2})_{fd}$ with respect to $\mu_{c}$ and triangle equations for $\mu_{c}, \epsilon_{c}$, we will be left with 

\begin{center}
	%
.
\end{center}
Similarly, doing the same pasting to the right hand side we will obtain 

\begin{center}
	%
.
\end{center}
so these two diagrams must be equal. However, we've already seen in the proof of the \textbf{Lemma \ref{lem:forgetful_from_coherent_fully_dual_pairs_essentially_surjective}} that these composites uniquely determine $\kappa$, so that also have $\kappa_{1} = \kappa_{2}$. Since constraint $\kappa$ and $\tau$ uniquely determine each other via $\psi$ or $\phi$-naturality, it follows that also $\tau_{1} = \tau_{2}$. This ends the proof of our "little claim", we will now show how to use it to obtain the statement of the lemma.

We want to show that $\pi: \mathcal{F}ullyDualPair(\mathbb{M}) \rightarrow \mathbb{K}(\mathbb{M}^{fd})$ is locally bijective on $2$-cells, so let  $(s_{1}, t_{1}) _{fd}, (s_{2}, t_{2}) _{fd}: \langle L, R \rangle _{fd} \rightarrow \langle L^\prime, R^\prime \rangle _{fd}$ be equivalences of fully dual pairs and suppose that we have an isomorphism $\Gamma _{L}: L \rightarrow L^\prime$. We want to lift it to an isomorphism to maps of fully dual pairs, which amounts to defining its component $\Gamma _{R}$ on the other object, see the proof of  \textbf{Lemma \ref{lem:bijectivity_on_isomorphisms_between_equivalences}}.

By the lemma we just mentioned, there is a unique way to define the lift $\Gamma = (\Gamma_{L}, \Gamma_{R})$  such that it will be natural with respect to the constraint isomorphisms $\gamma, \delta$, that is, the ones coming from the structure of a dual pair itself. This is enough to define an isomorphism $\Gamma: (s_{1}, t_{1}) \rightarrow (s_{2}, t_{2})$ of maps of $G_{fd}^{\prime \prime}$-shapes, where $G_{fd} ^{\prime \prime}$ is the $0$-truncation of $G_{fd}$. 

We can use promotion of isomorphisms between such maps, which is \textbf{Lemma \ref{lem:promotion_of_isomorphisms_between_maps_of_gpp_shapes}}, to transport the structure of a map of $G_{fd}$-shapes from $(s_{1}, t_{1}) _{fd}$ to $(s_{2}, t_{2})$ along the isomorphism $\Gamma$ to obtain a new map of $G_{fd}$-shapes. Since it is isomorphic to one, it is also a map of fully dual pairs and we will denote it by $\widetilde{(s_{2}, t_{2})} _{fd}$. It agrees with $(s_{2}, t_{2}) _{fd}$ on its components on objects by construction.

The isomorphism $\Gamma$ was already natural with respect to $\gamma, \delta$ and so the equivalences $\widetilde{(s_{2}, t_{2})} _{fd}$ and $(s_{2}, t_{2}) _{fd}$ must also agree on these two constraint isomorphisms. By our "little claim" above, they must in fact be equal and we deduce that $\Gamma$ was a map of equivalences of fully dual pairs to begin with, ending the proof of the lemma.
\end{proof}

\section{Bordism bicategories}

In this chapter we start the geometric part of current work, giving a brief introduction to bordism bicategories. We will highlight some issues related to their construction and describe the program we will follow to obtain a presentation of the framed bordism bicategory. Nothing here is new, the definitions are due to Christopher Schommer-Pries. We don't include proofs and our exposition is rather informal, a detailed reference for the material is \cite{chrisphd}. 

The idea of axiomatizing a \emph{topological field theory} in terms of a symmetric monoidal homomorphism from a bordism category into vector spaces is due to Atiyah, see \cite{atiyah_tqfts}. The adjective \emph{topological} is used to emphasize that the value of the theory depends only on the underlying topology of the manifold, one can also consider variants involving some geometrical structures.

One source of intuition is as follows. Let $Z$ be an $n$-dimensional topological field theory. Observe that a closed manifold $M^{n}$ can be understood as a bordism from the empty set to itself, since the latter is a monoidal unit in the category of $n$-bordisms, $Z$ will assign to $M$ an element $Z(M) \in Z(\emptyset) \simeq \mathbb{C}$. Seen this way, a topological field theory is a family of numerical invariants of $n$-manifolds, the composition laws for monoidal homomorphism then imply that the value of $Z$ at any given $M$ can be computed by cutting it along codimension $1$ submanifolds.

In their landmark paper \cite{baezdolanhigheralgebraandtqfts} John Baez and James Dolan introduced the idea of using higher categories to formalize the idea of cutting along submanifolds of higher codimension, this at the same time allows a simplficiation of the process of computing a value of a field theory and on the other hand enforces its locality. 

There are many technicalities involved in the construction of higher bordism categories, or higher categories in general, but by now many adequate definitions are known. The one used by Jacob Lurie in his proof of the Cobordism Hypothesis is based on complete $n$-fold Segal spaces, see \cite{lurietqfts}.

The work of Christopher Schommer-Pries, which we extend, is based on symmetric monoidal bicategories, this effectively allows one to consider submanifolds up to codimension $2$. The advantage of this approach is that it is fully algebraic, which means that all the composites one expects to exist in a higher category are always explicitly defined. It is also well studied by category theorists, with many examples of interesting targets for topological field theories. 

We start by giving an informal definition of the unoriented bordism bicategory.

\begin{defin}[Informal]
The bordism bicategory $\mathcal{B}ord_{2}$ is the symmetric monoidal bicategory described as follows.

\begin{itemize}
\item Its objects are $0$-manifolds.
\item Its morphisms are given by bordisms, ie. $1$-manifolds with a suitable decomposition of a boundary.
\item Its $2$-cells are given by diffeomorphism classes of $2$-bordisms relative to the boundary. 
\end{itemize}
Horizontal and vertical composition are given by gluing manifolds along boundaries, the monoidal structure is induced by disjoint union.
\end{defin}
The notions of $0$-manifolds and bordisms between these are familiar, $2$-bordisms are more mysterious. The intuition is that a $2$-bordism $\alpha: w_{1} \rightarrow w_{2}$, where $w_{i}$ are $1$-bordisms with the same source $S$ and target $T$, is a surface with boundary identified with some $1$-manifold glued from $w_{1}, w_{2}$. 

Since the source and target of $w_{i}$ coincide, one could simply glue them along these, but this is not what we want, as it would make horizontal and vertical composition difficult. The usual thing to do is to consider surfaces with boundary modeled on the square, with $w_{1}$ being at the top, $w_{2}$ being at the bottom and the vertical segments being trivial bordisms over the source and target, ie. $S \times I, T \times I$, where $I$ is the interval. 

\begin{figure}[htbp!]
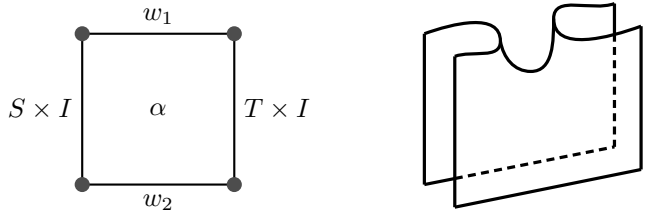


\caption{Schematic for the boundary of a $2$-bordism on the left, an example on the right}
\label{fig:schematic_for_a_2_bordism}
\end{figure}

A schematic for the boundary of a $2$-bordism is given in \textbf{Figure \ref{fig:schematic_for_a_2_bordism}}\footnote{We draw $1$-bordisms from left to right, that is, on the left we have always the source and on the other side the target of a $1$-bordism. Likewise, we draw $2$-bordisms from top to bottom.} on the left, on the right side we draw a classic example. Observe that in general the underlying manifold of a $2$-bordism is a manifold with corners. 

This informal definition of the unoriented bordism bicategory can be made precise, but there are some difficulties. 

One of the problems arises from the fact that if $M, N$ are smooth manifolds and we are given a diffeomorphism $f: \partial M \simeq \partial N$, then the gluing $M \cup _{f} N$ is a well-defined topological manifold, but the smooth structure on it is only unique up to a non-unique diffeomorphism. Hence, there is some ambiguity in defining composition of $1$-bordisms.

This is not an issue in the usual bordism categories, as therein the morphisms are given only by \emph{diffeomorphism classes} of bordisms and so the composition given by gluing is well-defined. We cannot do that, as to have an acceptable notion of a $2$-bordism, one should work with actual $1$-bordisms and not just diffeomorphisms classes of them. 

A standard approach, and the one taken in \cite{chrisphd}, is to consider manifolds equipped with some additional structure that facilitates well-defined composition. For example, a gluing of two manifolds with boundary is essentially unique as soon as we choose collars on both sides. 

In our case, this additional data attached to a manifold will take an intuitive form of a "infinitesimal neighbourhood of a $2$-manifold". This has an added advantage of making it much easier to discuss manifolds with some extra structure, like orientations or framings, as this topological structure can always be taken to live in that imaginary neighbourhood and so there are no problems with restricting it to submanifolds of smaller dimension. 

\begin{rem}
Another possible approach would be to adapt a non-algebraic definition of a higher category, so that the gluing need not be well-defined, as long as it is "unique up to a coherent isomorphism". This is the approach taken in \cite{lurietqfts}.
\end{rem}

\subsection{Construction of bordism bicategories} 

In this section we will give a sketch of construction of bordism bicategories, including the notion of an "infinitesimal neighbourhood" or a \emph{halo}. The right notion to formalize this idea is that of a pro-manifold.

\begin{defin}
A \textbf{pro-manifold} $X$ is a diagram $D \rightarrow \mathcal{M}an$ in the category of not necessarily compact manifolds, where $D$ is a small cofiltered category. Its \textbf{associated copresheaf} $LX$ is the colimit of $X$ taken in the category of copresheaves over manifolds. A \textbf{map} $X \rightarrow Y$ of pro-manifolds is a map of associated copresheaves $LY \rightarrow LX$. 
\end{defin}

Observe that in particular every manifold gives a pro-manifold indexed on the terminal category, moreover the inclusion $\mathcal{M}an \hookrightarrow \text{pro-} \mathcal{M}an$ is full and faithful. This is basically the content of co-Yoneda lemma.

\begin{rem}
If $X: D_{1} \rightarrow \mathcal{M}an, Y: D_{2} \rightarrow \mathcal{M}an$ are pro-manifolds, then one can also compute the set of maps $X \rightarrow Y$ directly as 

\begin{center}
$\text{pro-}\mathcal{M}an(X, Y) = lim _{d_{2} \in D_{2}} colim _{d_{1} \in D_{1}} \ \mathcal{M}an(X(d_{1}), Y(d_{2}))$
\end{center}
with both limit and colimit taken in sets. In particular, if we test a pro-manifold $X$ using ordinary manifolds, it behaves like the limit of the diagram it represents, even if the actual limit doesn't exist. Mapping out of a pro-manifold is, however, more interesting. 
\end{rem}

Before introducing haloes, we introduce their components, which are termed \emph{halations}. The intuition is that a halation over a manifold $X$ is what is left from the embedding $X \rightarrow Y$ into an interior of some other manifold $Y$ if we remember only the infinitesimal neighbourhood of $X \subseteq Y$. A halo will then constitute a sequence of halations of varying dimensions.

\begin{defin}
Let $i: X \hookrightarrow Y$ be an embedding of manifolds, possibly with boundary or corners, such that $X$ is contained in the interior of $Y$. Then the \textbf{halation associated to $i$} is the inclusion of pro-manifolds $X \hookrightarrow \hat{X}$, where $\hat{X}$ is the inverse system of open submanifolds $Z \subseteq Y$, without boundary and containing $X$, ordered by inclusion.

We define the category of \textbf{halations over a manifold $X$} to be the full replete subcategory of the category of pro-manifolds over $X$ spanned on halations associated to embeddings as above. In other words, an arrow $X \hookrightarrow \hat{X}$ is a halation if and only if it is isomorphic to one associated to some embedding of manifolds as above.
\end{defin}

By definition any halation of a manifold $X$ is isomorphic to one coming from an embedding of $X$ into some other manifold, however, for two halations coming from $X \subseteq Y_{1}$ and $X \subseteq Y_{2}$ to be isomorphic it is enough for some arbitrarily small neighbourhoods of $X$ to be diffeomorphic, the manifolds $Y_{1}, Y_{2}$ can be wildly different. This formalizes the intution of an "infinitesimal neighbourhood". It can be a little simpler to work with than $X$ itself, because it always consists of an inverse system of manifolds without boundary.

\begin{defin}
The \textbf{dimension} of a halation associated to an embedding $i: X \rightarrow Y$ is the dimension of $Y$. It is an intrinsic invariant of a halation.
\end{defin}

A principal example of a halation over $X$, where $X$ is without boundary, is a one coming from the total space of a vector bundle over $X$. One can prove using the tubular neighbourhood theorem that any halation over such manifold is isomorphic to a one obtained from a vector bundle. This  fact can be generalized to manifolds with corners, but a little bit more data than just the vector bundle is necessary.

If $X \hookrightarrow \hat{X}$ is a halation over $X$, then it restricts naturally to \emph{any} submanifold $X^\prime$ by choosing any representative embedding $X \hookrightarrow Y$ and restricting it to $Z$, this allows one to naturally define bordisms with halation between manifolds with halation. If we consider $n$-bordisms between $n-1$-manifolds, then the additional data of a dimension $n$-halation does facilitate well-defined composition, essentially because gluing of smooth manifolds along diffeomorphisms of open subsets \emph{is} well-defined.

Hence, in our bordism bicategory we will consider manifolds which are equipped with halations of dimensions $1$ and $2$, this gives a well-defined composition of $1$-bordisms and $2$-bordisms. However, there is some additional trivialization data that needs to be taken into account. 

To see it is necessary, assume for simplicity that $X^{n}$ is a manifold with boundary. Then not every $n$-dimensional halation over $\partial X$ comes from one defined over $X$, as the latter ones can always be trivialized using the inward-pointing normal vector. It follows that a manifold equipped with a non-trivial codimension $1$ halation cannot be a source or a target of any $n$-bordism. 

Thus, our halations need to be trivial, that is,  isomorphic to ones coming from a trivial vector bundle. Additionally, we should also remember the trivializations, so that two haloes over the same manifold will be canonically isomorphic in a suitable sense, with the space of such canonical isomorphisms being contractible. The technical tool to trivialize halations is that of a \emph{coorientation}.

\begin{defin}
Let $X \hookrightarrow \hat{X}$ be a halation which we identify with a halation coming from an embedding $X \hookrightarrow Y$, so that $\hat{X}$ runs through all the open $Z \subseteq Y$ without boundary, containing $X$. A \textbf{coorientation} of $X \hookrightarrow \hat{X}$ is a compatible choice of orientations for the normal bundles $N_{X \subseteq Z}$. 
\end{defin} 
Observe that a choice of an orientation can be understood as a weak form of trivializing a vector bundle. In particular, if $X \rightarrow \hat{X}$ is of codimension $1$, then a choice of coorientation is the same as a choice of a trivialization up to isotopy. 

One can show that a coorientation is in fact a structure intrinsic to the halation $X \hookrightarrow \hat{X}$, independant on the choice of identification with a halation of some embedding. By abuse of language, we assume that a codimension $0$ halation is always cooriented in a unique possible way. 

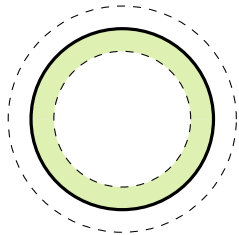
\begin{figure}[htbp!]
	\begin{tikzpicture}[scale=0.6]
		\fill [halo_green!30] (1.5, 0) to (2, 0) arc (0 : 180 : 2) to (-1.5, 0) arc (180 : 0 : 1.5);
		\fill [halo_green!30] (1.5, 0) to (2, 0) arc (0 : -180 : 2) to (-1.5, 0) arc (180 : 360 : 1.5);

		\draw [dashed] (0,0) circle (2.5cm);
		\draw [very thick] (0,0) circle (2cm);
		\draw [dashed] (0,0) circle (1.5cm);
	\end{tikzpicture}

\caption{A codimension $1$ cooriented halation on a circle}
\label{fig:codimension_1_cooriented_halation_on_a_circle}
\end{figure}

In \textbf{Figure \ref{fig:codimension_1_cooriented_halation_on_a_circle}} we present a circle equipped with a halation, which we draw as a small neighbourhood. The distinguished normal direction coming from coorientation is denoted with a light green shading. 

The precise definition of a manifold equipped with a suitably trivialized germ of neighbourhoods is then the following.

\begin{defin}
Let $X^{k}$ be a manifold, possibly with corners. If $k < 2$, then a  \textbf{$2$-halo} over $X$ is a sequence of inclusions 

\begin{center}
$X \subseteq \hat{X}_{1} \subseteq \hat{X}_{2}$ 
\end{center}
such that $X \subseteq \hat{X}_{1}$, $X \subseteq \hat{X}_{2}$ have a structure of cooriented halations of dimensions $1$ and $2$, respectively. If $X$ is a surface, then a \textbf{$2$-halo} over $X$ is simply a cooriented halation $X \subseteq \hat{X}_{2}$ of dimension $2$. 
\end{defin}
When defining haloes of higher codimension, the sequence of cooriented halations included as part of the data would get suitably longer. An interesting phenomena occurs with respect to coorientations, as in general a cooriented halation need not be trivial. Here, however, all of the coorientations \emph{taken together} will force each of the halations occuring to be trivial, more or less because it equips the normal space to the halation at each point of $X$ with a complete flag of relatively oriented subspaces. 

\begin{figure}[htbp!]
	\begin{tikzpicture}[scale=0.7]
		\draw [dashed ] (0,0) circle (1.5cm);
		\fill [halo_green!30] (-1.5, 0) to (1.5,0) arc (0 : 180 : 1.5);
		
		\draw [very thick, dashed] (-1.5, 0) to (0,0);
		\draw [very thick] (0,0) to (1.5, 0);
		\node [circle, fill=black, inner sep=1.5pt] at (0,0) {};
	\end{tikzpicture}
\caption{Haloed point}
\label{fig:haloed_point}
\end{figure}
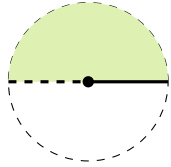

We will now explain how to make the informal definition of the bordism bicategory $\mathbb{B}ord_{2}$ precise. Its objects are given by $2$-haloed, or simply \emph{haloed}, $0$-manifolds. An example of a haloed $0$-manifold whose underlying manifold is a point is drawn in \textbf{Figure \ref{fig:haloed_point}}.

Its first halation $\hat{pt}_{1}$ is drawn as a small open interval, its coorientation is denoted using a solid line. The second halation $\hat{pt}_{2}$ is drawn as a small neighbourhood of an open surface, the coorientation is depicted using shading, which denotes the the orientation of the normal bundle of the $1$-manifold inside the surface. This is the same data as a coorientation $pt \subseteq \hat{pt}_{2}$, as the normal bundle of $pt \subseteq \hat{pt}_{1}$ is also oriented.

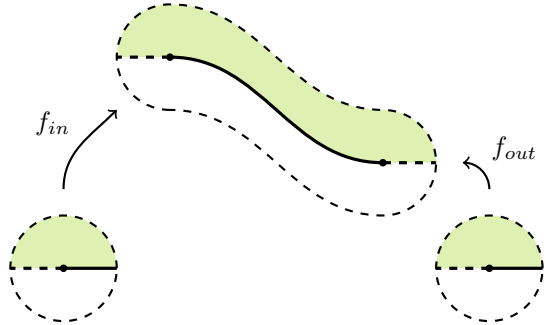
\begin{figure}[htbp!]
	\begin{tikzpicture}[thick, scale=0.7]
		\node (A) at (2,4) {};
		\node (B) at (2, 3) {};
		\node (X) at (6, 3) {};
		\node (Y) at (6,2) {};
		\fill [halo_green!30] (A.center) -- +(-1, 0) arc (180 : 90 : 1cm) to [ out = 0, in = 180] (X.center)
			arc (90 : 0 : 1cm) -- (Y.center) to [out = 180, in = 0] (A.center);
		\draw [dashed] (B.center) arc (270 : 90: 1cm)  to [ out = 0, in = 180] (X.center)
			arc (90: -90: 1cm) to [out = 180, in = 0] (B.center);
		\draw [very thick] (A.center) to [out = 0, in = 180]  (Y.center);
		\draw [dashed, very thick] (A.center) to +(-1,0) (A.center);
		\draw [dashed, very thick] (Y.center) to +(1,0) (Y.center);

		\node [circle, fill=black,inner sep=1pt] at (A.center) {};
		\node [circle, fill=black,inner sep=1pt] at (Y.center) {};

		\node (C) at (0,0) {};
		\node (D) at (8, 0) {};
		\fill [halo_green!30] (C.center) -- +(-1, 0) arc (180 : 0: 1cm) -- (C.center);
		\fill [halo_green!30] (D.center) -- +(-1, 0) arc (180 : 0: 1cm) -- (D.center);
		\draw [dashed] (C.center) circle (1cm);
		\draw [dashed] (D.center) circle (1cm);
		\draw [very thick] (C.center) -- +(1,0) (D.center) -- +(1,0);
		\draw [very thick, dashed] (C.center) -- +(-1,0) (D.center) -- +(-1,0);
		\node [circle, fill=black,inner sep=1pt] at (C.center) {};
		\node [circle, fill=black,inner sep=1pt] at (D.center) {};
	
		\draw [->] (0, 1.5) to [out = 90, in = 225] node [above left] {$f_{in}$} (1, 3);
		\draw [->] (8, 1.5) to [out = 90, in = 0] node [above right] {$f_{out}$} (7.5, 2);
	\end{tikzpicture}	
\caption{A haloed 1-bordism}
\label{fig:haloed_1_bordism}
\end{figure}

\begin{defin}
Let $A, B$ be haloed $0$-manifolds. A \textbf{haloed $1$-bordism} $w: A \rightarrow B$ consists of a haloed $1$-manifold together with a decomposition $\partial w \simeq \partial _{in} w \sqcup \partial _{out} w$ and isomorphisms 

\begin{center}
$f_{in}: (\partial _{in} w,  \hat{w}_{1} |_{\partial _{in}}, \hat{w} _{2} |_{\partial {in}}) \simeq (A, \hat{A}_{1}, \hat{A}_{2})$ \\
$f_{out}: (\partial _{out} w, \hat{w} _{1} |_{\partial _{out}}, \hat{w} _{2} |_{\partial {out}}) \simeq (B, \hat{B}_{1}, \hat{B}_{2})$
\end{center}
of haloed manifolds, where $\hat{w} _{1} |_{\partial _{in}}$ is equipped with the coorientation using the inward pointing normal vector and $\hat{w} _{1} |_{\partial _{out}}$ using the outward pointing normal vector.
\end{defin}
An example of a haloed $1$-bordism from a point to itself is given below in \textbf{Figure \ref{fig:haloed_1_bordism}}. The inclusions are the obvious ones. Observe how the coorientation of the $0$-manifold in its germ of $1$-manifold agrees with the inward direction of the $1$-bordism in the case of the domain and disagrees in the case of the codomain, this is one of the compatibility conditions we impose.

\begin{defin}
Let $A, B$ be haloed $0$-manifolds and $w, v: A \rightarrow B$ haloed $1$-bordisms between them. A \textbf{haloed 2-bordism} $\alpha: v \rightarrow w$ consists of 

\begin{itemize}
\item a haloed $2$-manifold with faces, which we also denote by $\alpha$
\item distinguished faces $\partial _{o}, \partial _{m}$ of $\alpha$ that satisfy $\partial _{o} \cup \partial _{m} = \partial \alpha$ and $\partial _{0} \cap \partial _{m} = \partial \partial _{o} = \partial \partial _{m}$
\item an isomorphism of cooriented halations
\[ (\partial _{m}, \hat{\alpha} |_{\partial _{m}}) \simeq (w, \hat{w}_{2}) \sqcup (v, \hat{v}_{2}) \]
where the coorientation of $w \subseteq \hat{w}_{2}$ agrees with the inward pointing normal direction in $\alpha$ and the coorientation of $v \subseteq \hat{v}_{2}$ agrees with the outward pointing normal direction. 
\item an isomorphism 
\[ (\partial _{o}, \hat{\alpha} _{2}) \simeq (A \times \{ 0 \} \times I, \widehat{A \times \mathbb{R}^{2}}) \sqcup (B \times \{ 0 \} \times I, \widehat{A \times \mathbb{R}^{2}}) \]
of cooriented halations, where the coorientation on $A \times I \subseteq \widehat{A \times \mathbb{R}^{2}}$ agrees with the inward pointing normal direction in $\alpha$ and the coorientation on $B \times I \subseteq \widehat{B \times \mathbb{R}^{2}}$ agrees with the outward orientation
\end{itemize}
Additionally, we require that the induced isomorphisms
\begin{gather*}
(A , \hat{A}_{1}, \hat{A}_{2}) \simeq (A \times \{ 0 \} \times \{ 1 \}, A_{\mathbb{R}}, A _{\mathbb{R}^{2}}) \simeq (A \times \{ 0 \} \times \{ 0 \}, A_{\mathbb{R}}, A _{\mathbb{R}^{2}}) \simeq (A, \hat{A}_{1}, \hat{A}_{2}) \\
(B , \hat{B}_{1}, \hat{B}_{2}) \simeq (B \times \{ 0 \} \times \{ 1 \}, B_{\mathbb{R}}, B _{\mathbb{R}^{2}}) \simeq (B \times \{ 0 \} \times \{ 0 \}, B_{\mathbb{R}}, B _{\mathbb{R}^{2}}) \simeq (B, \hat{B}_{1}, \hat{B}_{2})
\end{gather*}
of haloed $0$-manifolds are identities. 

An \textbf{isomorphism} between $2$-bordisms $\alpha, \beta: A \rightarrow B$ is an isomorphism of underlying haloed manifolds respecting all of the above structural identifications. 
\end{defin}

\begin{rem}
The compatibility conditions might seem slightly mysterious at first, but they stem from the fact that the vertical part of the boundary of bordisms is modeled on either $A \times \{ 0 \} \times I$ or $B \times \{ 0 \} \times I$, in particular it has its own halo which has nothing to do with the intrinsic haloes of $A, B$. The end result is that the top or bottom boundary of $A \times \{ 0 \} \times I$ is diffeomorphic to $A$ as a manifold, but the relevant halations are isomorphic only non-canonically, and likewise for $B$. 

The needed isomorphisms of haloed manifolds are induced by the $2$-bordism $\alpha$, for definiteness let's focus on the domain $A$, but the same can be said about the codomain $B$. Using $\alpha$, one can identify the top boundary of $A \times \{ 0 \} \times I$ with the domain of the source $1$-bordism $w$, which is exactly $A$. Likewise, the bottom boundary of $A \times \{ 0 \} \times I$ can be identified with the domain of the target $1$-bordism $v$, which is $A$ again. 

These two isomorphisms of haloed manifolds could be in principle different and so we ask them to coincide, otherwise what we would get is a notion of a $2$-bordism where the vertical bordisms are not trivial, but are possibly "twisted" by some automorphism of $A$ and $B$. Such "twisted" $2$-bordisms are, as we will explain, used in the construction of the bordism bicategory, but are not allowed in $\mathbb{B}ord_{2}$ itself.
\end{rem}
One can show that using the notions presented it is possible to give well-defined vertical and horizontal composition on isomorphism classes of $2$-bordisms, it is then possible to construct $\mathbb{B}ord_{2}$ by hand as a bicategory. The monoidal structure, however, is more problematic. 

There are no issues in the underlying geometry, as disjoint union is an extremely well-behaved operation, but the structure of a symmetric monoidal bicategory is designed to be very weak in the categorical sense and so it contains large amounts of coherence data. To properly define such a bicategory one needs to give it all explicitly.

It is perhaps slightly surprising that describing a symmetric monoidal structure induced by disjoint union could be complicated, as disjoint union of manifolds is associative and unital up to a coherent system of diffeomorphisms. However, observe that it is impossible to directly express these properties as a symmetric monoidal structure on $\mathbb{B}ord_{2}$. An example of the latter will always formally consist of multiplication on objects which is associative and commutative up to an invertible $1$-bordism, where invertibility is not a condition, but another structure, namely a family of $2$-bordisms satisfying some equations on their own.

The route taken in \cite{chrisphd} is to construct the symmetric monoidal bicategory $\mathbb{B}ord_{2}$ from a so called \emph{symmetric monoidal pseudo double category}. The latter is a structure not much more complicated than a pair of two symmetric monoidal categories related by a triple of functors. This method of construction of symmetric monoidal bicategories was originally devised by Mike Shulman, see \cite{shulman_constructing_symmetric_monoidal_bicategories}. 

In the case of the bordism bicategory, the symmetric monoidal categories used to construct it are essentially the category $(\mathbb{B}ord_{2}) _{(0)}$ of haloed $0$-manifolds and isomorphisms and the category $(\mathbb{B}ord_{2}) _{(1)}$ of haloed $1$-bordisms and isomorphisms classes of haloed $2$-bordisms between them. Here we allow 2-bordisms "twisted" by some isomorphism of haloed $0$-manifolds, as (co)domains of source and target $1$-bordisms will not in general coincide.

The two symmetric monoidal homomorphisms $S, T: (\mathbb{B}ord_{2}) _{(1)} \rightarrow (\mathbb{B}ord_{2}) _{(0)}$ are induced by taking the source and target of a $1$-bordism, the functor $U: (\mathbb{B}ord_{2}) _{(0)} \rightarrow (\mathbb{B}ord_{2})_{(1)}$ constructs identity bordisms. Construction of a symmetric monoidal bicategory out of this data leads to the following theorem.

\begin{thm}[Schommer-Pries]
There exists a symmetric monoidal bicategory $\mathbb{B}ord_{2}$ with objects given by haloed $0$-manifold, morphisms given by haloed $1$-bordisms and $2$-cells given by isomorphisms classes of $2$-bordisms. Composition is given by uniquely defined gluing along boundaries and the symmetric monoidal structure is induced by disjoint union of manifolds. 
\end{thm}

\begin{proof}
This is \cite[Theorem 3.39]{chrisphd}.
\end{proof}

\begin{rem}
There is an alternative construction of the bordism bicategory which does not use haloes. Instead, one defines composition of $1$-bordisms by choosing for each possible gluing \emph{some} smooth structure on it and then choosing \emph{some} associativity diffeomorphisms for each triple of composable $1$-bordisms. These diffeomorphisms are not canonical and so in general they will fail to satisfy the pentagon axiom, but the isomorphism class of an invertible $2$-bordism constructed from them \emph{does} turn out to be canonical. This data can then be assembled into associativity isomorphisms for composition.

We do not follow that way of constructing the bordism bicategory, as haloes are still very conveniant when one wants to discuss manifolds with additional structure, and it is the \emph{framed} bordism bicategory which will be of interest to us. 
\end{rem}

\begin{rem}
It is not necessary to start at points. With the obvious shift in dimension, the same definitions and methods allow one to define a bordism bicategory with objects haloed $(d-2)$-dimensional manifolds, morphisms haloed $(d-1)$-dimensional bordisms and $2$-cells given by isomorphism classes of $d$-dimensional $2$-bordisms between these. 
\end{rem}

Because of the use of haloes, it is not difficult to adapt the definitons to manifolds with additional structure, like orientations, $G$-structures or framings. In essence, one attaches a given structure not to a manifold itself, but to its top-dimensional halation. The structure isomorphisms of $2$-bordisms and $1$-bordisms which identify their boundary with domain and codomain are then required to respect that structure. Care must be taken, however, to ensure that we do end up with a bicategory in the end. 

Observe that when we glue $2$-bordisms horizontally, the object part of the boundary of gluing is most naturally modeled on a product of the source/target $0$-manifold with the "glued intervals" $I \cup _{pt } I$, and not a product of $0$-manifold with the interval itself. Hence, to define the composite $2$-bordism, we choose some diffeomorphism $I \cup _{pt} I \simeq I$ respecting end-points and use it to define structural isomorphisms of the composite $2$-bordism. Since we only work with isomorphisms classes of $2$-bordisms, the result does not depend on the choice and this prescription gives composition which is both associative and unital.

In the case of bordisms with structure, it might not even be true in the first place that there is an isomorphism $I \cup _{pt} I \simeq I$ respecting that additional structure. Intuitively, to define bordism bicategories, one can only allow structures which are "topological" or "feeble enough". A non-example would be the structure of a metric, where one can easily define some numerical invariants of $2$-bordisms, like volume, and show that it is not possible to have a unital composition. For a possible axiomatization of a "topological structure", giving precise properties one has to assume, see \cite[Definition 3.42]{chrisphd}.

\subsection{Obtaining presentations of bordism bicategories}

In this section we will informally discuss methods developed by Schommer-Pries to obtain presentations of the bordism bicategory $\mathbb{B}ord_{2}$ and its oriented variant $\mathbb{B}ord_{2}^{or}$. It is in some sense a warm-up to the later chapters, where we will extend these results to the framed case. 

We will use the theory of freely generated symmetric monoidal bicategories, for precise definitions see \cite[Section 2.7]{chrisphd}. As an emergency reference, one can also use \textbf{Appendix \ref{appendix_free_monoidal_bicategories}}, where we present a formally analogous theory for non-symmetric monoidal bicategories. We have used both notions extensively in the categorical part of current work.

A \emph{presentation} of a symmetric monoidal bicategory $\mathbb{M}$ consists first and foremost of a \emph{generating datum} for a symmetric monoidal bicategory $(G, \mathcal{R})$.  Such a data is given by a triple of sets $(G_{0}, G_{1}, G_{2})$ of \emph{generating $i$-cells}, together with a set of \emph{relations} $\mathcal{R}$. The generating $1$-cells and $2$-cells are allowed to have domains and codomains which are not required to be generating cells themselves, they can only be their consequences, the relations are only allowed to live at the level of $2$-cells. 

Given a generating datum, one can construct a free symmetric monoidal bicategory $\mathbb{F}(G, \mathcal{R})$ generated by it. Additionally, as part of the presentation, we require a strict homomorphism $\phi: \mathbb{F}(G, \mathcal{R}) \rightarrow \mathbb{M}$, which is an equivalence of symmetric monoidal bicategories. It is easy to specify $\phi$ using the universal property of a freely generated symmetric monoidal bicategory, which identifies strict homomorphisms out of it with certain diagrams, which we call \emph{$(G, \mathcal{R})$-shapes with values in $\mathbb{M}$}. Such a shape consists of an association of an $i$-cell $\phi_{x}$ of $\mathbb{M}$ to each generating $i$-cell in $x \in G_{i}$ satisfying some compatibility conditions.

To say that the homomorphism $\phi$ is an equivalence is to say that it satisfies the following three conditions. 

\begin{enumerate}
\item It must be essentialy surjective on objects, that is, surjective on their equivalence classes. This conditions just says that the tensor products of objects $\phi_{X}$ for $X \in G_{0}$ should run through all equivalence classes of objects in $\mathbb{M}$. 
\item It should be essentially full on morphisms, ie. locally essentially surjective. This on the other hand says that $G_{1}$ is "big enough" so that its elements generate locally - that is, in $\mathbb{M}(\phi(A), \phi(B))$ for all pairs $A, B \in \mathbb{F}(G, \mathcal{R})$ - all morphisms of $\mathbb{M}$, considered up to invertible $2$-cells.
\item It should be locally bijective on $2$-cells. This at the same time says that we have enough generating $2$-cells to locally generate all $2$-cells in $\mathbb{M}$ and also that we have $\alpha \sim _{c(\mathcal{R})} \beta$ if and only if $\phi(\alpha) = \phi(\beta)$, where $c(\mathcal{R})$ is the closure of the relation $\mathcal{R}$. 
\end{enumerate}

Of course, all of these conditions can be satisfied trivially by making everything "as large as possible", or more precisely, by taking $G_{i}$ to consists of all $i$-cells in $\mathbb{M}$ and by taking the property from $(3)$ as the definition of $\mathcal{R}$. Proceeding this way we would obtain a cofibrant replacement of $\mathbb{M}$, whose existence is theoretically useful, but to write down a homomorphism from such a replacement is not easier at all than giving one defined on $\mathbb{M}$ to begin with. Usually, one is interested in possibly small presentations. 

In practice, when $\mathbb{M}$ is some bordism bicategory, it is condition $(3)$ which is difficult, and $(1)$ and $(2)$ are easy. This stems partly from the fact that the first two conditions only ask for some sort of surjectivity, without asking for injectivity, and partly from the fact that in bordism bicategories objects and morphisms correspond to lower dimensional manifolds, which are often much easier.  

\begin{figure}[htbp!]
	\begin{tikzpicture}[thick, scale=0.6, every node/.style={scale=0.6}]
	\begin{scope}[xshift= 7cm]
	
	\node[style={scale=1.6}] at (-2.6, 0) {$ pt_{+} := $};
	\node [circle, minimum width = 3cm] (A) at (0,0) {};
	\node (B) at (A.0) {};
	\node (C) at (A.180) {};

	\fill [halo_green!30]  (C.center) arc (180 : 360 : 1.5cm) to [out = 180, in = 0] (A.center) to [out = 180, in = 0] (C.center); 
	\draw [very thick, black] (A.center) to [out = 0, in = 180] (B.center);
	\draw [ultra thick, dotted, black] (C.center) to [out = 0, in = 180] (A.center);
	\draw [dashed] (A.center) circle (1.5cm);
	\node [circle, fill=black,inner sep=1.5pt] at (A.center) {};
	
	\foreach \x / \y in {0.5/0.7,-0.85/0.7, 0.5/-0.6, -0.85/-0.6}
		{
		\begin{scope}[xshift=\x cm, yshift=\y cm]
			\draw[->] (0.2,-0.3) arc (-90 : 180 : 0.2);
		\end{scope}
		}
	\end{scope}

	\node[style={scale=1.6}]  at (-2.6, 0) {$ pt_{-} := $};
	\node [circle, minimum width = 3cm] (A) at (0,0) {};
	\node (B) at (A.0) {};
	\node (C) at (A.180) {};

	\fill [halo_green!30]  (C.center) arc (180 : 360 : 1.5cm) to [out = 180, in = 0] (A.center) to [out = 180, in = 0] (C.center);
	\draw [very thick, black] (A.center) to [out = 0, in = 180] (B.center);
	\draw [ultra thick, dotted, black] (C.center) to [out = 0, in = 180] (A.center);
	\draw [dashed] (A.center) circle (1.5cm);
	\node [circle, fill=black,inner sep=1.5pt] at (A.center) {};
	
	\foreach \x / \y in {0.8/0.7,-0.4/0.7, 0.8/-0.6, -0.4/-0.6}
		{
		\begin{scope}[xshift=\x cm, yshift=\y cm]
			\draw[->] (-0.2,0.1) arc (90 : -180 : 0.2);
		\end{scope}
		}
		
	\end{tikzpicture}
\caption{Standard positive and negative points}
\label{fig:standard_oriented_positive_and_negative_points}
\end{figure}
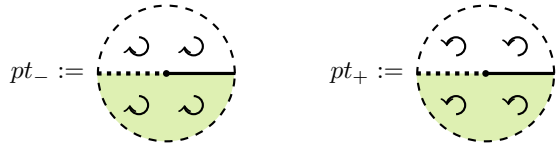

For example, any object of $\mathbb{B}ord_{2}$ is a $0$-manifold and so up to equivalence any object can be written as an interated tensor product of the standard point $pt$ with itself. The situation in the oriented bordism bicategory $\mathbb{B}ord_{2}^{or}$ is slightly more interesting, as objects are again $0$-manifolds, but equipped with an orientation on their $2$-dimensional halation. Two possible orientations allow one to distinguish between positive and negative points, one then sees that under tensor products all objects are generated by the two given below in \textbf{Figure \ref{fig:standard_oriented_positive_and_negative_points}}.

Finding the generating morphisms is not difficult, either. We have given a similar argument in the proof of \textbf{Theorem \ref{thm:presentation_of_the_framed_bordism_category}}, but for completeness, let us repeat it. 

If $w: A \rightarrow B$ is a $1$-bordism, then there exists a Morse function $h$ from the underlying manifold of $w$ into $[0, 1]$ such that the preimage $h^{-1}(0)$ is exactly $A$ and similarly $h^{-1}(1) = B$. Using such a Morse function, one can present $w$ as a composition of elementary bordisms, one for each critical point of $h$. Since we're in dimension $1$, there are only two kinds of handles, there is also no problem with the way they are glued, as these handles are attached along a $0$-manifold and thus the attaching data is essentially combinatorial. 

\begin{figure}[htbp!]
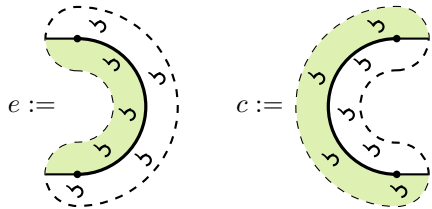

	
\caption{Generating 1-bordisms of the oriented bordism bicategory}
\label{fig:generator_1_bordisms_of_the_oriented_bordism_bicategory}
\end{figure}

The end result is that any $1$-bordism can be presented as a composition of left and right elbows, which we draw in the oriented case in \textbf{Figure \ref{fig:generator_1_bordisms_of_the_oriented_bordism_bicategory}}. In the unoriented case they look the same, except there is no orientation. Note that their domains and codomains are respectively $e: pt_{+} \otimes pt_{-} \rightarrow \emptyset$ and $c: \emptyset \rightarrow pt_{-} \otimes pt_{+}$, to verify that it is important to pay attention to coorientations of haloes. 

The real work needs to be done at the level of $2$-cells, which in the case of $\mathbb{B}ord_{2}$ and $\mathbb{B}ord_{2}^{or}$ correspond to diffeomorphism classes of respectively unoriented and oriented $2$-bordisms. Since at this level one seeks local bijectivity, one is essentially tasked with a form of a classification problem for surfaces. 

Classification of oriented surfaces dates back more than a hundred years, but here we are interested not so much in finding all classes of surfaces, but more so in finding the \emph{pieces} from which any surfaces can be glued, which will give generating set of $2$-cells $G_{2}$, and the \emph{relations between them}, which will give $\mathcal{R}$. 

This is a problem of a different sort and heavily depends on the categorical framework one works in. For example, it is true that any surface can be obtained by gluing disks along boundaries, as one sees by choosing a triangulation, but it is not \emph{true} that any $2$-bordism can be glued from a square $2$-bordism just by disjoint unions and horizontal and vertical composition. 

One can show that ordinary bordism categories are generated by elementary bordisms by mapping any given bordism into $\mathbb{R}$. If the map is sufficiently generic, the preimages of small open disks will give the desired decomposition into simple pieces. In the case of bordism bicategories, one is essentially allowed compositions in two different directions and so the natural idea is to generically map into $\mathbb{R}^{2}$. This is the program followed in \cite{chrisphd}, which we will now try to briefly explain. 

\begin{figure}[htbp!]
\begin{tikzpicture}[thick, yscale=0.8]
	\draw [very thick, red, rounded corners] (2.5, 4) to [out = -90, in = 45] (0.5, 2) to [out = 225, in = 90]  (0,1) parabola bend (1.5, 0) (3,1) parabola bend (3.5, 2) (4,1);
	\draw [very thick, red, rounded corners] (2.5, 4) to [out = -90, in = 135] (4.5, 2) to [out = -45, in = 90] (5,1) parabola bend (3.5, 0) (2, 1) parabola bend (1.5, 2) (1,1);
	\draw [very thick, red, rounded corners] (1,1) to [out = 90, in = -45] (0.5, 2) to [out = 135, in = -90] (0, 3) parabola bend (2.5, 5) (5, 3) to [out = -90, in = 45] (4.5, 2) to [out = 225, in = 90] (4, 1);
\end{tikzpicture}
\caption{A possible graphic of $\mathbb{RP}^{2}$}
\label{fig:a_possible_graphic_of_the_real_projective_plane}
\end{figure}
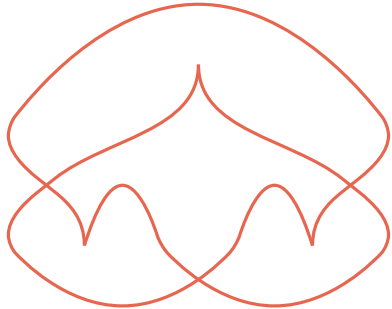

For simplicity, let us first focus on the case of closed surfaces, so let $W$ be one. Given a generic map $f: W \rightarrow \mathbb{R}^{2}$, where here "generic" means a multitude of transversality conditions with respect to a  carefully chosen collection of submanifolds of jet and multijet spaces, the image $\Psi = f(C)$ of the set of critical points will consist of a finite number of closed embedded curves in $\mathbb{R}^{2}$ and a finite number of singular points. We call the image the \textbf{graphic of critical values} of the map $f$, in \textbf{Figure \ref{fig:a_possible_graphic_of_the_real_projective_plane}} we give one possible graphic of the real projective plane. 

We call the embedded curves the \emph{arcs} of the graphic, the projection $\mathbb{R}^{2} \rightarrow \mathbb{R}$ onto a second coordinate restricts to a diffeomorphism on each of them. The idea is that one can recover the surface from its graphic together with some amount of combinatorial data, we will now describe this process. 

\begin{reading_graphic*}
The behaviour of the map $f$ in the preimage of small neighbourhood of any given point in $x \in \mathbb{R}^{2}$ can be completely classified in terms of the graphic and a small amount of additional \emph{labelling data}. There are essentially three possible cases. 
 
\textsl{(A point not lying on the graphic)} - If a point $x \in \mathbb{R}^{2}$ does not lie on the graphic, then the map $f$ induces an isomorphism on tangent spaces at all its preimages. It follows that there is just a finite number of them and that each of them has a small neighbourhood that is taken diffeomorphically to a neighbourhood of $x$. 

Hence, locally around $x$, the map $f$ has a form of a trivial covering. Thus, the labelling data needed to recover $f$ is given by the \textbf{set of sheets} of that trivial covering. 

\textsl{(A point lying on an arc of the graphic)} - If $x$ lies on an arc of $\Psi$, then it has precisely one preimage which is critical and possibly many where the map $f$ restricts to a local diffeomorphism. In the neighbourhood of the critical preimage, the map $f$ takes the form of a projection from a fold surface, pictured below in \textbf{Figure \ref{behaviour_of_a_generic_map_around_a_fold_arc}}. 

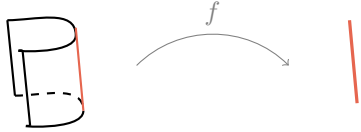
\begin{figure}[htbp!]
	\begin{tikzpicture}
		\draw[thick] (-0.5, 0) to [out=0, in=90] (0.4, -0.2) to [out=-90, in=180] (-0.3, -0.4);
		\draw[thick, dashed] (-0.4, -1) to [out=0, in=90] (0.5, -1.2);
		\draw[thick] (0.5, -1.2) to [out=-90, in=180] (-0.2, -1.4);
		\draw[thick] (-0.5, 0) to (-0.4, -1);
		\draw[thick, red] (0.4, -0.1) to (0.5, -1.2);
		\draw[thick] (-0.3,-0.4) to (-0.2, -1.4);
		
		\draw [->, gray] (1.2, -0.6) [out=45, in=135] to node[auto] {$ f $} (3.2, -0.6);
		
		\draw[very thick, red] (4, 0) to (4.1, -1.1);
		
	\end{tikzpicture}
\caption{Behaviour of a generic map in the preimage of an arc of the graphic}
\label{behaviour_of_a_generic_map_around_a_fold_arc}
\end{figure}
The arcs come in two types, as the fold surface can be either a cyllinder on a right or a left elbow. The needed labelling data consists of a distinguished pair of sheets over either the left or right side of the arc, these are the sheets "folded together", and a bijection between the rest of them. 

\textsl{(A codimension 2 singularity of the graphic)} - If $x$ is one of the singular points of the graphic, then its preimage consists of a finite number of points, at most two of which are critical. There is a finite number of possibilites of the behaviour of the map around such critical points. They are all given in \textbf{Figure \ref{fig:behaviour_of_a_generic_map_around_codimension_2_singularity}}, the map $f$ should be understood as a projection onto the plane from the given surface. The set of critical points is drawn using red colour, its projection describes the behaviour of the graphic around this type of singularity.

\begin{figure}[htbp!]
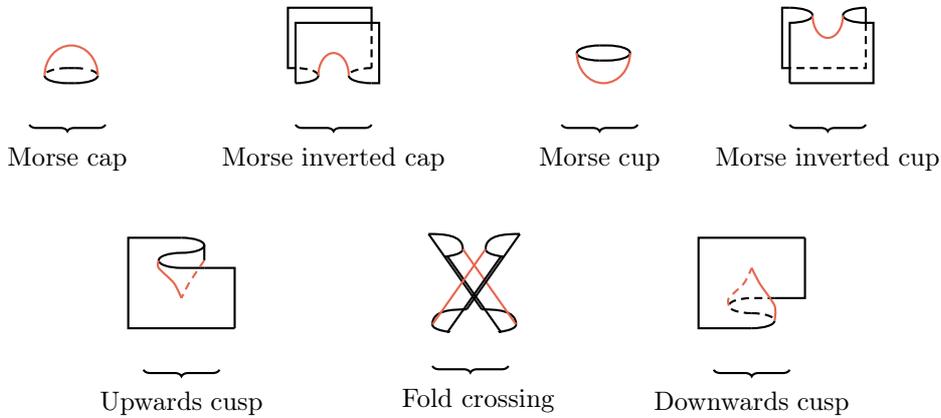


\caption{Behaviour of a generic map in the preimage of a codimension $2$ singularity of the graphic}
\label{fig:behaviour_of_a_generic_map_around_codimension_2_singularity}
\end{figure}
The needed labelling data consists of a singularity type together with some compatibility conditions between bijections being part of data over the arcs of the graphic leaving the given singularity. 
\end{reading_graphic*}

One can use these methods to decompose any surface into the promised simple pieces. Indeed, to do so, choose a generic map $f: W \rightarrow \mathbb{R}^{2}$. If one covers $\mathbb{R}^{2}$ with a suitably fine and well-behaved covering $(U_{i})$, for example by small squares, then the preimages of these squares will be disjoint unions of a finite number of $2$-bordisms. They will either be squares themselves, fold surfaces -  which one can interpret as identity $2$-cells of the left and right elbows, or one of the finite types of surfaces corresponding to singularities of the graphic, all glued together nicely along their boundaries. This shows that $2$-bordisms from \textbf{Figure \ref{fig:behaviour_of_a_generic_map_around_codimension_2_singularity}} are a set of generators for $\mathbb{B}ord_{2}$. 

This decomposition of $W$ depends heavily on the choice of the generic function $f$ and so to understand how different decompositions are related one studies generic functions $W \times I \rightarrow \mathbb{R}^{2} \times I$, which can be understood as "generalized paths" of maps into the plane. Eventually, it is possible to derive a finite list of relations between the generating $2$-bordisms. 

Some adaptations need to be made when $W$ is not a closed surface, but a general $2$-bordism, instead of mapping into $\mathbb{R}^{2}$, one instead maps generically into $I^{2}$, with prescribed behaviour at the boundary. However, these are minor changes and the methods developed give a good grasp of behaviour of $2$-bordisms. 

These Morse-theoretic results allow one to write down an intuitive presentation, however, it takes a little bit more work to show that the induced map $\phi: \mathbb{F}(G, \mathcal{R}) \rightarrow \mathbb{B}ord_{2}^{or}$ is indeed an equivalence. This stems from the fact that to do so one also needs to establish some control over the $2$-cells in $\mathbb{F}(G, \mathcal{R})$. In principle this should not be difficult, as it is a freely generated symmetric monoidal bicategory and so admits an explicit description, but the sheer amount of coherence data in this type of structure is overwhelming and obscures the view. 

Hence, the actual way of showing that $\phi$ is an equivalence is more indirect. In a sense, one pushes $\mathbb{F}(G, \mathcal{R})$ to be more geometric and $\mathbb{B}ord_{2}$ to be more algebraic, constructing a zig-zag of equivalences between them. 

To do the latter, a notion of a \textbf{planar diagram} is introduced, which codifies data needed to recover a surface over $\mathbb{R}^{2}$ in a combinatorial way. The definition is slightly technical, but the three main components of a planar diagram are as follows.

\begin{enumerate}
\item a \textbf{graphic} $\Psi$, which is a type of diagram embedded in the square which behaves like a set of critical points of a generic map from a $2$-bordism into $I^{2}$
\item a \textbf{chambering graph} $\Gamma$, which is a particular type of a graph embedded in $I^{2}$
\end{enumerate}
The purpose of the chambering graph is to divide $I^{2}$ into smaller regions. One calls the connected components of the complement of the graphic and the chambering graph in $I^{2}$ the \textbf{chambers} of the diagram, they play a role of a "nice open covering". 

\begin{enumerate}
\setcounter{enumi}{2}
\item \textbf{labelling data} attached to each chamber, to each arc and singularity of $\Psi$ and to each edge of $\Gamma$, satisfying a number of compatibility conditions
\end{enumerate}
Morse-theoretic methods presented can be used to design a combinatorial equivalence relation on planar diagrams which has the property that $2$-bordisms represented by diagrams some $P_{1}, P_{2}$ are diffeomorphic if and only if the diagrams themselves are equivalent. The promised "algebraicization" of $\mathbb{B}ord_{2}$ proceeds by construction of a symmetric monoidal bicategory $\mathbb{B}^{PD}$, where $2$-cells are given by equivalences of planar diagrams, one can then relate it to $\mathbb{B}ord_{2}$ via a span of equivalences. 

To establish control over $\mathbb{F}(G, \mathcal{R})$, one introduces a notion of unbiased semistrict symmetric monoidal bicategory. The latter is a partially strict variant of a symmetric monoidal bicategory, which is easier to work with in practice, similarly to how a \textbf{Gray}-monoid is a partially strict variant of a monoidal bicategory. It is also slightly more geometric, admitting a variant of string calculus.

Using the same generating datum, one can form the \emph{free unbiased semistrict symmetric monoidal bicategory} $\mathbb{U}(G, \mathcal{R})$. It is a particular type of a symmetric monoidal bicategory and so the universal property of $\mathbb{F}(G, \mathcal{R})$ induces a strict homomorphism $\mathbb{F}(G, \mathcal{R}) \rightarrow \mathbb{U}(G, \mathcal{R})$, which was shown by Schommer-Pries to be an equivalence for \emph{any} generating datum. This is a difficult result in coherence theory. 

Since unbiased semistrict symmetric monoidal bicategories admit string diagram calculus, $2$-cells in $\mathbb{U}(G, \mathcal{R})$ can be represented by equivalence classes of such string diagrams, with labellings coming from the generating sets $G_{i}$ and a number of additional relations coming from $\mathcal{R}$. The presentation is chosen in such a way that these diagrams closely resemble planar diagrams representing $2$-bordisms, one then shows that $\mathbb{B}^{PD}$ and $\mathbb{U}(G, \mathcal{R})$ are not only equivalent, but in fact isomorphic. This implies that $\mathbb{F}(G, \mathcal{R})$ and $\mathbb{B}ord_{2}$ are equivalent too, and a quick run through the relevant diagrams shows that it is the homomorphism $\phi$ which establishes this equivalence.

To adapt this proof to the case of bordism bicategories of manifolds with additional structure, one has to devise a correct notion of an "enriched planar diagram". Such a diagram would carry extra data with the property that using it one can recover not only a $2$-bordism, but also the additional structure over it, one also has to find an appropriate equivalence relation on enriched diagrams. 

If all of this is done in sufficiently combinatorial fashion, one will be able to use these "enriched planar diagrams" to compare the given bordism bicategory with some freely generated unbiased semistrict symmetric monoidal bicateory, establishing a presentation of the former. It is explained in \cite{chrisphd} how to do in the case of $\mathbb{B}ord_{2}^{or}$ through what we will call \emph{oriented planar diagrams}. These are very similar to their unoriented counterparts, differing only in labelling data, where to each sheet we additionally attach a plus or a minus sign, denoting its orientation relative to the square.

We will be later concerned with finding the correct notion of a \emph{framed planar diagram}, modelling diffeomorphism-isotopy classes of framed surfaces, as our goal is to find a presentation of the framed bordism bicategory $\mathbb{B}ord_{2}^{fr}$. 

\section{Framed bordism bicategory}

In this chapter we will construct the framed bordism bicategory $\mathbb{B}ord_{2}^{fr}$ and perform a few simple calculations with framings on generating $1$-bordisms and $2$-bordisms from the oriented presentation. We will try to highlight a few issues related to framings which are perhaps less visible in the language of framed planar diagrams, which we will introduce later. 

\begin{defin}
Let $W$ be a two-dimensional manifold. A \textbf{framing} $v$ on $W$ is a trivialization of its tangent bundle. In other words, it is a linear isomorphism $T_{W, w} \simeq \mathbb{R}^{2}$ which varies smoothly with $w \in W$ or, more compactly, a section of the $GL(2, \mathbb{R})$-principal bundle $Iso(T_{W}, \mathbb{R}^{2})$. An \textbf{isotopy} $s: v \rightarrow w$ between framings $v, w$ is a homotopy of such sections. 
\end{defin}

To construct the framed bordism bicategory $\mathbb{B}ord_{2}^{fr}$ as a symmetric monoidal bicategory, we will adapt the methods used to construct the unoriented bordism, essentially by carefully replacing all surfaces without boundary that appear in the definition with framed surfaces. Most likely, one of the ways to do so would be to use the language of topological $\mathcal{F}$-structures introduced in \cite{chrisphd}. We do not follow this route, as in our case it would introduce some slight complications, and instead build up our bicategory in more simple terms. 

We start by discussing framings on halations, but before doing so, let us first make some conventions that will hold till the end of this chapter. We assume that all haloes are implicitly two-dimensional, in particular, all haloed manifolds we consider are themselves of dimension $k \leq 2$. We also agree that by the word \emph{surface} we will always mean a two-manifold without boundary, possibly non-compact. 

By a \textbf{pro-surface} we mean a pro-object in the category of surfaces and their embeddings, it is a particular kind of a pro-manifold. Observe that when defining a two-dimensional halation, we may as well assume that it is a pro-surface in this sense, as the halations $X \hookrightarrow \hat{Y}$ associated to an embedding $X \hookrightarrow Y$ into an interior of the surface are always of this form. Likewise, we will now call two halations isomorphic if they are isomorphic as pro-surfaces, that is, admit an isomorphism whose representatives are embeddings. The advantage of this slightly more restrictive notion is that framings can be easily extended to this case.

\begin{defin}
Let $X: I \rightarrow \mathcal{M}an^{2}$ be a pro-surface. A \textbf{framing} on $X$ is a framing on one of the $X_{i}$, two such framings $v_{i}$ on $X_{i}$ and $v_{j}$ on $X_{j}$ are considered to be the same is some index $k \leq i, j$ such that $v_{i} |_{X_{k}} = v_{j} |_{X_{k}}$. 

More concisely, the space of framings over $X$ is precisely $colim \ Fr(X_{i})$, where $Fr(X_{i})$ is the space of framings over $X_{i}$. An \textbf{isotopy} of framings on $X$ is a path in this space.
\end{defin}
Observe that when $f: \hat{Z} \rightarrow \hat{X}$ is a map of pro-surfaces which admits a system of representatives given by local diffeomorphisms, for example when it is an embedding, then any framing $v$ on $\hat{X}$ induces one a framing $f^{*}v$ on $\hat{Z}$ just by pulling it back any representative. If $\hat{Z}, \hat{X}$ are both framed by $v, w$ respectively, then we will say that the map $f$ \textbf{preserves the framing} if $f^{*}v = w$. 

In particular, if $X \subseteq \hat{X}$ is a framed halation and $Z \subseteq X$ is a submanifold, then the induced halation $Z \subseteq \hat{Z}$ admits a unique framing such that the embedding $\hat{Z} \subseteq \hat{X}$ preserves framings on both sides. 

The notion of a framing on a pro-surface then leads to a notion of a framed halo. 

\begin{defin}
Let $X$ be a manifold. A \textbf{framed halo} $\hat{X}$ over $X$ is a halo equipped with a framing on the top-dimensional halation $\hat{X}_{2}$. 
\end{defin}
We will construct the framed bordism bicategory by replacing haloes throughout all the relevant definitions by framed haloes, we will then need to quotient out the $2$-bordisms by a suitable notion of isotopy of framings. In particular, the objects of $\mathbb{B}ord_{2}^{fr}$ will be framed haloed $0$-manifolds. Before going further, let's gain a little bit of visual intution by giving a few pictures.

We draw haloes themselves the same way we did before, by drawing small neighbourhoods which represent them, coorientations are denoted using shadings. An isomorphism of a finite-dimensional vector space with $\mathbb{R}^{n}$ is the same thing as a choice of an ordered basis and so framings can be understood - and are often defined - as smoothly varying ordered bases of the relevant tangent spaces. We will draw framings as sparse collections of such ordered bases on different points of the small $2$-dimensional neighbourhood representing the top-dimensional halation, a reader should visualize the obvious continous extension. The leading vector will be distinguished by an arrow. 

\begin{figure}[htbp!]
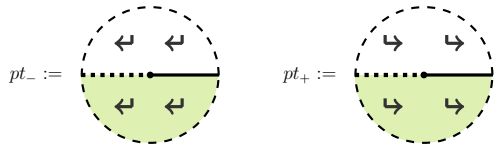


\caption{Standard positively and negatively framed points}
\label{fig:standard_positively_and_negatively_framed_points}
\end{figure}

In \textbf{Figure \ref{fig:standard_positively_and_negatively_framed_points}} we give an example of two framed haloed $0$-manifolds, which are framed analogues of the positive and negative points from the oriented bordism bicategory. These are the \textbf{standard positively and negatively framed points}, we denote them by $pt_{+}, pt_{-}$. 

The haloes can be a little inconveniant when one tries to focus on the underlying manifolds, but as we will see in the case of framings one can at least partially ignore them. We first state a basic property of the space of framings.

\begin{prop}
Let $(M, v)$ be a framed 2-manifold and $A \subseteq M$ a possibly empty submanifold. The space of framings $Fr(M, v |_{A})$ that pointwise agree with $v$ at all points of $A$ is a canonically homeomorphic to the space $C^{\infty}((M, A), (GL(\mathbb{R}, 2), id))$ of all smooth maps that take $A$ to the identity.
\end{prop}

\begin{proof}
This is true for any principal $G$-bundle, but for completeness let us give the argument. The space of maps acts on the space of framings $Fr(M)$ via change of framing and between any two framings there is a unique smooth map that takes one to the other. The canonical homeomorphism is given by evaluating this action at $v$. 
\end{proof}
Note that because of smooth approximation, up to homotopy we can replace the space of smooth maps by the space of all continous maps as the inclusion
\begin{gather*}
C^{\infty}((M, A), (GL(\mathbb{R}, 2), id)) \hookrightarrow Map((M, A), (GL(\mathbb{R}, 2), id))
\end{gather*}
into the latter is a homotopy equivalence. The latter space then admits an even smaller description as the space of base-point preserving maps $Map_{*}(M / A, O(2))$, as we also have a homotopy equivalence $GL(2, \mathbb{R}) \simeq O(2)$ given by Gram-Schmidt orthogonalization. This gives the following corollary. 

\begin{cor}
Let $M$ be a frameable surface, let $A \subseteq M$ be a submanifold such that the inclusion $A \hookrightarrow M$ is a homotopy equivalence. Then, the restriction map $Fr(M) \rightarrow Fr(A)$ is a homotopy equivalence. 

Consequently, if $A$ is a framed haloed manifold, then up to a unique homotopy class of isotopies the framing on the halo is determined by its values at the points of $A$. 
\end{cor}
Note that the framings considered are always $2$-dimensional, so that by $Fr(A)$ we mean the space of sections of isomorphisms of $T_{M} |_{A}$ with $\mathbb{R}^{2}$.
\begin{proof}
The first part follows directly from our up to homotopy description of the space of framings as the space of maps into $O(2)$. 

To see the second part, observe that the two-dimensional halation $\hat{A} _{2}$ can be represented by a pro-surface $(A_{i})$ indexed on natural numbers, such that each $A_{i+1} \subseteq A_{i}$ is a homotopy equivalence, each $A \subseteq A_{i}$ is a homotopy equivalence and additionally $\bigcap A_{i} = A$.  For example, this can be done by identifying the halation with one induced by the normal bundle of $A \subseteq \hat{A}_{2}$, and taking $A_{i}$ to be disk bundles of radii $\frac{1}{i}$. 

The space of framings on the halo $\hat{A}$ can then be identified with $colim \ Fr(A_{i})$, which is homotopy equivalent to $Fr(A)$. Passing to fundamental groupoids we get the desired result.
\end{proof}
Using the corollary, we can easily show that any framed haloed $0$-manifold with underlying manifold a single point is isomorphic to either the standard positively framed or standard negatively framed point. This is intuitively obvious, but technically does require some understanding of what exactly is a difference between a framing defined on an "infinitesimal neighbourhood" and a framing defined only on the manifold itself.

\begin{prop}
Let $A$ be a framed haloed $0$-manifold whose underlying manifold is a single point. Then either we have an isomorphism $A \simeq pt_{+}$ of haloed manifolds or an isomorphism $A \simeq pt_{-}$, preserving the framings up to isotopy. 
\end{prop}

\begin{cor}
\label{cor:any_framed_0_manifold_isomorphic_to_tensor_product_of_standard_points_up_to_isotopy}
Under disjoint unions, the standard positively and negatively framed points $pt_{+}, pt_{-}$ generate all framed haloed $0$-manifolds up to isomorphism preserving the framings up to isotopy. 
\end{cor}

\begin{proof}
Observe that $A$ is non-canonically isomorphic to both of them, we just have to ensure that there exists an isomorphism which preserves the framings up to isotopy.

By our results above, the framing on the two-dimensional halation $\hat{A}_{2}$ is determined by its restriction to the manifold itself, which consists of a single point. We can identify the tangent space $T_{\hat{A}_{2}, A}$ with the normal space of $A$ inside $\hat{A}$, as $A$ is zero-dimensional. It follows that the $2$-dimensional tangent space $T_{\hat{A}, A}$ inherits an orientation from the coorientations attached to a halo. 

The space of framings on $A$ is the same as the space of ordered bases of $T_{\hat{A}_{2}, A}$, it has exactly two path-components, classified by orientations induced. If the orientation induced by the framing agrees with the one coming from the structure of a halo, then any isomorphism $A \simeq pt_{+}$ will preserve the framings up to isotopy, as the same is true for the positively framed point. Conversely, if the orientations don't agree, then any isomorphism $A \simeq pt_{-}$ will do. 
\end{proof}

Having established these basic facts about framings on halations, we proceed with our definition of the framed bordism bicategory. We start by defining $1$-bordisms, which is rather easy to do.

\begin{defin}
Let $A, B$ be framed haloed $0$-manifolds. A \textbf{framed 1-bordism} $w: A \rightarrow B$ consists of an ordinary $1$-bordism $w$ together with a framing on its top-dimensional halation. Additionally, we require the boundary identifications 
\begin{gather*}
f_{in}: (\partial _{in} w,  \hat{w}_{1} |_{\partial _{in}}, \hat{w} _{2} |_{\partial {in}}) \simeq (A, \hat{A}_{1}, \hat{A}_{2}) ,\\
f_{out}: (\partial _{out} w, \hat{w} _{1} |_{\partial _{out}}, \hat{w} _{2} |_{\partial {out}}) \simeq (B, \hat{B}_{1}, \hat{B}_{2})
\end{gather*}
that are part of the structure of the bordism to preserve framings.
\end{defin}
Observe that we require that the framings of $w$ coincide with the framing of $A, B$ not only at the $0$-manifolds itself, but also on the infinitesimal neighbourhood provided by the halo. Additionally, we require them to agree pointwise and not only up to some kind of isotopy. This strong form of identification would be a little more difficult to phrase if we tried to consistently use topological $\mathcal{F}$-structures, which is one of the reasons we decided not to pursue this approach. 

\begin{figure}[htbp!]
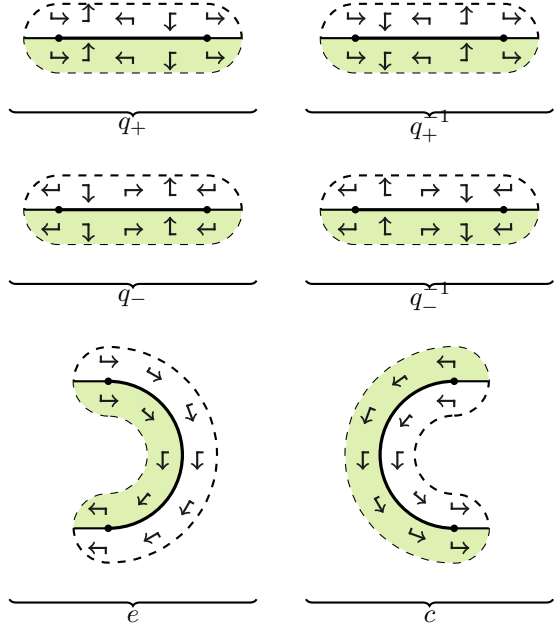

	
\caption{Generating 1-bordisms of the framed bordism bicategory}
\label{fig:generator_1-bordisms_of_the_framed_bordism_bicategory}
\end{figure}

We give examples of framed $1$-bordisms in \textbf{Figure \ref{fig:generator_1-bordisms_of_the_framed_bordism_bicategory}}. Their domains and codomains should be understood to be suitable tensor product of positively or negatively framed points, the structural identifications are either the obvious inclusions or the obvious inclusions composed with a reflection along the one-dimensional part of the halo. 

In precise terms, we have $e: pt_{+} \otimes pt_{-} \rightarrow \emptyset$, $c: \emptyset \rightarrow pt_{-} \otimes pt_{+}$ like in the oriented case, $q_{+}, q_{+}^{-1}$ are endomorphisms of the positively framed point and $q_{-}, q_{-}^{-1}$ are endomorphisms of the negatively framed point. We will see later that the latter four are in fact Serre autoequivalences or their pseudoinverses of their respective domains, it is interesting to observe that their underlying oriented bordisms are in fact isomorphic to the identity, which is not true for the framed ones.

We will now proceed with defining $2$-bordisms, where we need to be a little bit more careful. We will first define framings on $2$-bordisms, this will be particular kinds of framings defined on the halation satisfying some boundary conditions. The actual $2$-cells in $\mathbb{B}ord_{2}^{fr}$ will then be isomorphism-isotopy classes of framed $2$-bordisms, that is, we will consider two framed $2$-bordisms $\alpha, \beta$ to represent the same $2$-cell if they differ by an isomorphism such that under this isomorphism the framing of $\alpha$ is taken to a framing isotopic to the one given on $\beta$. 

\begin{defin}
Let $A, B$ be framed haloed $0$-manifolds and let $w, v: A \rightarrow B$ be framed $1$-bordisms between them. If $\alpha: w \rightarrow v$ is a $2$-bordism, then a \textbf{framing} on $\alpha$ consists of a framing on its halo satisfying the following compatibility conditions. 

\begin{enumerate}
\item Under the identification 
\[ (\partial _{m}, \hat{\alpha} |_{\partial _{m}}) \simeq (w, \hat{w}_{2}) \sqcup (v, \hat{v}_{2}) \]
the framing of $\alpha$ restricts to the framings of $w, v$ at all the points of the underlying $1$-manifolds.
\item Under the identification 
\[ (\partial _{o}, \hat{\alpha} _{2}) \simeq (A \times \{ 0 \} \times I, \widehat{A \times \mathbb{R}^{2}}) \sqcup (B \times \{ 0 \} \times I, \widehat{A \times \mathbb{R}^{2}}) \]
the framing of $\alpha$ restricts to the constant framing on $A \times \{ 0 \} \times I$ inherited from either its top or botton boundary and similarly for $B \times \{ 0 \} \times I$ at all the points of the underlying $1$-manifolds.
\end{enumerate}
An \textbf{isotopy} of framings on $\alpha$ consists of an isotopy of framings on its halo such that at all times of the isotopy the above conditions are satisfied. 
\end{defin}
There are a few delicate issues in our definition of a framing on a $2$-bordism. First, unlike in the case of framed $1$-bordisms, a framing on a $2$-bordism is required to satisfy conditions only at the points lying directly on the boundary, not in any infinitesimal neighbourhood. This is not problematic, as we will only consider framings on $2$-bordisms up to isotopy, and so this stronger condition can in fact always be enforced if one wished so.

We also have to make sense of the \emph{constant framing} on $A \times \{ 0 \} \times I$ and $B \times \{ 0 \} \times I$. For definiteness let us focus on the former, the other case is analogous. Observe that we have an isomorphism of haloed manifolds $(A, \hat{A} ^{1}, \hat{A}^{2}) \simeq (A \times \{ 0 \} \times \{ 1 \}, A _{\mathbb{R}}, A_{\mathbb{R}^{2}})$ induced by $\alpha$, this induces a framing on the top boundary of $A \times \{ 0 \} \times I$. In more detail, the top boundary $A \times \{ 0 \} \times I$ is identified with the domain of $w$ and this gives the desired isomorphism, one then sees that the induced framing is exactly the restricted framing of $w$. 

We can then extend this framing in a constant manner over all of $A \times \{ 0 \} \times I$, this is possible as the tangent space to the halo is at each point canonically isomorphic to $\mathbb{R}^{2}$. One could also extend the framing induced on $A \times \{ 0 \} \times \{ 0 \}$ from the domain of $v$, these two possibly  different framings on $A \times \{ 0 \} \times I$ actually coincide by the compatibility conditions on $2$-bordisms. We call this framing the \emph{constant framing} on $A \times I$, observe the subtle fact that this framing is only defined in the presence of $\alpha$. 

\begin{defin}
Let $w, v: A \rightarrow B$ be framed $1$-bordisms. Two framed $2$-bordisms $\alpha, \beta: w \rightarrow v$ are considered \textbf{isomorphic-isotopic} if there exists an isomorphism $\alpha \simeq \beta$ of $2$-bordisms, ie. an isomorphism of haloed manifolds commuting with structural identifications of the boundaries, such that under this isomorphism the framing of $\alpha$ is taken to a framing isotopic to the given framing on $\beta$. 
\end{defin}
Here we mean the isotopy of framings on $2$-bordisms in the sense we defined it above, in other words, it is required to be relative to the boundary. Observe that isomorphism-isotopy is clearly na equivalence relation. 

We are now ready to prove the existence of our main object of study in the later chapters.

\begin{thm}
There exists a symmetric monoidal bicategory, the \textbf{framed bordism bicategory} $\mathbb{B}ord_{2}^{fr}$, with objects framed haloed $0$-manifolds, morphisms framed $1$-bordisms and $2$-cells given by isomorphism-isotopy classes of framed $2$-bordisms. Horizontal and vertical composition are given by gluing along boundaries, the monoidal structure is induced from disjoint union of manifolds.
\end{thm}

\begin{proof}
Analogously to the unoriented case, we use the method of constructing symmetric monoidal bicategories from symmetric monoidal pseudo double categories, see \cite{shulman_constructing_symmetric_monoidal_bicategories}. The relevant symmetric monoidal categories are the category $(\mathbb{B}ord_{2}^{fr}) _{(0)}$ of framed haloed $0$-manifolds and isomorphisms that preserve the framing and the category $(\mathbb{B}ord_{2}^{fr}) _{(1)}$ of framed haloed $1$-bordisms and isomorphism-isotopy classes of framed $2$-bordisms between them. The former of these two categories is self-explanatory, so let's focus on the latter.

Similarly to the unoriented case, one has to allow "possibly twisted" framed $2$-bordisms between $1$-bordisms whose (co)domains do not necessarily coincide. More precisely, if $w: A \rightarrow B$, $v: A^\prime \rightarrow B^\prime$ are framed haloed $1$-bordisms, then a "possibly twisted" framed $2$-bordism between them is a tuple $(h_{in}, h_{out}, \alpha)$, where $h_{in}: A \simeq A^\prime$, $h_{out}: B \simeq B^\prime$ are isomorphisms and $\alpha: v \rightarrow w$ is a framed $2$-bordism, where we identify $w$ with a framed bordism of type $A \rightarrow B$ by composing its structural identifications with $h_{in}, h_{out}$.

Composition in $(\mathbb{B}ord_{2}^{fr})_{(1)}$ is induced by vertical gluing of $2$-bordisms. Observe that to have an induced \emph{smooth} framing on the gluing $\beta \circ \alpha$ one has to deform framings on both of $\alpha, \beta$ to agree with the framings on the boundary not only pointwise, but also is some infinitesimal neighbourhood. Since we work up to isotopy, this is always possible and the result is up to isotopy independant of a choice of such deformation. However, this new smooth framing does not in general make $\beta \circ \alpha$ into a framed $2$-bordism, it needs to be slightly modified in a way which we now describe.

Observe that the vertical boundary of $\beta \circ \alpha$ is modeled on $I \cup _{pt} I$, to make it into a $2$-bordism one composes with some rescaling $I \cup _{pt} I \simeq I$. However, the constant framing on $I \cup _{pt} I$ will in general not coincide with the constant framing on $I$, the difference between them is essentially the same as between the framings the two diffeomorphic manifolds $[0, 1] \times \{ 0 \}$ and $[0, 2] \times \{ 0 \}$ inherit from $\mathbb{R}^{2}$. 

Intuitively, the constant framing on $I \cup _{pt} I$ should differ from the one on $I$ "not by much", where we compare the two framings using some chosen rescaling $I \cup _{pt} I$. This is indeed the case. The two framings differ by tilting by a function $f: I \rightarrow \mathbb{R}^{+} \subseteq GL(2, \mathbb{R})$, where we identify $\mathbb{R}^{+}$ with the group of transformations $g_{r}$ of the form $g_{r}(e_{1}) = e_{1}, g_{r}(e_{2}) = r e_{2}$, where $e_{i}$ is the standard basis of $\mathbb{R}^{2}$ and $r \in \mathbb{R}^{+}$. The map $f$ in question is the derivative of the rescaling map $I \rightarrow I \cup _{pt} I \simeq [0, 2]$ composed with $x \mapsto \frac{1}{x}$. This is "not much" as the group $\mathbb{R}^{+}$ is contractible.

We define $\beta \circ \alpha$ as a framed $2$-bordism in the following way. One first chooses a rescaling $I \simeq I \cup _{pt} I$ that restricts to the identity near the boundary of $I$, observe that all such rescalings are isotopic. This makes $\beta \circ \alpha$ into an oriented $2$-bordism, it is, as we explained above, also equipped with a smooth framing on its halo, we only have to deform the framing to make it satisfy the conditions of a framed bordism. The difference between the framing on $\partial (\beta \circ \alpha)$ and the one required on the boundary of a framed $2$-bordism of this type can be expressed as a function $f: \partial (\beta \circ \alpha) \rightarrow GL(2, \mathbb{R})$. 

By the discussion above, $f$ has image contained in $\mathbb{R}^{+}$ and moreover takes value $1$ in some neighbourhood of the horizontal boundary. As $\mathbb{R}^{+}$ is contractible, we can define an extension $\widetilde{f}: \beta \circ \alpha \rightarrow \mathbb{R}^{+}$, moreover this extension can also be taken to take value $1$ in some neighbourhood of the horizontal boundary. By tilting the framing of $\beta \circ \alpha$ using this extension $\widetilde{f}$ we obtain the needed structure of a framed $2$-bordism. Note that all such $\widetilde{f}$ are necessarily homotopic and so the structure of a framed $2$-bordism does not depend on the choice of this extension.

To complete the two categories $(\mathbb{B}ord_{2}^{fr}) _{(0)}$, $(\mathbb{B}ord_{2}^{fr}) _{(1)}$ into a symmetric monoidal pseudo double category we also need to supply the needed functors. The homomorphisms $S, T: (\mathbb{B}ord ^{fr} _{2}) _{(1)} \rightarrow (\mathbb{B}ord ^{fr} _{2}) _{(0)}$ are induced by taking the source and target of a $1$-bordism. The functor $U: (\mathbb{B}ord ^{fr}_{2}) _{(0)} \rightarrow (\mathbb{B}ord ^{fr}_{2})_{(1)}$ is more complicated, it constructs identity bordisms by taking a product with the interval and endowing it with a "horizontally constant" framing. Such a framing is not unique, but it is unique up to canonical homotopy class of isotopies relative to the boundary and so all the resulting choices are canonically isomorphic in $(\mathbb{B}ord_{2}^{fr}) _{(1)}$. 

The proof that the resulting symmetric monoidal pseudo double category is fibrant, and so leads to a symmetric monoidal bicategory, is analogous to the unoriented case.
\end{proof}

\begin{rem}
The definition of vertical composition of framed $2$-bordisms we have given above is slightly messy, but this complication seems unavoidable if one wants to use our more rigid definition of a framed $2$-bordism. 

Alternatively, when defining the latter one could only require that the framing on the vertical part of the boundary is constant up to some tilting function $f: I \rightarrow \mathbb{R}^{+}$ and allow isotopy of framings that restricts to isotopy of such tilting functions on the boundary. The symmetric monoidal bicategory one would obtain in this way is isomorphic to the one we have constructed. 
\end{rem}

We will finish this chapter by making a few simple computations related to the possible framings on the generator $2$-bordisms of the oriented bicategory. These results are in some sense starting data that we will need to feed later into the theory of framed planar diagrams. We start with a simple result which greatly simplifies construction of $2$-cells in $\mathbb{B}ord_{2}^{fr}$. 

\begin{lem}
\label{lem:uniqueness_of_framings_on_generators}
Let $\alpha: w_{1} \rightarrow w_{2}$ be one of the generating $2$-bordisms of $\mathbb{B}ord_{2}$, that is, a cusp or a generator of Morse type. Suppose we are given a compatible pair of framings on $w_{1}, w_{2}$ respectively, where by compatible we mean that the obtained framed $1$-bordisms have the same domain and codomain. If a framing of $\alpha$ making it into a framed $2$-bordism exists, then it is unique up to isotopy. In particular, once domain and codomain are fixed, there is at most one $2$-cell in $\mathbb{B}ord_{2}^{fr}$ with underlying oriented $2$-bordism $\alpha$. 
\end{lem} 

\begin{proof}
Let $v, v^\prime$ be two different framings making $\alpha$ into a framed $2$-bordism. The difference between $v, v^\prime$ is an element of $d \in Map(\alpha, GL(\mathbb{R}^{2}))$, to say that the two framings are isotopic is to say that he map $d$ is smootly homotopic to the constant map relative to the boundary. However, we have a homeomorphism $(\alpha, \partial \alpha) \simeq (D^{2}, S^{1})$ of pairs and so a continous such homotopy can be found by an elementary calculation, as $GL(\mathbb{R}^{2})$ is a one-type. Applying smooth approximation, we get the desired result. 
\end{proof}

\begin{rem}
Observe that in the argument above it is not at all important that we are working with a generator $2$-bordism, all that matters is that it is homeomorphic to a disk.
\end{rem}

\begin{prop}
Up to isomorphism-isotopy, there exist unique framed $2$-bordisms $\alpha, \beta$ lifting the usual cusp generators and expressing $pt_{+}, pt_{-}$ as a dual pair with (co)evaluation maps given by $e, c$. We call the resulting $2$-cells in $\mathbb{B}ord_{2}^{fr}$ the \textbf{framed cusp generators}. 
\end{prop}

\begin{proof}
We only have to verify that there exists some framing on $\alpha$ extending the given framing on its boundary, uniqueness will then follow from the lemma above. 

Observe that the domain of $\alpha$ is $(e \otimes pt_{+}) \circ (pt_{+} \otimes c)$, the codomain is $id_{pt_{+}}$. These two $1$-bordisms are diffeomorphic and $\alpha$, as an unoriented $2$-bordism, is the $2$-bordism associated to that diffeomorphism. It follows that the framing on the boundary of $\alpha$ can be extended to the whole $2$-bordism precisely when the diffeomorphism between the domain and codomain preserves framings up to isotopy. 

However, both the framing on $(e \otimes pt_{+}) \circ (pt_{+} \otimes c)$ and $id_{pt_{+}}$ have the property of being \emph{tangential} in the sense that the leading vector is at all times tangent to the underlying $1$-manifold. Since they also agree near the boundary, having the same domain and codomain, it follows that the two framings are isotopic. 

The argument for $\beta$ and the inverses $\alpha^{-1}, \beta^{-1}$ is the same. The composites $\alpha \circ \alpha^{-1}, \alpha^{-1} \circ \alpha, \beta \circ \beta^{-1}, \beta^{-1} \circ \beta$ are all identities again by uniqueness of framings, as they have correct (co)domains and are homeomorphic to disks, and it follows that $\alpha, \beta$ present $pt_{+}, pt_{-}$ as a dual pair.
\end{proof}

\begin{rem}
\label{rem:coherency_of_a_dual_pair_of_positive_and_negative_points}
The dual pair $\langle pt_{+}, pt_{-} \rangle_{d}$ is in fact coherent, as is not difficult to verify using the lemma, as the relevant composites that need to be compared will be diffeomorphic to disks and have the same framings on the boundary. 
\end{rem}

\begin{prop}
There exists an adjunction $e \circ (q_{+} \otimes pt_{-}) \circ b \dashv c$ with (co)unit maps $\mu_{c}$ with underlying unoriented $2$-bordism the Morse cup and $\epsilon_{c}$ with underlying inverted cap. Likewise, we have an adjunction $b \circ (pt_{-} \otimes q_{+}^{-1}) \circ c \dashv e$ with (co)unit maps $\mu_{e}$ with underlying inverted cup and $\epsilon_{e}$ with underlying cap. We call the resulting $2$-cells the \textbf{framed Morse generators}. 
\end{prop}

\begin{cor}
The framed $1$-bordism $q_{+}$ is the Serre autoequivalence of the positively framed point.
\end{cor}
Observe that there any many other framed $2$-bordisms with underlying unoriented bordisms the Morse generators, as there are many other choices of framings on domain and codomain that allow a lift to a framed $2$-bordism. Here our purpose is only to define a particular class of them, this will be useful in the next chapter. 

\begin{proof}
Again, it's enough to prove that the framing on the boundary of a given $2$-bordism can be extended to the whole manifold. This is a matter of counting "twists" the framing performs as one travels along the boundary.

Consider the counit $\mu_{c}$ of $e \circ (q_{+} \otimes pt_{-}) \circ b \dashv c$. Here, the whole boundary of the $2$-bordism is concentrated in the domain, which we picture below.

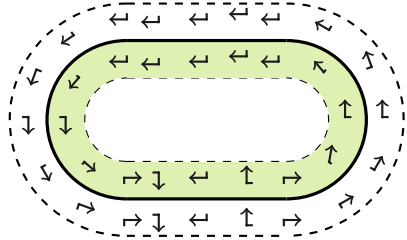
\begin{center}

\end{center}
The Morse cap $2$-bordism can be visualized as a disk attached to the inner part of the circle, verification that the framing on the boundary will extend to it amounts to observing that the framing drawn can be extended to the inner region bounded by the $1$-bordism, which is easy to see. All the other cases are analogous, except a little more difficult to draw in a "flat" way. 
\end{proof}

\begin{rem}
\label{rem:coherency_of_a_fully_dual_pair_of_positive_and_negative_points}
One easily verifies that the adjunctions $e \circ (q_{+} \otimes pt_{-}) \circ b \dashv c$ and $b \circ (pt_{-} \otimes q_{+}^{-1}) \circ c \dashv e$ we constructed equip the pair of standard positively and negatively framed points with a structure of a coherent \emph{fully dual pair}. Again, all the relevant composites will be diffeomorphic to disks and share the same framing on the boundary. 
\end{rem}

\section{Calculus of framed surfaces}

In this chapter we will describe how to enrich the oriented planar diagrams of Christopher Schommer-Pries with additional data that will also remember the isotopy class of a framing on the reconstructed surface. Our main application if this new variant of diagrammatic calculus will be to obtain a presentation of the framed bordism bicategory, which we do in the next chapter, although some of the results and arguments given here could be of independant interest. 

We first describe how given an oriented surface and a generic function into $\mathbb{R}^{2}$, one can enrich the induced oriented planar diagram by making it into what we call a \emph{framed planar diagram}. Conversely, given a framed planar diagram one can construct a framed surface mapping into $\mathbb{R}^{2}$, the two processes are inverse to each other up to isotopy of framings. 

We then present two different equivalence relations on framed planar diagrams, that of \emph{strong equivalence} and simply \emph{equivalence}. The first one applies only to framed planar diagrams with some chosen underlying oriented planar diagram, strong equivalence classes of such diagrams correspond to isotopy classes of framings on the reconstructed oriented surface. 

The second relation, that of equivalence, applies to all framed planar diagrams at once and describes diffeomorphism-isotopy classes of surfaces. It is the latter which we will use when deriving the presentation of the framed bordism bicategory. 

In this chapter we will work with closed two-manifolds and we make a convention of referring to them simply by the word \emph{surface}. This restriction to boundary-less manifolds is done for greater clarity, the more general case of $2$-bordisms, which are modeled by \emph{relative} framed planar diagrams, will be developed in the next chapter. There are no real problems in passing to this relative case and a more adventurous reader is invited to think in terms of surfaces with boundary from the get-go.

\subsection{Framed graphics}

In this section we will describe the "visual part" of a framed planar diagram, the \emph{framed graphic}. 

The general idea is as follows. We fix an oriented surface $W$ together with a generic map into $\mathbb{R}^{2}$, this induces a graphic of critical values $\Psi$. We use this data to define a canonical isotopy class of \emph{ambient framings} of $W$, defined on the complement of the preimage of the singularities of the graphic. On this complement, any other framing $v$ of $W$ compatible with the given orientation can be represented as a function into $GL_{+}(2, \mathbb{R})$ by comparing it to the ambient framing.

This function is well-defined up to homotopy, as all ambient framings are isotopic. Since we have a homotopy equivalence $GL_{+}(2, \mathbb{R}) \simeq SO(2)$ and the latter is a circle, the homotopy type of this function can be represented by a preimage of any regular value of this map, which is a codimension $1$ submanifold. Essentially, the framed graphic will consists of the usual graphic $\Psi$ together with the image of this submanifold under $f$. For this to be sufficiently well-behaved, we will need to consider "particularly nice" or what we call \emph{diagrammable} framings of $W$, we will also show that any framing is isotopic to a diagrammable one.

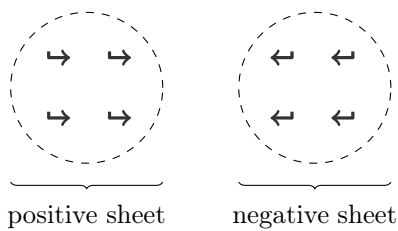
\begin{figure}[htbp!]
	\begin{tikzpicture}
		\begin{scope}
		\draw[dashed] (1,0) arc (0 : 360 : 1);

		\foreach \x / \y in {0.3/0.4,-0.5/0.4, 0.3/-0.4, -0.5/-0.4}
			{
			\begin{scope}[xshift=\x cm, yshift=\y cm]
				\node [circle, minimum width = 0.8cm] (Lcir) at (0,0) {};
				\node [circle, minimum width = 0.4cm] (Scir) at (0,0) {};
				\node (Arrow) at (Lcir.0) {};
				\node (Dash) at (Scir.90) {};
			
				\node [circle, fill= black!80,inner sep=0.4pt] at (Lcir.center) {};
				\draw [->, very thick, black!80] (Lcir.center) to (Arrow);
				\draw [very thick, black!80] (Scir.center) to (Dash);
			\end{scope}
			}
		
		\draw [
    			decoration={
    	    		brace,mirror,
      	  		raise=.55cm
    		},
  		decorate
		] (-1,-0.7) -- (1,-0.7);
		\node at (0, -1.7)  {positive sheet};
		\end{scope}
		
		\begin{scope}[xshift=3cm]
		\draw[dashed] (1,0) arc (0 : 360 : 1);

		\foreach \x / \y in {0.5/0.4,-0.3/0.4, 0.5/-0.4, -0.3/-0.4}
			{
			\begin{scope}[xshift=\x cm, yshift=\y cm]
				\node [circle, minimum width = 0.8cm] (Lcir) at (0,0) {};
				\node [circle, minimum width = 0.4cm] (Scir) at (0,0) {};
				\node (Arrow) at (Lcir.180) {};
				\node (Dash) at (Scir.90) {};
			
				\node [circle, fill= black!80,inner sep=0.4pt] at (Lcir.center) {};
				\draw [->, very thick, black!80] (Lcir.center) to (Arrow);
				\draw [very thick, black!80] (Scir.center) to (Dash);
			\end{scope}
			}
		
		\draw [
    			decoration={
    	    		brace,mirror,
      	  		raise=.55cm
    		},
  		decorate
		] (-1,-0.7) -- (1,-0.7);
		\node at (0, -1.7)  {negative sheet};
		\end{scope}
	\end{tikzpicture}
\caption{Ambient framing on a region with no singularities}
\label{fig:ambient_framing_on_a_region_with_no_singularities}
\end{figure}

\begin{defin}
\label{defin:ambient_framing}
Let $W$ be an oriented surface and $f: W \rightarrow \mathbb{R}^{2}$ be a generic function inducing a graphic $\Psi \subseteq \mathbb{R}^{2}$. Let $W_{0} = f^{-1}(\Psi_{0})$ be the preimage of singularities of $\Psi$. Observe that the complement $W_{> 0} = W \setminus W_{0}$ is naturally foliated by the preimages of the horizontal lines $L_{y} = \mathbb{R} \times \{ y \}$, as the map $f$ is transversal to this folation when restricted to $W_{> 0}$. 

A framing on $W_{> 0}$ is called \textbf{ambient} if it is compatible with the chosen orientation, its leading vector is always tangent to the foliation and at all the non-critical points $x \in W$, the image of the second vector under $df_{x}$ is pointing upwards in $\mathbb{R}^{2}$. 
\end{defin}
Observe that since an ambient framing always induces the given orientation of $W$, the second condition is equivalent to saying that at non-critical points the image of the leading vector under $df$ points to the right if $x$ lies on a positive sheet and to the left otherwise. Hence, the ambient framing on the set of non-critical points of $f$ can be pictured as the framing taken diffeomorphically to one of the framings on $\mathbb{R}^{2}$ given in \textbf{Figure \ref{fig:ambient_framing_on_a_region_with_no_singularities}}, depending on the sheet.

\begin{figure}
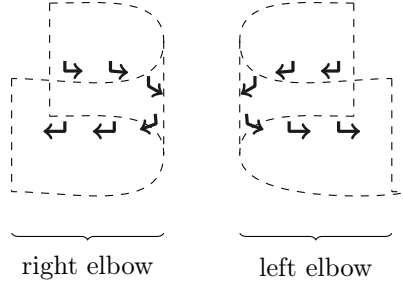


\caption{Ambient framing around a fold singularity}
\label{fig:ambient_framing_around_a_fold_singularity}
\end{figure}

One similarly sees that the behaviour of the ambient framing around a fold singularity is the one given in \textbf{Figure \ref{fig:ambient_framing_around_a_fold_singularity}}. Observe that the map $f$ induces two different diffeomorphisms of a neighbourhood of a fold singularity with the surfaces pictures, as the two folded sheets can be interchanged, but only under one of these diffeomorphisms the pictured framing will be taken to one compatible with the chosen orientation.

\begin{prop}
Let $W$ be an oriented surface and $f: W \rightarrow \mathbb{R}^{2}$ generic map. An ambient framing on $W _{> 0}$ always exists and all such framings are isotopic.
\end{prop}

\begin{proof}
The property of being ambient is local and so framings of $W _{> 0}$ satisfying it form a sheaf. Observe that pointwise, the leading vectors of any two ambient framings have the same leading vector up to a positive scalar. Since they also induce the same orientation, it follows that they differ only by an element of $G \subseteq GL(2, \mathbb{R})$ of the subgroup of transformation that preserve the complete flag $0 \subseteq \mathbb{R} \subseteq \mathbb{R}^{2}$ of subspaces together with their orientation. 

One sees immediately that ambient framings form in fact a sheaf of sections of a principal $G$-bundle. Since $G$ is contractible, so is the space of global sections and both claims follow at once.
\end{proof}
Because of the proposition above, we will generally - by abuse of langauge - speak of \emph{the} ambient framing, as it is well-defined up to a canonical class of isotopies.

The first property that distinguishes generic framings from diagrammable ones will be that the latter have a very controlled behaviour in some neighbourhood of the preimage $f^{-1}(\Psi_{0})$ of codimension $2$ singularities of $\Psi$. This is important, as these preimages consists of exactly these points where the ambient framing is not defined and so where our identification of framings with functions into $GL_{+}(2, \mathbb{R})$ breaks down. As we have discussed previously, there are essentially three types of such singularities: cusps, singularities of Morse type and fold crossings.

The preimage of any small neighbourhood of such singularity will consist of a disjoint union of disks, geometry of which we can control by the work of Schommer-Pries. More precisely, this preimage will consist of a disjoint union of trivial sheets taken diffeomorphically onto the chosen small neighbourhood and of one non-trivial sheet, which the map $f$ allows us to identify with one of the generator $2$-bordisms, that is, cup/cap, saddle, cusp or a fold crossing. 

One observes these small preimages are component-wise contractible and so in particular frameable. Moreover, up to isotopy they admit a unique framing compatible with the given orientation. Hence, our "nice framings" can have any prescribed behaviour on the preimages of a small neighbourhood of such singularity, the homotopy extension property will imply that any framing is isotopic to a "nice" one.

\begin{defin}
Let $W$ be an oriented surface and let $f: W \rightarrow \mathbb{R}^{2}$ be a generic function inducing a graphic $\Psi \subseteq \mathbb{R}^{2}$. Let $x \in \Psi_{0}$ be a codimension $2$ singularity of the graphic and let $V_{x}$ be a rectangular neighbourhood of it which is small in the sense that it does not contain any other codimension $2$ singularities. We will say that the framing of $W$ is \textbf{standard in} $V_{x}$ if the framing on the preimage $U_{x} = f^{-1}(V_{x})$

\begin{itemize}
\item is ambient on the trivial sheets, where the map $f$ restricts to a local diffeomorphism,
\item in the case of a fold crossing, is ambient on each of the two components of the non-trivial sheet, which the map $f$ identifies with fold surfaces
\item in the case of Morse or cusp singularity, on the non-trivial sheet coincides with the framing of one of the framed $2$-bordism generators described in the previous chapter, according to its singularity type and \textbf{Table \ref{tab:standard_framings_on_singularities}}.
\end{itemize}
\end{defin}

\begin{table}[!htbp]
\caption{Standard framings on singularities}
\begin{tabular}{l || c}
	\textbf{Singularity type} & \textbf{Framing} \\

	\hline
	Morse cap & $\mu_{e}$ \\
	\hline
	Morse inverted cap & $\mu_{c}$ \\
	\hline 
	Morse cup & $\epsilon_{e}$ \\
	\hline
	Morse inverted cup & $\epsilon_{c}$  \\
	\hline
	Upwards cusp on a positive sheet & $\alpha$ \\
	\hline
	Upwards cusp on a negative sheet & $\beta$ \\
	\hline
	Downwards cusp on a positive sheet & $\alpha^{-1}$ \\
	\hline
	Downwards cusp on a negative sheet & $\beta^{-1}$ 
\end{tabular}
\label{tab:standard_framings_on_singularities}
\end{table}

\begin{rem}
It is not exactly important what exactly framing we choose around codimension $2$-singularities, as long as we control the framing on the boundary of the non-trivial sheet. We need the latter because these singularities will later give the generators of the framed bordism bicategory and we would like to have a possibly small amount of them to work with. The framing on the boundary corresponds to the domain and codomain of the generator, so it is very important.

On the other hand, what exactly choice we make in the interior is not really relevant. Since non-trivial sheets are component-twise diffeomorphic to disks, any two framings that coincide on the boundary will be isotopic relative to it. 
\end{rem}
We now give the definition of diagrammable framings, it is unfortunately extremely technical. The idea is that we assume that our framing is standard on some neighbourhoods of the points on the surface $W$ surface that become codimension $2$ singularities and only try to describe it on the complement, which we term $\widetilde{W}$. There, it is possible to identify any framing with a function into $GL_{+}(2, \mathbb{R}) \simeq SO(2) \simeq S^{1}$ by comparing it to the ambient framing.

Such a function is, up to homotopy, described by the preimage of any regular point and orientation of its normal bundle, this follows from the Thom-Pontryagin construction. We assume that this function has $-id \in SO(2)$ as a generic point and call the preimage the \emph{Serre submanifold}. It lies on the surface $\widetilde{W}$, which is not ideal, as we would rather remember the framing via some data embedded in $\mathbb{R}^{2}$. Hence, a string of technical conditions on this submanifold follow, their purpose is to make the image somehow well-behaved. We term this image the \emph{twist graphic}, it will be a part of a structure of a framed planar diagram.

\begin{defin}
Let $W$ be an oriented surface and let $f: W \rightarrow \mathbb{R}^{2}$ be a generic function inducing a graphic $\Psi \subseteq  \mathbb{R}^{2}$. We will say a framing $v$ of $W$, compatible with the orientation, is \textbf{diagrammable} with respect to a system of small rectangular neighbourhoods $V_{x}$ of codimension $2$ singularities $x \in \Psi_{0}$ if it is standard in the preimages $U_{x} = f^{-1}(V_{x})$ of these neighbourhoods and such that on the complement $\widetilde{W} = W \setminus (\bigcup \ U_{x})$ it has the following properties.

\begin{enumerate}
\item The function $\widetilde{V}$ is transversal to $-id \in SO(2)$, where $\widetilde{V}: \widetilde{W} \rightarrow SO(2)$ is the function comparing $v$ to the ambient framing composed with the retraction $p: GL_{+}(2, \mathbb{R}) \rightarrow SO(2)$ given by Gram-Schmidt orthogonalization.
\end{enumerate}
This implies that the preimage $\widetilde{V}^{-1}(-id) \subseteq \widetilde{W}$ is a codimension $1$ submanifold, we call it the \textbf{Serre submanifold}.

\begin{enumerate}
\setcounter{enumi}{1}

\item The submanifold $S$ is transversal to the set of critical points $C = \{ w \in \widetilde{W} \ | \ rk(df _{w}) < 2 \}$.
\end{enumerate}
The statement makes sense, as $C$ is a submanifold since $\widetilde{W}$ does not contain any points that become codimension $2$ singularities in the graphic. We call the finite set of intersections $C \cap S$ the set of \textbf{crossing points}.

\begin{enumerate}
\setcounter{enumi}{2}
\item The tangent space to $S$ at the crossings points is not contained in the kernel of the map $f$. 
\item The function $f$ immerses $S$ into $\mathbb{R}^{2}$ with at most double transversal self-intersections.
\end{enumerate}
We call the image $f(S)$ the \textbf{twist graphic} and denote it by $\Upsilon \subseteq \mathbb{R}^{2}$, by the condition above it is an immersed codimension $1$ submanifold of the plane. It intersects $\Psi$ in finitely many points, which are exactly the images of the crossing points. We call the complement of this image the \textbf{free part} of $\Upsilon$.

\begin{enumerate}
\setcounter{enumi}{4}
\item The double self-intersections of the twist graphic do not lie on the graphic $\Psi$.
\item The free part of $\Upsilon$ is transversal to the graphic $\Psi$. 
\item The "height" function $\pi_{2} f$, where $\pi_{2}: \mathbb{R}^{2} \rightarrow \mathbb{R}$ is the projection to the second coordinate, is Morse when restricted to the Serre submanifold.
\end{enumerate}
The last property implies that the height function is Morse on the twist graphic and is in fact equivalent to it, as we assumed $f: S \rightarrow \Upsilon$ to be an immersion.

\begin{enumerate}
\setcounter{enumi}{7}
\item The critical points of the height function on $\Upsilon$ do not occur at intersections with the graphic $\Psi$.
\end{enumerate}
We call the tuple $(\Psi, \Upsilon)$ the \textbf{framed graphic} of $W$. 
\end{defin}

\begin{rem}
We will in general identify $GL_{+}(2, \mathbb{R})$ with $SO(2)$, since the inclusion $SO(2) \rightarrow GL_{+}(2, \mathbb{R})$ is a homotopy equivalence and a group homomorphism with an explicit left inverse given by Gram-Schmidt orthonormalization. It follows that for our purposes the two groups can be virtually interchanged at all times.

However, using $SO(2)$ is technically more conveniant for two reasons. One is that it is simply a smaller space, being homeomorphic to a circle. The second reason is that the "comparison with the ambient framing" function $W_{> 0} \rightarrow GL_{+}(2, \mathbb{R})$ \emph{does} depend on the choice of an ambient framing, however the composite with Gram-Schmidt \emph{does not}, as one easily verifies.
\end{rem}

\begin{rem}
The above definition makes it clear that the notion of a diagrammable framing depends on the choice of small rectangular neighbourhoods of codimension $2$ singularities of the graphic of critical values. These are assumed to be implicitly chosen whenever we fix a generic map $f: W \rightarrow \mathbb{R}^{2}$ inducing a given oriented planar diagram, we will later see the equivalence class of a framed planar diagram being induced by a given framing on $W$ does not depend on these choices.
\end{rem}

The next two remarks are rather technical and concern why we chose exactly the set of genericity conditions presented above. We would advise the reader to skip them on the first reading, as it is better to first gain intuition about the behaviour of framed graphics, which we do later with the aid of pictures, and only then verify that these are exactly the conditions which enforce this behaviour.

\begin{rem}
We do not ask the double self-intersections of the graphic $\Psi$ not to lie on $\Upsilon$, as we treated these double intersections as singularities and cut them out of $W$ when defining $\widetilde{W}$. One could also imagine an alternative definition where we do not cut out these singularities, but then require the twist graphic to be disjoint from them. The definition of framed planar diagram one would then obtain would be probably completely the same. 
\end{rem}

\begin{rem}
We only ask the free part of $\Upsilon$ to be transversal to $\Psi$, as this condition cannot possibly be satisfied at the images of the crossing points. Indeed, if $x \in S \cap C \subseteq \widetilde{W}$ is a crossing point, then the tangent space to $\Upsilon$ (resp. to $\Psi$) is the image of the tangent space $TS_{x}$ (resp. $TC_{x}$) under $df$ and both of these images are contained in $im(TW _{x}) \subseteq T \mathbb{R}^{2} _{f(x)}$. The latter is one-dimensional, as $x$ is a critical point, and this shows that $\Upsilon$ and $\Psi$ \emph{always} intersect non-transversally at the images of crossing points.

However, if it is the free part of $\Upsilon$ that intersects $\Psi$, then this means that we have points $x \in S$, $y \in C$ with $f(x) = f(y)$, but $x \neq y$. Hence, there is no obstruction to transversality like in the case above, as the map $f$ will be a local diffeomorphism around $x$. We will study both of these situations in more detail when discussing singularities of framed graphics and we will then give figures illustrating them.
\end{rem}
The key result is that any framing can be deformed via an isotopy to a diagrammable one. However, the definition given was very technical and one cannot expect a proof of existence of such framings to be any different. The intuition here is that a generic framing will have all of these properties, as they can all be formulated as some form of transversality. This makes the argument rather standard, if slighly unenlightening.

\begin{prop}[Approximation by diagrammable framings]
Let $W$ be an oriented surface, let $f: W \rightarrow \mathbb{R}^{2}$ be a generic function inducing a graphic $\Psi$. Then, any framing on $W$ inducing the given orientation is isotopic to a diagrammable one.
\end{prop}
Here, when we say diagrammable, we implicitly assume that a system of small rectangular neighbourhoods $V_{x}$ of codimension $2$ singularities of the graphic has been chosen. The statement then says that for any such choice, any framing can be deformed to a diagrammable one.

\begin{proof}
Let $v$ be a framing on $W$, let $U_{x} = f^{-1}(V_{x})$ be the system of preimages of small rectangular neighbourhoods. Observe that each of the preimages consists of a disjoint union of disks by the description of the map $f |_{U_{x}}$ one can read off from the graphic. More precisely, each of them will consists of a disjoint union of trivial sheets where the map $f$ restricts to a diffeomorphism and one non-trivial sheet where the geometry of $f$ is described by the singularity type of $x$. 

What is important is that the framing on a disk inducing a given orientation is essentially unique up to homotopy and the same is true for disjoint unions of disks. Using the homotopy extension property for the pair $\bigcup U_{x} \hookrightarrow W$ we may assume that $v$ is standard with respect to $V_{x}$.

We will now work with the complement $\widetilde{W} = W \setminus (\bigcup \ U_{x})$, which is a manifold with boundary. We can identify $v |_{\widetilde{W}}$ with a function $V: \widetilde{W} \rightarrow GL_{+}(2, \mathbb{R})$ via the comparison with ambient framing, which is well-defined on $\widetilde{W}$.  Under this identification, isotopy of framings corresponds to homotopy of such functions and so we just want to show that $V$ is homotopic with a map having all the needed properties. This will be done, as we already announced, by expressing these properties as transversality with respect to submanifolds of (multi)jet spaces, jet transversality will then give the needed statements.

Note that the conditions are in fact genericity conditions on the composition $pV: \widetilde{W} \rightarrow SO(2)$, but the map $p$ is a submersion and hence so are all the maps induced by it on jet spaces we work with, any such conditions are then also generic on $V$ itself. Moreover, if $pV$ is homotopic to a map having the needed properties, $V$ will also be, as $p$ is a homotopy equivalence, too. Hence, by abuse of langauge, we will give the proof as if $V$ itself was a function into $SO(2) \subseteq GL_{+}(2, \mathbb{R})$, although it need not be.

A minor problem in applying genericity statements is that we need to keep the function $V$ fixed on the boundary, so that together with the standard framing on $U_{x}$ it will still give a global framing of $W$. Jet transversality is rarely stated in a form where only functions fixed on the boundary are considered, but it \emph{does} apply in this case, as long as the function fixed on the boundary admits at least one extension satisfying all the needed conditions in some neighbourhood. This is the case we are in, as one readily verifies, as the framing is assumed to be standard on the boundary and so is well-behaved. One possible reference for this type of jet transversality is \cite[Theorem 1.7]{chrisphd}. 

We will now list the relevant submanifolds of jet spaces.

\begin{enumerate}
\item This is just ordinary transversality with respect to $A_{1} = \{ -id \} \subseteq SO(2)$.

\item To enforce this property we enforce transversality to the submanifold $A_{2} \subseteq J_{1}(\widetilde{W}, SO(2))$ lying in the first jet space, points of which consists of tuples $(w, z, \phi)$ where $w \in \widetilde{W}, z \in SO(2)$ and $\phi: T_{w, \widetilde{W}} \rightarrow T_{z, SO(2)}$ is a linear map on tangent spaces. The relevant subset is given by

\begin{center}
$A_{2} = \{ (c, -id, \phi) \ | \ c \in C, ker(\phi) = T_{c,C} \}$,
\end{center}
in other words by the points lying on the set of critical points $C$ with value $-id$ and the kernel of the differential contained in the tangent space of $C$. It's a submanifold, as maps $\phi$ satisfying $ker(\phi) = T_{c, C}$ are the total space of the homomorphism vector bundle

\begin{center}
$\text{Hom}(T_{\widetilde{W}} / T_{C}, T_{-id, SO(2)})$
\end{center}
over the submanifold $C \times \{ -id \} \subseteq J_{0}(\widetilde{W}, SO(2))$. By a simple count, the codimension of $A_{2}$ in the jet space is $3$. It follows that the first jet of a generic function $\widetilde{W} \rightarrow SO(2)$ does not meet $A_{2}$ at all.

\item This is similar to the one above, we have to use the subset $A_{3} \subseteq J_{1}(\widetilde{W}, SO(2))$ given by 

\begin{center}
$A_{3} = \{ (c, -id, \phi) \ | c \in C, ker(\phi) = ker(f) _{c} \}$.
\end{center}
It's a codimension $3$ submanifold for the very same reason, as $ker(f) |_{C}$ is a vector bundle over $C \subseteq \{ -id \}$, being a kernel of a rank $1$ map of vector bundles.

\item Observe that $f$ being an immersion of $S$ already follows from the properties given above. Indeed, away from the crossing points $f$ is a local diffeomorphism and at the crossing points the tangent space to the Serre submanifold is also not contained in $ker(df)$ by $(3)$. Thus, we only have to enforce the "at most double intersections" property. We first rule out triple intersections by using a submanifold $A_{4}$ of the third $J^{(3)}_{0}$ multijet space consisting of tuples of three preimages of $-id$ that have the same value under $f$. This is a codimension $7$ property, since specifying the value of $-id$ cuts out $3 \cdot dim(SO(2)) = 3$ dimensions and coincidence of values of a function $f$ cuts out $(3 - 1) \cdot dim(\mathbb{R}^{2})$. Since the dimension of $\widetilde{W} ^{(3)}$ is only six, generically such singularities do not occur at all.

We then use the submanifold $A^{\prime} _{4} \subseteq J^{(2)} _{1}(\widetilde{W}, SO(2))$ that consists of those pairs of preimages of $-id \in SO(2)$ that happen to have the same value under $f$ and such that the kernels of differential coincide in the tangent space of the image. It is of codimension $5$ and hence transversality with respect to this submanifold implies that the self-intersections of the arcs of the twist graphic are all transversal.

\item This is similar to the one above, we use $A _{5} \subseteq J^{(2)} _{1}(\widetilde{W}, SO(2))$ that consists of those pairs of preimages of $-id \in SO(2)$ that happen to have the same value under $f$ and moreover this value is critical, ie. lies on the graphic $\Psi$. This is again of codimension $5$ and hence generically such singularities do not occur.

\item Observe that by saying \emph{arcs} of the twist graphic $\Upsilon ^{1}$ we specifically exclude the images of crossing points, so when $\Upsilon^{1}$ does cross the graphic $\Psi$, the corresponding arc of the Serre submanifold and arc of critical points of $f$ do not cross at all. Hence, it's enough to ask for transversality with respect to $A_{6} \subseteq J_{1}(\widetilde{W}, SO(2))$ with respect to the submanifold consisting of tuples of points $w \in \widetilde{W}$ which are not critical, but their image lies on the graphic $\Psi$, taking value $-id \in SO(2)$ with the kernel of the differential coinciding under $f$ with the tangent space of the graphic.

\item Away from critical points, this can be rephrased using the implicit function theorem into saying that at all preimages of $-id \in SO(2)$ at which the partial differential $\partial _{x} V$ vanishes, where here we lift the tangent vector $\partial _{x}$ from $\mathbb{R}^{2}$ to $\widetilde{W}$ using $df$, the second differential $\partial^{2} _{x} V$ doesn't. This follows from the fact that when $\partial _{x}$ vanishes at a preimage of $-id \in SO(2)$, then $\partial _{y}$ necessarily can't, as otherwise we would lose the required transversality to $-id$. That condition defines a submanifold $A_{7} \subseteq J_{2}(\widetilde{W}, SO(2))$ of the second jet space, it is of codimension $3$, as we have one condition on the value, one on the first differential and one on the second differential. Hence, generically such singularities do not occur.

\item This is similar to $(7)$, we have to enforce transversality to the submanifold $A_{8} \subseteq J_{1}(\widetilde{W}, SO(2))$ consisting of tuples of these points in $\widetilde{W}$ that take value $-id \in SO(2)$, $\partial _{x} V$ vanishes at these points and additionally their image lies on the graphic $\Psi$. This is a codimension $3$ condition and hence can't occur.

\end{enumerate}
\end{proof}
Before discussing the possibly behaviour of the framed graphic $(\Psi, \Upsilon)$, we first agree on a convention on how to draw it.

\begin{figure}[htbp!]
	\begin{tikzpicture}
	
	\begin{scope}
		\draw [twist, serre_blue] (0,0) to (1, -1);
	
			\draw [
	    			decoration={
	    	    		brace,mirror,
	      	  		raise=.55cm
	    		},
	  		decorate
			] (-0.5,-0.7) -- (1.5,-0.7);
			\node at (0.5, -1.7) {Positive sheet};
	\end{scope}
	
	\begin{scope}[xshift=4cm]
		\draw [twist, serre_green] (0,0) to (1, -1);
	
			\draw [
	    			decoration={
	    	    		brace,mirror,
	      	  		raise=.55cm
	    		},
	  		decorate
			] (-0.5,-0.7) -- (1.5,-0.7);
			\node at (0.5, -1.7) {Negative sheet};
	\end{scope}
	\end{tikzpicture}
\caption{Twist graphic representing the Serre submanifold lying on either a positive or a negative sheet}
\label{fig:colours_of_twist_graphic}
\end{figure}
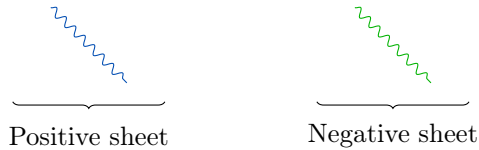

If $x \in \Upsilon$ lies on the free part of the twist graphic and is not a double intersection, then the preimage of its small neighbourhood under $f$ will be disjoint union of sheets taken diffeomorphically to this neighbourhood. Moreover, the Serre submanifold will lie on only one of those sheets and will also be taken diffeomorphically to the twist graphic. We agree to draw the the corresponding part of $\Upsilon$ using a wiggly line of blue colour if that distinguished sheet on which the Serre submanifold lies is positive and using a wiggly line of green colour if it is negative. This is pictured in \textbf{Figure \ref{fig:colours_of_twist_graphic}}.

We will later discuss the necessary labellings needed to recover a framed surface from its framed graphic and in particular all arcs of the free part of $\Upsilon$ will be labeled with an associated sheet. Hence, this distinction in colour is just a convention which is not needed from the mathematical point of view. We found, however, that it makes the diagrams a little easier to read and visualize. 

\begin{rem}
The colours are chosen in such a way that there is still a slight difference between them in grayscale, with green being slightly brighter, however, they are only supposed to be a minor help and nothing we write in dependant on this way of drawing things.
\end{rem}

Let us also make a convention on drawing framed graphics in a neighbourhood of a singularity of the graphic of critical values $\Psi$. Observe that the way we defined it, the Serre submanifold is contained in the complement of preimages of small rectangular neighbourhoods of these singularities. This implies that the twist graphic $\Upsilon$ lies on the complement of these rectangular neighbourhoods, possibly intersecting only their boundary, and then ends abruptly. 

This situation is not particularly conveniant, especially since it puts too much highlight on the choice of these small neighbourhoods, which, as we will later show, will be ultimately inconsequential. However, one can easily correct this minor problem.

Since the framing of $W$ is standard in the interior of these neighbourhoods, one observes that the intersection of the twist graphic with their boundary cannot be arbitrary and is in fact controlled by the singularity type. This allows us to artificially extend the twist graphic also into the interior of these neighbourhoods in a way which we now describe.

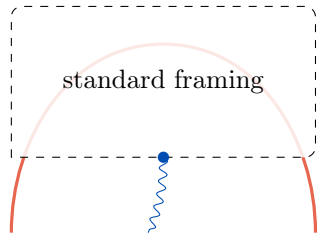
\begin{figure}[htbp!]
	\begin{tikzpicture}
		\draw[very thick, red] (0,0) [out=90, in=180] to (2,2.5) [out=0, in=90] to (4,0);
		\fill[white!50, opacity=0.85, rounded corners] (0,1) to (0,3) to (4,3) to (4,1) to (0,1);
		\node at (2,2) { standard framing };
		\draw[dashed, rounded corners] (0, 1) to (0,3) to (4,3) to (4,1) to (0,1);
		
		\node [circle, fill=serre_blue, inner sep=1.5pt] at (2,1) {};
		\draw [twist, serre_blue] (2,1) [out=-90, in=70] to (1.8, 0) {};
	\end{tikzpicture}
\caption{Behaviour of the twist graphic around Morse cap singularity}
\label{fig:unextended_twist_graphic_around_morse_cap}
\end{figure}

As a concrete example, suppose we are dealing with the Morse cap singularity and moreover the layer further from the reader is the positive one, so in terms of oriented generators it corresponds to a morphism $id_{\emptyset} \rightarrow e \widetilde{c}$. The standard framing on singularity of this type is given by the unit $\mu _{e}$ of the adjunction $ (q^{-1}_{+} \otimes id_{pt_{-}}) \circ \widetilde{c} \dashv e$, this fixes the framing on the total boundary of our generator, which in this case is just given by the codomain. In particular, when compared to the ambient framing, the framing on the boundary is transversal to $-id \in SO(2)$, as we claimed, and moreover the preimage of $id \in SO(2)$ consists of a precisely one point, a common feature of all of our choices for standard framings of Morse generators. It follows that in some small neighbourhood the configuration of the framed graphic is as pictured in \textbf{Figure \ref{fig:unextended_twist_graphic_around_morse_cap}}.

The same will happen for any other of the geometric codimension $2$-singularities, arcs of the twist graphic will join this small standard neighbourhood either from top to the bottom in a manner controlled by the domain and codomain of our standard framed generators. In this particular example, we had an occurance of $q_{+}$ in the codomain of this Morse cap generator and this is why we had one blue twist arc joined from the bottom.

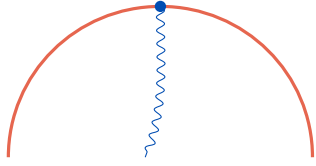
\begin{figure}[htbp!]
	\begin{tikzpicture}
		\draw[very thick, red] (0,0.5) [out=90, in=180] to (2,2.5) [out=0, in=90] to (4,0.5);
		
		\node [circle, fill=serre_blue,inner sep=1.5pt] at (2,2.5) {};
		\draw [twist, serre_blue] (2,2.5) [out=-90, in=70] to (1.8, 0.5) {};
	\end{tikzpicture}
\caption{Extended graphic representing a Morse cap singularity}
\label{fig:extended_twist_graphic_around_more_cap}
\end{figure}

We can, however, directly connect the point of intersection with the Morse signularity,  so that we would be left with the graphic given below in \textbf{Figure \ref{fig:extended_twist_graphic_around_more_cap}}.

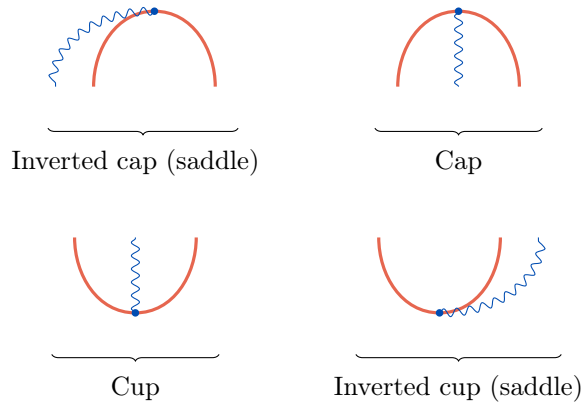
\begin{figure}[htbp!]
\begin{center}
	\begin{tikzpicture}	

	\begin{scope}
		\draw [very thick, red] (0.2,0) [out=90, in=180] to (1,1) [out=0, in=90] to (1.8,0);
		\node [circle, fill=serre_blue, inner sep=1pt] at (1,1) {};
		\draw [twist, serre_blue] (1,1) to [out=180, in=90] (-0.3,0);
		
		\draw [
    			decoration={
    	    		brace,mirror,
      	  		raise=.55cm
    		},
  		decorate
		] (-0.4,-0) -- (2.1,-0);
		\node at (0.75, -1)  { Inverted cap (saddle) };
	\end{scope}
	
	\begin{scope}[xshift=4cm]
		\draw [very thick, red] (0.2,0) [out=90, in=180] to (1,1) [out=0, in=90] to (1.8,0);
		\node [circle, fill=serre_blue, inner sep=1pt] at (1,1) {};
		\draw [twist, serre_blue] (1,0) to (1,1);
		
		\draw [
    			decoration={
    	    		brace,mirror,
      	  		raise=.55cm
    		},
  		decorate
		] (-0.1,-0) -- (2.1,-0);
		\node at (1, -1)  { Cap };
	\end{scope}
	
	\begin{scope}[xshift=-0.25cm, yshift=-3cm]
		\draw [very thick, red] (0.2,1) [out=-90, in=180] to (1,0) [out=0, in=-90] to (1.8,1);
		\node [circle, fill=serre_blue, inner sep=1pt] at (1,0) {};
		\draw [twist, serre_blue] (1,0) to (1,1);
		
		\draw [
    			decoration={
    	    		brace,mirror,
      	  		raise=.55cm
    		},
  		decorate
		] (-0.1,-0) -- (2.1,-0);
		\node at (1, -1)  { Cup };	
	\end{scope}
	
	\begin{scope}[yshift=-3cm, xshift=3.75cm]
		\draw [very thick, red] (0.2,1) [out=-90, in=180] to (1,0) [out=0, in=-90] to (1.8,1);
		\node [circle, fill=serre_blue, inner sep=1pt] at (1,0) {};
		\draw [twist, serre_blue] (1,0) to [out=0, in=-90] (2.3,1);
		
		\draw [
    			decoration={
    	    		brace,mirror,
      	  		raise=.55cm
    		},
  		decorate
		] (-0.1,-0) -- (2.6,-0);
		\node at (1.25, -1)  { Inverted cup (saddle) };
	\end{scope}
	
	\end{tikzpicture}
\end{center}
\caption{Extension of the twist graphic to the interior of standard neighbourhoods of Morse type singularities}
\label{fig:twist_graphic_around_morse_singularities}
\end{figure}

This way, we extend the twist graphic to the interior of the chosen small rectangular neighbourhoods to obtain the \textbf{extended twist graphic}. We show how to do this for all other framed generators of Morse type in \textbf{Figure \ref{fig:twist_graphic_around_morse_singularities}}. In the case of cusp generators and fold crossing, there is nothing to be done, as these always have neighbourhoods disjoint from the twist graphic.

By \emph{twist graphic} without any adjectives we will refer to the extended version, although often a system of small rectangular neighbourhoods of codimension $2$-singularities will be implicitly chosen. In such a case one can pass from the non-extended version to the extended version by filling out the interiors of small rectangular neighbourhoods in the way described above and the other way around by cutting out the intersection of the extended twist graphic with the interior of these neighbourhoods. This ends our discussion of conventions related to drawings.

We will now describe the possible singularities of a framed graphic. The word "singularity" is used here in a rather broad sense, unlike singularities of the usual graphic $\Psi$, they do not necessarily arise as images of some singular points of a generic map $W \rightarrow \mathbb{R}^{2}$. Another fitting word would perhaps be "generators", as together with the usual codimension $2$ singularities of graphics, they will form pieces into which any framed graphic can be decomposed into. 

We are only trying to describe the "new" singularities, we are not interested in those occuring already in graphics of unoriented surfaces, which are singularities of Morse type, cusps and fold crossings. In other words, we are only describing the behaviour of a framed graphic outside of small rectangular neighbourhoods of these singularities, we know that in their interior the behaviour will be as in \textbf{Figure \ref{fig:twist_graphic_around_morse_singularities}}, directly from the way we defined it. 

The first two singularities arise when the fold graphic $\Psi$ intersects the twist graphic $\Upsilon$. There are two ways in which this can happen, the first one is when the intersection is \emph{not} an image of a crossing point of $W$. This happens precisely when the corresponding arc of the Serre submanifold does not actually intersect the set of critical points, as they lie on different sheets. Hence, this intersection is in some sense "apparent".

\begin{figure}[htbp!]
	\begin{tikzpicture}
		\node (MarginNode) at (0,0.2) {};
		\draw[very thick, red] (0,0) to (2, -2);
		\draw[twist, serre_blue] (0, -2) to (2,0);
	\end{tikzpicture}
\caption{Fold-twist interchange singularity}
\label{fig:fold-twist_interchange_singularity}
\end{figure}
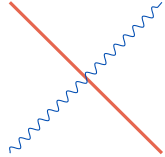
By our assumption on diagrammable framings the intersection will then be transversal, we call this singularity the \textbf{fold-twist interchange}. It is pictured in \textbf{Figure \ref{fig:fold-twist_interchange_singularity}}. It will not correspond to any generators in the presentation of the framed bordism bicategory, since in unbiased semistrict string diagram calculus in can be simply interpreted as the interchange isomorphism.

The second way the fold graphic $\Psi$ and the twist graphic $\Upsilon$ can intersect is when the intersection happens actually in the surface $W$, in other words when the intersection point is an image of a crossing point. Such an intersection cannot be transversal by a remark we made earlier, as crossing points are in particular critical for the map $f$ and so the whole image of the tangent space of $W$ in $\mathbb{R}^{2}$ is only one-dimensional.

\begin{figure}[htbp!]
	\begin{tikzpicture}
	
	\begin{scope}
		\draw [very thick, dashed] (-0.5, 0) to [out=0, in=90] (1.5, -0.25) to [out=-90, in=180] (0, -0.5);
		\draw [very thick, dashed] (0.5, -3) to [out=0, in=90] (2.5, -3.25) to [out=-90, in=180] (1, -3.5);
		\draw [very thick, dashed] (-0.5,0) to (0.5, -3);
		\draw [very thick, dashed, red] (1.5, -0.25) to (2.5, -3.25);
		\draw [very thick, dashed] (0, -0.5) to (1, -3.5);
		
		\draw [twist, serre_orange] (1, -2.5) to [out=20, in=-130] (2, -1.75);
		\draw [twist, serre_orange] (0.75, -1.25) to [out=-20, in=130] (2, -1.75);
		
		\draw [->, gray] (3, -1.8) [out=45, in=135] to node[auto] {$ f $} (4, -1.8);
	\end{scope}
	
	\begin{scope}[xshift=6cm]
		\draw[very thick, red] (0,-0.5) to (0, -2.5);
		\draw[twist, serre_blue] (-1, -2.5) [out=90, in=-90] to (0, -1.5);
		\draw[twist, serre_green] (0, -1.5) [out=90, in=-90] to (-1, -0.5);
	\end{scope}
	\end{tikzpicture}
\caption{Fold-twist crossing singularity}
\label{fig:fold_twist_crossing_singularity}
\end{figure}
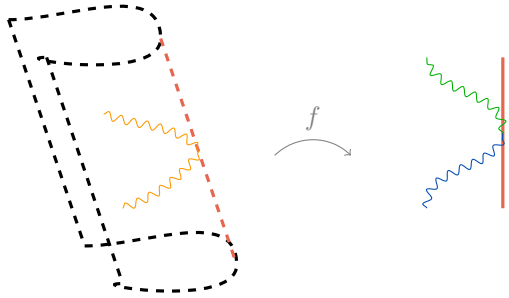

By assumption, this is happening away from codimension $2$ singularities of the map $f$, so such a point has a small neighbourhood such that its preimage will consist of a disjoint union of trivial sheets and one in which the map $f$ is a projection from a fold surface. The Serre submanifold crosses the set of critical points and by our assumptions on diagrammable framings it does so transversally. We picture the situation, together with the corresponding graphic, below in \textbf{Figure \ref{fig:fold_twist_crossing_singularity}}.

On the left we have the fold surface, the set of critical points is drawn in red and the Serre submanifold crossing it in orange. On the right we have the graphic, the Serre submanifold crosses the set of critical values in a tangential way. Here the convention on blue and green colours reminds us that the Serre submanifold "crossed to the other side of the fold surface", in the pictured situation it starts at a negative sheet and as it moves downwards it goes to the positive sheet.

This is the \textbf{fold-twist crossing} singularity and it indeed will give us a non-trivial generator of the framed bordism bicategory. Under the interpretation of $\mathbb{B}ord_{2}^{fr}$ as the free symmetric monoidal bicategory on a fully dualizable object, this $2$-cell corresponds to an isomorphism between one of the (co)evaluation maps composed with Serre autoequivalence of the left dual and the same map composed with Serre autoequivalence of the right dual. 

The next singularity is analogous to the first one and it corresponds to the double crossing of the twist graphic with itself. Again, by the properties of diagrammable framings, such a crossing will be transversal.
\begin{figure}[htbp!]
	\begin{tikzpicture}
		\node (MarginNode) at (0,0.1) {};
		\draw[twist, serre_blue] (0,0) to (2, -2);
		\draw[twist, serre_blue] (0, -2) to (2,0);
	\end{tikzpicture}
\caption{Twist-twist interchange singularity}
\label{fig:twist-twist_interchange_singularity}
\end{figure}
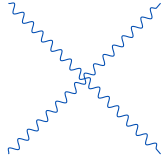
This is the \textbf{twist-twist interchange} singularity, we picture it in \textbf{Figure \ref{fig:twist-twist_interchange_singularity}}. The associated arcs of the Serre submanifold must lie on different sheets, as otherwise it would fail to be a submanifold at all and so this intersection is also "apparent". Similarly to the fold-twist interchange, this singularity will not correspond to any generators of the framed boridsm bicategory, as in the language of unbiased semistrict string diagrams it is expressible as the interchange isomorphism.

The last singularity happens at the critical points of the height function restricted to the twist graphic. We assumed that this restriction is Morse and hence the situation is one of the two pictured below.

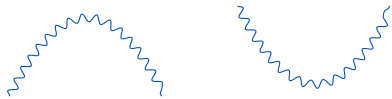
\begin{figure}[htbp!]
	\begin{tikzpicture}
	\node (MarginNode) at (0,1.2) {};
	\draw[twist, serre_blue] (0,-0.2) parabola bend (1, 0.8) (2,-0.2);
	\draw[twist, serre_blue] (3, 1) parabola bend (4, 0) (5, 1);
	\end{tikzpicture}
\caption{Twist birth-death singularity}
\label{fig:twist_birth_death_singularity_graphic}
\end{figure}
This is the \textbf{twist birth-death} singularity, we picture it in \textbf{Figure \ref{fig:twist_birth_death_singularity_graphic}}. Similarly to the case of the Morse singularities of the graphic of critical values it is not really a "singularity" of the twist graphic, as the latter remains smooth there. It is, however, important to distinguish it for the purpose of interpreting framed planar diagrams as certain unbiased semistrict string diagrams, which we will do alter. 

This singularity corresponds to a non-trivial generator of the framed bordism bicategory. We will see that under the interpretation given by the Cobordism Hypothesis, it corresponds to (co)unit isomorphisms for the adjoint equivalence between Serre autoequivalence and its pseudoinverse.

This basically ends our description of a graphical langauge one can use to draw framings on oriented surfaces, we will do a brief summary in the form of a theorem. For completness, we also describe the usual graphic $\Psi$ of a surface, partly restating \cite[Thm. 1.28, Definition 1.29]{chrisphd}, although only the results where we describe the twist graphic are new to the current work.

\begin{defin}
An (abstract) \textbf{framed graphic} consists of a tuple $(\Psi, \Upsilon)$ of diagrams embedded in the plane. The diagram $\Psi$, the \textbf{graphic of critical values}, consists of a finite number of embedded closed curves in $\mathbb{R}^{2}$ such that the height function restricts to a diffeomorphism on each of them and of two finite collections of points, corresponding to singularities of Morse type and cusps.  The diagram $\Upsilon$, the \textbf{twist graphic}, consists of a finite number of embedded closed curves in $\mathbb{R}^{2}$ such that the height function restricts to a diffeomorphism on each of them and of two finite collections of points, corresponding to crossing points and points of twist birth-death. Moreover, we require that 

\begin{itemize}
\item each Morse point has a \textbf{small rectangular neighbourhood} where the embedded curves restrict as in \textbf{Figure \ref{fig:twist_graphic_around_morse_singularities}}
\item each cusp point has a \textbf{small rectangular neighbourhood} disjoint from the twist graphic and where $\Psi$ restricts to the cusp graphic
\item each crossing point has a neighbourhood where the embedded curves restrict as in \textbf{Figure \ref{fig:fold_twist_crossing_singularity}}
\item each twist birth-death point has a neighbourhood disjoint from the graphic of critical values and where $\Upsilon$ restricts as in \textbf{Figure \ref{fig:twist_birth_death_singularity_graphic}}
\end{itemize}
\end{defin}
We emphasized the notion of small rectangular neighbourhoods in the case of singularities of the graphic of critical values to make it clear that whenever we will use this term, we will implicitly assume that the given neighbourhoods have the properties we required above. Observe that if a framed graphic is induced by a diagrammably framed surface, then the small rectangular neighbourhoods one chose to obtain a notion of a diagrammable framing will have these properties. 

Using the above definition, our analysis up to this point can be compactly stated as the following result.

\begin{thm}
Let $W$ be an oriented surface, let $f: W \rightarrow \mathbb{R}^{2}$ be a generic function and $v$ a diagrammable framing of $W$. Then the induced framed graphic $(\Psi, \Upsilon)$ is an abstract framed graphic.
\end{thm}

\subsection{Framed planar diagrams}

In this section we will describe framed planar diagrams, a combinatorial kind of diagram which can be used to model framed surfaces.

Analogously to the unoriented and oriented cases, a framed planar diagram will consists of a framed graphic, a chambering graph and some amount of sheet data. 

The chambering graph plays a role of a "nice open covering" of a graphic, as together with the graphic of critical values it divides the complement in $\mathbb{R}^{2}$ into connected components, which we call chambers. The sheet data will make it possible to recover a framed surface from its graphic, this is very similar to the oriented case, which we will briefly recall.

Let $W$ be an oriented surface, $f: W \rightarrow \mathbb{R}^{2}$ a generic function inducing a graphic $\Psi$ and let $\Gamma$ be a chambering graph compatible with $\Psi$. Then, over each chamber of $(\Psi, \Gamma)$ the map $f$ restricts to a local diffeomorphism and so the preimage is a disjoint union of sheets taken diffeomorphically to the chamber. As $W$ is oriented, over each sheet $f$ either preserves or reverses the orientation, this divides the set of sheets into the positive and negative parts.

We can recover $W$ by appropriately gluing these sheets together and the labelling data tells us precisely how to do this. Thus, over each arc of the chambering graph $\Gamma$ we have a bijection between the set of sheets over the chambers lying to the left and right of the arc. Similarly, over each arc of the graphic of critical values $\Psi$ we have a bijection between the set of sheets over a chamber on one side and the set of sheets over a chamber on the other side with the exception of two distinguished sheets, which are being "folded together". Additionally, one has more complicated data attached to singularities of $\Psi$ which tells us how to glue sheets around these points.

In this way, given $(\Psi, \Gamma)$ and the labelling data one can recover the oriented surface up to oriented diffeomorphism over $\mathbb{R}^{2}$, this is the oriented variant of Planar Decomposition Theorem,  \cite[Thm 1.52]{chrisphd}.

Suppose now that $W$ also comes with a framing diagrammable with respect to $f$, as we have seen, this enriches the graphic of critical values into the framed graphic $(\Psi, \Upsilon)$. 

If we were to recover the framing on $W$ from the framed graphic and some combinatorial data, we can start by defining our "recovered" framing to be standard in the preimage of small rectangular neighbourhoods of singularities of $\Psi$. Indeed, the framing of $W$ was assumed to be diagrammable, so it certainly had this property. This leaves us with recovering the framing outside of the preimages of these small neighbourhoods.

Observe that it is not possible to recover the Serre submanifold just from the twist graphic, as if we take a small arc of the latter, then it is not possible to determine on which sheet the corresponding arc of the Serre submanifold was lying. However, given this data we can recover it except for a finite number of crossing points, fold-twist interchange points and points of twist birth-death singularity, simply because the map from a sheet to the chamber under it is a diffeomorphism. This is enough, the whole Serre submanifold is then given by the closure of that "lifted" subset.

Thus, one part of sheet data for framed planar diagrams will be an assignment for each chamber, for each arc of twist graphic lying in that chamber, a distinguished sheet on which the corresponding arc of the Serre submanifold lies. 

The data has to satisfy some obvious compatibility conditions. Namely, that when an arc of the twist graphic crosses and edge of the chambering graph, then the sheets attached to that arc on either side of the edge are related by the bijection associated to it. This ensures that the "lift" of the twist graphic to the reconstructed surface is "continous", that is, doesn't jump at the places where we glue. To make this well-behaved with respect to later applications we will additionally require that the intersection of the chambering graph and the twist graphic consists only of transversal intersections between edges and arcs.

We require the same condition on the fold-twist interchange singularities, in particular at such singularities the sheet associated to the arc of the twist graphic cannot be either of the sheets being "folded together", as they do not participate in this bijection. This conforms with our discussion that the fold-twist interchange singularity happens precisely when the twist graphic crosses the graphic of critical values, but the corresponding of the Serre submanifold lies on a different sheet than the one being folded. 

At the crossing point singularity two arcs of the twist graphic meet. We require that they both lie on that side of the arc of the graphic of critical values where we have two distinguished sheets being folded together and that one of the arcs is attached to one of these sheets and the second to the other. This is, again, consistent with the geometry of the crossing point singularity which we discussed before.

At twist birth-death two arcs of the twist graphic meet, we only have to require that they are both attached to the same sheet. We can summarize our discussion in the following proposition, where we also state a compatibility condition on the chambering graph. 

\begin{defin}
We will say that the chambering graph $\Gamma$ is \textbf{compatible with the framed graphic} $(\Psi, \Upsilon)$ if it is compatible with $\Psi$ and additionally the intersection $\Gamma \cap \Upsilon$ consists only of transversal intersections of edges and arcs.

A \textbf{partly-framed sheet data} for $(\Psi, \Upsilon, \Gamma)$, where $\Gamma$ is assumed to be compatible with the framed graphic, consists of oriented sheet data and additionally the following assignment. To each chamber of $(\Psi, \Gamma)$, to each arc of $\Upsilon$ lying on it, we assign one of the sheets lying over the given chamber. This has to satisfy the compatibility conditions discussed above on fold-twist interchanges, intersections of the chambering graph and the twist graphic, crossing point singularities and twist birth-deaths. 
\end{defin}

\begin{prop}
\label{prop:partial_reconstruction_of_framings}
For any diagrammably framed surface with respect to a generic map $f: W \rightarrow \mathbb{R}^{2}$ it is possible to choose a chambering graph $\Gamma$ compatible with the induced framed graphic $(\Psi, \Upsilon)$ and any such choice gives rise to partly-framed sheet data for $(\Psi, \Upsilon, \Gamma)$. 

Conversely, given a framed graphic $(\Psi, \Upsilon, \Gamma)$ together with a choice of a compatible chambering graph and partly-framed sheet data it is possible to to construct an oriented surface $W^\prime$ mapping into $\mathbb{R}^{2}$, inducing $\Psi$, equipped with a codimension $1$ submanifold $S^\prime \subseteq \widetilde{W}^\prime$ of the complement of the preimage of small rectangular neighbourhoods of singularities of the graphic $\Psi$. 

Moreover, if $(\Psi, \Upsilon, \Gamma)$ and the partly-framed sheet data is itself induced by a diagrammably framed surface $W$ mapping into $\mathbb{R}^{2}$, then $W^\prime$  and $W^\prime$ will be diffeomorphic over $\mathbb{R}^{2}$ and this diffeomorphism will take the Serre submanifold $S \subseteq W$ to $S^\prime \subseteq \widetilde{W}^{\prime}$. 
\end{prop}

\begin{rem}
Here it doesn't matter if we work with extended or non-extended twist graphics, since we implicitly assume that some system of small rectangular neighbourhoods around  singularities of the graphic of critical values $\Psi$ has been chosen. 

In the third part of the statement, we also assume that we made the same choice for $W$ and $W^\prime$. This makes sense, as they both share the same graphic of critical values by construction of $W^\prime$. It is then clear that the diffeomorphism constructed will take the complement of the preimage of standard neighbourhoods $\widetilde{W} \subseteq W$ to $\widetilde{W}^\prime \subseteq W$, as it is a diffeomorphism over $\mathbb{R}^{2}$. 

These difficulties with certain non-canonicity will disappear soon, when we will recover a framed surface and not only a surface with a "candidate Serre submanifold", like we do above.
\end{rem}

\begin{proof}
There are few parts to this proposition. First, we claim that one can choose a compatible chambering graph, this is \cite[Proposition 1.46]{chrisphd} in the unoriented case, however, in the same way as one can choose the chambering graph to have edges transversal to the arcs of the graphic of critical values one can make it transversal to arcs of the twist graphic, so the proof also works here.

That one can choose an oriented sheet data and that one can then reconstruct the surface up to an oriented diffeomorphism over $\mathbb{R}^{2}$ is \cite[Thm 1.52]{chrisphd} in the unoriented case, the necessary changes in the oriented case are summarised in the discussion leading to \cite[Theorem 3.50]{chrisphd}.

The partly-framed sheet data additionally contains the data needed to obtain a codimension $1$ submanifold $S^\prime \subseteq \widetilde{W}^\prime$, its intersection with any sheet is given by a union of all twist graphics attached to that sheet, lying on the corresponding chamber of $(\Psi, \Gamma)$ and it can be recovered as the closure of the sum of all these intersections. The compatibility data ensure that it is indeed a closed submanifold with its boundary contained in $\partial \widetilde{W}^\prime$.

The diffeomorphism will take $S$ to $S^\prime$ as it is over $\mathbb{R}^{2}$, over any given chamber $S$ and $S^\prime$ are taken diffeomorphically to the arcs of the twist graphic and their arcs lie on sheets corresponding to each other under this diffeomorphism by construction.
\end{proof}

Hence, from the framed graphic, the chambering graph and partly-framed sheet data one can recover the surface together with a map to $\mathbb{R}^{2}$ and its Serre submanifold. This is almost the data needed to recover the framing of $W$, which we will do using the Thom-Pontryagin construction, the missing piece is the normal framing of the Serre submanifold.

In fact, we can do with a little less, as it is enough to be given the normal framing up to isotopy and a framing of a $1$-dimensional bundle considered up to isotopy is the same thing as an orientation. Thus, to complete our sheet data, we will attach to each chamber, to each arc of the twist graphic lying on it, an orientation of the normal bundle of the arc in the chamber.

One of the properties of arcs of the twist graphic is that they are "never horizontal", in the sense that height map restricts to a diffeomorphism on each of them. This allows us to have a standard convention on naming the relevant normal orientation, it is then easy to phrase compatibility conditions in terms of it.

More precisely, using the metric on $\mathbb{R}^{2}$ one can identify the normal bundle of any arc as the bundle orthogonal to the tangent space of the arc. As the latter is never horizontal, any vector lying in the former will always be pointing either to the left or right. That is, if we write a non-zero vector in the normal space as $v = a \partial _{x} + b \partial _{y}$, where $(x, y)$ are coordinates in $\mathbb{R}^{2}$, then we always have either $a > 0$, that is, $v$ is a right-pointing vector, or $a < 0$,  $v$ being a left-pointing vector.

We will call the orientation consisting of right-pointing vectors \emph{straightforward} and the orientation consisting of left-pointing vectors \emph{inverted}. A more natural terminology would perhaps speak of positive and negative orientations, but this could potentially collide with our language of positive and negative sheets. 

\emph{Framed sheet data} will consist of partly-framed data together with a choice of straightforward or inverted orientation for each arc of the twist graphic. This choice has to satisfy some compatibility conditions which are there to ensure that these orientations will glue to a global orientation of the normal bundle of Serre submanifold in the reconstructed surface.  

The conditions are as follows. First, we require that at places where a twist arc crosses a chamber, which are fold-twist interchange singularities and crossings with an edge of the chambering graph, both ends of the arc leaving from the crossing are either straightforward or inverted. In other words, the normal orientation is a property of the whole arc of the twist graphic, not depending on the chamber.

Moreover, there are some compatbility conditions at the singularities delimiting the arcs. At twist birth-death, two arcs meet and we require one of them to be inverted and the other to be straightforward. Similarly, two arcs meet at the crossing point singularity and we require that one of them is inverted and the other straighforward. This is consistent with the geometry of these singularities, in the latter case we need to have a "switch" of normal orientation because the orientation of the sheets to which these arcs are associated changes, as one of them is always positive and the other always negative. 

Lastly, we require that any arc of the twist graphic leaving an inverted Morse cap singularity is straightforward, any arc leaving Morse cap is inverted, any arc joining Morse cup is straightforward and any arc joining inverted Morse cup is inverted. That this is necessary can be read off from the framing on the boundary of the framed bordism generators that these singularities correspond to. In each case, we have one occurance of the Serre autoequivalence or its pseudoinverse in either the domain and codomain, one can verify that $q_{+}$ corresponds to the straghtforward blue arc and $q_{+}^{-1}$ to an inverted one. 

\begin{rem}
In fact, as we will see later, under the comparison of framed planar diagrams with unbiased semistrict string diagrams on the presentation of $\mathbb{B}ord_{2}^{fr}$, the straightforward blue arcs correspond to paths of $q_{+}$ and the inverted ones to $q_{+}^{-1}$. Similarly, straightforward and inverted green arcs will correspond to paths of respectively $q_{-}$ and $q_{-} ^{-1}$.
\end{rem}

\begin{rem}
The passage to extended framed graphic does not add any more arcs to the twist graphic and so for the issues of normal orientation it is not relevant whether one works with extended or non-extended framed graphics.
\end{rem}
In the case of framed graphic being induced by a generic map $f: W \rightarrow \mathbb{R}^{2}$ on a diagrammably framed surface, this additional data of the orientation of normal bundle of arcs descends from the orientation of the normal bundle of the Serre submanifold $S \subseteq \widetilde{W}$. The latter has this structure, as it is defined as a preimage of $-id \in SO(2)$ and this point is canonically oriented by the global orientation of $SO(2) \simeq S^{1}$. 

We can now summarize our discussion leading up to this point and show that given all the data we talked about one is indeed able to reconstruct the framed surface up to, in particular, diffeomorphism-isotopy. 

\begin{defin}
\label{defin:framed_planar_diagram}
A \textbf{framed sheet data} for a tuple $(\Psi, \Upsilon, \Gamma)$ consisting of a framed graphic together with a choice of a compatible chambering graph $\Gamma$ consists of a partly-framed sheet data and a choice of either straightforward or inverted normal orientation for each arc of the twist graphic, satisfying the compatibility conditions outlined above. A \textbf{framed planar diagram} $\aleph$ consists of a such a tuple together with a choice of framed sheet data.
\end{defin}

\begin{notation*}
We will usually denote framed planar diagrams using the letter $\aleph$, but sometimes, by abuse of notation, we might refer to it using just the tuple $(\Psi, \Upsilon, \Gamma)$, with sheet data being implicit. We chose the Hebrew $\aleph$, as it is severely underused and also to distinguish it from Greek letters, with which we have denoted the graphics and the chambering graph, which are all things embedded in $\mathbb{R}^{2}$. 
\end{notation*}

\begin{rem}
We already know that for any diagrammably framed surface $W$ together with a choice of a generic map $f: W \rightarrow \mathbb{R}^{2}$ it is possible to choose a chambering graph $\Gamma$ compatible with the induced framed graphic $(\Psi, \Upsilon)$. Any such choice leads to framed sheet data, giving us a framed planar diagram.
\end{rem}

\begin{thm}[Reconstruction of framed surfaces]
\label{thm:reconstruction_of_framed_surfaces}
Given a framed planar diagram $\aleph$ and a choice of small rectangular neighbourhoods of singularities of its graphic of critical values, it is possible to construct a framed surface $W^\prime$ mapping into $\mathbb{R}^{2}$ and inducing the given diagram. Moreover, the framing is well-defined up to isotopy and in this sense does not depend on the choice of these neighbourhoods.

If the framed planar diagram $\aleph$ is itself induced by some diagrammably framed surface $W$ and a generic map $f: W \rightarrow \mathbb{R}^{2}$, then there is a diffeomorphism $W \simeq W^\prime$ over $\mathbb{R}^{2}$ and under this diffeomorphisms the framings of $W, W^\prime$ correspond to each other up to isotopy.
\end{thm}

\begin{proof}
We first apply \textbf{Proposition \ref{prop:partial_reconstruction_of_framings}} to obtain an oriented surface $W^\prime$ mapping into $\mathbb{R}^{2}$ together with a codimension $1$ submanifold $S^\prime \subseteq \widetilde{W}^\prime$ of the complement of the preimage of small rectangular neighbourhoods of singularities of the graphic of critical values. As $S^\prime$ was constructed from the twist graphic by lifting the arcs of the twist graphic to the attached sheets, the orientations of the normal bundle of these arcs give orientations of the normal bundle of corresponding "subarcs" of $S^\prime \subseteq \widetilde{W}^{\prime}$, these orientations glue to a global normal orientation of $S^\prime$ by the compatibility conditions we imposed.

Choose a framing of $S^\prime \subseteq \widetilde{W}^\prime$ compatible with the given normal orientation, we can now apply the Thom-Pontryagin construction to construct a map $\widetilde{W}^{\prime} \rightarrow S^{1}$ with $S^{\prime}$ as the preimage of $-id \in SO(2) \simeq S^{1}$, moreover this map can have the required behaviour on the boundary. By "tilting" the ambient framing of $\widetilde{W}^{\prime}$ using this function we obtain a framing of $\widetilde{W}^{\prime}$. One can the extend this framing to a framing of whole of $W^\prime$, this is possible as we have the correct framing on the boundary of $\partial \widetilde{W} ^\prime$ by the compatibility conditions we imposed on the behaviour of the twist graphic in small rectangular neighbourhoods, namely that the Serre submanifold intersects it the correct number of times.

To prove that the framing on $W^\prime$ does not depend, up to isotopy, on the choice of small rectangular neighbourhoods, we just show that it doesn't change if we use strictly smaller ones. This is enough, as in any two rectangular neighbourhoods of a point in $\mathbb{R}^{2}$ one can find a third rectangular neighbourhood, strictly contained in both of them.

Hence, suppose we have two systems of small rectangular neighbourhoods of singularities of the graphic of critical values and let  $V_{x}$ be the larger of the two, with $U_{x} = w^{-1}(V_{x})$ being the preimage under the constructed map $w: W^{\prime} \rightarrow \mathbb{R}^{2}$. Observe that both framings may be taken to agree on $W^\prime \setminus (\bigcup U _{x})$, as there they both come by definition from the Thom-Pontryagin construction applied the same normally oriented codimension $1$ submanifold. But they also agree on $\bigcup U_{x}$, at least up to isotopy, as the latter is a disjoint union of disks and so is componentwise contractible. This isotopy may be taken to be constant on the boundary and so defines a global isotopy, showing that the framing on $W^\prime$ indeed doesn't depend on the rectangular neighbourhoods.

We have a diffeomorphism $W \simeq W^\prime$ again by $\textbf{Proposition \ref{prop:partial_reconstruction_of_framings}}$ and it will take the framing of $W$ to a possibly different framing of $W^\prime$, but they will agree on the preimage of the chosen small rectangular neighbourhoods, as they are both standard there. Moreover, both framings restricted to the complement $\widetilde{W}^\prime$ induce the same Serre submanifold, seen an a normal oriented codimension $1$ submanifold. By the Thom-Pontryagin theorem, this means that the associated functions into $S^{1} \simeq SO(2)$ are homotopic and so the framings are isotopic. As they agree on $\partial \widetilde{W}^\prime$, the homotopy can be taken to be constant there and so gives a global isotopy of framings on $W^\prime$.
\end{proof}

\subsection{Relations on framed planar diagrams}

In the previous section we have defined framed planar diagrams. These are combinatorial objects one can attach to an oriented surface $W$, a generic map into $\mathbb{R}^{2}$ and a diagrammable raming. They are enough to recover the surface with its framing up to isotopy. 

Our goal in this section will be to describe two equivalence relations on the set of framed planar diagrams. The more important one is that of \emph{equivalence}, as we will need it to derive the presentation of the framed bordism bicategory. This relation models diffeomorphism-isotopy classes of surfaces and so has the following properties.

\begin{enumerate}
\item If $\aleph _{0}$ and $\aleph _{1}$ are equivalent, then the framed surfaces one can construct from them are diffeomorphic and this diffeomorphism can be chosen to preserve the framing up to isotopy. 
\end{enumerate}
This implies that the "surface reconstructed from a framed planar diagram" gives a well-defined map from the set of equivalence classes of framed planar diagrams into the set of diffeomorphism-isotopy classes of framed surfaces.

\begin{enumerate}
\setcounter{enumi}{1}
\item If $W$ is an oriented surface equipped with a generic map $f: W \rightarrow \mathbb{R}^{2}$, then isotopic diagrammable framings lead to equivalent framed planar diagrams.
\end{enumerate}
This implies that there is a well-defined map from framed surfaces equipped with a generic map into $\mathbb{R}^{2}$ into equivalence classes of framed planar diagrams given by "deform the framing to be diagrammable and take the associated diagram". 

\begin{enumerate}
\setcounter{enumi}{2}
\item If $W$ is a framed surface, then the equivalence class of the associated diagram does not depend on the choice of a generic map into $\mathbb{R}^{2}$.
\end{enumerate}
The last property, together with our reconstruction theorem, implies that there is a bijection between equivalence classes of framed planar diagrams and diffeomorphism-isotopy classes of framed surfaces.

As in the case of unoriented planar diagrams, the equivalence relations we will describe are generated by certain modifications one can perform on framed planar diagrams, these modifications come in essentially two types. The first type is that of isotopy, where two framed planar diagrams are considered to be equivalent if they are isotopic as diagrams embedded in $\mathbb{R}^{2}$. That is, like in the case of ordinary planar diagrams, the isotopy can move the singular points of the graphic of critical values, the singular points of the twist graphic, the vertices of the chambering graph and corresponding arcs and edges freely in a smooth way, however, we require that at each time we still have a framed planar diagram.  

The second type of modification is that of a "local move", where one is allowed to replace a small part of the framed graphic by a different one, without changing it anywhere else. This can usually only be performed not only when the framed graphic has a particular form but also when the sheet data satisfies some compatibility conditions.

\begin{defin}
Two framed planar diagrams $\aleph_{0}, \aleph_{1}$ are \textbf{equivalent} if they differ by an isotopy and possibly a finite number of applications of round twist local moves, moves of crossing point type, twist homotopy, change of function or change in chambering graph moves.
\end{defin}
The meaning of the above definition will become clear as we go, as we will spend the rest of the section deriving it and describing the relevant local moves. We will start by first describing moves related to property $(2)$, that is, invariance of diagrams under the isotopy class of framings.

Our way to proceed will be to fix an oriented surface $W$ together with a generic map $f: W \rightarrow \mathbb{R}^{2}$, this gives us a diffeomorphism of $W$ with the surface reconstructed from the oriented planar diagram associated to $f$. Any diagrammable framing of $W$ will now induce a framed planar diagram with the given underlying oriented planar diagram, conversely, any such framed planar diagram describes a framing on the reconstructed surface by \textbf{Theorem \ref{thm:reconstruction_of_framed_surfaces}} and so, through the diffeomorphism, also a framing on $W$. Our goal will be to understand this relation better and since $f: W \rightarrow \mathbb{R}^{2}$ will be fixed, this is equivalent to understanding the relation between framed planar diagrams with a given underlying oriented planar diagram and the framings on the oriented surface reconstructed from it. 

To do so, we will describe an auxiliary equivalence relation on framed planar diagrams with a fixed underlying oriented planar diagram which we will call \emph{strong equivalence}. We will then have, for any oriented planar diagram, a bijection between strong equivalence classes of framed planar diagrams living over it  and isotopy classes of framings on the surface reconstructed from it. 

\begin{defin}
Two framed planar diagrams $\aleph_{0}$, $\aleph _{1}$ extending a given oriented planar diagram are \textbf{strongly equivalent} if they differ by an isotopy of the twist graphic, round twist, crossing point moves or changes in the chambering graph.
\end{defin}
This equivalence relation, as the name suggests, will be strictly stronger than usual equivalence of framed planar diagrams, as we defined it above. The latter models diffeomorphism-isotopy classes of surfaces and since it is perfectly possible for the same surface with two non-isotopic framings to still be diffeomorphic-isotopic to itself, one can thus find two equivalent framed planar diagrams with the same underlying oriented planar diagram which are not strongly equivalent. 

This additional result on strong equivalence is not needed for the presentation of the framed bordism bicategory, but we would need to undertake a similar study anyway, so there is no reason not to perform it in this generality. The conclusion is that framed planar diagrams under this finer equivalence relation can be used to model framings on a fixed surface, which we find interesting in its own right.

A minor nuisance is that the chambering graph, which allows one to combinatorially describe sheet data, is not very useful for the purposes of the strong equivalence relation. The reason is that we are now really working over a fixed oriented surface mapping into $\mathbb{R}^{2}$ and any two oriented planar diagrams describing it will be "essentially the same" in the sense that their graphics coincide and the sets of sheets living over any point in the plane not lying on the graphic can be canonically identified with each other. However, the diagrams themselves can still be wildly different if one chooses different chambering graphs. 

The most elegant way to deal with this issue would be to redefine an oriented planar diagram to consist of a graphic together with a certain "sheaf of sheets" on the complement of it, together with some identifications along the arcs and singularities of the graphic. We would then say that two oriented planar diagrams are \emph{the same} if their graphics coincide and one is given an identification of their sheaves, this is the same as asking for an oriented diffeomorphism over $\mathbb{R}^{2}$ of the surfaces they describe. 

This would take us too far astray and we instead proceed by ignoring the issue altogether. In the derivation of the strong equivalence relation given below we - sometimes implicitly - assume that the chambering graph is allowed to change at all times and we will identify two planar diagrams if they have the same graphic and sheets over the relevant open subsets of $\mathbb{R}^{2}$ have been identified. In particular, we allowed the chambering graph to change in the definition of strong equivalence given above. We hope we will be forgiven for this abuse of notation and language.

We will now proceed to derive the strong equivalence relation, it will relate different framed planar diagrams over a given oriented planar diagram, mirroring the isotopy relation on framings on the surface described. 

Since a framed planar diagram always comes from a diagrammable framing, our task would be the easiest if not only any framing was isotopic to a diagrammable one, which we proved, but also any diagrammable framings that were isotopic would be isotopic through a family of other diagrammable framings. We could then simply conclude that two isotopic diagrammable framings give framed planar diagrams which differ at most by an isotopy of the twist graphic. This is not the case. However, to some extent we will be able to measure a failure of an isotopy of framings to be homotopic to a pointwise diagrammable one.

Recall that a framing was diagrammable if it had a standard form in the preimage of small rectangular neighbourhoods of singularities of the graphic of critical values and was sufficiently generic outside of it. A part of the proof that any framing is isotopic to diagrammable one was to observe that the preimages of these rectangular neighbourhoods were disjoint unions of disks. We then used that any two framings on a disk are isotopic, this follows from the fact that the space of all such framings can be identified with

\begin{center}
$Map(D^{2}, SO(2)) \simeq SO(2)$,
\end{center}
which happens to be path-connected. However, the same computation shows that the space of self-isotopies of any framing on a disk can be identified with 

\begin{center}
$\Omega SO(2) \simeq \mathbb{Z}$,
\end{center}
which is \emph{not} path-connected. It follows that there are some self-isotopies of the standard framing on a preimage of a rectangular neighbourhood that cannot be deformed to be constant. Hence, any isotopy between diagrammable framings that restricts to such a non-trivial isotopy on the preimage cannot be deformed to an isotopy of diagrammable framings. We can measure this failure precisely by introducing the notion of \emph{degree}.

\begin{defin}
Let $W$ be an oriented surface together with a generic map $f: W \rightarrow \mathbb{R}^{2}$, let $h: v \rightarrow v^\prime$ be an isotopy of diagrammable framings of $W$. If $x$ is a singularity of either Morse or cusp type in the induced graphic $\Psi$, then we define the \textbf{degree} of $h$ along $x$ to be its isotopy class when restricted to $y$, considered as an element of 
\begin{center}
$h |_{\{ y \}} \in \pi_{0} (\Omega Map(\{ y \}, SO(2))) \simeq \pi_{1}(SO(2)) \simeq \mathbb{Z}$,
\end{center}
where $y$ is the unique preimage of $x$ lying on the non-trivial sheet. If $x$ is a fold crossing singularity, we define the degree to be its isotopy class when restricted to $y$ and $y^\prime$, considered as an element of

\begin{center}
$h |_{\{ y, y^\prime \}} \in \pi_{0}(\Omega Map(\{ y, y^\prime \}, SO(2))) \simeq \pi_{1}(SO(2) \times SO(2)) \simeq \mathbb{Z} \oplus \mathbb{Z}$,
\end{center}
where $y, y^\prime$ are the two preimages of $x$ lying on the two components of the non-trivial sheet.
\end{defin}

\begin{rem}
Observe that any self-isotopy of framings on a disk can be deformed to one that is constant on the boundary. This follows from the homotopy extension property applied to the pair $\partial D^{2} \subseteq D^{2}$ and the fact that there are no interesting self-isotopies in the space of framings over $\partial D^{2} \simeq S^{1}$, as the latter is homeomorphic to $Map(S^{1}, SO(2))$ and so is componentwise contractible.
\end{rem}
This idea that we can measure the obstruction for an isotopy to be homotopic to one constant along some subspace is at the crux of the relations we introduce. Before discussing them in detail, let us give a general idea of the derivation. Recall that our goal right now is to find relations under which two framed planar diagrams with the same underlying oriented planar diagram describe isotopic framings on the reconstructed surface.

We will first introduce a certain class of relations which we will call \emph{round twist}. We will show that using just these moves one can for any framed planar diagram and any prescription of degrees for each singularity of its graphic of critical values, find a strongly equivalent framed planar diagram, such that the two framings they describe are isotopic by an isotopy realizing the prescribed degrees. 

This is a step forward, as an argument about \emph{cancellation} of degrees - they lie in an abelian group, after all - will show that after introducing these it's enough to describe relations that will imply that any two diagrammable framings which are related by an isotopy constant along preimages of singularities lead to equivalent framed planar diagrams.

We then similarly study the behaviour of isotopies along arcs of critical values lying on the surface. We again introduce local moves, which we call \emph{crossing point relations}, that ensure that any such isotopy - up to homotopy of isotopies - can be realized as an isotopy relating two strongly equivalent framed planar diagrams. A similar cancellation argument then allows us to reduce to the case of diagrammable framings related by an isotopy constant along some neighbourhood of all critical values of the fixed generic map. This is a case we attack directly using isotopy and so called moves of \emph{twist homotopy}.

The \textbf{round twist} local move relations come in two kinds for each of the Morse and cusp singularities and four kinds for the fold crossing, the kinds correspond to generators of the group of self-isotopies of the standard framing on the non-trivial sheet associated to the given singularity. This group is $\mathbb{Z}$ for Morse and cusp singularities and $\mathbb{Z} \oplus \mathbb{Z}$ for the fold crossing, as we have computed when defining the degree, so as a group it is generated by one, resp. two elements, but we are interested in their generators as a monoid, which for us will be $\pm 1$ and $(\pm 1, \pm 1)$. 

The relations are given in \textbf{Figure \ref{fig:round_twist_relations_for_framed_planar_diagrams}} for Morse singularities, in \textbf{Figure \ref{fig:cusp_round_twist_relations_for_framed_planar_diagrams}} for cusps and in \textbf{Figure \ref{fig:fold_crossing_round_twist_relations_for_framed_planar_diagrams}} for fold crossings. Note that the correspondence between generators and these pictures is not one-to-one, as some pictures can be interpreted in two different ways depending on the labellings, we explain how to do this below.

\begin{rem}
A diligent reader might observe that $\mathbb{Z} \oplus \mathbb{Z}$ needs only three generators as a monoid, for example $(1, 0), (0, 1)$ and $(-1, -1)$ would suffice. This means that we could exchange the four relations we present for the fold crossing for three different ones. However, our goal at this point is not necessarily to minimize the number of relations in itself, but rather minimize the overhead needed to remember them.
\end{rem}

\begin{figure}[htbp!]
	\begin{tikzpicture}
		\node (MarginNode) at (0, 2) {};
		
		\draw [dashed] (0.2,0) to (0.3, 1.8) to (2.8, 1.8) to (2.8,0) to (0.2,0);
		\draw [very thick, red] (0.5,0) [out=90, in=180] to (1.5,1.5) [out=0, in=90] to (2.5,0);
		\node [circle, fill=serre_blue, inner sep=1.5pt] at (1.5,1.5) {};
		\draw [twist, serre_blue] (1.5,1.5) to (1.5, 0);
		
		\node at (3.5,0.9) {$ \leftrightarrow $};
		
		\draw [dashed] (4.2,0) to (4.2, 1.8) to (6.8,1.8) to (6.8, 0) to (4.2, 0);
		\draw [very thick, red] (4.5,0) [out=90, in =180] to (5.5,1.5) [out=0, in=90] to (6.5,0);
		\node [circle, fill=serre_blue, inner sep=1.5pt] at (5.5,1.5) {};
		\draw [twist, serre_blue] (5.5,1.5) to [out=-90, in=60] (4.8, 1.1);
		\draw [twist, serre_green] (4.8, 1.1) to [out=-120, in=90] (6.4, 0.4);
		\draw [twist, serre_blue] (6.4, 0.4) to [out=-90, in=90] (5.5, 0);
	\end{tikzpicture}
\caption{One of the round twist relations for the Morse cap singularity}
\label{fig:twist_around_relation_for_morse_cap_singularity}
\end{figure}
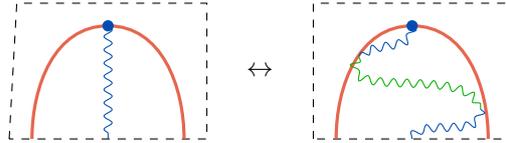

\begin{figure}[htbp!]
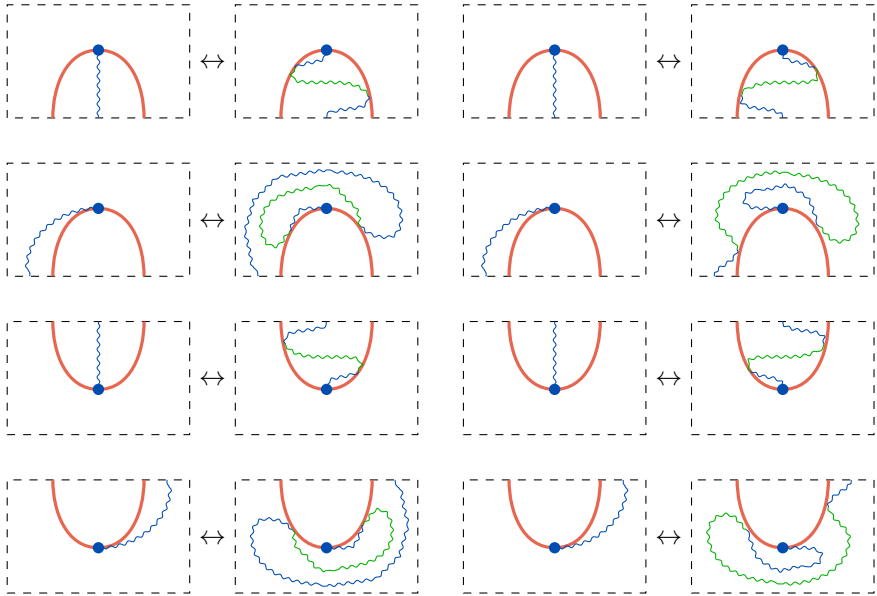


\caption{Round twist moves for Morse singularities}
\label{fig:round_twist_relations_for_framed_planar_diagrams}
\end{figure}

\begin{figure}[htbp!]
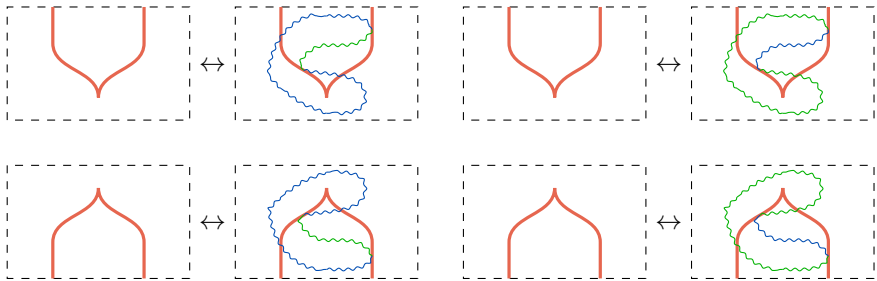

	%
\caption{Round twist moves for cusp singularities}
\label{fig:cusp_round_twist_relations_for_framed_planar_diagrams}
\end{figure}

\begin{figure}[htbp!]
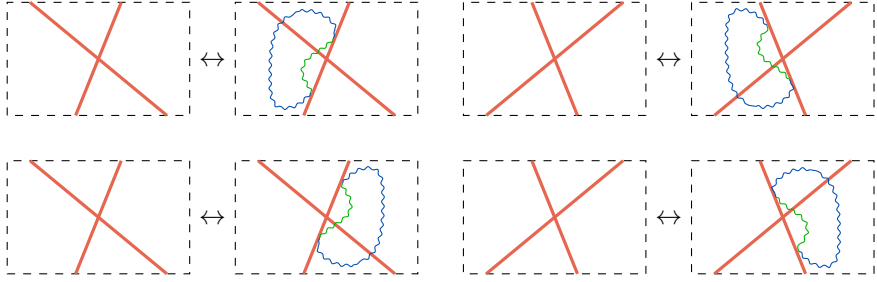

	%
\caption{Round twist moves for fold crossing singularities}
\label{fig:fold_crossing_round_twist_relations_for_framed_planar_diagrams}
\end{figure}

The relations are presented as local moves, they are analogous to the standard local moves one can perform on unoriented planar diagrams, as explained in \cite{chrisphd}, although we describe how to interpret them below. For concreteness, consider one of the round twist relations associated to Morse singularity, pictured below in \textbf{Figure \ref{fig:twist_around_relation_for_morse_cap_singularity}}.

A local move relation should be understood in the following way: Whenever we have a framed planar diagram and a fragment of its graphic which looks like either of the sides, then one can modify the diagram by replacing this fragment by the other side of the relation, without changing anything else outside the dashed region, and obtain an equivalent diagram. 

Hence, a sanity check for this type of relation is whether both sides of the move look the same "on the boundary", otherwise after performing it we might wind up with something which is no longer a framed graphic. In round twist for Morse cap, for example, on both sides we have a fragment which from the bottom ends with two arcs of the graphic of critical values and one arc of the twist graphic. Similarly, on both sides there are no arcs joining from either left, right or the top.

It could be that some changes to sheet data are necessary when performimg such a local move, but the relations we describe now are part of the strong equivalence, so the underlying graphic of critical values and oriented sheet data is assumed not to change.

Note that a fragment of twist graphic coming from a framed planar diagram will also come with some labellings as part of framed sheet data, for example each arc has an associated sheet and is normally oriented in either straightforwad or inverted way. When presenting the relations as local moves, we make a convention not to draw these labellings. The relations presented are then assumed to hold for each possible type of labelling which is not ruled out by some compatibility conditions present already between the pieces of the fragment.

For example, on the left hand side of the presented round twist relation for Morse cap, the two arcs of critical values are either left or right elbows and the arc of twist graphic can be either straightforward or inverted. However, this type of singularity of a framed planar diagram, when a downward cap of graphic of critical values is joined from the bottom by an arc of the twist graphic, is only possible if the arcs are labeled as to make a Morse cap and the arc of twist graphic is straightforward.

We can use this information to figure out what are the necessary labellings on the right hand side, the blue twist arc at the bottom has to be straightforward too, as the two sides of the relation agree on the boundary. Hence, by the compatibility conditions of a crossing point singularity the green arc has to be inverted and the top blue arc straightforward, as normal orientation of an arc of twist graphic changes when we pass a crossing point. 

On the other hand, any drawing from \textbf{Figure \ref{fig:cusp_round_twist_relations_for_framed_planar_diagrams}} or \textbf{Figure \ref{fig:fold_crossing_round_twist_relations_for_framed_planar_diagrams}} can be interpreted in two different ways, depending on whether the blue twist graphic is taken to be straightforward or inverted. Both moves are valid.

The reason we introduce precisely these relations is because of the following proposition, to which we aluded before.

\begin{prop}
\label{prop:existence_of_isotopies_with_prescribed_degrees}
Let $\aleph = (\Psi, \Upsilon, \Gamma)$ be a framed planar diagram, let $d_{x}$ be a family of possible degrees of an isotopy between diagrammable framings indexed by singularities of the graphic of critical values $\Psi$. 

Then, there exists a framed planar diagram $\aleph ^\prime$ differing from $\aleph$ only by round twist moves, hence in particular describing the same oriented surface and strongly equivalent, such that the two framings on the reconstructed surface corresponding to these diagrams are isotopic by an isotopy realizing the prescribed degrees.
\end{prop}

\begin{proof}
As degree is additive with respect to concatenation of isotopies, it is enough to prove this in the situation where the prescribed degree vanishes for all but one singularity $x \in \Psi _{0}$ and $d_{x}$ is a generator.

For each singularity type and for each generator, we have a round twist relation with only this type of singularity on the left hand side and the same singularity on the right hand side, but with a different framing, differing from the original by an isotopy representing the chosen generator. 

We obtain $\aleph ^\prime$ from $\aleph$ by applying the corresponding relation, from the left hand side to the right. This is always possible, as a small neighbourhood of any such codimension $2$ singularity must have a neighbourhood looking exactly like the left hand side.

This does not alter the graphic of critical values or the oriented sheet data, one also does not need to change the chambering graph. The isotopy from the framing represented by $\aleph$ to the one represented by $\aleph ^\prime$ is given by an isotopy representing the needed generator on the non-trivial sheet associated the chosen singularity. It can be chosen to be constant on the boundary and so it extends to an isotopy defined on the whole surface.
\end{proof}

\begin{rem}
The proof explains why we needed to have round twist relations indexed by generators of the corresponding group \emph{as a monoid}. If passing from the left hand side to the right hand side of a round twist relation realizes some generator degree, then passing from the right hand side to the left realizes its negation. However, the left hand side is special in that one can always find a small neighbourhood of singularity with exactly that configuration of the graphic, this is not necessarily true for the right hand side. 
\end{rem}

The \textbf{crossing points} moves come in three types, that of \emph{color shift}, \emph{birth} and \emph{death}, but have variations depending on whether we have a left or right elbow and on the color of the arcs of the twist graphic. We present all of them in \textbf{Figure \ref{fig:crossing_points_moves_color_shift}} and \textbf{Figure \ref{fig:crossing_points_moves_birth_death}} below. 

\begin{figure}[htbp!]
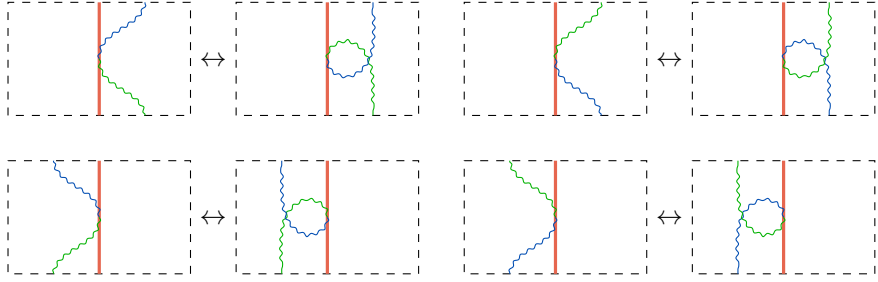


\caption{Crossing points moves, color shift}
\label{fig:crossing_points_moves_color_shift}
\end{figure}

\begin{figure}[htbp!]
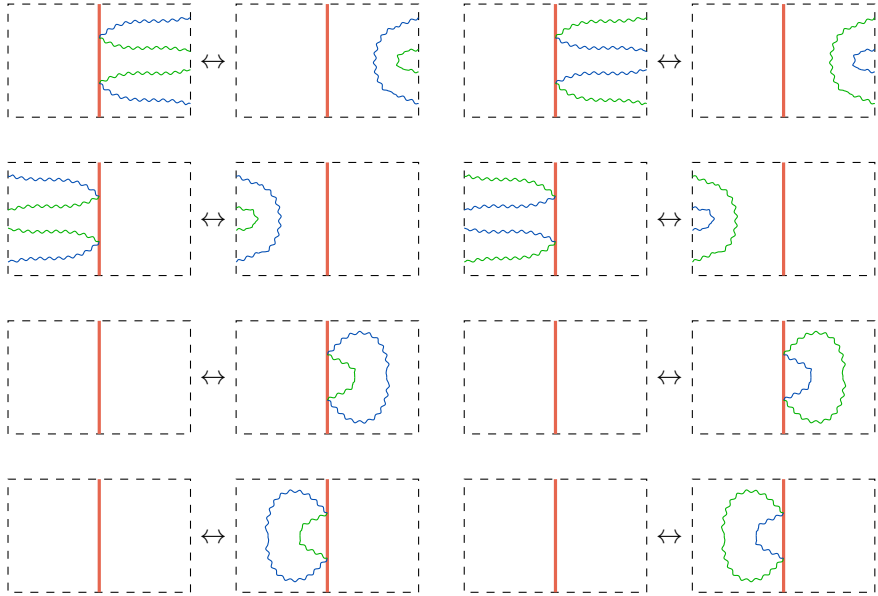

	%
\caption{Crossing points moves, birth-death}
\label{fig:crossing_points_moves_birth_death}
\end{figure}
The idea that leads one to introduce them is the study of isotopies of framings restricted to a single arc of a graphic of critical values, which corresponds to an arc of critical points in the oriented surface under discussion. 

The ends of the arc are "grounded" in some singularities of the graphic of critical values, if one chooses small rectangular neighbourhoods around them, then the complement of these neighbourhoods in the arc can be parametrized by the interval $I$. Any framing on the parametrized part of the arc can now be described as a function $I \rightarrow SO(2)$ by comparing it to ambient framing. If the framing in question is described by the framed graphic, then the function can be even taken to take $\partial I$ to $id \in SO(2)$, as under the parametrization ends of the interval lie on an a part of an arc disjoint from the twist graphic, by the required configuration of the graphic in the chosen small rectangular neighbourhood.

Then, isotopies of such framings which are constant on the parametrized ends of the arc can be identified with homotopies of maps $(I, \partial I) \rightarrow (SO(2), \{ id \})$. Such maps, up to homotopy, are parametrized by an integer. More precisely, if we start with a function transversal to $-id \in SO(2)$ - which we may assume we do, as this is true for framings obtained by restricting a diagrammable framing from an oriented surface - then the homotopy class of a map $(I, \partial I) \rightarrow (SO(2), \{ id \})$ is given by the number of preimages of $-id$, counted with multiplicity.

It's easy to read off such preimages from the framed planar diagram in which the arc lies, these are exactly the crossing points. We have one preimage for each crossing point, but its multiplicity depends also on the normal orientation of the twist arc crossing the fold, that is, on whether the twist arc is straightforward or inverted. We formalize this in the following definition.

\begin{defin}
In a given framed planar diagram $\aleph$, a crossing point is called \textbf{positive} if of the two arcs meeting at it, the arc pointing downwards is straightforward and the arc pointing upwards is inverted. It is \textbf{negative} otherwise. A \textbf{crossing degree} $cd(A)$ of an arc $A$ of critical values in $\aleph$ is given by 

\begin{center}
$cd(A) = p_{A} - n_{A}$,
\end{center}
the number of positive crossing points lying on $A$ minus the number of the negative ones.
\end{defin}

\begin{rem}
Observe that in our convention on drawing framed planar diagrams, on either the left or right elbow there are still two types of positive (resp. negative) crossing points, as the downwards arc could be either blue or green, that is, it could lie on either a positive or a negative sheet. This, as it turns out, does not affect the behaviour of the framing on the arc itself.
\end{rem}

\begin{rem}
Using the ideas we presented, it is now relatively easy to write down two framed planar diagrams describing the same oriented surface such that the associated diagrammable framings on that surface are not isotopic through an isotopy of diagrammable framings. Namely, it's enough to have an arc that has a different crossing degree in each of the two framed planar diagrams. An isotopy of diagrammable framings would be constant in some small rectangular neighbourhoods of the singularities and hence its restriction to the arc would give a homotopy of maps $(I, \partial I) \rightarrow (SO(2), \{ id \})$ by the discussion above. However, if the crossing degrees are different, the two maps will not be homotopic in the first place.
\end{rem}

\begin{prop}
\label{prop:deformation_of_framed_planar_diagrams_in_neighbourhood_of_an_arc}
Let $\aleph$ and $\aleph ^\prime$ be two framed planar diagrams with the same underlying oriented planar diagram. If $A$ is an arc of the graphic of critical values which has the same crossing degree in both of the diagrams, then one can modify $\aleph$ using only isotopy of the twist graphic, crossing points moves and possibly changes in the chambering graph, such that both diagrams will agree in some neighbourhood of $A$.
\end{prop}

\begin{proof}
After parametrizing the arc, we can identify it with the interval $I$ and the crossing points on it with a configuration of positive and negative points. Each positive or negative crossing point can additionally come in two types, but we can freely change from one type to the other using color shift moves. 

Since the total number of crossing points, counted with multiplicity, is the same in both diagrams, we can go from one to the other by isotopy or simultaneous creation or annihilation of pairs of opposite points. This can be done using birth-death crossing points moves, by passing from the left hand side to the right.
\end{proof}

The last type of relations needed to derive the strong equivalence relation will be that of moves of \textbf{twist homotopy}. These can be further divided into two kinds, the \textbf{topological} and \textbf{positional} ones.

They all describe changes that the twist graphic can undergo under a general isotopy, that is, one that is not necessarily part of an isotopy of framed planar diagrams. The topological ones describe a situation where the topology of the twist graphic, seen as a diagram, changes, these lead to actual relations in the framed bordism bicategory. 

The latter kind, the positional relations, describe a change of the position of the twist graphic relative to itself, the graphic of critical values or the chambering graph. This type of relations is very numerous, but rather straightforward in principle. The possible variants of positional moves are as follows.

\begin{itemize}
\item an arc of the twist graphic passes through a singularity of the graphic of critical values, the twist arc lies on a different sheet than the sheets being folded together over all the arcs of the graphic of critical values meeting at the given singularity
\item a twist birth-death singularity passes through an arc of critical values, the two sheets being folded together over it are different than the sheets on which the two twist arcs lie
\item a crossing point singularity passes through an arc of critical values, the two pairs of sheets being folded together over the two relevant arcs are disjoint from each other
\item an arc of the twist graphic passes through an arc of critical values, forming two transversal intersections, it lies on a different sheet than the sheets being folded together over that arc
\item an arc of the twist graphic passes through a vertex of the chambering graph
\item a twist birth-death singularity passes through an edge of the chambering graph
\item a crossing point singularity passes through an edge of the chambering graph
\item an arc of the twist graphic passes through an edge of the chambering graph, forming two transversal intersections
\item an arc of the twist graphic passes through a twist birth-death singularity, it lies on a different sheet then both of the arcs involved in the singularity
\item an arc of the twist graphic passes through a crossing point singularity, it lies on a different sheet then the sheets being folded together over the relevant arc of critical values
\item an arc of the twist graphic passes through a twist arc lying on a different sheet, forming two transversal intersections
\end{itemize}
Note that a general rule is just that things are allowed to pass through each other in a generic way exactly when they describe phenomena happening on disjoint sets of sheets. An exception is given by the moves involving the chambering graph, which is always allowed to pass through the twist graphic, as the latter only describes something that happens locally on one sheet and so no compatibility conditions are needed.

Geometrically, as we will see when we prove that these relations are sufficient for our purposes, they describe a situation where the Serre submanifold undergoes an isotopy and its image, the twist graphic, passes through itself, the graphic of critical values or the chambering graph, even though the Serre submanifold did not cross anything, as its arcs were moving on different sheets. This "different sheet" condition is usually implicitly present in our graphical representation of the local moves, enforced by compatibility properties of framed planar diagrams.

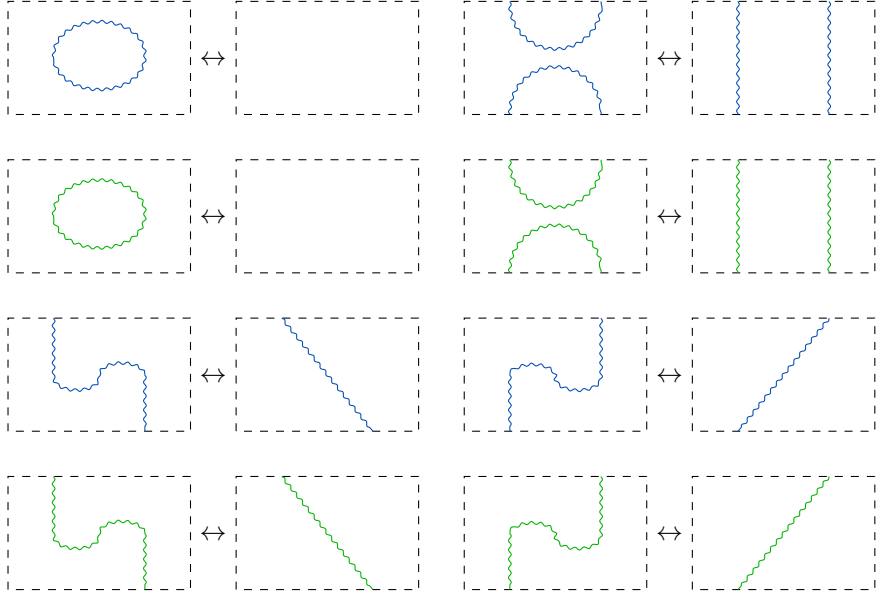
\begin{figure}[htbp!]
	\begin{tikzpicture}[scale=0.6]
	
		\begin{scope}
			\draw[dashed] (0, 0) to (0,2.5) to (4,2.5) to (4,0) to (0,0);
			
			\draw[small_twist, serre_blue] (1, 1.25) [out=90, in=180] to (2, 2) [out=0, in=90] to (3, 1.25) [out=-90, in=0] to (2, 0.5) [out=180, in=-90] to (1, 1.25);
			
			\node at (4.5, 1.2) {$ \leftrightarrow $};
			
			\draw[dashed] (5, 0) to (5,2.5) to (9,2.5) to (9,0) to (5,0);
		\end{scope}
		
		\begin{scope}[xshift=10cm]
			\draw[dashed] (0, 0) to (0,2.5) to (4,2.5) to (4,0) to (0,0);
			
			\draw[small_twist, serre_blue] (1, 2.5) [out=-90, in=180] to (2, 1.5) [out=0, in=-90] to (3, 2.5);
			\draw[small_twist, serre_blue] (1, 0) [out=90, in=180] to (2, 1) [out=0, in=90] to (3, 0); 
			
			\node at (4.5, 1.2) {$ \leftrightarrow $};
			
			\draw[dashed] (5, 0) to (5,2.5) to (9,2.5) to (9,0) to (5,0);
			\draw[small_twist, serre_blue] (6, 2.5) to (6, 0);
			\draw[small_twist, serre_blue] (8, 2.5) to (8,0);
		\end{scope}
	
		\begin{scope}[yshift=-3.5cm]
			\draw[dashed] (0, 0) to (0,2.5) to (4,2.5) to (4,0) to (0,0);
			
			\draw[small_twist, serre_green] (1, 1.25) [out=90, in=180] to (2, 2) [out=0, in=90] to (3, 1.25) [out=-90, in=0] to (2, 0.5) [out=180, in=-90] to (1, 1.25);
			
			\node at (4.5, 1.2) {$ \leftrightarrow $};
			
			\draw[dashed] (5, 0) to (5,2.5) to (9,2.5) to (9,0) to (5,0);
		\end{scope}
		
		\begin{scope}[xshift=10cm, yshift=-3.5cm]
			\draw[dashed] (0, 0) to (0,2.5) to (4,2.5) to (4,0) to (0,0);
			
			\draw[small_twist, serre_green] (1, 2.5) [out=-90, in=180] to (2, 1.5) [out=0, in=-90] to (3, 2.5);
			\draw[small_twist, serre_green] (1, 0) [out=90, in=180] to (2, 1) [out=0, in=90] to (3, 0); 
			
			\node at (4.5, 1.2) {$ \leftrightarrow $};
			
			\draw[dashed] (5, 0) to (5,2.5) to (9,2.5) to (9,0) to (5,0);
			\draw[small_twist, serre_green] (6, 2.5) to (6, 0);
			\draw[small_twist, serre_green] (8, 2.5) to (8,0);
		\end{scope}
		
		\begin{scope}[yshift=-7cm]
			\draw[dashed] (0, 0) to (0,2.5) to (4,2.5) to (4,0) to (0,0);
			
			\draw[small_twist, serre_blue] (1, 2.5) to (1, 1.25) [out=-90, in=-90] to (2, 1.25) [out=90, in=90] to (3, 1.25);
			\draw[small_twist, serre_blue] (3, 1.25) to (3, 0); 
			
			\node at (4.5, 1.2) {$ \leftrightarrow $};
			
			\draw[dashed] (5, 0) to (5,2.5) to (9,2.5) to (9,0) to (5,0);
			\draw[small_twist, serre_blue] (6, 2.5) to (8, 0);
		\end{scope}
		
		\begin{scope}[yshift=-7cm, xshift=10cm]
			\draw[dashed] (0, 0) to (0,2.5) to (4,2.5) to (4,0) to (0,0);
			
			\draw[small_twist, serre_blue] (3, 2.5) to (3, 1.25) [out=-90, in=-90] to (2, 1.25) [out=90, in=90] to (1, 1.25);
			\draw[small_twist, serre_blue] (1, 1.25) to (1, 0); 
			
			\node at (4.5, 1.2) {$ \leftrightarrow $};
			
			\draw[dashed] (5, 0) to (5,2.5) to (9,2.5) to (9,0) to (5,0);
			\draw[small_twist, serre_blue] (8, 2.5) to (6, 0);
		\end{scope}
		
		\begin{scope}[yshift=-10.5cm]
			\draw[dashed] (0, 0) to (0,2.5) to (4,2.5) to (4,0) to (0,0);
			
			\draw[small_twist, serre_green] (1, 2.5) to (1, 1.25) [out=-90, in=-90] to (2, 1.25) [out=90, in=90] to (3, 1.25);
			\draw[small_twist, serre_green] (3, 1.25) to (3, 0); 
			
			\node at (4.5, 1.2) {$ \leftrightarrow $};
			
			\draw[dashed] (5, 0) to (5,2.5) to (9,2.5) to (9,0) to (5,0);
			\draw[small_twist, serre_green] (6, 2.5) to (8, 0);
		\end{scope}
		
		\begin{scope}[yshift=-10.5cm, xshift=10cm]
			\draw[dashed] (0, 0) to (0,2.5) to (4,2.5) to (4,0) to (0,0);
			
			\draw[small_twist, serre_green] (3, 2.5) to (3, 1.25) [out=-90, in=-90] to (2, 1.25) [out=90, in=90] to (1, 1.25);
			\draw[small_twist, serre_green] (1, 1.25) to (1, 0); 
			
			\node at (4.5, 1.2) {$ \leftrightarrow $};
			
			\draw[dashed] (5, 0) to (5,2.5) to (9,2.5) to (9,0) to (5,0);
			\draw[small_twist, serre_green] (8, 2.5) to (6, 0);
		\end{scope}
	\end{tikzpicture}
\caption{Twist homotopy moves, topological}
\label{fig:twist_homotopy_moves_same_topological}
\end{figure}

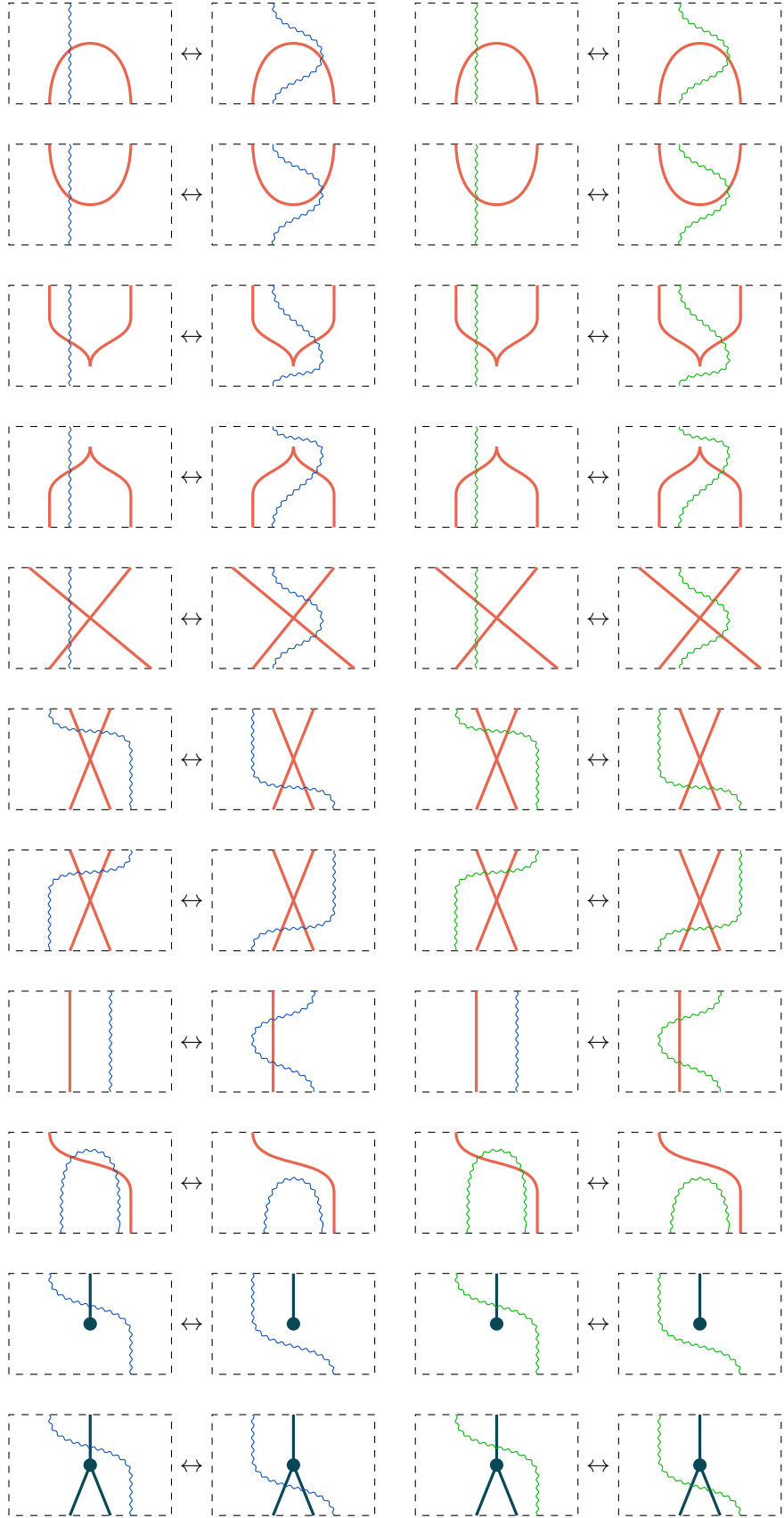
\begin{figure}[htbp!]
	\begin{tikzpicture}[scale=0.6]
	
		\begin{scope}
			\draw[dashed] (0, 0) to (0,2.5) to (4,2.5) to (4,0) to (0,0);
			\draw[very thick, red] (1, 0) [out=90, in=180] to (2, 1.5) [out=0, in=90] to (3, 0);
			\draw[small_twist, serre_blue] (1.5, 2.5) to (1.5, 0);
			
			\node at (4.5, 1.2) {$ \leftrightarrow $};
			
			\draw[dashed] (5, 0) to (5,2.5) to (9,2.5) to (9,0) to (5,0);
			\draw[very thick, red] (6, 0) [out=90, in=180] to (7, 1.5) [out=0, in=90] to (8, 0);
			\draw[small_twist, serre_blue] (6.5, 2.5) [out=-90, in=90] to (7.7, 1.25) [out=-90, in=90] to (6.5, 0);
		\end{scope}
		
		\begin{scope}[xshift=10cm]
			\draw[dashed] (0, 0) to (0,2.5) to (4,2.5) to (4,0) to (0,0);
			\draw[very thick, red] (1, 0) [out=90, in=180] to (2, 1.5) [out=0, in=90] to (3, 0);
			\draw[small_twist, serre_green] (1.5, 2.5) to (1.5, 0);
			
			\node at (4.5, 1.2) {$ \leftrightarrow $};
			
			\draw[dashed] (5, 0) to (5,2.5) to (9,2.5) to (9,0) to (5,0);
			\draw[very thick, red] (6, 0) [out=90, in=180] to (7, 1.5) [out=0, in=90] to (8, 0);
			\draw[small_twist, serre_green] (6.5, 2.5) [out=-90, in=90] to (7.7, 1.25) [out=-90, in=90] to (6.5, 0);
		\end{scope}

		\begin{scope}[yshift=-3.5cm]
			\begin{scope}[yshift=2.5cm, yscale=-1]
			\draw[dashed] (0, 0) to (0,2.5) to (4,2.5) to (4,0) to (0,0);
			\draw[very thick, red] (1, 0) [out=90, in=180] to (2, 1.5) [out=0, in=90] to (3, 0);
			\draw[small_twist, serre_blue] (1.5, 2.5) to (1.5, 0);
			\end{scope}
			
			\node at (4.5, 1.2) {$ \leftrightarrow $};
			
			\begin{scope}[yshift=2.5cm, yscale=-1]
			\draw[dashed] (5, 0) to (5,2.5) to (9,2.5) to (9,0) to (5,0);
			\draw[very thick, red] (6, 0) [out=90, in=180] to (7, 1.5) [out=0, in=90] to (8, 0);
			\draw[small_twist, serre_blue] (6.5, 2.5) [out=-90, in=90] to (7.7, 1.25) [out=-90, in=90] to (6.5, 0);
			\end{scope}
		\end{scope}
		
		\begin{scope}[xshift=10cm, yshift=-3.5cm]
			\begin{scope}[yshift=2.5cm, yscale=-1]
			\draw[dashed] (0, 0) to (0,2.5) to (4,2.5) to (4,0) to (0,0);
			\draw[very thick, red] (1, 0) [out=90, in=180] to (2, 1.5) [out=0, in=90] to (3, 0);
			\draw[small_twist, serre_green] (1.5, 2.5) to (1.5, 0);
			\end{scope}
			
			\node at (4.5, 1.2) {$ \leftrightarrow $};
			
			\begin{scope}[yshift=2.5cm, yscale=-1]
			\draw[dashed] (5, 0) to (5,2.5) to (9,2.5) to (9,0) to (5,0);
			\draw[very thick, red] (6, 0) [out=90, in=180] to (7, 1.5) [out=0, in=90] to (8, 0);
			\draw[small_twist, serre_green] (6.5, 2.5) [out=-90, in=90] to (7.7, 1.25) [out=-90, in=90] to (6.5, 0);
			\end{scope}
		\end{scope}
		
		\begin{scope}[yshift=-7cm]
			\draw[dashed] (0, 0) to (0,2.5) to (4,2.5) to (4,0) to (0,0);
			\draw[very thick, red] (1, 2.5) to (1, 1.7) [out=-90, in=90] to (2, 0.5) [out=90, in=-90] to (3, 1.7) to (3, 2.5);
			\draw[small_twist, serre_blue] (1.5, 2.5) to (1.5, 0);
			
			\node at (4.5, 1.2) {$ \leftrightarrow $};
			
			\draw[dashed] (5, 0) to (5,2.5) to (9,2.5) to (9,0) to (5,0);
			\draw[very thick, red] (6, 2.5) to (6, 1.7) [out=-90, in=90] to (7, 0.5) [out=90, in=-90] to (8, 1.7) to (8, 2.5);
			\draw[small_twist, serre_blue] (6.5, 2.5) [out=-90, in=90] to (7.7, 0.7) [out=-90, in=90] to (6.5, 0);
		\end{scope}
		
		\begin{scope}[xshift=10cm, yshift=-7cm]
			\draw[dashed] (0, 0) to (0,2.5) to (4,2.5) to (4,0) to (0,0);
			\draw[very thick, red] (1, 2.5) to (1, 1.7) [out=-90, in=90] to (2, 0.5) [out=90, in=-90] to (3, 1.7) to (3, 2.5);
			\draw[small_twist, serre_green] (1.5, 2.5) to (1.5, 0);
			
			\node at (4.5, 1.2) {$ \leftrightarrow $};
			
			\draw[dashed] (5, 0) to (5,2.5) to (9,2.5) to (9,0) to (5,0);
			\draw[very thick, red] (6, 2.5) to (6, 1.7) [out=-90, in=90] to (7, 0.5) [out=90, in=-90] to (8, 1.7) to (8, 2.5);
			\draw[small_twist, serre_green] (6.5, 2.5) [out=-90, in=90] to (7.7, 0.7) [out=-90, in=90] to (6.5, 0);
		\end{scope}
		
		\begin{scope}[yshift=-10.5cm]
			\begin{scope}[yshift=2.5cm, yscale=-1]
			\draw[dashed] (0, 0) to (0,2.5) to (4,2.5) to (4,0) to (0,0);
			\draw[very thick, red] (1, 2.5) to (1, 1.7) [out=-90, in=90] to (2, 0.5) [out=90, in=-90] to (3, 1.7) to (3, 2.5);
			\draw[small_twist, serre_blue] (1.5, 2.5) to (1.5, 0);
			\end{scope}
			
			\node at (4.5, 1.2) {$ \leftrightarrow $};
			
			\begin{scope}[yshift=2.5cm, yscale=-1]
			\draw[dashed] (5, 0) to (5,2.5) to (9,2.5) to (9,0) to (5,0);
			\draw[very thick, red] (6, 2.5) to (6, 1.7) [out=-90, in=90] to (7, 0.5) [out=90, in=-90] to (8, 1.7) to (8, 2.5);
			\draw[small_twist, serre_blue] (6.5, 2.5) [out=-90, in=90] to (7.7, 0.7) [out=-90, in=90] to (6.5, 0);
			\end{scope}
		\end{scope}
		
		\begin{scope}[xshift=10cm, yshift=-10.5cm]
			\begin{scope}[yshift=2.5cm, yscale=-1]
			\draw[dashed] (0, 0) to (0,2.5) to (4,2.5) to (4,0) to (0,0);
			\draw[very thick, red] (1, 2.5) to (1, 1.7) [out=-90, in=90] to (2, 0.5) [out=90, in=-90] to (3, 1.7) to (3, 2.5);
			\draw[small_twist, serre_green] (1.5, 2.5) to (1.5, 0);
			\end{scope}
			
			\node at (4.5, 1.2) {$ \leftrightarrow $};
			
			\begin{scope}[yshift=2.5cm, yscale=-1]
			\draw[dashed] (5, 0) to (5,2.5) to (9,2.5) to (9,0) to (5,0);
			\draw[very thick, red] (6, 2.5) to (6, 1.7) [out=-90, in=90] to (7, 0.5) [out=90, in=-90] to (8, 1.7) to (8, 2.5);
			\draw[small_twist, serre_green] (6.5, 2.5) [out=-90, in=90] to (7.7, 0.7) [out=-90, in=90] to (6.5, 0);
			\end{scope}
		\end{scope}
		
		\begin{scope}[yshift=-14cm]
			\draw[dashed] (0, 0) to (0,2.5) to (4,2.5) to (4,0) to (0,0);
			\draw[very thick, red] (0.5, 2.5) to (3.5, 0);
			\draw[very thick, red] (3, 2.5) to (1, 0);
			\draw[small_twist, serre_blue] (1.5, 2.5) to (1.5, 0);
			
			\node at (4.5, 1.2) {$ \leftrightarrow $};
			
			\draw[dashed] (5, 0) to (5,2.5) to (9,2.5) to (9,0) to (5,0);
			\draw[very thick, red] (5.5, 2.5) to (8.5, 0);
			\draw[very thick, red] (8, 2.5) to (6, 0);
			\draw[small_twist, serre_blue] (6.5, 2.5) [out=-90, in=90] to (7.7, 1.25) [out=-90, in=90] to (6.5, 0);
		\end{scope}	
		
		\begin{scope}[xshift=10cm, yshift=-14cm]
			\draw[dashed] (0, 0) to (0,2.5) to (4,2.5) to (4,0) to (0,0);
			\draw[very thick, red] (0.5, 2.5) to (3.5, 0);
			\draw[very thick, red] (3, 2.5) to (1, 0);
			\draw[small_twist, serre_green] (1.5, 2.5) to (1.5, 0);
			
			\node at (4.5, 1.2) {$ \leftrightarrow $};
			
			\draw[dashed] (5, 0) to (5,2.5) to (9,2.5) to (9,0) to (5,0);
			\draw[very thick, red] (5.5, 2.5) to (8.5, 0);
			\draw[very thick, red] (8, 2.5) to (6, 0);
			\draw[small_twist, serre_green] (6.5, 2.5) [out=-90, in=90] to (7.7, 1.25) [out=-90, in=90] to (6.5, 0);
		\end{scope}
	
		\begin{scope}[yshift=-17.5cm]
			\draw[dashed] (0, 0) to (0,2.5) to (4,2.5) to (4,0) to (0,0);
			\draw[very thick, red] (1.5, 2.5) to (2.5, 0);
			\draw[very thick, red] (2.5, 2.5) to (1.5, 0);
			\draw[small_twist, serre_blue] (1, 2.5) [out=-90, in=90] to (3, 1.5) [out=-90, in=90] to (3, 0);
			
			\node at (4.5, 1.2) {$ \leftrightarrow $};
			
			\draw[dashed] (5, 0) to (5,2.5) to (9,2.5) to (9,0) to (5,0);
			\draw[very thick, red] (6.5, 2.5) to (7.5, 0);
			\draw[very thick, red] (7.5, 2.5) to (6.5, 0);
			\draw[small_twist, serre_blue] (8,0) [out=90, in=-90] to (6, 1) [out=90, in=-90] to (6, 2.5);
		\end{scope}	
		
		\begin{scope}[xshift=10cm, yshift=-17.5cm]
			\draw[dashed] (0, 0) to (0,2.5) to (4,2.5) to (4,0) to (0,0);
			\draw[very thick, red] (1.5, 2.5) to (2.5, 0);
			\draw[very thick, red] (2.5, 2.5) to (1.5, 0);
			\draw[small_twist, serre_green] (1, 2.5) [out=-90, in=90] to (3, 1.5) [out=-90, in=90] to (3, 0);
			
			\node at (4.5, 1.2) {$ \leftrightarrow $};
			
			\draw[dashed] (5, 0) to (5,2.5) to (9,2.5) to (9,0) to (5,0);
			\draw[very thick, red] (6.5, 2.5) to (7.5, 0);
			\draw[very thick, red] (7.5, 2.5) to (6.5, 0);
			\draw[small_twist, serre_green] (8,0) [out=90, in=-90] to (6, 1) [out=90, in=-90] to (6, 2.5);
		\end{scope}		
		
		\begin{scope}[yshift=-21cm]
			\draw[dashed] (0, 0) to (0,2.5) to (4,2.5) to (4,0) to (0,0);
			\draw[very thick, red] (1.5, 2.5) to (2.5, 0);
			\draw[very thick, red] (2.5, 2.5) to (1.5, 0);
			\draw[small_twist, serre_blue] (3, 2.5) [out=-90, in=90] to (1, 1.5) [out=-90, in=90] to (1, 0);
			
			\node at (4.5, 1.2) {$ \leftrightarrow $};
			
			\draw[dashed] (5, 0) to (5,2.5) to (9,2.5) to (9,0) to (5,0);
			\draw[very thick, red] (6.5, 2.5) to (7.5, 0);
			\draw[very thick, red] (7.5, 2.5) to (6.5, 0);
			\draw[small_twist, serre_blue] (6,0) [out=90, in=-90] to (8, 1) [out=90, in=-90] to (8, 2.5);
		\end{scope}	
		
		\begin{scope}[xshift=10cm, yshift=-21cm]
			\draw[dashed] (0, 0) to (0,2.5) to (4,2.5) to (4,0) to (0,0);
			\draw[very thick, red] (1.5, 2.5) to (2.5, 0);
			\draw[very thick, red] (2.5, 2.5) to (1.5, 0);
			\draw[small_twist, serre_green] (3, 2.5) [out=-90, in=90] to (1, 1.5) [out=-90, in=90] to (1, 0);
			
			\node at (4.5, 1.2) {$ \leftrightarrow $};
			
			\draw[dashed] (5, 0) to (5,2.5) to (9,2.5) to (9,0) to (5,0);
			\draw[very thick, red] (6.5, 2.5) to (7.5, 0);
			\draw[very thick, red] (7.5, 2.5) to (6.5, 0);
			\draw[small_twist, serre_green] (6,0) [out=90, in=-90] to (8, 1) [out=90, in=-90] to (8, 2.5);
		\end{scope}	
		
		\begin{scope}[yshift=-24.5cm]
			\draw[dashed] (0, 0) to (0,2.5) to (4,2.5) to (4,0) to (0,0);
			\draw[very thick, red] (1.5, 2.5) to (1.5, 0);
			\draw[small_twist, serre_blue] (2.5, 2.5) to (2.5, 0);
			
			\node at (4.5, 1.2) {$ \leftrightarrow $};
			
			\draw[dashed] (5, 0) to (5,2.5) to (9,2.5) to (9,0) to (5,0);
			\draw[very thick, red] (6.5, 2.5) to (6.5, 0);
			\draw[small_twist, serre_blue] (7.5, 2.5) [out=-90, in=90] to (6, 1.35) [out=-90, in=90] to (7.5, 0);
		\end{scope}
		
		\begin{scope}[xshift=10cm, yshift=-24.5cm]
			\draw[dashed] (0, 0) to (0,2.5) to (4,2.5) to (4,0) to (0,0);
			\draw[very thick, red] (1.5, 2.5) to (1.5, 0);
			\draw[small_twist, serre_blue] (2.5, 2.5) to (2.5, 0);
			
			\node at (4.5, 1.2) {$ \leftrightarrow $};
			
			\draw[dashed] (5, 0) to (5,2.5) to (9,2.5) to (9,0) to (5,0);
			\draw[very thick, red] (6.5, 2.5) to (6.5, 0);
			\draw[small_twist, serre_green] (7.5, 2.5) [out=-90, in=90] to (6, 1.35) [out=-90, in=90] to (7.5, 0);
		\end{scope}
		
		\begin{scope}[yshift=-28cm]
			\draw[dashed] (0, 0) to (0,2.5) to (4,2.5) to (4,0) to (0,0);
			\draw[very thick, red] (1, 2.5) [out=-90, in=90] to (3, 1) to (3, 0);
			\draw[small_twist, serre_blue] (1.3, 0) [out=90, in=180] to (2, 2) [out=0, in=90] to (2.7, 0);
			
			\node at (4.5, 1.2) {$ \leftrightarrow $};
			
			\draw[dashed] (5, 0) to (5,2.5) to (9,2.5) to (9,0) to (5,0);
			\draw[very thick, red] (6, 2.5) [out=-90, in=90] to (8, 1) to (8, 0);
			\draw[small_twist, serre_blue] (6.3, 0) [out=90, in=180] to (7, 1.3) [out=0, in=90] to (7.7, 0);
		\end{scope}
		
		\begin{scope}[xshift=10cm, yshift=-28cm]
			\draw[dashed] (0, 0) to (0,2.5) to (4,2.5) to (4,0) to (0,0);
			\draw[very thick, red] (1, 2.5) [out=-90, in=90] to (3, 1) to (3, 0);
			\draw[small_twist, serre_green] (1.3, 0) [out=90, in=180] to (2, 2) [out=0, in=90] to (2.7, 0);
			
			\node at (4.5, 1.2) {$ \leftrightarrow $};
			
			\draw[dashed] (5, 0) to (5,2.5) to (9,2.5) to (9,0) to (5,0);
			\draw[very thick, red] (6, 2.5) [out=-90, in=90] to (8, 1) to (8, 0);
			\draw[small_twist, serre_green] (6.3, 0) [out=90, in=180] to (7, 1.3) [out=0, in=90] to (7.7, 0);
		\end{scope}
		
		\begin{scope}[yshift=-31.5cm]
			\draw[dashed] (0, 0) to (0,2.5) to (4,2.5) to (4,0) to (0,0);
			\draw[very thick, chambering_blue] (2, 2.5) to (2, 1.25);
			\node [circle, fill=chambering_blue, inner sep=2pt] at (2, 1.25) {};
			\draw[small_twist, serre_blue] (1, 2.5) [out=-90, in=90] to (3, 1) to (3, 0);
			
			\node at (4.5, 1.2) {$ \leftrightarrow $};
			
			\draw[dashed] (5, 0) to (5,2.5) to (9,2.5) to (9,0) to (5,0);
			\draw[very thick, chambering_blue] (7, 2.5) to (7, 1.25);
			\node [circle, fill=chambering_blue, inner sep=2pt] at (7, 1.25) {};
			\draw[small_twist, serre_blue] (6, 2.5) to (6, 1.5) [out=-90, in=90] to (8, 0);
		\end{scope}
		
		\begin{scope}[xshift=10cm, yshift=-31.5cm]
			\draw[dashed] (0, 0) to (0,2.5) to (4,2.5) to (4,0) to (0,0);
			\draw[very thick, chambering_blue] (2, 2.5) to (2, 1.25);
			\node [circle, fill=chambering_blue, inner sep=2pt] at (2, 1.25) {};
			\draw[small_twist, serre_green] (1, 2.5) [out=-90, in=90] to (3, 1) to (3, 0);
			
			\node at (4.5, 1.2) {$ \leftrightarrow $};
			
			\draw[dashed] (5, 0) to (5,2.5) to (9,2.5) to (9,0) to (5,0);
			\draw[very thick, chambering_blue] (7, 2.5) to (7, 1.25);
			\node [circle, fill=chambering_blue, inner sep=2pt] at (7, 1.25) {};
			\draw[small_twist, serre_green] (6, 2.5) to (6, 1.5) [out=-90, in=90] to (8, 0);
		\end{scope}
		
		\begin{scope}[yshift=-35cm]
			\draw[dashed] (0, 0) to (0,2.5) to (4,2.5) to (4,0) to (0,0);
			\draw[very thick, chambering_blue] (2, 2.5) to (2, 1.25);
			\node [circle, fill=chambering_blue, inner sep=2pt] at (2, 1.25) {};
			\draw[very thick, chambering_blue] (2, 1.25) to (1.5, 0);
			\draw[very thick, chambering_blue] (2, 1.25) to (2.5, 0);
			\draw[small_twist, serre_blue] (1, 2.5) [out=-90, in=90] to (3, 1) to (3, 0);
			
			\node at (4.5, 1.2) {$ \leftrightarrow $};
			
			\draw[dashed] (5, 0) to (5,2.5) to (9,2.5) to (9,0) to (5,0);
			\draw[very thick, chambering_blue] (7, 2.5) to (7, 1.25);
			\draw[very thick, chambering_blue] (7, 1.25) to (6.5, 0);
			\draw[very thick, chambering_blue] (7, 1.25) to (7.5, 0);
			\node [circle, fill=chambering_blue, inner sep=2pt] at (7, 1.25) {};
			\draw[small_twist, serre_blue] (6, 2.5) to (6, 1.5) [out=-90, in=90] to (8, 0);
		\end{scope}
		
		\begin{scope}[xshift=10cm, yshift=-35cm]
			\draw[dashed] (0, 0) to (0,2.5) to (4,2.5) to (4,0) to (0,0);
			\draw[very thick, chambering_blue] (2, 2.5) to (2, 1.25);
			\node [circle, fill=chambering_blue, inner sep=2pt] at (2, 1.25) {};
			\draw[very thick, chambering_blue] (2, 1.25) to (1.5, 0);
			\draw[very thick, chambering_blue] (2, 1.25) to (2.5, 0);
			\draw[small_twist, serre_green] (1, 2.5) [out=-90, in=90] to (3, 1) to (3, 0);
			
			\node at (4.5, 1.2) {$ \leftrightarrow $};
			
			\draw[dashed] (5, 0) to (5,2.5) to (9,2.5) to (9,0) to (5,0);
			\draw[very thick, chambering_blue] (7, 2.5) to (7, 1.25);
			\draw[very thick, chambering_blue] (7, 1.25) to (6.5, 0);
			\draw[very thick, chambering_blue] (7, 1.25) to (7.5, 0);
			\node [circle, fill=chambering_blue, inner sep=2pt] at (7, 1.25) {};
			\draw[small_twist, serre_green] (6, 2.5) to (6, 1.5) [out=-90, in=90] to (8, 0);
		\end{scope}

	\end{tikzpicture}
\caption{Some of the twist homotopy moves, positional}
\label{fig:twist_homotopy_moves_positional}
\end{figure}

The topological moves are presented in \textbf{Figure \ref{fig:twist_homotopy_moves_same_topological}}. The positional ones admit so many variations differing little from each other, that we instead present only some of them in \textbf{Figure \ref{fig:twist_homotopy_moves_positional}} and trust that the reader can derive all the others from the description we have given above. 

\begin{rem}
The moves presented in this section are not implied to be independant from each other and in fact it's fairly easy to find examples of moves which are consequences of others. We did not remove even the "most obvious redundancies", instead opting to be complete. This is another reason why the list of twist homotopy moves is very long. 
\end{rem}

\begin{thm}[Calculus of framings on an oriented surface]
\label{thm:description_of_framings_on_surface}
Let $W$ be an oriented surface together with a fixed generic map $f: W \rightarrow \mathbb{R}^{2}$. Then, two isotopic diagrammable framings of $W$ give framed planar diagrams which are strongly equivalent, in other words differ at most by an isotopy of the twist graphic, round twist, crossing points and twist homotopy moves and possibly changes in the chambering graph. 

Conversely, two strongly equivalent framed planar diagrams with the underlying oriented diagram the one induced by $f$ describe isotopic framings of $W$. In particular, there is a bijection between strong equivalence classes of framed planar diagrams with a fixed underlying oriented planar diagram and isotopy classes of framings on the reconstructed oriented surface.
\end{thm}

\begin{cor}
If $W$ is a framed surface equipped with a generic map into $\mathbb{R}^{2}$, then one can define the \textbf{associated framed planar diagram} by deforming the framing to be diagrammable and then applying the usual machinery. Up to strong equivalence of framed planar diagrams, so in particular up to equivalence, the result does not depend on the choice of the diagrammable deformation.
\end{cor}

\begin{proof}
We start with the second part, that the described relation preserves the isotopy class of described framing, as it is much easier.

Isotopy of the twist graphic will only change the embedding of the Serre submanifold by an isotopy. In particular, the Serre submanifolds corresponding to the beginning and ending configuration of that isotopy will be bordant in $\widetilde{W} \times I$, the bordism given by the isotopy. Hence, by the Thom-Pontryagin theorem, the associated functions into $S^{1}$ are homotopic.

The local moves described all replace a framing on a part of a surface by an isotopic one, this can be easily verified by hand. Alternatively, one can also observe that the moves only change the behaviour of the framing either on a disk or on a disjoint union of disks, so the framings on both sides of each move are isotopic by necessity. The isotopy can be taken to be constant on the boundary and so extends to a global isotopy of framings on $W$.

We now go to the interesting part. Suppose we have two isotopic framings $v_{0}, v_{1}$ of the oriented surface $W$, let $\aleph_{0}, \aleph_{1}$ be the associated framed planar diagrams, both with the same underlying oriented planar diagram by construction, so in particular agreeing on the graphic of critical values $\Psi$. 

Let $h: v_{0} \rightarrow v_{1}$ be the isotopy, let $d_{x}$ be its degrees, where $x \in \Psi_{0}$ runs through the singularities of the graphic of critical values. By \textbf{Proposition \ref{prop:existence_of_isotopies_with_prescribed_degrees}}, there exists some other framed planar diagram $\aleph _{1} ^\prime$ differing from $\aleph_{1}$ only by round twist moves and such that the framings described by these two diagrams differ by an isotopy realizing the family $-d_{x}$. Since the framing is additive under concatenation of isotopies, the framings described by $\aleph_{0}$ and $\aleph_{1} ^\prime$ differ by an isotopy of vanishing degree. Replacing  $\aleph_{1}$ by $\aleph_{1} ^\prime$ we may assume that the isotopy $h$ had this property to begin with.

Using the homotopy extension property we can now replace the isotopy $h$ by one which is constant on the non-trivial sheets over a small rectangular neighbourhood of each singularity. Its existence implies that for each arc $A$ of the underlying oriented planar diagram, its crossing degrees in $\aleph_{0}$ and $\aleph_{1}$ are equal. By repeated application of \textbf{Proposition \ref{prop:deformation_of_framed_planar_diagrams_in_neighbourhood_of_an_arc}} for each arc, we can replace $\aleph_{0}$ by $\aleph_{0} ^\prime$, agreeing with $\aleph_{1}$ on the neighbourhood of each arc, using only isotopy of the twist graphic, crossing points moves and possibly some changes in the chambering graph. Replacing $\aleph_{0}$ by $\aleph_{0} ^\prime$, we may assume that it agreed with $\aleph_{1}$ in some neighbourhood of all the arcs to begin with. They still differ only by an isotopy of vanishing degree, as $\aleph_{0}$ and $\aleph_{0} ^\prime$ did. 

Again, choose an isotopy $h$ which is constant in some neighbourhood of the critical preimages of each singularity of the graphic of critical values, it is then also up to homotopy constant on each arc. As $\aleph_{0}, \aleph_{1}$ agree in some neighbourhood of the arcs of critical values, we may assume that the two framings described agree in the neighbourhood of critical points. Shrinking it if necessary and again using the homotopy extension property, we can assume that the chosen homotopy $h$ is constant in this neighbourhood.

Let $\hat{W} \subseteq \widetilde{W}$ be the complement of some neighbourhood of critical points of the generic map $f$ in which $h$ is constant, small enough so that $h$ is constant even in some neighbourhood of $\partial \hat{W}$. By construction, $\hat{W}$ consists of a disjoint union of path-connected manifolds with boundary, whose interior is taken diffeomorphically to open subsets of $\mathbb{R}^{2}$ by $f$.

Since $\hat{W} \subseteq \widetilde{W}$, the framing on it can be identified with a function into $SO(2) \simeq S^{1}$ by comparing it to the ambient framing, $h$ then corresponds to a homotopy of such functions. Replacing this isotopy by a sufficiently generic one, we may assume it is transversal to $-id \in SO(2)$. Taking preimages, we then get a bordism inside $\hat{W} \times I$ between $\hat{S}_{1} = \hat{W} \cap S_{1}$, $\hat{S}_{2} = \hat{W} \cap S_{2}$, where $S_{1}, S_{2}$ are the Serre submanifolds of $v_{1}$ and $v_{2}$. 

By results of Schommer-Pries, any such $2$-dimensional bordism can be decomposed into cusps and Morse singularities. Its effect on the Serre submanifold is then either that of isotopy or of one of the topological twist homotopy moves.
 
The isotopy of $\hat{S} \subseteq \hat{W}$ induces an "isotopy" of the twist graphic, we put the latter in quotation marks, as it can be wild, that is, not a part of an isotopy of framed planar diagrams. For example, as the arcs of $\hat{S}$ move, their images in $\mathbb{R}^{2}$ could start forming new intersections, ones that are not allowed in framed planar diagrams, or their images could intersect some singularities of the underlying oriented planar diagram. This is possible, as over any point of $\mathbb{R}^{2}$ we can have many components of $\hat{W}$ mapping into it.

However, making $h$ sufficiently generic we can arrange it so that it induces an isotopy of the twist graphic which is part of an isotopy of framed planar diagrams except for a finite number of times, where one of the twist homotopy moves happens. The precise meaning of "sufficiently generic" would then be transversality conditions for that bordism with respect to preimages of edges of the chambering graph, preimages of the graphic of critical values, together with multijet transversality to ensure that the images of arcs of $\hat{S}$ living on different components of $\hat{W}$ are sufficiently generic with respect to each other. The derivation of such conditions together with a list of "possible moves" is formally analogous to the results leading to moves in \cite[Table 1.9, Figure 1.28]{chrisphd}, hence we do not give the details.

This gives a way to pass from the framed graphic of $\aleph_{0}$ to the one of $\aleph_{1}$ using only isotopy of the twist graphic and twist homotopy moves. The normal orientation of the arcs will coincide, as it is induced by the normal orientation of the Serre submanifold, which is preserved by the moves we have used. Then, the two framed planar diagrams we end up with can differ at most in the chambering graph, ending the proof. 
\end{proof}
The theorem gives us a graphical calculus with which one can describe framings on a fixed oriented surface, considered up to isotopy. This is an interesting result on its own right, but is not sufficient for the applications to the presentation of the framed bordism bicategory. For the latter, we need to be able to describe diffeomorphism-isotopy classes of surfaces and it is perfectly possible for the same surface equipped with two non-isotopic framings to be still diffeomorphic-isotopic to itself, as we have already remarked. 

In the language of expected properties of equivalence of diagrams we have given in the beginning of the section, the theorem above implies that is has property $(2)$, namely that the equivalence class of the associated framed planar diagram depends only on the isotopy class of the framing used to define it. We will now try to establish property $(3)$, which is that the equivalence type of the diagram is also independant of the choice of generic function $f$, which was up till now kept fixed.

Finding such a combinatorial equivalence relation is of course not easy even in the case of unoriented surfaces, where one would be instead asking for a diagrammatic calculus modelling diffeomorphism classes of surfaces. The relevant result, \cite[Theorem 1.52 (Planar Decomposition Theorem)]{chrisphd}, is a culmination of over seventy pages of careful study involving Cerf theory. We will see, however, that it is reasonably straightforward to build upon the techniques of Schommer-Pries and extend them to the framed case.

\begin{figure}[htbp!]
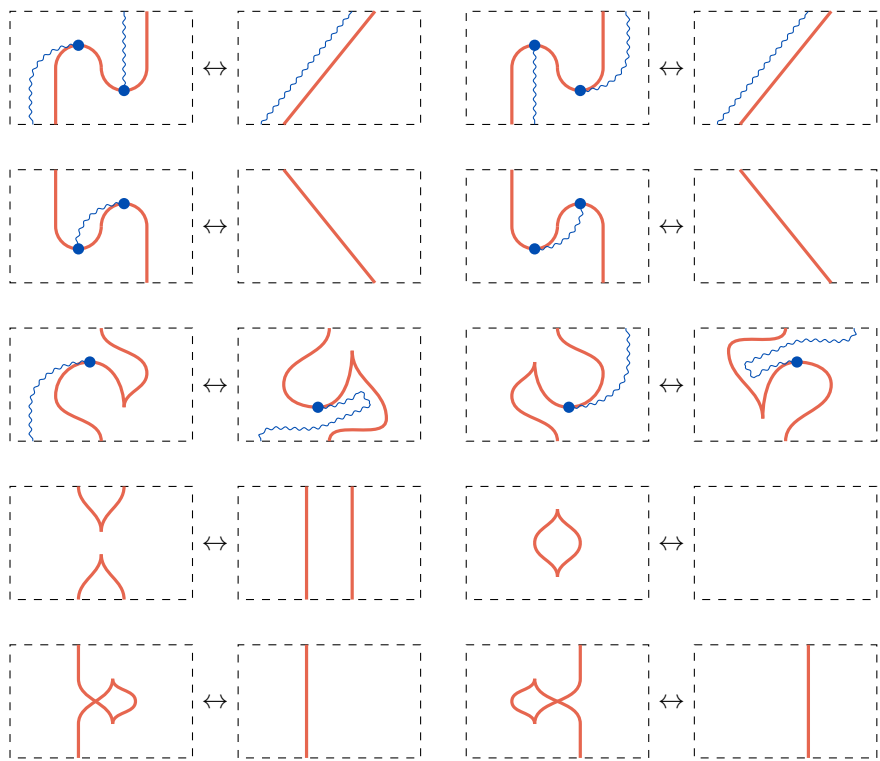


\caption{Change of function moves, topological}
\label{fig:change_of_function_moves_topological}
\end{figure}

\begin{figure}[htbp!]
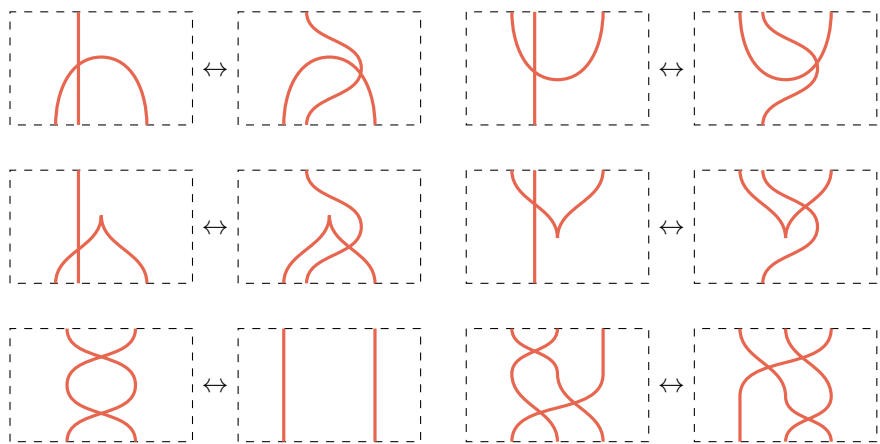

	%
\caption{Some of the change of function moves, positional}
\label{fig:change_of_function_moves_positional}
\end{figure}

\begin{figure}[htbp!]
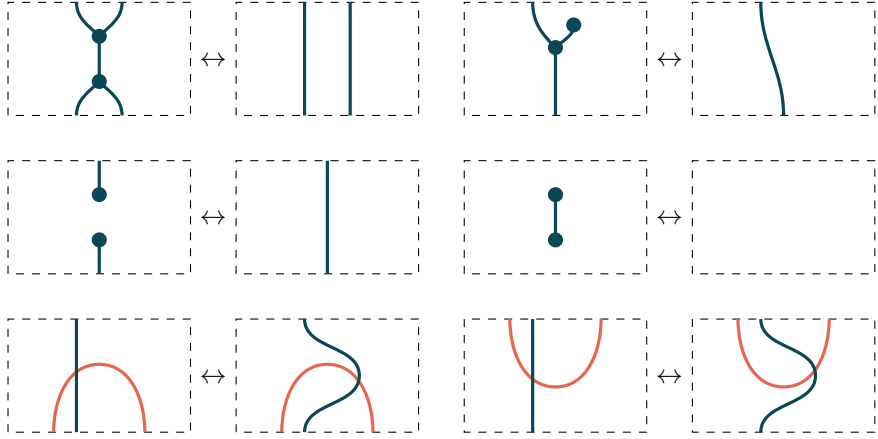

	%
\caption{Some of the change in chambering graph moves}
\label{fig:change_in_chambering_graph_moves}
\end{figure}

The new classes of moves we will introduce now are the \emph{change of function} moves, which one can again divide into \emph{topological} and \emph{positional}, and \emph{change in chambering graph} moves, presented in \textbf{Figure \ref{fig:change_of_function_moves_topological}}, \textbf{Figure \ref{fig:change_of_function_moves_positional}} and \textbf{Figure \ref{fig:change_in_chambering_graph_moves}}. Similarly to the case of twist homotopy, we present only some of the positional moves and only some of the moves allowed in the chambering graph, all the others are available in the original reference \cite{chrisphd}.

A fitting name for both of the groups at once would be \emph{moves of unoriented planar diagrams}, as what we do is we take the moves from \cite{chrisphd} and "lift" them to the framed case by adding the twist graphic where needed. A formal argument will then show us that these moves are sufficient.

Note that the moves, especially the changes in chambering graph, do actually require some compatibility conditions on sheet data. We do not discuss this in detail here, as all of them are already well explained in \cite{chrisphd} and the adaptation to the framed case is immediate.

Recall that two framed planar diagrams diagrams are equivalent if they differ by an isotopy and possibly a finite application of round twist, crossing points, twist homotopy, change of function or change in chambering graph moves. We have now given complete meaning to this definition and thus we are ready to establish the main theorems of the section. We start by showing that our relation has the property $(1)$ we discussed in the beginning, as it is a little easier and also sets up ground for the proof of property $(3)$. 

\begin{prop}[Invariance of reconstructed framed surfaces]
\label{prop:invariance_of_reconstructed_framed_surfaces}
If $\aleph_{0}$ and $\aleph_{1}$ are equivalent framed planar diagrams, then the framed surfaces one can reconstruct from them are diffeomorphic-isotopic. 
\end{prop}

\begin{rem}
Note that from a framed planar diagram we can reconstruct only an oriented surface together with an isotopy class of framings. This is of course enough to compare two such surfaces by a diffeomorphism-isotopy relation.
\end{rem}

\begin{proof}
Since $\aleph _{0}$ and $\aleph_{1}$ are equivalent, they differ by an isotopy and a finite application of local moves. Hence, it is enough to show the statement in the case where the given framed planar diagrams differ by only an isotopy and in the case where they differ by precisely one of the moves.

The latter is easier. The round twist, crossing points and twist homotopy moves do not affect the reconstructed oriented surface at all and only change the framing to an isotopic one, that was part of \textbf{Theorem \ref{thm:description_of_framings_on_surface}} on a calculus of framings on a fixed surfaces. Similarly, the change in chambering graph moves do not affect the reconstructed surface at all and are only visible in the diagram. 

The change of function moves do not alter the diffeomorphism type of the reconstructed surface but only modify its given map into $\mathbb{R}^{2}$. One can see by direct inspection that this modification of the generic function is always performed in some region that one can cover by disks, see \cite[Figures 1.10-1.13]{chrisphd}.  It follows that surfaces reconstructed from framed planar diagrams differing only by such a move are diffeomorphic and their framing can be taken to coincide outside of this region. As the region in contractible, up to isotopy the framing also coincides inside and so the two surfaces are diffeomorphic-isotopic.

Suppose now that $\aleph_{0}$ and $\aleph_{1}$ are isotopic through a family of framed planar diagrams $\aleph _{t}$. The totality of their graphics 

\begin{center}
$\Phi = \bigcup _{t \in [0, 1]} \Psi_{t} \times \{ t \} \subseteq \mathbb{R}^{2} \times I$
\end{center}
assembles to give a \emph{spatial graphic}, the oriented sheet data then induces a structure of an \emph{oriented spatial diagram}. The latter is a three-dimensional analogue of an oriented planar diagram, from $\Phi$ one can reconstruct a three-manifold mapping into $\mathbb{R}^{2} \times I$ in the same way as one reconstruct a surface mapping into $\mathbb{R}^{2}$ from a planar diagram, see \cite[Lemma 1.37]{chrisphd}. This three-manifold is diffeomorphic to a product $W \times I$, where $W$ is an oriented surface, such that the map $W \times I \rightarrow \mathbb{R}^{2} \times I$ is sufficiently generic in a suitable sense and induces the given spatial diagram. 

In our case, since the spatial diagram comes from an isotopy, we can assume that the constructed map $F: W \times I \rightarrow \mathbb{R}^{2} \times I$ is simply a \emph{path} of generic maps $f_{t}: W \rightarrow \mathbb{R}^{2}$, where $f_{t} = F(-, t)$. The underlying oriented planar diagram of $\aleph _{t}$ is then simply the diagram induced by $f_{t}$, although we are yet to establish any control over the framings.

(By saying it is a path of functions $W \rightarrow \mathbb{R}^{2}$ we mean that $\pi_{2} F(-, t) = t$ identically, in the general case of a spatial diagram which is perhaps not an isotopy of planar diagrams one can only assume that $\partial _{t} \pi_{2} F(-, t)$ is always strictly positive.)

This already establishes that the two surfaces reconstructed from $\aleph _{0}, \aleph_{1}$ are diffeomorphic, we now have to compare the framings. The idea is that we can suitably enrich the spatial diagram $\Psi$ by assembling all of the twist graphics of $\aleph_{t}$ together, a suitable preimage of that should give us a "Serre submanifold" of $W \times I$ which will now be a normally oriented surface. We can then apply Thom-Pontryagin to obtain a function $W \times I \rightarrow SO(2)$ and tilt the ambient framing using this function to obtain some smoothly varying family of framings on $W_{t} = W \times \{ t \}$, ie. an isotopy. Since the framings on $W_{0}, W_{1}$ will be by construction the framings associated to $\aleph_{0}, \aleph_{1}$, this shows they are in fact isotopic, ending the proof. 

In some sense, we are performing a parametrized version of the reconstruction of framed surfaces \textbf{Theorem \ref{thm:reconstruction_of_framed_surfaces}}, the only real difficulty is that we need to ensure that things such as ambient framings and small rectangular neighbourhoods are chosen in a smoothly-varying way. 

We start by choosing a family of small rectangular neighbourhoods of singularities of the graphic of critical values. Observe that the singularities of the graphics $\Psi_{t}$ can all be identified with each other, as the graphics only differ by a chosen isotopy. Fix some real numbers $a, b > 0$. We can define small rectangular neighbourhoods of the codimension $2$ singularities of $\Psi _{t}$ by declaring that their centre is always a given singularity and the lenghts of their sides are the chosen real numbers, let $U_{t}$ be the union of all these neighbourhoods. 

If $a, b$ are chosen sufficiently close to $0$, these rectangles will be disjoint at all times and also small in the needed sense, that is, the singularity they are centred at is the only one they enclose. This way we obtain a family $U_{t}$ of neighbourhoods which is smooth in the sense that it consists of images of an isotopy of embeddings of a single smooth manifold, in particular their totality 

\begin{center}
$U = \bigcup _{t \times [0,1]} U_{t} \times \{ t \} \subseteq \mathbb{R}^{2} \times I$ 
\end{center}
is a smooth open submanifold. Let $\widetilde{W \times I}$ be $F^{-1}(\mathbb{R}^{2} \times I \setminus U)$. Observe that $\widetilde{W \times I} \cap \{ t \}$ is the familiar manifold $\widetilde{W_{t}}$, it consists of $W_{t} \simeq W$ with the preimages of small rectangular neighbourhoods $U_{t}$ along the generic function $f_{t}$ cut out. In particular, on the complements $W_{t} \setminus \widetilde{W_{t}}$ we have well-defined standard framings. The smoothly-varying choice of neighbourhoods and functions allowed us to assemble $W_{t}$ into one three-manifold.

On each of $W_{t}$ we have an ambient framing or more precisely a contractible space of these. Recall that it was defined in terms of the foliation defined as the preimage of the horizontal foliation of $\mathbb{R}^{2}$. By transversality, $\widetilde{W \times I}$ is also foliated by the preimage of this horizontal foliation, this allows one to choose a smoothly varying time-wise ambient family of framings on $\widetilde{W \times I}$. More precisely, such a family gives an ambient framing on each of $\widetilde{W} _{t}$ and this framing varies smoothly with $t$ under the canonical identification of $\widetilde{W} _{t} \subseteq W_{t} \simeq W$ with a submanifold of $W$. 

Since $\aleph _{t}$ was an isotopy of framed planar diagrams, we also have an isotopy of twist graphics $\Upsilon _{t}$. This family can be time-wise lifted using sheet data to obtain a family of Serre submanifolds $S_{t} \subseteq \widetilde{W_{t}}$, as the twist graphic varied smoothly their totality 

\begin{center}
$\hat{S} = S_{t} \times \{ t \} \subseteq \widetilde{W \times I}$
\end{center}
is again a codimension $1$, normally oriented submanifold. We can use the Thom-Pontryagin construction to obtain from $\hat{S}$ a smooth function $\widetilde{W \times I} \rightarrow SO(2)$, tilting the ambient framings using it we obtain a family $v_{t}$ of framings on $\widetilde{W}_{t}$. We can assume that it agrees with the standard framing on each of $\widetilde{W}_{t}$ and so extends to give a family of framings on $W_{t}$, this family can be understood as an isotopy of framings of $W$ under the identification $W_{t} \simeq W$. 

However, the framing on $W_{0}$ is exactly the framing reconstructed from $\aleph _{0}$ as to define it we used the Serre submanifold $S_{0}$ lifted from $\Upsilon _{0}$, similarly the framing on $W_{1}$ is the one reconstructed from $\aleph _{1}$. They are isotopic and so $W_{0}, W_{1}$ are diffeomorphic-isotopic, which is what we needed to prove.
\end{proof}

\begin{thm}[Framed Planar Decomposition]
\label{thm:framed_planar_decomposition}
If $W$ is a framed surface, then the equivalence class of the associated framed planar diagram does not depend on the choice of generic function used to obtain it. 

Consequently, the associated framed planar diagram and the reconstruction of a framed surface from a diagram establish a bijection between equivalence classes of framed planar diagrams and diffeomorphism-isotopy classes of framed surfaces.
\end{thm}

\begin{proof}
Let $W$ be a framed surface and let $f_{0}, f_{1}: W \rightarrow \mathbb{R}^{2}$ be generic. Choose diagrammable framings $v_{i}$ for $i = 0, 1$ compatible with the respective functions, let $\aleph_{0}$ and $\aleph_{1}$ be the associated framed planar diagrams. We have to show they are equivalent. 

By \textbf{Theorem \ref{thm:description_of_framings_on_surface}}, it's enough to show that $\aleph_{0}$ is equivalent to some framed planar diagram $\aleph_{0} ^\prime$ with the same underlying oriented planar diagram as $\aleph_{1}$ and describing the same framing, up to isotopy, on the reconstructed surface, which one can identify with $W$ using $f_{1}$. Indeed, then $\aleph_{0} ^\prime$ and $\aleph_{1}$ are even strongly equivalent and so differ by an isotopy, some of the moves we presented and possibly changes in the chambering graph. By \cite[Lemma 1.44]{chrisphd}, one can modify any chambering graph into any other using only isotopy, positional change of function, positional twist homotopy and change in chambering graph moves we presented and it follows that $\aleph _{0} ^\prime$ and $\aleph _{1}$ are equivalent, too, ending the argument.

Again, we will basically ignore chambering graphs in what follows, as they are not essential, and instead think of framed planar diagrams as describing an oriented surface over $\mathbb{R}^{2}$ equipped with an isotopy class of framings. We can do so by the theorem we just mentioned.

By results of Schommer-Pries, specifically \cite[Corollary 1.38]{chrisphd}, there exists a generic function $F: W \times I \rightarrow \mathbb{R}^{2} \times I$, restricting to $f_{0}, f_{1}$ at the ends of the interval. Similarly how one associates to a generic function into the plane a planar diagram,  associated to $F$ there is a \emph{spatial diagram} $\Phi \subseteq \mathbb{R}^{2} \times I$, a structure we already used in the proof of \textbf{Proposition \ref{prop:invariance_of_reconstructed_framed_surfaces}}.

The diagram $\Phi$ restricts to the planar diagram associated to $f_{0}$ (resp. $f_{1}$) on $\mathbb{R}^{2} \times \{ 0 \}$ (resp. $\mathbb{R}^{2} \times \{ 1 \}$) by construction, its points in-between can be understood as a recipe on how to transform one diagram into the other. More precisely, at almost all times $t \in I$, the restriction $\Phi \cap \mathbb{R}^{2} \times \{ t \}$ is a planar diagram, too, and generically as we vary $t$ continously, this restriction changes by an isotopy of diagrams. At a finite number of critical times, where the restriction can fail to be a planar diagram, some more drastic change occurs. These changes are classified by classifying singularities of generic maps into $\mathbb{R}^{2} \times I$, one obtains this way a finite list of "moves" that can occur. Incorporating these moves into the equivalence relation on planar diagrams one shows that the equivalence class of a planar diagram associated to a surface does not depend on the choice of a generic function, this is the main idea of the corresponding statement for oriented and unoriented surfaces.

We claim that we can lift the changes that the underlying oriented planar diagram of $\aleph_{0}$ undergoes as $t$ goes from $0$ to $1$ into analogous changes of framed planar diagrams, ie. isotopy with possibly a finite number of moves. This way, at time $t = 1$ we will have the promised equivalent framed planar diagram $\aleph_{0} ^\prime$ with the underlying oriented planar diagram coinciding with the one associated to $f_{1}$, which is the same as the underlying oriented diagram of $\aleph_{1}$. We will then show that $\aleph_{0} ^\prime$ describes the framing of $W$ up to isotopy, using the diffeomorphism of the reconstructed surface with $W$ induced by $f_{1}$. This will imply that $\aleph_{0} ^\prime$ and $\aleph_{1}$ are strongly equivalent, finishing the argument.

We can show that such a lift is possible by induction on the structure of the spatial diagram. In other words, it's enough to show how to lift isotopy and each of the moves separately, and that in each case the framings on both sides really do coincide. Any two generic functions from $W$ into the plane can then be connected by a sequence of spatial diagrams corresponding to these changes of diagrams.

Again, it is the case of the isotopy which is more difficult. To show that the two needed framings will be equivalent we will use the same argument as in \textbf{Proposition \ref{prop:invariance_of_reconstructed_framed_surfaces}}, where we lifted some isotopy of twist graphics to a family of Serre submanifolds and using it constructed an isotopy of framings. However, in that case we started with an isotopy of twist graphics, as we had an isotopy of framed planar diagrams to begin with, here we will need to construct it using bare hands.

On the other hand, analogous argument for local moves will be easy simply because the moves are "local enough", as we've observed before, in that the change in function they represent can be contained to some componentwise contractible region of $W$. There is then no trouble with representation of framings, as they need to coincide by necessity. 

We proceed with the proof. Let $f_{0}, f_{1}$ be connected by a spatial diagram representing an isotopy of oriented planar diagrams. By slightly deforming the generic function leading into the diagram we can actually assume that these are restrictions to times $0, 1$ of a path of generic functions $F: W \times I \rightarrow \mathbb{R}^{2} \times I$.

We will show that the isotopy of oriented planar diagrams induced by $F$ can be lifted to an \emph{almost-isotopy} of framed planar diagrams. By an almost-isotopy we mean that twist graphic changes by an isotopy, too, but at a finite number of times the resulting diagram can fail to be a framed planar diagram, but only due to the fact that the height function might stop being Morse when restricted to the free part of the twist graphic. In fact, what will happen is that at such times a Morse birth-death singularity occurs for the height function restricted to the twist graphic. This singularity is one of the topological twist homotopy moves and so the two end-point diagrams $\aleph _{0}$ and $\aleph _{0} ^{\prime}$ will still be equivalent. 

A graphic of critical values $\Psi$ can be understood as a particular type of graph embedded in $\mathbb{R}^{2}$, whose vertices are given by Morse, cusp and fold crossing singularities and whose edges are the arcs. By a \emph{directed edge} of $\Psi$ we mean an arc together with a choice of \emph{source} and \emph{target} vertices among the two vertices it connects. A sequence of directed edges $e_{i}$ of $\Psi$ is called a \emph{delimiting curve} if the target of $e_{i}$ coincides with the source of $e_{i+1}$, the source of the first edge coincides with the target of the last one and additionally each vertex of $\Psi$ appears as a source or target of $e_{i}$ at most twice. These are in some sense the non-self intersecting directed closed curves inside $\Psi$. 

Each delimiting curve bounds an open region of $\mathbb{R}^{2}$ which we call a \emph{tile}, it is homeomorphic to a disk by the Jordan curve theorem. However, since the arcs are always horizontal and their behaviour around vertices is controlled, in this case we do not need to use the theorem, as an explicit coordinatization of the interior as a union of disks glued along the boundary is possible. This coordinatization is presented in \textbf{Figure \ref{fig:coordinatization_of_tiles_as_a_union_of_disks}} below to give a general idea, we do not go into details.

\begin{figure}[htbp!]
	\begin{tikzpicture}
		\begin{scope}
			\draw [very thick, red] (0, -1.5) [out=90, in=-90] to (-1.5, 0) [out=90, in=-90] to (0, 1.5);
			\draw [very thick, red] (0, -1.5) [out=90, in=-90] to (1.5, 0) [out=90, in=-90] to (0, 1.5);
			
			\draw [<->] (1.8, 1.5) to node[auto] {$ x $} (1.8, -1.5);
			\draw [<->] (-1.5, 0) to node[auto] {$ y  $} (1.5, 0);
			\draw [<->] (-0.6, 0.8) to (0.6, 0.8);
			\draw [<->] (-1, -0.6) to (1, -0.6);
			
			\node at (0.1, -2) {Disk-like coordinates on a tile,};
			\node at (0.1, -2.45) { here $x, y \in [0, 1]$ };
		\end{scope}
	
		\begin{scope}[xshift=5cm]
			\draw [very thick, red] (0,0) [out=90, in=180] to (2, 1) [out=0, in=90] to (4, 0);
			\draw [very thick, red] (0,0) [out=-90, in=180] to (1, -1) [out=0, in=180] to (2, 0) [out=0, in=180] to (3, -1) [out=0, in=-90] to (4, 0);
			
			\draw [very thick, dashed] (0,0) to (4,0);
			\node at (2, 0.5) { \scriptsize Disk (1)};
			\node at (0.85, -0.5) { \scriptsize Disk (2)};
			\node at (3.15, -0.5) { \scriptsize Disk (3)};
	
			\node at (2, -2) {A tile coordinatized as a union of three disks};
		\end{scope}
	\end{tikzpicture}
\caption{Coordinatization of tiles as unions of disks}
\label{fig:coordinatization_of_tiles_as_a_union_of_disks}
\end{figure}
This explicit coordinatization is important, as one sees from the definition that it varies smoothly as the delimiting curve undergoes a smooth isotopy. This allows to extend this isotopy of the delimiting curve to an isotopy of the whole tile, by declaring that during the isotopy the position of any given point in the interior of the tile stays constant under the coordinate system we presented. The image in $\mathbb{R}^{2}$ will, however, vary, as the coordinates depend on the embedding of the delimiting curve, which can change throughout the isotopy. 

The fact that we only consider isotopy of graphics ensures that the coordinate system cannot "break". For example, the two closed curves from \textbf{Figure \ref{fig:coordinatization_of_tiles_as_a_union_of_disks}} are isotopic in $\mathbb{R}^{2}$, but there is no isotopy of graphics between them, which is important, as they are coordinatized in a different way. 

Our idea is that we can extend the isotopy of the graphic of critical values to an isotopy of the twist graphic by simply declaring that we use these induced isotopies on the interior of all the minimal tiles. Here, we call a tile \emph{minimal} if it cannot be presented as a union of different tiles.

This method needs a little tweaking, however, as minimal tiles are not in general disjoint. By minimality, their intersection must contain at least one of them and it follows that the only way two minimal tiles can intersect is when one is completely contained in the other. A general situation is presented in \textbf{Figure \ref{fig:tiles_embedded_one_in_another}}, where we have two small minimal tiles embedded into a bigger one.

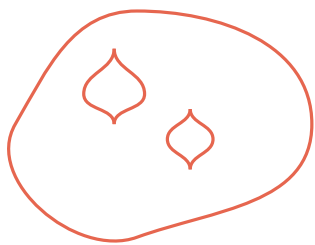
\begin{figure}[htbp!]
	\begin{tikzpicture}
		\draw [very thick, red] (0, -1.5) [out=200, in=-120] to (-1.6, 0) [out=60, in=180] to (0, 1.5) [out=0, in=90] to (2.3, 0) [out=-90, in=20] to (0, -1.5);
		\draw [very thick, red] (0.7, -0.6) [out=90, in=-90] to (0.4, -0.2) [out=90, in=-90] to (0.7, 0.2);
		\draw [very thick, red] (0.7, -0.6) [out=90, in=-90] to (1, -0.2) [out=90, in=-90] to (0.7, 0.2);
		
		\draw [very thick, red] (-0.3, 0) [out=90, in=-90] to (-0.7, 0.4) [out=90, in=-90] to (-0.3, 1);
		\draw [very thick, red] (-0.3, 0) [out=90, in=-90] to (0.1, 0.4) [out=90, in=-90] to (-0.3, 1);
	\end{tikzpicture}
\caption{Minimal tiles embedded one in another}
\label{fig:tiles_embedded_one_in_another}
\end{figure}
We claim that it's enough to solve the problem of extending the isotopy of the graphic of critical values to the twist graphic only in the case outlined above, where we have one big tile with some smaller one embedded in its interior, with additional assumptions that the delimiting curve of the big tile does not move during the isotopy and the extension is already defined in the interior of the smaller tiles. 

Indeed, if we can solve this problem, then we can also solve it without the additional assumption that the delimiting curve of the big tile does not move. To do so we use the disk-like coordinates, as in them any isotopy of the big tile looks as if its boundary did not move. More precisely, we can compose the isotopy of the big tile with some embedding into $\mathbb{R}^{2}$ induced by disk-like coordinates, solve the extension problem there and pull the extension back. This composition sends framed planar diagrams into framed planar diagrams, as the transformation between ambient coordinates on the big tile coming from the embedding into $\mathbb{R}^{2}$ and the disk-like coordinates preserves both vertical and horizonal directions.

Once we can solve the extension of isotopy problem in the interior of any tile, assuming that it was already solved in the interiors of all the smaller tiles it contains, the case of the whole graphic of critical values follows. We first extend the isotopy into tiles which do not contain any other tiles in their interior, then into tiles containing only tiles with this property and so on. This extends it into interior of all tiles contained in the graphic of critical values $\Psi$ and all of the twist graphic is contained in such interiors. 

(As if we had twist graphic in some region not contained in any tile, then we could find a path between it and the infinity of $\mathbb{R}^{2}$. As twist graphic cannot lie on a region with no sheets, this would mean that in the given oriented planar diagram there is an unbounded region of $\mathbb{R}^{2}$ containing some sheets over it. This contradicts the compactness of the surface.)

We now outline how to extend the isotopy from the graphic of critical values to the whole framed graphic in this simplified case from above, ie. we have a big tile $T$ containing some smaller tiles $S_{1}, \ldots, S_{n}$, the extension was already performed in the interior of smaller tiles and additionally the boundary of the big tile stays constant during the isotopy. Here we assume that the list $S_{i}$ contains only \emph{maximal} small tiles, that is, ones that are not contained in any other tile which is not $T$ itself. 

If two small tiles share an edge, then there is some maximal tile containing both of them. Hence, the maximal tiles $S_{i}$ are disjoint in the sense that neither their interiors or their their delimiting curves meet. It follows that there is some minimal distance $d > 0$ such that two points lying in $S_{i}, S_{j}$ for $i \neq j$ are never in a distance smaller than $d$ througout the whole time of the isotopy. 

The extension of the isotopy to the whole vector field will be directed by a flow of a vector field $v_{t}$ for $t \in [0, 1]$, defined in the complement of the interiors of small tiles. Thus, at least in the complement, the isotopy is in fact ambient. 

We want it to extend the given isotopy $h$, which is already defined on the boundary of all the small tiles, this is asking for the vector field to satisfy

\begin{center}
$v_{t_{0}}(h_{t_{0}}(x)) = \partial _{t} h_{t}(x)$
\end{center}
for all $x \in \partial S_{i}$. In words, this means that the flow of a vector field at any given time $t_{0}$ of a point on a boundary of one of the small tiles $S_{i}$ coincides with its movement which is already prescribed by the isotopy $h$. 

This defines $v$ at all times at the boundary of the small tiles, we then take $v$ to be any extension of this to the complement of interior of the small tiles which is also identically zero on the boundary of the big tile $\partial T$. One can define such an extension by finding some smoothly varying family $U_{i}$ of neighbourhoods of $S_{i}$ and define $v$ inside them by interpolating between the given prescription on $\partial S_{i}$ and zero, extending outside of $U_{i}$ by putting zero identically. This is well-defined as long as the radius of $U_{i}$ is smaller than $\frac{d}{2}$, as then these neighbourhoods have throughout the whole isotopy disjoint closures.

Intuitively, we create small "force fields" around the smaller tiles and their effect is that whenever the twist graphic becomes too close to one of $S_{i}$, it is "pushed back". This intuitive description is pictured in \textbf{Figure \ref{fig:small_tile_pushes_back_the_twist_graphic}}.

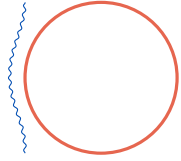
\begin{figure}[htbp!]
	\begin{tikzpicture}[scale=0.5]
		\draw [very thick, red] (2,0) circle (2cm);
		\draw [small_twist, serre_blue] (0, -2) [out=90, in=-90] to (-0.4, 0) [out=90, in=-90] to (0, 2);
	\end{tikzpicture}
\caption{One of the small tiles pushes back the twist graphic}
\label{fig:small_tile_pushes_back_the_twist_graphic}
\end{figure}

The construction above gives us some extension of the isotopy to the twist graphic, however it won't in general satisfy the conditions needed to make it an isotopy of framed planar diagrams, as the height function can quickly stop being Morse when restricted to the twist graphic. One reason for this is that our control of the process is too skeletal, but in fact such an extension is simply not always possible. 

As a concrete example, consider the situation given in \textbf{Figure \ref{fig:an_isotopy_of_graphic_of_critical_values_which_cannot_be_extended_to_framed_graphic}}. The presented fragment of a framed graphic can be completed to a framed planar diagram by extending it downwards, we assume that the arc of the twist graphic is "grounded" at the bottom.

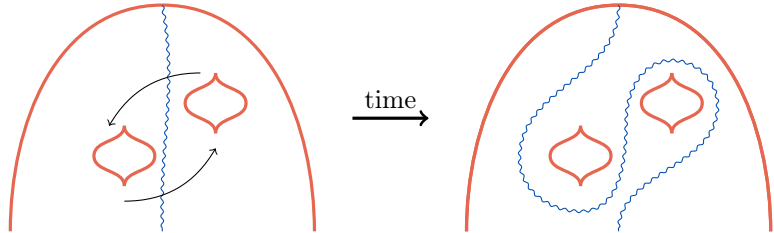
\begin{figure}[htbp!]
	\begin{tikzpicture}
		\begin{scope}
			\draw [very thick, red] (-2, -1) [out=90, in=180] to (0, 2) [out=0, in=90] to (2, -1);
			\draw [very thick, red] (-0.5, -0.4) [out=90, in=-90] to (-0.9, 0) [out=90, in=-90] to (-0.5, 0.4);
			\draw [very thick, red] (-0.5, -0.4) [out=90, in=-90] to (-0.1, 0) [out=90, in=-90] to (-0.5, 0.4);
			\draw [very thick, red] (0.7, 0.3) [out=90, in=-90] to (0.3, 0.7) [out=90, in=-90] to (0.7, 1.1);
			\draw [very thick, red] (0.7, 0.3) [out=90, in=-90] to (1.1, 0.7) [out=90, in=-90] to (0.7, 1.1);
			
			\draw [small_twist, serre_blue] (0, 2) [out=-80, in=93] to (0, -1); 
			
			\draw [->] (0.5, 1.1) [out=180, in=60] to (-0.7, 0.4);
			\draw [->] (-0.5, -0.6) [out=0, in=-120] to (0.7, 0.1);
			
			\draw [very thick, ->] (2.5, 0.5) to node[auto] { time } (3.5, 0.5);
		\end{scope} 
		
		\begin{scope}[xshift=6cm]
			\draw [very thick, red] (-2, -1) [out=90, in=180] to (0, 2) [out=0, in=90] to (2, -1);
			\draw [very thick, red] (-2, -1) [out=90, in=180] to (0, 2) [out=0, in=90] to (2, -1);
			\draw [very thick, red] (-0.5, -0.4) [out=90, in=-90] to (-0.9, 0) [out=90, in=-90] to (-0.5, 0.4);
			\draw [very thick, red] (-0.5, -0.4) [out=90, in=-90] to (-0.1, 0) [out=90, in=-90] to (-0.5, 0.4);
			\draw [very thick, red] (0.7, 0.3) [out=90, in=-90] to (0.3, 0.7) [out=90, in=-90] to (0.7, 1.1);
			\draw [very thick, red] (0.7, 0.3) [out=90, in=-90] to (1.1, 0.7) [out=90, in=-90] to (0.7, 1.1);
			
			\draw [small_twist, serre_blue] (0,2) [out=-90, in=90] to (-1.3, 0) [out=-90, in=180] to (-0.5, -0.7) [out=0, in=-90] to (0.1, 0.5) [out=90, in=180] to (0.7, 1.3) [out=0, in=90] to (1.3, 0.7) [out=-90, in=93] to (0, -1);
		\end{scope}		
	\end{tikzpicture}
\caption{An isotopy of graphic of critical values which cannot be extended to an isotopy of the whole framed graphic}
\label{fig:an_isotopy_of_graphic_of_critical_values_which_cannot_be_extended_to_framed_graphic}
\end{figure}
In the figure we have one big tile and two smaller ones within it. As the smaller tiles move throughout the isotopy, they switch places by passing around each other. It is then simply impossible to extend this isotopy to one of framed planar diagrams, as the arc of the twist graphic cannot stay horizontal throughout the whole time of the isotopy. In the extension we constructed above, the two small tiles would push this arc as they move, winding it around both of them, giving us the situation pictured on the right.

However, the isotopy of the twist graphic we constructed can be slightly deformed and by standard Cerf theory a generic such deformation will have the property that the height function is a Morse function on the twist graphic at almost all times and at times when it doesn't, a Morse birth-death singularity occurs. This ends the construction of the needed almost-isotopy of framed planar diagrams.

The restriction of this isotopy to time $t = 1$ gives us the needed framed planar diagram $\aleph _{0} ^{\prime}$. It is equivalent to $\aleph _{0}$, as they differ only by an honest isotopy of framed planar diagrams and some applications of the topological twist homotopy moves, representing these Morse birth-death singularities, and has the same underlying oriented planar diagram as $\aleph_{1}$ by construction.

We can now show that $\aleph _{0} ^\prime$ describes the framing of $W$ under the diffeomorphism of its reconstructed surface induced by $f_{1}$. This can be done in the same way as we proved \textbf{Proposition \ref{prop:invariance_of_reconstructed_framed_surfaces}}, by lifting the isotopy of the twist graphics to obtain a family of framings on $W \times I$, which on $W \times \{ 0 \}$ restricts to the framing described by $\aleph _{0}$ and on $W \times \{ 1 \}$ to the framing described by $\aleph _{1}$. As the first one is the original framing of $W$, this shows they are isotopic as needed. Note that in the proof of \textbf{Proposition \ref{prop:invariance_of_reconstructed_framed_surfaces}} we lifted an honest isotopy of twist graphics, here we have what we termed an almost-isotopy. This makes no difference, as the latter also can be lifted to a family of Serre submanifolds. This ends the argument in the case where $f_{0}, f_{1}$ differ by an isotopy of generic functions.

We are then left with the case of local moves, which, as we promised, is substantially easier. In terms of functions $f_{0}, f_{1}$, this means that they differ only in some small region of $W$ mapping into a fragment of graphics where the move is being performed. Our change of function moves lift the corresponding moves from \cite{chrisphd} and so we can simply apply it to $\aleph _{0}$ to get the needed $\aleph _{0} ^{\prime}$ with the same underlying oriented planar diagram as $\aleph _{1}$. The framings described by $\aleph _{0} ^\prime $ and $\aleph _{1}$ on the reconstructed surface are the same up to isotopy, as the two diagrams agree outside of this small fragment where the move was performed, this means that the two described framings can be taken to agree outside of this small region where the functions $f_{0}, f_{1}$ differ. This region can always be covered by a disjoint union of disks and hence the framings of $\aleph_{0} ^\prime, \aleph _{1}$ agree up to isotopy by necessity. 

Note that it is possible that to directly apply a change of function move to $\aleph _{0}$ to obtain $\aleph _{0} ^\prime$ we will first need to replace the former by a strongly equivalent framed planar diagram. This is because the configuration of the twist graphic in this small fragment of $\aleph _{0}$ might not be the one that is presented in the local move. 

However, up to isotopy, the framings described by the actual configuration of graphic in $\aleph _{0}$ and the configuration we need  do agree by necessity, as the region can be taken to be a disjoint union of disks. It then follows from \textbf{Theorem \ref{thm:description_of_framings_on_surface}} that it is possible to perform the needed change in a given strong equivalence class of framed planar diagrams with the same underlying oriented diagram as $\aleph _{0}$. This ends the proof. 
\end{proof}

\section{Classification of two-dimensional framed topological field theories}

In this chapter we will use our diagrammatical calculus of framed planar diagrams to obtain a presentation of the framed bordism bicategory $\mathbb{B}ord_{2}^{fr}$. That is, we will construct an equivalence between the framed bordism bicategory and a freely generated symmetric monoidal bicategory $\mathbb{F}(G^{fr}_{bord}, \mathcal{R}_{bord}^{fr})$ on a generating datum we describe explicitly. 

This is analogous to the main result of \cite{chrisphd}, where presentations of the unoriented and oriented bordism bicategories are derived. Together with the theory of shapes modeled on a generating datum, which allows one to describe symmetric monoidal homomorphisms out of a freely generated symmetric monoidal bicategories in a compact way, the presentation gives an explicit description of the bicategory of $2$-dimensional framed topological field theories with arbitrary target. Thus, the result can be understood as providing a complete classification of this type of field theories.

\subsection{Description of the presentation}

In this section we will describe the presentation of the framed bordism bicategory and the way it should be interpreted. Additionally, we give a rough outline of our proof, which we give in detail in the next section.

\begin{figure}[htbp!]
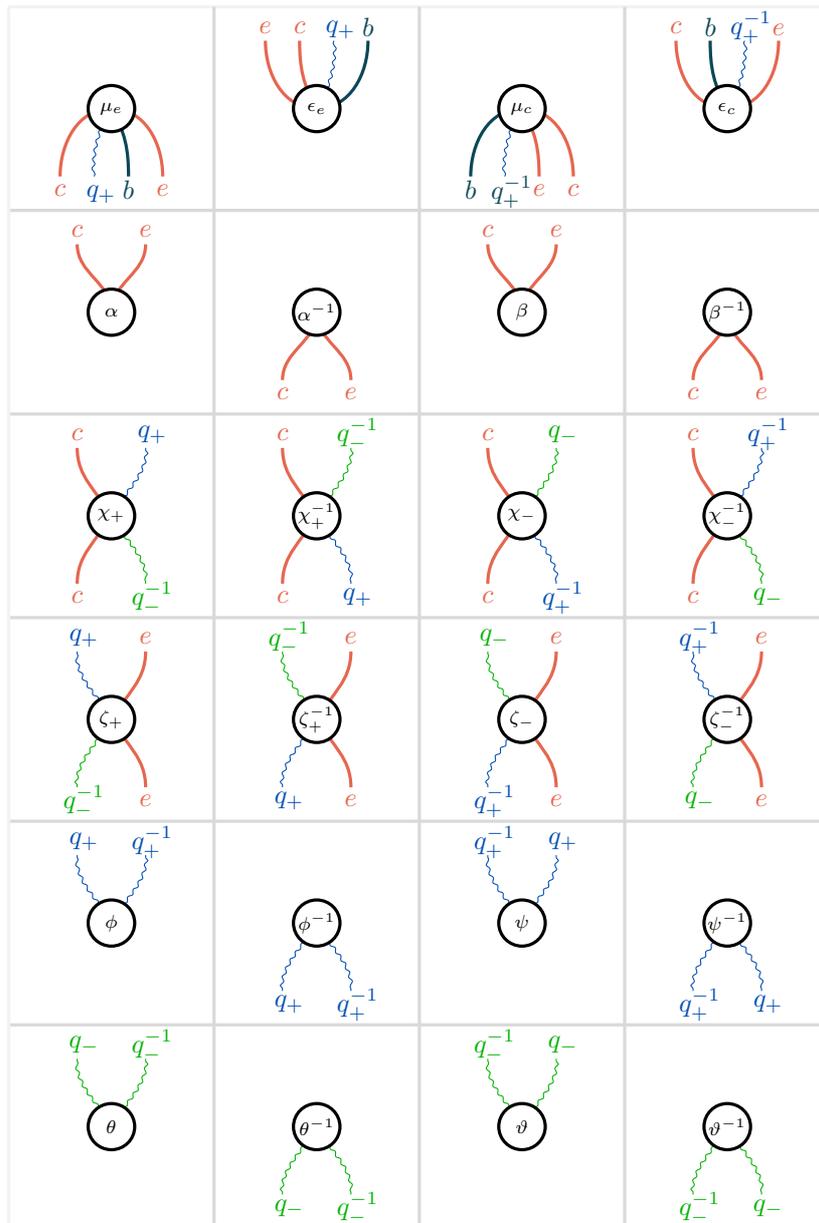


	
	\caption{Generating $2$-cells of the framed bordism bicategory, drawn using string diagrams}
	\label{fig:generating_2_cells_of_the_framed_bordism_bicategory}
\end{figure}

\begin{figure}[htbp!]
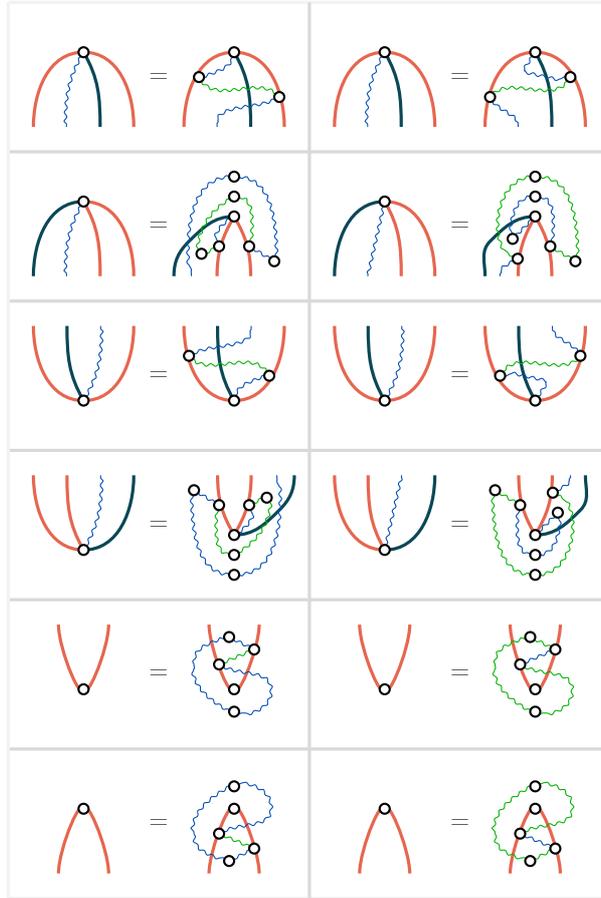

	%
\caption{Round twist relations in the framed bordism bicategory}
\label{fig:round_twist_relations_in_the_framed_bordism_bicategory}
\end{figure}

\begin{figure}[htbp!]
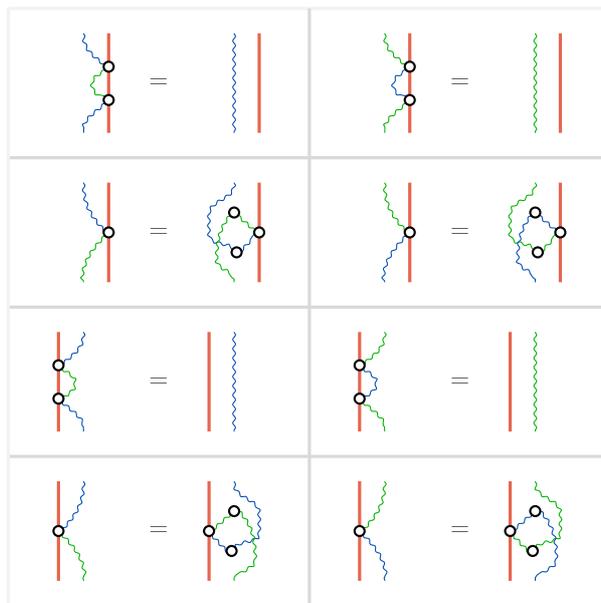

	%
\caption{Crossing point relations in the framed bordism bicategory}
\label{fig:crossing_point_relations_in_the_framed_bordism_bicategory}
\end{figure}

\begin{figure}[htbp!]
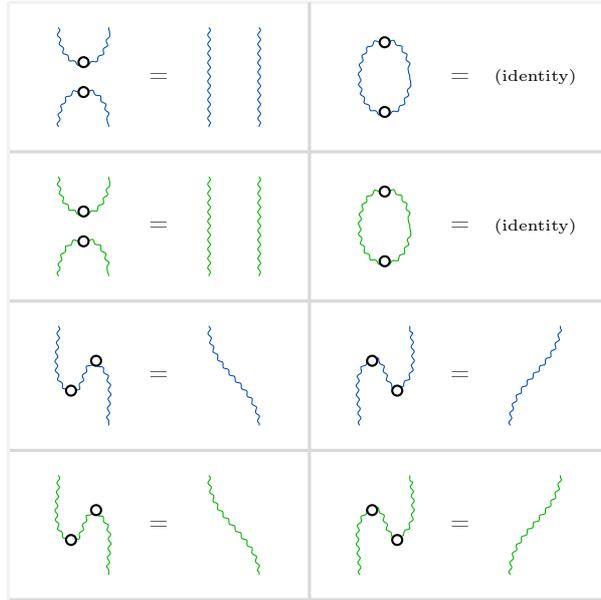

	%
\caption{Twist homotopy relations in the framed bordism bicategory}
\label{fig:twist_homotopy_relations_in_the_framed_bordism_bicategory}
\end{figure}

\begin{figure}[htbp!]
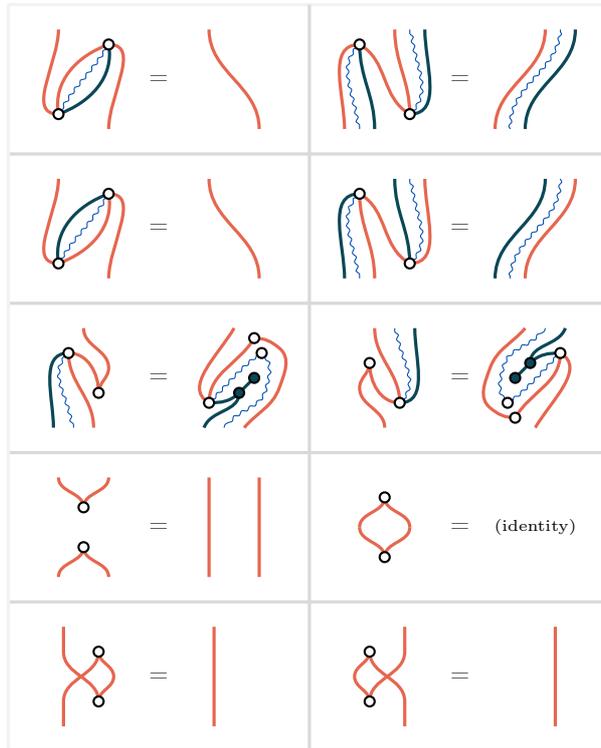

	%
\caption{Relations of oriented bordism in the framed bordism bicategory}
\label{fig:relations_of_oriented_bordism_in_the_framed_bordism_bicategory}
\end{figure}

\begin{thm}[Presentation of the Framed Bordism Bicategory]
\label{thm:presentation_of_the_framed_bordism_bicategory}
The framed bordism bicategory $\mathbb{B}ord_{2}^{fr}$ admits the following presentation.

\begin{center}
\textbf{Generators}
\begin{itemize}
\item Two generating objects $pt_{+}$ and $pt_{-}$. 
\item Six generating morphisms $e: pt_{+} \otimes pt_{-} \rightarrow I$, $c: I \rightarrow pt_{-} \otimes pt_{+}$, $q_{+}, q_{+}^{-1}: pt_{+} \rightarrow pt_{+}$ and $q_{-}, q_{-}^{-1}: pt_{-} \rightarrow pt_{-}$.
\item Twenty-four generating $2$-cells, as given in \textbf{Figure \ref{fig:generating_2_cells_of_the_framed_bordism_bicategory}}.
\end{itemize}
\end{center}

\begin{center}
\textbf{Relations}
\begin{itemize}
\item Round twist relations, as presented in \textbf{Figure \ref{fig:round_twist_relations_in_the_framed_bordism_bicategory}}.
\item Crossing points relations, as presented in \textbf{Figure \ref{fig:crossing_point_relations_in_the_framed_bordism_bicategory}}.
\item Twist homotopy relations, as presented in \textbf{Figure \ref{fig:twist_homotopy_relations_in_the_framed_bordism_bicategory}}.
\item Relations of oriented bordism, as presented in \textbf{Figure \ref{fig:relations_of_oriented_bordism_in_the_framed_bordism_bicategory}}.
\end{itemize}
\end{center}
More precisely, it is equivalent to the free symmetric monoidal bicategory on the above generating datum.
\end{thm}

\begin{notation*}
We denote the generating datum for a symmetric monoidal bicategory given in the statement of the theorem by $(G_{bord} ^{fr}, \mathcal{R}_{bord}^{fr})$, where $G_{bord}^{fr}$ is the set of generating cells and $\mathcal{R}_{bord}^{fr}$ is the set of relations. The free symmetric monoidal bicategory on this datum will be denoted by $\mathbb{F}_{bord}^{fr}$, the free unbiased semistrict symmetric monoidal bicategory by $\mathbb{U}_{bord}^{fr}$.
\end{notation*}

\begin{reading_presentation*}
The presentation we give cannot be directly interpreted as a generating datum for a free symmetric monoidal bicategory, as we don't formally specify sources and targets of generating $2$-cells, we also don't use any bracketings. Additionally, to make sense of the relations one would need to include a large number of coherence cells. Trying to specify all this data would obfuscate the presentation to the point of being completely unreadable. Instead, the presentation above should be interpreted as a generating datum for an unbiased semistrict symmetric monoidal bicategory, which partially strict kind of symmetric monoidal bicategory introduced in \cite{chrisphd}. All unbiased semistrict symmetric monoidal bicategories admit a variant of string calculus, which we use to draw generating $2$-cells and relations between them. 

By making a consistent choice of bracketings and necessary inclusions of coherence cells, one can lift the presentation above to \emph{some} generating datum for a free symmetric monoidal bicategory. We assume such a choice has been made, this is the datum $(G_{bord}^{fr}, \mathcal{R}_{bord}^{fr})$ we referred to above. By a coherence theorem of Schommer-Pries, the induced homomorphism $\mathbb{F}_{bord}^{fr} \rightarrow \mathbb{U}_{bord}^{fr}$ between, respectively, the free symmetric monoidal and the free unbiased semistrict symmetric monoidal bicategories in this datum, will be an equivalence. Since the target does not depend on the choice of the lift of our data to a generating datum for a free symmetric monoidal bicategory, up to equivalence neither does the source. Thus, the ambiguity created in this way is neglible. 

Additionally, in our string diagrams of generating $2$-cells and relations between them we do not label the open regions with tensor products of the generating objects. The labellings are implicitly assumed to correspond to the minimal sheet data around the corresponding singularity or a "local move" of framed planar diagrams. Precisely, the generating $2$-cells of $\epsilon_{e}, \mu_{e}, \epsilon_{c}, \mu_{c}$ lie in the endomorphism category of the monoidal unit, generating $2$-cells $\alpha, \alpha^{-1}, \phi, \phi^{-1}, \psi, \psi^{-1}$ lie in the endomorphism category of $pt_{+}$, $\beta, \beta^{-1}, \theta, \theta^{-1}, \vartheta, \vartheta^{-1}$ in the endomorphism category of $pt_{-}$, $\chi_{+}, \chi_{+}^{-1}, \chi_{-}, \chi_{-}^{-1}$ in the category of morphisms $I \rightarrow pt_{-} \otimes pt_{+}$ and, finally, $\zeta_{+}, \zeta_{+}^{-1}, \zeta_{-}, \zeta_{-}^{-1}$ lie in the category of morphisms $pt_{+} \otimes pt_{-} \rightarrow I$.

We do not label $2$-cells and $1$-cells in our list of generating relations, they are assumed to hold for all labellings consistent with our notation.
\end{reading_presentation*}

Our approach to the proof of the above theorem will be completely analogous to the one taken in \cite{chrisphd} for unoriented and oriented bordism bicategories. A general outline is as follows. We first introduce a relative version of framed planar diagrams, one which is embeded in the square rather than the plane and which can model not only framed surfaces but also framed $2$-bordisms. 

Then, we construct auxiliary symmetric monoidal bicategories $\mathbb{B}ord / \mathbb{P}D_{2}^{fr}$ and $\mathbb{P}D _{2}^{fr}$. The former has "diagrammatic" objects and $1$-cells, but has diffeomorphism-isotopy classes of framed $2$-bordisms as $2$-cells, the latter has $2$-cells given by equivalence classes of relative framed planar diagrams. We construct a span of homomorphisms
\[ \mathbb{B}ord_{2}^{fr} \leftarrow \mathbb{B}ord / \mathbb{P}D_{2}^{fr} \rightarrow \mathbb{P}D_{2}^{fr}, \]
using the Framed Planar Decomposition \textbf{Theorem \ref{thm:framed_planar_decomposition}} we will show that both of maps are equivalences.

The last step will be to compare $\mathbb{P}D_{2}^{fr}$ to the free unbiased semistrict symmetric monoidal bicategory $\mathbb{U}_{bord}^{fr}$ on the generating datum given above. We will do so by directly identifying relative framed planar diagrams with string diagrams in generating $2$-cells using the string diagram calculus of unbiased semistrict symmetric monoidal bicategories, the two bicategories will not only be equivalent, but in fact isomorphic.  

By coherence for unbiased semistrict symmetric monoidal bicategories, $\mathbb{U}_{bord}^{fr}$ is equivalent to the free symmetric monoidal bicategory $\mathbb{F}_{bord}^{fr}$ on the same generating datum. Chasing through all the relevant homomorphisms, one sees that this establishes an equivalence $\mathbb{B}ord_{2}^{fr} \simeq \mathbb{F}_{bord}^{fr}$, which is exactly the result we are after.

\subsection{Proof of the presentation theorem}

In this section we will proceed with the proof of the presentation theorem. We have explained the general idea in the previous section and so will start by introducing the needed notions, which are readily adapted from their unoriented counterparts. 

Since we will be working with framed $2$-bordisms and not only closed surfaces, we need to modify the theory of framed planar diagrams a little bit. Instead of diagrams embedded in $\mathbb{R}^{2}$, we will consider diagrams embedded in the square $I^{2}$. All the codimension $2$ singularities will be assumed to lie in the interior, the intersection of graphics with the vertical boundary $\partial I \times I$ will be assumed to be empty. 

The graphics will be allowed to intersect the horizontal boundary $I \times \partial I$, but only transversally. This way, restriction to either the top or bottom boundary will yield what we call a \emph{framed linear diagram}.The latter is a particular kind of diagram embedded in $I$ that can be used to model framed $1$-bordisms. 

Similarly to the case of framed surfaces, any $1$-bordism equipped with a choice of a generic function into $\mathbb{R}$ and a sufficiently generic framing will induce a framed linear diagram. Conversely, given the diagram one can recover the $1$-bordism together with its map into $\mathbb{R}$ and its framing up to isotopy. 

\begin{defin}
A \textbf{framed linear diagram} $L$ consists of a tuple $(\psi, \upsilon, \gamma)$ of disjoint finite collections of points in interior of the interval $I$ together with labelling data. We call the connected components of $I \setminus (\psi \cup \gamma)$ its \textbf{chambers}, the required labelling data is as follows.
\begin{itemize}
\item to each chamber we associate a natural number $(n)$, which we identify with the set $\{ 0, 1, \ldots, n-1 \}$, its \textbf{set of sheets}, each sheet is labeled as either \emph{positive} or \emph{negative}
\item to each point of the \textbf{chambering set} $\gamma$ we associate a bijection between the sets of sheets over the chambers to its left and right, the bijection is require to take positive elements to positive elements and likewise for the negative ones
\item each point of the \textbf{critical set} $\psi$ is labeled as either a \emph{left elbow} or a \emph{right elbow}
\item to each left elbow we associate two consecutive sheets $i-1, i$ over the chamber to its right, with $i-1$ positive and $i$ negative, the set of sheets to its left is required to be the one obtained by deleting these two elements 
\item to each right elbow we associate two consecutive sheets $i-1, i$ over the chamber to its left, with $i-1$ negative and $i$ positive, the set of sheets to its right is required to be the one obtained by deleting these two elements 
\item to each point of the \textbf{twist set} $\upsilon$ we associate a sheet over its chamber
\item each twist point is labeled as either \emph{straightforward} or \emph{inverted}
\end{itemize}
\end{defin}
A framed linear diagram should be understood as representing a morphism between disjoint unions of standard positively and negatively framed points taken according to the labellings of the sheets over the chambers next to the endpoints of the interval. We make a convention to read linear diagrams from left to right, so that the left-most chamber represents the domain and the right-most chamber the codomain. This is consistent with our way of drawing bordisms, although not consistent with some other appearances in the literature. \footnote{In fact, in their appearance in \cite{chrisphd}, on which our approach is based, linear diagrams and bordisms are read right to left. The tension, of course, stems from the fact that we draw functions as arrows going from left to right, yet we compose them the other way around.}

\begin{figure}[htbp!]
	\begin{tikzpicture}
	
		\draw [|-] (-1.5,0) to node[above] {$ (+, -, +) $} (1,0);
		
		\node at (0.6, -0.4) {$ + $};
		\node at (0.6, -0.8) {$ - $}; 
		\node at (0.6, -1.2) {$ + $};
		
		\draw [->, black!40] (0.8, -0.4) to (1.2, -0.4);
		\draw [->, very thick, serre_blue!70] (0.8, -0.8) to (1.2, -0.8);
		\draw [->, black!40] (0.8, -1.2) to (1.2, -1.2);
		\node[serre_blue] at (1, -1.5) { \tiny straightforward };
		
		\node at (1.4, -0.4) {$ + $};
		\node at (1.4, -0.8) {$ - $};
		\node at (1.4, -1.2) {$ + $};
		
		\draw (1,0) to node[above] {$ (+, -, +) $} (4, 0);
		
		\node at (3.6, -0.4) {$ + $};
		\node at (3.6, -0.8) {$ - $}; 
		\node at (3.6, -1.2) {$ + $};
		
		\draw [->, chambering_blue!70] (3.8, -0.4) to (4.2, -0.4);
		\draw [->, chambering_blue!70] (3.8, -0.8) [out=0, in=180] to (4.2, -1.2);
		\draw [->, chambering_blue!70] (3.8, -1.2) [out=0, in=180] to (4.2, -0.8);
		
		\node at (4.4, -0.4) {$ + $};
		\node at (4.4, -0.8) {$ + $};
		\node at (4.4, -1.2) {$ - $};
		
		\draw (4,0) to node[above] {$ (+, +, -) $} (7, 0); 
		
		\node at (6.6, -0.4) {$ + $};
		\node at (6.6, -0.8) {$ + $}; 
		\node at (6.6, -1.2) {$ - $};
		
		\draw [->, black!40] (6.8, -0.4) [out=0, in=180] to (7.2, -0.8);
		\draw [red] (6.8, -0.8) [out=0, in=90] to (7, -1) [out=-90, in=0] to (6.8, -1.2);
		
		\node at (7.4, -0.8) {$ + $};	
		
		\draw [-|] (7,0) to node[above] {$ (+) $} (8.5, 0);
		
		\node [circle, fill=serre_blue, inner sep=2pt] at (1,0) {};
		\node [circle, fill=chambering_blue, inner sep=2pt] at (4,0) {};
		\node [circle, fill=red, inner sep=2pt] at (7,0) {};
	\end{tikzpicture}
\caption{A framed linear diagram}
\label{fig:framed_linear_diagram}
\end{figure}
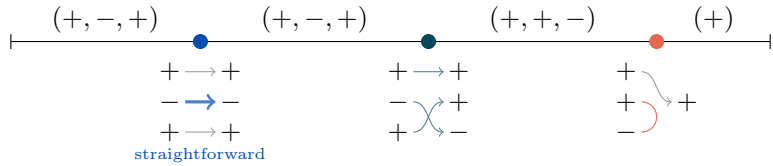

We give an example of a framed linear diagram, drawn with all the necessary labellings, in \textbf{Figure \ref{fig:framed_linear_diagram}}. In the presented diagram, each of the twist set, critical set and chambering set consist of precisely one point. The domain of the map corresponds to the disjoint union $pt_{+} \otimes pt_{-} \otimes pt_{+}$, the codomain to $pt_{+}$. The critical point is a right elbow.

Another possible source of  intuition is that a framed linear diagram is exactly what one would obtain by restricting a framed planar diagram to a sufficiently generic horizontal slice. The labellings over the chambers correspond to positive and negative sheets over a given open region, critical points are restrictions of a fold surface and so represent elbows, chambering points are restrictions of an edge of the chambering graph and so come with a permutation of the sheets. The two consecutive attached to an elbow are the ones that are being "folded together". Lastly, twist points correspond to places where the framing makes a full twist, their labellings corresponds to the direction of the twist.

A new feature is that the set of sheets over any chamber is always identified with $(n)$, which is the same data as a linear ordering. This was not the case for sets of sheets in framed planar diagrams, this change should be understood as being related to the progressing algebraicization of our bicategory as we pass from $\mathbb{B}ord_{2}^{fr}$, through $\mathbb{B}ord / \mathbb{P}D _{2}^{fr}$, to $\mathbb{P}D_{2}^{fr}$. After all, our goal is to compare the framed bordism bicategory to some freely generated symmetric monoidal bicategory and the definition of the latter is purely algebraic.

\begin{prop}[Framed Linear Decomposition]
Let $A= pt_{s_{1}} \otimes \ldots \otimes pt_{s_{k}}$, $B = pt_{t_{1}} \otimes \ldots \otimes pt_{t_{l}}$ be disjoint unions of standard positively and negatively framed points, where $s, t$ are signs. Given an oriented $1$-bordism $w: A \rightarrow B$, a Morse function $f: w \rightarrow \mathbb{R}$ with distinct critical values and a sufficiently generic framing of $w$, one can define a linear diagram $L_{w}$ from $(s_{1}, \ldots, s_{k})$ to $(t_{1}, \ldots, t_{l})$. 

Conversely, given a linear diagram diagram $L$ one can construct an oriented $1$-bordism $w(L)$ mapping into $I$, equipped with a canonical isotopy class of generic framings and inducing $L$ through the above construction, so that $L_{w(L)} = L$ on the nose. Conversely, $w(L_{w}), w$ are diffeomorphic over $I$ as oriented $1$-bordisms and under this diffeomorphism their framings differ at most by an isotopy relative to the boundary.
\end{prop}

\begin{cor}
\label{cor:up_to_iso_any_framed_1_bordism_arises_from_a_framed_linear_diagram}
Up to isomomorphism, any framed $1$-bordism in $\mathbb{B}ord_{2}^{fr}(A, B)$ arises as $w(L)$ for some framed linear diagram $L$.
\end{cor}

\begin{proof}
The construction of the associated framed linear diagram proceeds as follows. The critical values of $f$ become the critical set of of $L_{w}$, on connected components of their complement the map $f$ restricts to a trivial covering and so over each connected component we have a set of positive and negative sheets, depending on whether $f$ preserves or reverses the orientation of the underlying oriented $1$-manifold of $w$. 

We then order the sheets over the chambers in an arbitrary way. It is possible that initially the conditions that the critical points are only allowed to "fold" two consecutive sheets and that the sheets over the endpoints agree with $A, B$ need not be satisfied. However, inserting chambering graph as needed we can correct the orderings to make $L_{w}$ into an oriented linear diagram. 

To define the twist set, we mimic the construction used for framed planar diagrams. Let an \emph{ambient framing} of $w$ be any framing compatible with the coorientation of its halo that has the property that its leading vector is tangent at all times to the underlying $1$-manifold and induces its intrinsic orientation. Note that this framing is defined only on the underlying $1$-manifold, not on the whole halo, and that it is unique up to isotopy. Comparing the given framing of $w$ with any chosen ambient framing we obtain a function into $GL_{+}(2, \mathbb{R})$, composing with Gram-Schmidt orthogonalization we obtain a map $v: w \rightarrow SO(2)$. Again, one readily verifies that this composite does \emph{not} depend on the choice of the ambient framing.

We will say that the framing of $w$ is \emph{sufficiently generic} if the map $v$ is transversal to $-id \in SO(2)$ and the images under $f$ of the preimages $v^{-1}(-id)$ are all distinct from each other and the critical and chambering sets. In this case we can enrich the oriented diagram we already constructed by declaring the twist set to consists of these images, labelling them as straightforward or inverted according to whether the map $v$ preserves or reverses the orientation in their neighbourhood. 

Conversely, if $L_{w}$ is a framed diagram, one can construct an oriented $1$-bordism $w: A \rightarrow B$ together with a map $f: w \rightarrow I$ by gluing the sheets over the given chambers in the way prescribed by the diagram and adding a trivial halo. One can then get a framing on $w$ inducing the given twist set by tilting the ambient framing of $w$ by inserting a full twist in some small neighbourhood of the preimage of the twist point on the prescribed sheet. 

It's clear that $w(L_{w})$, $w$ can be identified as oriented $1$-bordisms mapping into $I$, their framings are isotopic relative to the boundary by Thom-Pontryagin, as up to homotopy they describe the same map into $SO(2)$, having the same preimage of $-id \in SO(2)$. 

To deduce the corollary, observe that if two framed $1$-bordisms in $\mathbb{B}ord_{2}^{fr}$ are isomorphic as oriented bordisms through an isomorphism preserving the framing up to isotopy relative to the boundary, then they are isomorphic as maps in the framed bordism bicategory, as the isotopy can be "spread out" along an invertible $2$-bordism. 

Hence, given a framed $1$-bordism $w: A \rightarrow B$ choose a Morse map $f: w \rightarrow I$ with distinct values and then deform the framing of $w$ by an isotopy to obtain an isomorphic $w^\prime$ whose framing is already sufficiently generic with respect to $f$. Since $w^\prime, w(L_{w^\prime})$ are also isomorphic through an isomorphism preserving the framings up to isotopy, the claim follows.
\end{proof}

\begin{rem}
Note that in the proof we only construct diffeomorphisms of underlying $1$-manifolds of oriented bordisms, we don't construct isomorphisms of haloed manifolds. However, any such diffeomorphism can be lifted to an isomorphism of haloed manifolds, though non-canonically.
\end{rem}

\begin{rem}
Note that the $1$-bordism $w(L)$ associated to a linear diagram $L$ is a well-defined oriented $1$-bordism, but it only carries an \emph{isotopy class} of framings, as there were mild choices involved in the construction given above.  By abuse of language, we will think of them as framed $1$-bordisms, assuming that for each $w(L)$ an appropriate choice has been made once and for all.
\end{rem}

We will now proceed with defining the promised relative version of framed planar diagrams, which can be used to model framed $2$-bordisms. A good idea might be to first become familiar with the corresponding notion for unoriented bordisms, which is \cite[Definition 3.46]{chrisphd} and with the notion of an (absolute) framed planar diagram, which is \textbf{Definition \ref{defin:framed_planar_diagram}}.

The biggest change we will have to make is to consider a slightly different kind of chambering graphs, which was perhaps already visible in the way we defined the generating $2$-cells of the presentation. The problem is that we would like the sets of sheets over given chambers in a relative framed planar diagram to be ordered, which in itself is not difficult.  However, we would also like to assume that the left elbow always folds together two consecutive sheets, of which the larger one is negative and likewise the right elbow folds two consecutive sheets of which the larger one is positive. 

This is, unfortunately, problematic as one sees that around a Morse type singularity there is no way to order sheets so that both of these conditions would be satisfied. Reversing our convention for one of the elbows does not solve the problem, as then this incompatibility disappears, but on the other hand the cusp singularities become impossible to label. Note that this slight inconveniance is not related to framings and appears in the same way in the oriented case, as discussed in \cite{chrisphd}. 

There are two solutions, one of which is to consider two kinds of left and right elbows, corresponding to whether the larger sheet is positive or negative. This leads to multitude of additional generators in the presentation of $\mathbb{B}ord_{2}^{fr}$, but has the advantage that one can follow the whole proof as we present it without any changes to the behaviour of the chambering graph. Once this larger presentation is established, one can use purely category-theoretical methods to obtain the one we described in the theorem, simply by showing that the additional generators are redundant, being isomorphic to the usual left and right elbows composed with symmetry. 

The advantage of that first approach is that it is very simple-minded, the disadvantage is that it forces one to work with a presentation which is essentially double the size of the one we present, even though it models the same phenomena. 

\begin{figure}[htbp!]
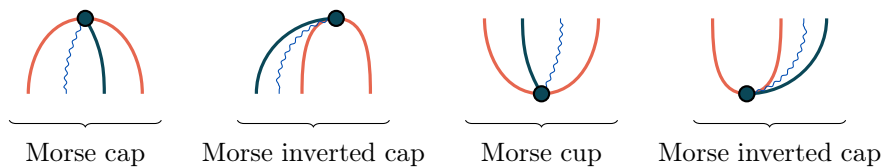


\caption{Special vertices of a chambering graph in relative framed planar diagrams}
\label{fig:special_vertices_of_a_chambering_graph_in_a_relative_framed_planar_diagram}
\end{figure}

The other approach is to consider a modified kind of a chambering graph, which is now not only allowed, but \emph{forced} to have vertices at the Morse singularities of the graphic of critical values. From each such special vertex, there is exactly one edge of the chambering graph leaving it, and we require that the configuration is a small neighbourhood is as given in \textbf{Figure \ref{fig:special_vertices_of_a_chambering_graph_in_a_relative_framed_planar_diagram}}, depending on the singularity type. 

The permutation associated to an edge leaving a special vertex is always required to be the one interchanging the two sheets being folded together. This ensures that the sheets living over chambers can be ordered in a compatible way also in a neighbourhood of a Morse type singularity. The proof of \cite[Proposition 1.46]{chrisphd}, which states that given a graphic a compatible chambering graph can always be chosen, also goes through in the case of such modified graphs. 

The other changes needed when one passes from absolute framed planar diagrams to relative ones are basically the same as in the unoriented case, leading to the following definition. 
\begin{defin}
A \textbf{relative framed planar diagram} consists of a tuple $(\Psi, \Upsilon, \Gamma)$ together with associated sheet data. The graphic of critical values $\Psi \subseteq I^{2}$ and the twist graphic $\Psi \subseteq I^{2}$ are diagrams embedded in a square, their singularities are assumed to occur only in the interior of $I^{2}$. Their arcs do not intersect the vertical boundary $\partial I \times I$ at all, they are allowed to intersect $I \times \partial I$, but only transversally. 

The vertices of the chambering graph $\Gamma \subseteq I^{2}$ are only allowed to lie in the interior, the edges do not intersect the vertical boundary and may intersect horizontal boundary in a transversal way. At each singularity of Morse type of the graphic of critical values $\Gamma$ is required to have a special vertex and an edge leaving it.

The required labelling data is as follows.

\begin{itemize}
\item to each chamber, ie. a connected component of the complement in $I^{2}$ of the graphic of critical values and the chambering graph, we associate a natural number $n = \{ 0, \ldots, n-1 \}$, its set of sheets, together with labelling of its elements as either \emph{positive} or \emph{negative}\
\item analogously to the definition of a linear diagram, each arc of critical values is labeled as either a \emph{left elbow} or a \emph{right elbow} and has associated two consecutive sheets in the chamber to its right or left, which determines sheets over the other chamber
\item to each edge of the chambering graph associate a permutation between the natural numbers corresponding to the chambers to the left and right of the edge, the permutation is required to preserve positivity and negativity of elements, a permutation associated to an edge leaving a special vertex is required to be the permutation exchanging the two sheets being folded together
\item to each arc of the twist graphic, to each chamber it crosses, we associate a sheet over that chamber 
\item each arc of the twist graphic is labeled as either \emph{straightforward} or \emph{inverted}
\item to each codimension $2$ singularity of the graphic $\Psi$ and to each vertex of $\Gamma$ we associate suitable singularity sheet data
\end{itemize}
The data is assumed to satisfy exactly the same conditions as in \textbf{Definition \ref{defin:framed_planar_diagram}} of a non-relative framed planar diagram.
\end{defin}
Observe that the restriction of a relative framed planar diagram to its top and bottom boundary yields two framed linear diagrams $L_{in}, L_{out}$, which we think of as domain and codomain. As no arcs and edges are allowed to cross the vertical boundary, the sources and targets of these linear diagrams are exactly the same. Proceeding as in the non-relative case, by gluing the sheets over chambers, one obtains a framed $2$-bordism $\alpha _{P}: w(L_{in}) \rightarrow w(L_{out})$. More specifically, we obtain an oriented $2$-bordism with framing only well-defined up to isotopy, but this is enough to get a $2$-cell in $\mathbb{B}ord_{2}^{fr}$. 

To reverse this construction, that is, to associate a relative framed planar diagram to a framed $2$-bordism, one has to be a little bit more careful than in the unoriented case. The problem lies in the fact that given a framed $1$-bordism $w$ it is in general \emph{not} possible to find a linear diagram $L$ and an isomorphism $w \simeq w(L)$ of haloed $1$-manifolds, preserving the decomposition of the boundary into domain and codomain and strictly preserving the framing. This is in contrast to the unoriented and oriented cases, where a choice of a Morse function $w \rightarrow I$ would immediately yield the required linear diagram. Thus, for the reverse construction to be possible at all, one has to restrict to $2$-bordisms between $1$-bordisms admitting such isomorphisms. 

Once this is taken into a account, the general idea of association is the same as in the unoriented case. One maps a given $2$-bordism into $I^{2}$, the chosen generic map is assumed to have particular behaviour on the boundary of the $2$-bordism. More specifically, when restricted to the domain $1$-bordism it should give a chosen Morse function into $I \times \{ 1 \}$, on codomain a chosen Morse function into $\{ 0 \} \times I$, and it should map the parts of the boundary modeled on the objects componentwise diffeomorphically to $\{ 0 \} \times I$ for the source and $\{ 1 \} \times I$ for the target. We now describe this in more detail.

Suppose $A, B$ are framed haloed $0$-manifolds, let $w_{1}, w_{2}: A \rightarrow B$ be framed $1$-bordisms. Assume that we are given isomorphisms of framed $1$-bordisms $f_{1}: w_{1} \simeq w(L_{1})$, $f_{2}: w_{2} \simeq w(L_{2})$, by which we mean isomorphism of haloed manifolds, preserving decomposition of the boundary into domain and codomain and strictly preserving the framing. This also fixes Morse functions $w_{1} \rightarrow I$, $w_{2} \rightarrow I$, as $1$-bordisms reconstructed from linear diagrams always come with a map into $I$. 

Observe that the isomorphisms $f_{1}, f_{2}$ induce identifications of framed haloed $0$-manifolds
\begin{gather*}
A \simeq s(w(L_{1})), B \simeq  t(w(L_{1})), 
A \simeq s(w(L_{2})), B \simeq t(w(L_{2}))
\end{gather*}
between $A, B$ and (co)domains of framed $1$-bordisms reconstructed from $L_{1}, L_{2}$. As a compatibility condition, we additionally require that $s(L_{1}) = s(L_{2}), t(L_{1}) = t(L_{2})$ and that the two pairs of isomorphisms above coincide. 

Suppose now that $\alpha: w_{1} \rightarrow w_{2}$ is framed $2$-bordism. We can choose a generic map $\alpha \rightarrow I^{2}$ extending the given Morse maps on the top and bottom boundary and mapping the vertical boundary componentwise diffeomorphically into $\partial I \times I$, in the way we described above. Deforming the framing of $\alpha$ to a diagrammable one and applying the construction described in chapter on calculus of surfaces, we obtain a framed graphic $(\Psi, \Upsilon)$ contained in $I^{2}$. 

Note that the deformation can be taken to be relative to the boundary. Indeed, the framing in the neighbourhood of the domain and codomain already satisfies the needed conditions, as $w_{i}$ are isomorphic to $w(L_{i})$ and so carry a framing transversal to $-id \in SO(2)$ under the comparison with the ambient one. Similarly, on the object part of the boundary of $\alpha$ the given framing in fact agrees with the ambient framing, this follows from the fact that $A, B$ are isomorphic to disjoint unions of standard positively and negatively framed points. This also ensures that the twist graphic is disjoint from $\partial I \times I$, which is one of the conditions we enforce on relative framed planar diagrams.

We now have to choose choose a compatible chambering graph $\Gamma$. One such a choice is made, such a graph will divide $I^{2}$ into chambers and over each such chamber we have sets of sheets and all the usual labelling data specific to non-relative framed planar diagrams. However, to make it into a data of a relative diagram, we need to choose an ordering of the set of sheets over each chamber. 

The orderings are required to satisfy compatibility conditions at arcs, specifically that the two sheets folded together are always consecutive in the ordering, they should also coincide with the given orderings of sheets at the top and bottom boundary, coming from the isomorphisms $f_{1}: w_{1} \simeq w(L_{1})$, $f_{2}: w_{2} \simeq w(L_{2})$. Such orderings can always be chosen, around the vertical boundary $\partial I \times I$ the two orderings coming from the bottom and top linear diagrams coincide by the compatibility conditions we impose on $f_{1}, f_{2}$. On the other hand, in the interior of the $2$-bordism one can make the chambering graph sufficiently fine so that any codimension $2$-singularity is enclosed by it, once this is the case, the compatibility conditions for different arcs do not interfere with each other and one can make the needed choices. 

This way, one obtains a framed planar diagram $P_{\alpha}$, by construction we will also have an isomorphism of oriented haloed $2$-manifolds $\alpha \simeq \alpha _{P_{\alpha}}$, this isomorphism respects the decomposition of boundarites into object and morphism part, in fact induces the chosen isomorphisms $w_{i} \simeq w(L_{i})$ and $A \simeq s(L_{i}), B \simeq t(L_{i})$. It also preserves the framing up to isotopy relative to the boundary.

Our results on calculus of framed surfaces will imply a relative analogue. Let us start by defining an appropriate notion of equivalence for relative framed planar diagrams.

\begin{defin}
Let $P_{1}, P_{2}$ be two relative framed planar diagrams such that their restrictions to the top and bottom boundary are the same linear diagrams $L_{in}, L_{out}$. We will say that $P_{1}, P_{2}$ are \textbf{equivalent} if they differ by an isotopy of diagrams relative to the boundary or the local moves of round twist, crossing point, twist homotopy, change of function or change in chambering graph, applied only in the interior of $I^{2}$. 
\end{defin}

\begin{rem}
The moves involving the singularities of Morse type need to be slightly modified in form to take into account that special vertex of a chambering graph. It is straightforward to do so, the modified form can be seen as relations in our presentation of the framed bordism bicategory. 
\end{rem}
The same arguments as in the non-relative case yield the following theorem. Note that in the statement we omit the isomorphism $w_{i} \simeq w(L_{i})$ used in our analysis above and instead work with a $2$-bordism whose domain and codomain are $w(L_{i})$ to begin with. This is sufficient for our purposes, as the added generality given by these isomorphisms is very illusory. 

\begin{thm}[Relative Framed Planar Decomposition]
\label{thm:relative_framed_planar_decomposition}
Let $\alpha: w(L_{1}) \rightarrow w(L_{2})$ be a framed $2$-bordism. Then, one can find a relative framed planar diagram $P$ restricting to $L_{1}, L_{2}$ at the top and bottom and an isomorphism $\alpha \simeq \alpha_{P}$ of haloed $2$-bordisms which preserves framings up to isotopy relative to the boundary. 

Moreover, if $\alpha ^\prime$ is any other such framed $2$-bordism and $P^\prime$ a diagram presenting it, then $\alpha, \alpha ^\prime$ are isomorphic as haloed $2$-bordisms in a way that preserves framings up to isotopy relative to the boundary if and only if $P, P^\prime$ are equivalent.
\end{thm}

\begin{cor}
If $\alpha, \alpha^\prime: w(L_{1}) \rightarrow w(L_{2})$ are framed $2$-bordisms, then they represent the same $2$-cell in $\mathbb{B}ord_{2}^{fr}$ if and only if relative framed planar diagrams $P, P^\prime$ presenting them are equivalent.
\end{cor}
Having defined all the necessary notions, we are now ready to proceed with the proof of the presentation theorem. We will start by introducting the symmetric monoidal bicategories $\mathbb{B}ord / \mathbb{P}D_{2}^{fr}$ and $\mathbb{P}D_{2}^{fr}$.

\begin{defin}
The bicategory $\mathbb{B}ord / \mathbb{P}D_{2}^{fr}$ of \textbf{framed linear diagrams and bordisms} is a symmetric monoidal bicategory whose objects are words in symbols $pt_{+}$ and $pt_{-}$, whose morphisms are linear diagrams $L$ with the given source and target and whose $2$-cells between two linear diagrams $L_{i}$ are framed $2$-bordisms $\alpha: w(L_{1}) \rightarrow w(L_{2})$ considered up to isomorphism-isotopy. 
\end{defin}
Vertical and horizontal compositions of $2$-cells is given by vertical and horizontal composition of $2$-bordisms. Compositions of $1$-cells concatenates linear diagrams and rescales them. 

The tensor product is more difficult, on objects it is given by concatenation. On linear diagrams one first defines tensor product with the identity diagram by adding the relevant labels from either left or right, a tensor product of arbitrary linear diagrams is then \emph{defined} using the interchange law. 

To formally construct $\mathbb{B}ord / \mathbb{P}D_{2}^{fr}$ as a symmetric monoidal bicategory one derives it from a symmetric monoidal pseudo double category, similarly to $\mathbb{B}ord_{2}^{fr}$ itself. We do not give the details of the argument, as an analogous bicategory for unoriented bordism is described in detail in \cite[3.4.4. Auxiliary bordism bicategories]{chrisphd}, as are composition and tensor products of the relevant diagrams, which we have only briefly discussed above.

\begin{prop}
The symmetric monoidal homomorphism $p: \mathbb{B}ord / \mathbb{P}D_{2}^{fr} \rightarrow \mathbb{B}ord_{2}^{fr}$ which takes a word in $pt_{+}, pt_{-}$ to the corresponding disjoint union of the standard positively and negatively framed points, which takes a linear diagram $L: A \rightarrow B$ to the framed $1$-bordism $w(L)$ and which takes an equivalence class of framed $2$-bordisms to itself, is an equivalence of symmetric monoidal bicategories.
\end{prop}

\begin{proof}
We have to verify three conditions, essential surjectivity on objects, local essential surjectivity on morphisms and local bijectivity on $2$-cells.  

The first condition follows directly from \textbf{Corollary \ref{cor:any_framed_0_manifold_isomorphic_to_tensor_product_of_standard_points_up_to_isotopy}}. Indeed, let $A$ be any framed haloed $0$-manifold and let $P = pt_{s_{1}} \otimes \dots \otimes pt_{s_k}$ be a tensor product of the positively or negatively framed points such that there exists an isomorphism $A \simeq P$ of haloed $0$-manifolds, preserving framings up to isotopy. 

Such an isotopy can be spread out along a $1$-bordism $A \rightarrow P$ whose underlying manifold is simply $A \times I$, this framed $1$-bordism will be an equivalence with the inverse given by $1$-bordism constructed from the inverse isotopy. It follows that $A$ and $P$ are equivalent in $\mathbb{B}ord_{2}^{fr}$, since the latter is in the image of $p$, we are done. 

Local essential surjectivity of $p$ is precisely \textbf{Corollary \ref{cor:up_to_iso_any_framed_1_bordism_arises_from_a_framed_linear_diagram}} and local bijectivity on $2$-cells is clear by construction, as $2$-cells on both sides coincide by definition.
\end{proof}

\begin{defin}
The bicategory $\mathbb{P}D_{2}^{fr}$ of \textbf{diagrammatical framed bordism} is a symmetric monoidal bicategory whose objects are words in symbols $pt_{+}, pt_{-}$, whose morphisms are isotopy classes of linear diagrams and whose $2$-cells are equivalence classes of relative framed planar diagrams.
\end{defin}
Analogously to other symmetric monoidal bicategories we consider in this section, $\mathbb{P}D_{2}^{fr}$ can be formally constructed from a symmetric monoidal pseudo double category.

\begin{prop}
The symmetric monoidal homomorphism q: $\mathbb{B}ord / \mathbb{P}D_{2}^{fr} \rightarrow \mathbb{PD}_{2}^{fr}$ which takes a word in $pt_{+}, pt_{-}$ to itself, takes a linear diagram to its isotopy class and takes a framed $2$-bordism $\alpha: w(L_{1}) \rightarrow w(L_{2})$ to the equivalence class of relative framed planar diagrams representing it, is an equivalence of symmetric monoidal bicategories.
\end{prop}

\begin{proof}
The homomorphism is well-defined on $2$-cells by \textbf{Theorem \ref{thm:relative_framed_planar_decomposition}}, which implies that all relative framed planar diagrams presenting a given framed $2$-bordism $\alpha$ are equivalent to each other. The same result implies that $q$ is locally bijective on $2$-cells. Since it is surjective on objects and morphisms by construction, the claim follows.
\end{proof}

\begin{prop}
The bicategory $\mathbb{P}D_{2}^{fr}$ of diagrammatical framed bordism is isomorphic to the free unbiased semistrict symmetric monoidal bicategory $\mathbb{U}_{bord}^{fr}$. 
\end{prop}

\begin{proof}
The same argument as in the proof of \cite[Proposition 3.49]{chrisphd} works here with little adaptations.

One first compares objects and morphisms. The former is straightforward, as objects on both sides are given by words in symbols $pt_{+}, pt_{-}$. Observe that here we already use the fact that we work with the free unbiased semistrict symmetric monoidal bicategory $\mathbb{U}_{bord}^{fr}$ and not $\mathbb{F}_{bord}^{fr}$ itself, as objects in the latter are given not by words, but by words together with a binary bracketing. 

We compare morphisms by associating to each linear diagram a composite of the generators of $\mathbb{U}_{bord}^{fr}$ obtained by reading the diagram from left to right. Left elbow critical points correspond to $c$, right elbows correspond to $e$, straightforward twist points lying on positive sheets correspond to $q_{+}$, inverted ones to $q_{+}^{-1}$ and the twist points lying on negative sheets correspond to either $q_{-}$ or $q_{-}^{-1}$. The chambering points correspond to applications of the symmetry. 

Comparison at the level of $2$-cell is similar, except more complicated. Recall that unbiased semistrict symmetric monoidal bicategories admit a variant of string calculus, it then follows that $2$-cells in $\mathbb{U}_{bord}^{fr}$ can be taken to be equivalence classes of such string diagrams, labeled by the generatoring cells of the presentation. 

Our presentation $(G_{bord}^{fr}, \mathcal{R}_{bord}^{fr})$ is chosen exactly in a way such that given a relative framed planar diagram one can associate to it a string diagram by careful interpretation. More precisely, the ordered sets of sheets over chambers give the needed labelling of open regions of a string diagram. The edges of the chambering graph become distinguished strings corresponding to symmetries, the arcs of left or right elbows become strings on the generators $e, c$, arcs of the twist graphic become strings on $q_{+}, q_{+}^{-1}, q_{-}$ or $q_{-}^{-1}$, according to the labelling data, in the same way as on the level of morphisms.

Lastly, singularities of relative framed planar diagrams become "coupons" labeled by generating $2$-cells in $G_{bord}^{fr}$ in the evident way. More precisely, Morse singularities become $\mu_{e}, \epsilon_{e}, \mu_{c}$ or $\mu_{c}$, cusp singularities become $\alpha, \alpha^{-1}, \beta$ or $\beta^{-1}$, crossing points become either $\chi$ or $\zeta$ and twist birth-death becomes $\phi, \psi, \theta$ or $\vartheta$. In the same way one can given a string diagram interpret it as a relative framed planar diagram. 

To establish that equivalence classes on both sides correspond to each other, one verifies by hand that there is a correspondence between "local moves" allowed in relative framed planar diagrams and those allowed in unbiased semistrict string diagrams.

Some moves of relative framed planar diagrams, in particular those which we have called positional and the round twist moves for fold crossing, are automatically allowed as moves in string diagrams, as they witness some form of naturality present in any unbiased semistrict symmetric monoidal bicategory. 

Those which are not automatic included are explicitly given as relations in the presentation. In most cases, the relevant moves were included "as is", there is only a slight difference between crossing point moves of framed planar diagrams and the crossing point relations we impose in the presentation. One can verify that these two different types of local moves in fact imply each other under the assumption of all the others, this minor change was needed because the crossing point local moves of framed planar diagrams do not have a form of a relation between $2$-cells in a symmetric monoidal bicategory.
\end{proof}

\section{The Cobordism Hypothesis in dimension two}

In this chapter we will prove the Cobordism Hypothesis of Baez-Dolan in dimension two. We will work with the following precise statement, which is a straightforward adaptation of the formulation due to Jacob Lurie to the setting of symmetric monoidal bicategories. 

\begin{thm}[The Cobordism Hypothesis of Baez-Dolan]
\label{thm:cobordism_hypothesis_in_dimension_two}
Let $\mathbb{B}ord_{2}^{fr}$ be the framed bordism bicategory and let $\mathbb{M}$ be any other symmetric monoidal bicategory. Then, the evaluation at the positively framed point homomorphism induces an equivalence
\[ \textbf{SymMonBicat}(\mathbb{B}ord_{2}^{fr}, \mathbb{M}) \rightarrow \mathbb{K}(\mathbb{M}^{fd}) \]
between the bicategory of framed topological field theories with values in $\mathbb{M}$ and the groupoid of fully dualizable objects in $\mathbb{M}$.
\end{thm}
The Cobordism Hypothesis was first conjectured by John Baez and James Dolan in their landmark paper \cite{baezdolanhigheralgebraandtqfts}, where they also introduced the idea of using higher categories to model extended topological field theories. 

In principle, it is a statement about classification of topological field theories, although one could argue that since the answer ("an $n$-dimensional framed topological field theory is determined by its value at a point, which can be any fully dualizable object") turns out to be so surprisingly simple, the Cobordism Hypothesis in fact reveals a deep connection between topology of manifolds and the theory of higher categories. 

The Cobordism Hypothesis was proven in all dimensions by Jacob Lurie, one of his many insights was to describe the problem using the theory of $(\infty, n)$-categories and not only $n$-categories, this allows one to set up an inductive proof where the Cobordism Hypothesis in one dimension is deduced from dimensions below. An extensive sketch of the argument appears in \cite{lurietqfts}. To our knowledge, a detailed account has yet to appear in the literature. 

Our approach to the proof is not technically dependant on the work of Jacob Lurie, rather extending results of Christopher Schommer-Pries on classification of two-dimensional topological field theories from \cite{chrisphd}, which first appeared at roughly similar time as \cite{lurietqfts}. One should, however, stress that our proof has a computational flavour and it is only possible because a precise, correct statement of the Cobordism Hypothesis is already known. 

\subsection{Proof of the Cobordism Hypothesis}

This whole section is devoted to the proof of the Cobordism Hypothesis in dimension two. 

Observe that the conditions of the Cobordism Hypothesis are satisfied by the free symmetric monoidal bicategory $\mathbb{F}_{cfd}$ on a coherent fully dual pair, this is essentially the content of the coherence for fully dualizable objects theorem, which we proved in the current work as \textbf{Theorem \ref{thm:coherence_for_fully_dualizable_objects}}. Indeed, consider the commutative diagram 

\begin{center}
	\begin{tikzpicture}
		\node (MarginNode) at (0, 0.2) {};
		\node (StrHomomorphisms) at (0.5,0) {$ \mathcal{C}ohFullyDualPair(\mathbb{M}) $};
		\node (Homomorphisms) at (6.5,0) {$ \textbf{SymMonBicat}(\mathbb{F}_{cfd}, \mathbb{M})$};
		\node (Mfd) at (3.5,-1) {$ \mathbb{K}(\mathbb{M}^{fd}) $};
		
		\draw [->] (StrHomomorphisms) to (Homomorphisms);
		\draw [->] (StrHomomorphisms) to (Mfd);
		\draw [->] (Homomorphisms) to (Mfd);
	\end{tikzpicture}.
\end{center}
of bicategories and homomorphisms. The horizontal arrow is given by associating to a coherent fully dual pair the induced strict homomorphism from the free symmetric monoidal bicategory on a coherent fully dual pair, it is an equivalence by the cofibrancy theorem for freely generated symmetric monoidal bicategories. 

Both vertical arrows are given by evaluating at the generating object $L$ of $\mathbb{F}_{cfd}$. The left one is an equivalence by the coherence theorem for fully dualizable objects, it then follows that the right vertical map is also an equivalence for all symmetric monoidal $\mathbb{M}$. This is what we mean by saying that $\mathbb{F}_{cfd}$ satisfies the conditions of the Cobordism Hypothesis.

The property of satisfying these conditions is clearly invariant under equivalences of symmetric monoidal bicategories preserving the distinguished object at which we evaluate, thus, to prove the Cobordism Hypothesis it is enough to establish an equivalence $\mathbb{B}ord_{2}^{fr} \simeq \mathbb{F}_{cfd}$.

Our main result of the previous chapter was a derivation of the presentation of the framed bordism bicategory, which is an equivalence $\mathbb{B}ord_{2}^{fr} \simeq \mathbb{F} _{bord} ^{fr}$ with a freely generated symmetric monoidal bicategory whose generating datum we described explicitly. This further reduces the proof of the Cobordism Hypothesis to showing that the two symmetric monoidal bicategories $\mathbb{F}_{cfd}, \mathbb{F}_{bord}^{fr}$ are equivalent ot each other. As both admit finite algebraic descriptions, a verification of this fact is essentially a matter of computation. 

The comparison can be made a little easier by coherence for unbiased semistrict symmetric monoidal bicategories. Indeed, instead of comparing freely generated symmetric monoidal bicategories one can instead compare freely generated unbiased semistrict symmetric monoidal bicategories. The latter have the same equivalence type, but are much smaller and admit string diagram calculus, making calculations more manageable. 

The precise statement we will prove is the following, by what we explained above it is equivalent to the Cobordism Hypothesis in dimension two. 

\begin{thm}
\label{thm:inclusion_i_from_ucfd_to_u_bord_fr_is_an_equivalence}
Let $i: \mathbb{U}_{cfd} \rightarrow \mathbb{U} _{bord}^{fr} $ be the strict homomorphism classifying the coherent fully dual pair consisting of the positive point, negative point, left and right elbows, cusp generators as cusp isomorphisms and Morse generators as (co)units of the adjunctions $e^{L} \dashv e, c^{L} \dashv c$. Then $i$ is an equivalence of symmetric monoidal bicategories.
\end{thm}
Observe that this homomorphism is well-defined,  as we verified in \textbf{Remarks \ref{rem:coherency_of_a_dual_pair_of_positive_and_negative_points},  \ref{rem:coherency_of_a_fully_dual_pair_of_positive_and_negative_points}} that the relevant fully dual pair is indeed coherent, and so induces a strict homomorphism $\mathbb{F}_{cfd} \rightarrow \mathbb{U}_{bord}^{fr}$. This homomorphism will factorize through $\mathbb{U}_{cfd}$ because of the universal property of the latter with respect to homomorphisms into unbiased semistrict symmetric monoidal bicategories. 

It is not necessary to involve geometry of framed manifolds to construct this homomorphisms, like we did above, as a more direct description is available. Algebraically, $i$ corresponds to an inclusion of free generating data $G_{fd} \subseteq G_{bord}^{fr}$ on which the two symmetric monoidal bicategories are constructed. More specifically, this inclusion takes $L$ to $pt_{+}$, $R$ to $pt_{-}$, $e$ to $e$, $c$ to $c$, $q$ to $q_{+}$, $q^{-1}$ to $q_{+}^{-1}$ and follows labels we have given on both sides at the level of $2$-cells. That is, the $2$-cell generators $\alpha, \beta, \alpha^{-1}, \beta^{-1}, \mu_{e}, \epsilon_{e}, \mu_{c}, \epsilon_{c}, \psi, \psi^{-1}, \phi, \phi^{-1}$ of $G_{fd}$ are taken to the corresponding $2$-cells in $G_{bord}^{fr}$. 

From this point of view, $\mathbb{U}_{bord}^{fr}$ may be intuitively seen as a larger symmetric monoidal bicategory which envelopes $\mathbb{U}_{cfd}$. This suggests an idea of showing that the two bicategories are equivalent to each other by exhibiting a retraction of unbiased semistrict symmetric monoidal bicategories, which is exactly what we will do. Before proceeding with the construction of the inverse, however, we will cut down a little bit on the presentation of the framed bordism bicategory by showing that some of the relations are redundant. 

\begin{lem}
\label{lem:morse_round_twist_relations_are_derivative}
The Morse round twist relations in the presentation $(G_{bord}^{fr}, \mathcal{R}_{bord}^{fr})$ of the framed bordism bicategory are derivative. That is, if $\mathcal{R}_{bord}^{fr \prime}$ is the class of relations obtained from $\mathcal{R}_{bord}^{fr}$ by removing them, then the quotient homomorphism $\mathbb{U}_{bord}^{fr \prime} \rightarrow \mathbb{U}_{bord}^{fr}$ is an isomorphism of symmetric monoidal bicategories.
\end{lem}
We warn the reader that the proof given below can be considered slightly unsatisfactory, as it does not produce an explicit derivation of the relevant relations from the others, as one might expect. Instead, after an obvious simplification, we rely on the details of the Framed Planar Decomposition theorem, which itself extends Morse-theoretic methods developed by Christopher Schommer-Pries. The lack of explicitness stems from the fact that it can be difficult to connect two generic functions $W \rightarrow I^{2}$ from a surface to a square by a generic function $W \times I \rightarrow I^{2} \times I$, even when such functions exist in abundance by transversality arguments. 
\begin{proof}
We have to show that the Morse round twist relations can be derived from cusp round twist, crossing point, twist homotopy and relations of oriented bordism. We will do this in steps, first showing that it is enough to prove only some of them. 

Observe that out of the eight Morse round twist relations, four correspond to cups and caps and four to saddles. We claim that the saddle ones follow from the cup/cap ones and conversely.

Indeed, the right hand sides of Morse round twist relations can be easily grouped in pairs which satisfy triangle equations, that is, they form (co)unit pairs of adjunctions $e^{L} \dashv e$ or $c^{L} \dashv c$, in such a pair we always have one relation concerning a cup/cap and one concerning a saddle. These (co)unit pairs could in principle be different from the generating $\mu_{e}, \epsilon_{e}, \mu_{c}, \epsilon_{c}$, the Morse round twist relations say that in fact they are the same. As (co)units determine each other, it is enough to show that either the unit or a counit of these new adjunctions is standard. This proves that the two relations forming a pair are equivalent and hence the round twist for cups/caps is equivalent to round twist for saddles.

We will now prove the round twist relations for Morse cups and caps follow from cusp twist, twist homotopy, crossing point and relations of oriented bordism, this will end the proof of the proof of the lemma. For concreteness, let us proceed in the specific case of the relation

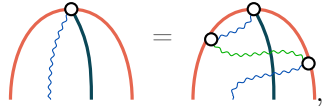
\begin{center}
,
\end{center}
the argument works analogously for all the other cups/caps. In fact, we will prove a stronger statement, namely that one can deduce the "unwinding relation" 

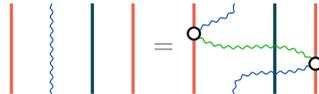
\begin{center}
%
\end{center}
This clearly implies all Morse round twists for cup/caps at once. In our argument we will exploit the details of the proof of Framed Planar Decomposition \textbf{Theorem \ref{thm:framed_planar_decomposition}}.

Observe that under geometrical interpretation as relative framed planar diagrams, both sides of the unwinding relation correspond to a cyllinder. If the two framings presented were isotopic, we could conclude that they two diagrams presented differ by an application of suitable moves, however, this is not the case. 

What is the case is that both framed cyllinders are diffeomorphic-isotopic, the diffeomorphism in question given by the Dehn twist. By Framed Planar Decomposition, the two diagrams are equivalent, that is, they differ by an isotopy and a sequence of applications of the relations in the framed bordism bicategory. However, this sequence will in general also involve relations of Morse round twist, we will show that with suitable care one can ensure that these do not occur. 

The idea is to instead consider the equivalent "modified unwinding relation" 

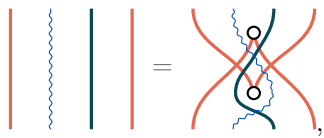
\begin{center}
\begin{tikzpicture}[scale=0.8]
			\begin{scope}
				\draw [very thick, red] (0,0) to (0,2);
				\draw [very thick, red] (2,0) to (2, 2);
				\draw [very thick, chambering_blue] (1.33, 0) to (1.33, 2);
				\draw [small_twist, serre_blue] (0.66, 0) to (0.66, 2);
			\end{scope}
			
			\node at (2.5, 1) { = };
			
			\begin{scope}[xshift=3cm]
			
				\draw [very thick, red] (0,0) [out=90, in=-90] to (1,1.7);
				\draw [very thick, red] (2,0) [out=90, in=-90] to (1,1.7);
				\draw [very thick, red] (0,2) [out=-90, in=90] to (1, 0.5);
				\draw [very thick, red] (2,2) [out=-90, in=90] to (1, 0.5);
				
				\draw [very thick, chambering_blue] (1.33, 2) [out=-90, in=60] to (1, 1.2);
				\draw [very thick, chambering_blue] (1, 1.2) [out=-120, in=90] to (0.7, 0.6);
				\draw [very thick, chambering_blue] (0.7, 0.6) [out=-90, in=90] to (1.33, 0);
				
				\draw [small_twist, serre_blue] (0.66, 2) [out=-90, in=120] to (1,1.2);
				\draw [small_twist, serre_blue] (1,1.2) [out=-60, in=90] to (1.3, 0.6);
				\draw [small_twist, serre_blue] (1.3, 0.6) [out=-90, in=90] to (0.66, 0);
				
				\draw [thick, fill=white] (1, 1.6) circle (0.10cm);
				\draw [thick, fill=white] (1, 0.6) circle (0.10cm);
			\end{scope}
\end{tikzpicture},
\end{center}
where the right hand side arises from the right hand side just above by an application of two Swallowtail and some crossing point relations. We will show that the proof of Framed Planar Decomposition theorem can be used to obtain a sequence of relations in framed bordism bicategory connecting these two diagrams, none of which is a Morse round twist, this will end the proof of the lemma.

Recall that in the proof fof Framed Planar Decomposition one shows that two framed planar diagrams $\aleph_{0}, \aleph_{1}$ are equivalent by first finding an oriented spatial diagram connecting the underlying oriented diagrams, this is possible by results from \cite{chrisphd}. One then separates this spatial diagram into a finite application of relations of oriented bordism and isotopy of oriented diagrams, this process is then suitably lifted to analogous changes in framed planar diagrams with diagram $\aleph_{0}$ as the starting point. The end result is a third framed planar diagram $\aleph_{0} ^ \prime$ which, as one shows, is strongly equivalent to $\aleph_{1}$, that is, has the same underlying oriented diagram and describes the same framing on the reconstructed oriented surface up to isotopy. 

In our case, if we let $\aleph_{0}$ and $\aleph_{1}$ be the left and right hand sides of the modified unwinding relation, then one can conclude that $\aleph_{0} ^\prime$ differs from $\aleph_{1}$ at most by an isotopy of twist graphics and a finite number of applications of cusp round twist, crossing point and twist homotopy relations. Indeed, no Morse round twist relations are needed as the underlying oriented diagram has no Morse singularities to begin with. 

Thus, to finish our argument, it is enough to show that $\aleph _{0}^{\prime}$ can be constructed from $\aleph_{0}$ without using Morse round twist moves. In fact, we will show that neither any type of round twist nor crossing point moves are needed.

Observe that the only moment in the process described in the proof of \textbf{Theorem \ref{thm:framed_planar_decomposition}} where one might apply one of round twist relations is when we have to lift an application of an oriented bordism relation of oriented diagrams to one of framed planar diagrams. Indeed, to perform the lift one needs the twist graphic to be of particular form, exactly as on one of the sides of the oriented bordism relation. 

As the region of the surface on which the relation takes place is componentwise contractible, all framings on it coincide up to isotopy and thus one invokes \textbf{Theorem \ref{thm:description_of_framings_on_surface}} on calculus of framings to conclude that the twist graphic can be put into particular form. However, one sees by examining both sides of oriented relations that if the framed planar diagram in question has no crossing points, then the twist graphic can be modified as needed using only twist homotopy.

Since $\aleph_{0}$ has no crossing points to begin with, it follows that the first lift can be performed using only relations of twist homotopy and oriented bordism. As neither isotopy nor this type of "careful" lift can introduce crossing points, the process can be continued without resorting to crossing point or round twist relations at any subsequent point. This shows that $\aleph_{0}$ and $\aleph_{0}^\prime$ can be assumed to differ by only the restricted set of moves, ending our proof of the lemma.
\end{proof}

\begin{rem}
To use the above proof to obtain an explicit derivation of Morse round twist from the other relations one would have to start by connecting the standard projection from the cyllinder to the square with the projection composed with the Dehn twist by a sufficiently generic map $(S^{1} \times I) \times I \rightarrow \mathbb{R}^{2} \times I$ and proceed from there. This might or might not be difficult.
\end{rem}
We will now begin the proof of \textbf{Theorem \ref{thm:inclusion_i_from_ucfd_to_u_bord_fr_is_an_equivalence}} which states that the strict homomorphism $i: \mathbb{U}_{cfd} \rightarrow \mathbb{U}_{bord}^{fr}$ classifying the coherent fully dual pair of standard positively and negatively framed points is an equivalence of symmetric monoidal bicategories, we have already explained how this entails the Cobordism Hypothesis in dimension two. 

We start by defining a strict homomorphism $r: \mathbb{U}_{bord}^{fr} \rightarrow \mathbb{U}_{cfd}$ which is an explicit pseudoinverse to $i$, we will proceed in such a way that it will be clear that $ri = id$ on the nose. We will then show that $ir$ is equivalent to the identity as an endomorphism of $\mathbb{U}_{bord}^{fr}$, in fact, we will prove that the needed natural equivalence can be chosen to have identity components on objects and so constitutes what in the case of bicategories is known as an \emph{icon}, see \cite{lack_icons}.

Since $r$ is supposed to be strict, it is enough to prescribe it on generators of $\mathbb{U}_{bord}^{fr}$, the definition is largely dictated by the requirement that $ri$ is the identity. Indeed, this immediately forces $r(pt_{+}) = L$ and $r(pt_{-}) = R$, as $i$ is bijective on objects. Similarly, we must have $r(e) = e$, $r(c) = c$, $r(q_{+}) = q$ and $r(q_{+}^{-1}) = q^{-1}$. 

The first place where something interesting happens is the case of generating morphisms $q_{-}, q_{-}^{-1}$, that is, Serre autoequivalences of the negatively framed point, as there are no corresponding generators of $\mathbb{U}_{cfd}$. As we want morphisms $q_{-}$ and $ir(q_{-})$ to be isomorphic, a little bit of thought shows that 
\begin{gather*} 
r(q_{-}) = (R \otimes e) \circ (R \otimes q_{+}^{-1} \otimes R) \circ (c \otimes R), \\
r(q_{-}^{-1}) = (R \otimes e) \circ (R \otimes q_{+} \otimes R) \circ (c \otimes R),
\end{gather*}
are sensible candidates. In the informal language of ordinary string diagrams the definition can be presented as 

\begin{center}
	\begin{tikzpicture}[scale=0.8]
		\begin{scope}
		\node at (-1.6, 0.8) {$ r(q_{-}) = $};
	
		
		\draw[thick, red] (0, 1.5) to (1.2, 1.5);
		\draw[thick, red] (0, 1.5) [out=180, in=180] to node[left] {$ c $} (0, 0.75);
		\draw[small_twist, serre_blue] (0, 0.75) to node[auto] {$ q^{-1} $} (1,0.75);
		\draw[thick, red] (1, 0.75) [out=0, in=0] to node[auto] {$ e $} (1, 0);
		\draw[thick, red] (-0.2, 0) to (1,0);
		\end{scope}
		
		\begin{scope}[xshift=4.9cm]
		\node at (-1.6, 0.8) {$ r(q_{-} ^{-1}) = $};
	
		
		\draw[thick, red] (0, 1.5) to (1.2, 1.5);
		\draw[thick, red] (0, 1.5) [out=180, in=180] to node[left] {$ c $} (0, 0.75);
		\draw[small_twist, serre_blue] (0, 0.75) to node[auto] {$ q $} (1,0.75);
		\draw[thick, red] (1, 0.75) [out=0, in=0] to node[auto] {$ e $} (1, 0);
		\draw[thick, red] (-0.2, 0) to (1,0);
		\end{scope}		
	\end{tikzpicture}.
\end{center}
The relevant composites are easily seen to be isomorphic to Serre autoequivalences of the negatively framed point, the needed isomorphisms given by a crossing point followed by a cusp.

We are left with defining $r$ on generating $2$-cells. Again, this is straightforward for those generators that are obviously the image of a generator of $\mathbb{U}_{cfd}$. These are precisely $\alpha, \beta, \alpha^{-1}, \beta^{-1}, \mu_{e}, \epsilon_{e}, \mu_{c}, \epsilon_{c}, \psi, \psi^{-1}, \phi, \phi^{-1}$ and in this case we prescribe $r(\Box) = \Box$, where $\Box$ is one of the symbols denoting aforementioned generators. 

Hence, we now only have to define $r$ on crossing point and twist birth-death for Serre autoequivalences of the negatively framed point. This is where things start to be a little tricky. 

As we work with unbiased semistrict symmetric monoidal bicategories, we will use unbiased semistrict string diagrams to give the needed definition. We will, however, supplement it with informal pictures of spaced-out ordinary string diagrams for symmetric monoidal categories, as they often turn out to be more readable. 

We start with a crossing generator where the Serre autoequivalence of the positive point passes through a coevaluation and becomes the Serre autoequivalence of the negative point, the definition is 

\begin{center}
	\begin{tikzpicture}
		\node at (-0.3, 0.5) {$ r( $};
	
		\draw[thick, red] (0,1) to (0,0);
		\draw[small_twist, serre_blue] (0.5, 1) to (0, 0.5);
		\draw[small_twist, serre_green] (0, 0.5) to (0.5, 0);
		
		\draw [thick, fill=white] (0, 0.5) circle (0.10cm);
		
		\node at (0.5, 0.5) {$ ) $};
		\node at (1, 0.5) {$ = $};
		
		\begin{scope}[xshift=1.5cm]
		\draw [small_twist, serre_blue] (0.25, 1.3) [out=-100, in=100] to (0.4, -0.3);
		\draw [thick, red] (0, 1.3) [out=-100, in=100] to (0.1, -0.3);
		\draw [thick, red] (0.5, 0.5) [out=-90, in=90] to (-0.2, -0.3);
		\draw [thick, red] (0.5, 0.5) [out=-90, in=90] to (0.7, -0.3);
		\draw [thick, fill=white] (0.5, 0.45) circle (0.10cm);
		
		\end{scope}
	\end{tikzpicture},
\end{center}
which in informal ordinary string diagrams takes the more accessible form 

\begin{center}
	\begin{tikzpicture}
	\begin{scope}
		\draw[thick, red] (0, 1) to (0.75, 1);
		\draw[thick, red] (0, 1) [out=180, in=180] to (0, 0.5);
		\draw[small_twist, serre_blue] (0, 0.5) to (0.75 ,0.5);
	\end{scope}
	
	\draw[thick, ->] (1.6, 0.7) [out=30, in=150] to node[auto] {\tiny cusp} (2.6, 0.7);
	
	\begin{scope}[xshift=3.5cm, yshift=0.5cm]
		\draw[thick, red] (0, 1) to (0.75, 1);
		\draw[thick, red] (0, 1) [out=180, in=180] to (0, 0.5);
		\draw[small_twist, serre_blue] (0, 0.5) to (0.75 ,0.5);
		\draw[thick, red] (0.75, 0.5) [out=0, in=0] to (0.75, 0);
		\draw[thick, red] (0, 0) to (0.75, 0);
		\draw[thick, red] (0,0) [out=180, in=180] to (0, -0.5);
		\draw[thick, red] (0, -0.5) to (0.75, -0.5);	
	\end{scope}
	\end{tikzpicture}.
\end{center}
Observe that like in the case of the presentation of the framed bordism bicategory, we omit the necessary labellings in unbiased semistrict string diagrams, as being diligent about them would make the pictures cluttered to the point of being incomprehensible. The needed data can be deduced from knowing the domain and codomain of the generating $2$-cell in question, it can also be easily read off the informal string diagrams, which is one of the reasons we include them. 

Additionally, already in the case above it is impossible to deduce whether the blue twist is $q$ or $q^{-1}$, in the first case we have just defined $r(\chi_{+})$, in the second, we have given $r(\chi _{-}^{-1})$. This ambiguity is introduced by design, the definition given above should be understood to apply in both possible cases. This is done mainly for brevity, as all of the generators of framed bordism except Morse ones, which we have already dealt with, are symmetric in the sense that each has its "mirror image" where all twists have opposite signs. To define $r$ on such generators separately would only make the situation seem more lengthy than it really is. 

The case of Serre autoequivalence of negative point passing through a coevaluation is similar, it is in fact dictated by the fact that this crossing point generator is inverse to the one we already considered. We prescribe

\begin{center}
,
\end{center}
which makes it clear that the crossing point inverse relation in $\mathbb{U}_{bord}^{fr}$ will be preserved by $r$. 

The definition of $r$ on the crossing point generators involving evaluation is analogous. In detail, we have 

\begin{center}
	%
.
\end{center}
Again, observe that the crossing point inverse relation will be trivially preserved. 

We are now left with defining $r: \mathbb{U}_{bord}^{fr} \rightarrow \mathbb{U}_{cfd}$ on the witnessing isomorphisms for the adjoint equivalence of $q_{-}$ and $q_{-}^{-1}$. Here the relevant unbiased semistrict string diagrams are a little more formidable, and the informal diagrams will be of great help. We put 

\begin{center}
	%
.
\end{center}
We note that these definitions are clearly mutually inverse as $2$-cells in $\mathbb{U}_{cfd}$, this is of course necessary for the homotopy twist relations to be preserved under $r$, as these state that twist birth-deaths are inverse to each other. This ends our description of $r$.

\begin{lem}
The above definition yields a well-defined strict symmetric monoidal homomorphism $r: \mathbb{U}_{bord}^{fr} \rightarrow \mathbb{U}_{cfd}$, moreover we have $ri = id$ on the nose. 
\end{lem}

\begin{proof}
The second part is clear by construction, to verify the first we have to check if the generating relations of framed bordism hold after applying $r$. 

There are some simplifications that can be made. By \textbf{Lemma \ref{lem:morse_round_twist_relations_are_derivative}}, we do not have to check Morse round twist relations, as they follow from the others. Moreover, twist homotopy relations concerning Serre autoequivalences of the positively framed point are taken by $r$ to generating relations of $\mathbb{U}_{cfd}$. We have also already seen that crossing point birth-death relations are preserved, as are the twist homotopy relations concerning Serre autoequivalences of the negatively framed point which state that twist birth and twist death are inverse to each other.

This leaves us with verifying crossing point colour shift, cusp round twist, twist homotopy relations concerning Serre autoequivalences of the negatively framed point which state that twist birth and twist death form a (co)unit pair and, finally, relations of oriented bordism. We will do these case by case. 

We start with colour shift relations. We will only verify that 

\begin{center}
,
\end{center}
the other variants of the relation are done in an analogous way. Expanding both sides according to our definition of $r$ yields

\begin{center}
	%
.
\end{center}
After changing the order of $2$-cells the left hand side becomes

\begin{center}
	%
,
\end{center}
which is exactly the right hand side after one application of Swallowtail relation and one application of twist homotopy. This ends the argument.

Observe that as much as the proof using unbiased semistrict diagrams is concise, since we do not use any form of labelling it might be rather difficult to see at first glance whether the change of order $2$-cells was allowed. We will hence present it once again using informal string diagrams, where the left hand side is 

\begin{center}
	%
.
\end{center}
On the other hand, the left hand side with changed order of $2$-cells we had drawn above using unbiased semistrict diagrams is 

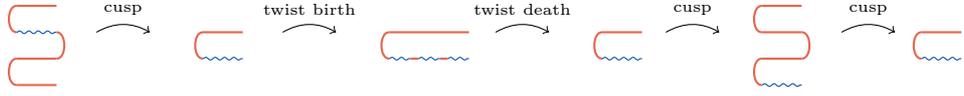
\begin{center}
	%
,
\end{center}
one sees it arises by changing the order of the first cusp with the twist birth and the order of the second cusp with the twist death. The Swallowtail relation then cancels out the last two cusps and twist homotopy cancels twist birth and death. 

We now move on to cusp round twist. Again, we will only show that the relation

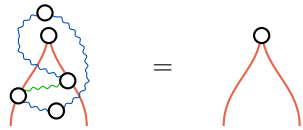
\begin{center}
	%
\end{center}
is preserved by $r$, the arguments for the other ones are analogous. As one sees by applying twist birth/death generators, which are invertible, to the top and bottom, the relation above is equivalent to the perhaps more simple

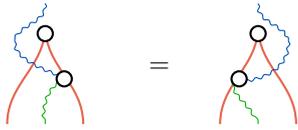
\begin{center}
	%
,
\end{center}
which after applying $r$ reads

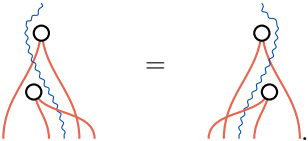
\begin{center}
	%
.
\end{center}
The two $2$-cells on both sides are equal, as the order of the two cusps can be interchanged. This is relatively easy to see in informal diagrams, where the relevant relation is
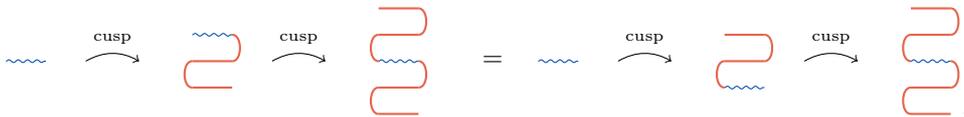
\begin{center}
	%
.
\end{center}

We will now show that $r$ preserves the triangle equations for witnessing isomorphisms between the adjoint equivalence $q_{-} \dashv q_{-}^{-1}$, these are part of twist homotopy relations. It is enough to verify that 

\begin{center}
	%
,
\end{center}
as in the case of adjoint equivalences any of the triangle equations implies the other. After expanding according to the definition, we obtain

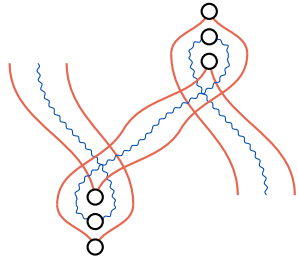
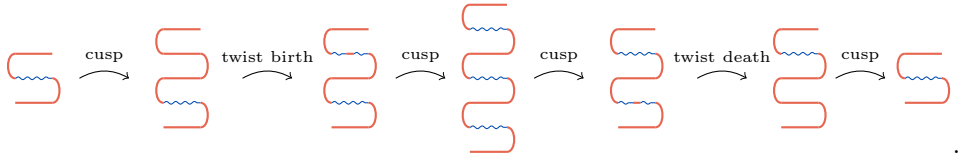
\begin{center}
	%
.
\end{center}
In this case we will just describe the proof and not give it in terms of unbiased semistrict diagrams, as even though the idea is clear, the sequence of modifications turns out to be rather long due to a large number of strings that need to passed through each other. 

Let us focus on the informal diagram given above. One first applies the Swallowtail relation to replace the last cusp by a cusp contracting the other (co)evaluation pair. Once this is done, the last cusp can then be pulled back to the front and canceled with the second cusp. Similarly, what is above the third cusp can be pulled back to the front and canceled with the first cusp with another Swallowtail. Once this is done, we are only left with twist birth and twist death, the relation then follows from triangle equations for the adjoint equivalence $q \dashv q^{-1}$, which is one of the generating relations of $\mathbb{U}_{cfd}$. 

We now move on to relations of oriented bordism, this is luckily the last type we have to consider. The case of Morse birth-death, cusp inverse and Swallowtail isomorphisms is obvious, as $r$ takes them to generating relations of $\mathbb{U}_{cfd}$. More precisely, their images correspond to, respectively, triangle equations for adjunctions $e^{L} \dashv e, c^{L} \dashv c$, equations that state that cusp isomorphisms are inverse to each other and the Swallowtail equations of a coherent dual pair. 

This leaves us with verifying that cusp flip equations also hold in $\mathbb{U}_{cfd}$, the main vehicle of proving this will be the cusp-counits equation defining a coherent fully dual pair. We will show that 

\begin{center}
,
\end{center}
the other cusp flip equation admits an analogous proof. After expanding out and rewriting informally, we need to verify that

\begin{center}
	%
\end{center}
is equal to 

\begin{center}
	%
.
\end{center}
One reasons as follows. Using the cusp-counit equation composed with one of the cusps, the second composite can be rewritten as 
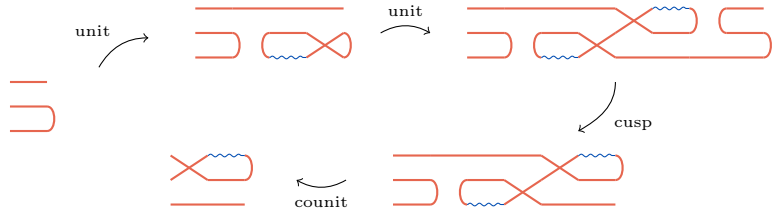
\begin{center}
	%
.
\end{center}
Once this is done, the last occuring counit can be pulled back and canceled with the first unit using triangle equations, we are then left with exactly the first composite. This ends the proof of the lemma.
\end{proof}
Since we have $ri = id$ on the nose by construction, to verify that $i, r$ are pseudoinverse equivalences of symmetric monoidal bicategories we only need to verify that $ir$ is naturally equivalent to the identity. We will do this in the following lemma, which ends our proof of \textbf{Theorem \ref{thm:inclusion_i_from_ucfd_to_u_bord_fr_is_an_equivalence}} and, thus, also of the Cobordism Hypothesis in dimension two. 

\begin{lem}
There is a natural equivalence $ir \simeq id$ of endohomomorphisms of $\mathbb{U}_{bord}^{fr}$ and it can be taken to have identity components on objects.
\end{lem}

\begin{proof}
Since $\mathbb{U}_{bord}^{fr}$ is freely generated and $ir, id$ are strict homomorphisms, to define a natural transformation between these it is enough to do so on generating objects and to give the relevant constraint isomorphisms for each generating $1$-cell. 

To see this, observe that we can compose both homomorphisms with the quotient map $\mathbb{F}_{bord}^{fr} \rightarrow \mathbb{U}_{bord}^{fr}$ from the free symmetric monoidal bicategory on the same generating datum. Since we are then given two strict homomorphisms out of a freely generated symmetric monoidal bicategory, the theory of $(G_{bord}^{fr}, \mathcal{R}_{bord}^{fr})$-shapes developed in \cite{chrisphd} and \textbf{Appendix \ref{appendix_free_monoidal_bicategories}} applies and shows that a semi-strict natural transformation between them is uniquely characterized by its values on generators. Since the target of the relevant homomorphisms is unbiased semistrict, the datum of a natural transformation will then descend to one giving a natural transformation of endohomomorphisms of $\mathbb{U}_{bord}^{fr}$.

We define a natural transformation $s: id \rightarrow ir$ to have identity components on objects, this leaves us with prescribing constraint $2$-cells for generators $e, c, q_{+}, q_{+}^{-1}, q_{-}$ and $q_{-}^{-1}$. 

Since $ir$ takes the first four of these generators to themselves, we define their constraint $2$-cells to be identities. In the case of $q_{-}, q_{-}^{-1}$, we let the constraints $c_{q_{-}}, c_{q_{-}^{-1}}$ be 

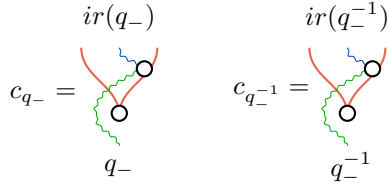
\begin{center}
	\begin{tikzpicture}
	
	\begin{scope}
		\node at (-0.5, 0.35) {$ c_{q_{-}} = $};
	
		\node at (0.5, 1.4) {$ ir(q_{-}) $};
	
		\draw[thick, red] (0,1) [out=-90, in=90] to (0.5, 0.1);
		\draw[thick, red] (1,1) [out=-90, in=90] to (0.5, 0.1);
		
		\draw[small_twist, serre_blue] (0.5, 1) [out=-90, in=45] to (0.8, 0.72);
		\draw[small_twist, serre_green] (0.8, 0.72) [out=-150, in=90] to (0.2, 0.1);
		\draw[small_twist, serre_green] (0.2, 0.1) [out=-90, in=90] to (0.5, -0.3);
		
		\draw [thick, fill=white] (0.83, 0.72) circle (0.10cm);
		\draw [thick, fill=white] (0.5, 0.13) circle (0.10cm);
		
		\node at (0.5, -0.6) {$ q_{-} $};
	\end{scope}
	
	\begin{scope}[xshift=3cm]
		\node at (-0.5, 0.35) {$ c_{q_{-}^{-1}} = $};
	
		\node at (0.5, 1.4) {$ ir(q_{-} ^{-1}) $};
	
		\draw[thick, red] (0,1) [out=-90, in=90] to (0.5, 0.1);
		\draw[thick, red] (1,1) [out=-90, in=90] to (0.5, 0.1);
		
		\draw[small_twist, serre_blue] (0.5, 1) [out=-90, in=45] to (0.8, 0.72);
		\draw[small_twist, serre_green] (0.8, 0.72) [out=-150, in=90] to (0.2, 0.1);
		\draw[small_twist, serre_green] (0.2, 0.1) [out=-90, in=90] to (0.5, -0.3);
		
		\draw [thick, fill=white] (0.83, 0.72) circle (0.10cm);
		\draw [thick, fill=white] (0.5, 0.13) circle (0.10cm);
		
		\node at (0.55, -0.6) {$ q_{-}^{-1} $};
	\end{scope}
	\end{tikzpicture}.
\end{center}
Observe that these are different $2$-cells with different domain and codomain, they just happen to look the same due to our convention of not labelling unbiased semistrict diagrams.

If this data gives rise to a natural transformation, it will clearly be a natural equivalence with identity components on objects. Hence, we only have to verify that this is indeed the case, that is, check that the data above is natural with respect to generating $2$-cells of $\mathbb{U}_{bord}^{fr}$. This is not difficult, as by presentation \textbf{Theorem \ref{thm:presentation_of_the_framed_bordism_bicategory}} the latter is equivalent to the framed bordism bicategory $\mathbb{B}ord_{2}^{fr}$ and so geometric arguments apply. 

Since our natural transformation will have identity components on objects by definition, naturality with respect to some generating $2$-cell $\omega: f_{1} \rightarrow f_{2}$ means that we have an equality

\begin{center}
	\begin{tikzpicture}
		\begin{scope}
			\node (A) at (0, 0) {$ A $};
			\node (B) at (3, 0) {$ B $};
		
			\draw[->] (A) [out=70, in=110] to node[above] {$ ir(f_{1}) $} (B);
			\draw[->] (A) to [out=-17, in=-163] node[above] {$ ir(f_{2}) $} (B);
			\draw[->] (A) [out=-70, in=-110] to node[below] {$ f_{2} $} (B);
		
			\node at (1.5, 0.57) {$ \Downarrow ir(\omega) $};
			\node at (1.5, -0.65) {$ \simeq c_{f_{2}} $};
		\end{scope}
		
		\node at (4, 0) {$ = $};
		
		\begin{scope}[xshift=5cm]
			\node (A) at (0, 0) {$ A $};
			\node (B) at (3, 0) {$ B $};
		
			\draw[->] (A) [out=70, in=110] to node[above] {$ ir(f_{1}) $} (B);
			\draw[->] (A) to [out=-17, in=-163] node[above] {$ f_{1} $} (B);
			\draw[->] (A) [out=-70, in=-110] to node[below] {$ f_{2} $} (B);
		
			\node at (1.5, 0.57) {$ \simeq c_{f_{1}} $};
			\node at (1.5, -0.65) {$ \Downarrow \omega $};
		\end{scope}
	\end{tikzpicture}.
\end{center}
We start with generating (co)units $\mu_{e}, \epsilon_{e}, \mu_{c}, \epsilon_{c}$ and generating cusp isomorphisms $\alpha, \alpha^{-1}, \beta, \beta^{-1}$. In this case, the equation above holds obviously, as the relevant generating $2$-cells do not have $q_{-}$ or $q_{-}^{-1}$ in their (co)domains and so the relevant constraint isomorphisms are identities by definition. 

We now move on to the rest of the generating $2$-cells of $\mathbb{U}_{bord}^{fr}$, we will see that the above equation can be proven already from the facts that these generators are invertible and that the bordisms corresponding to their domain and codomain always have an interval as the underlying manifold. This holds true for all of the relevant generators by a quick check. 

Indeed, observe that since the generators are invertible, both sides of the naturality equation are invertible aswell. It is hence enough to prove that either their domain or codomain has no non-trivial automorphisms, this will force both sides to be equal. As we have an equivalence $\mathbb{U}_{bord}^{fr} \simeq \mathbb{B}ord_{2}^{fr}$, we can prove this geometrically. 

Let $\zeta$ be an invertible framed $2$-autobordism of some $1$-bordism $f$ whose underlying $1$-manifold is the interval. As an oriented, compact surface whose boundary components we know, $\zeta$ must be homeomorphic to a genus $g$ oriented surface with one disk cut out. As it is invertible as a $2$-bordism, van Kampen theorem implies that it must in fact be a disk. 

It follows by uniqueness of framings on disks up to isotopy relative to the boundary, which we've proven as \textbf{Lemma \ref{lem:uniqueness_of_framings_on_generators}}, that $\zeta$ must be the identity. This ends the argument.
\end{proof}

\subsection{Future directions}

In this section we will highlight some possible extensions of the current work. It is not our goal to be thorough, for a beautiful exposition of possible applications of the Cobordism Hypothesis in general see \cite[Ch. 4]{lurietqfts}. Rather, we will focus on those results that seem particularly accessible. 

\begin{so2_action_on_fully_dualizable_objects*}

As we have explained in the chapter devoted to fully dualizable objects in symmetric monoidal bicategories, any such object $X \in \mathbb{M}$ admits a canonical up to isomorphism invertible endomorphism, the Serre autoequivalence. We claimed that one way to define it would be to use the Cobordism Hypothesis, by observing that we have an equivalence 

\begin{center}
$\textbf{SymMonBicat}(\mathbb{B}ord_{2}^{fr}, \mathbb{M}) \simeq \mathbb{K}(\mathbb{M}^{fd})$ 
\end{center}
between the space of fully dualizable objects and the bicategory of framed topological field theories with values in $\mathbb{M}$. As the group $SO(2)$ acts on framed bordism through change of framing, it also acts on the left hand side of the above equivalence and hence also on the right. The generator of $\pi_{1}(SO(2))$ is taken by this action to an invertible endohomomorphism of $\mathbb{K}(\mathbb{M}^{fd})$, its component at a given fully dualizable object $X$ is precisely its Serre autoequivalence.

We did not use this definition of the latter, instead giving a direct description in terms of string diagrams. We then verified some of its properties, like naturality, by hand. We believe it would be interesting to assemble all these components into an honest action of $SO(2)$ on the space $\mathbb{K}(\mathbb{M}^{fd})$, which we did not do. 

Note that since we have proven the Cobordism Hypothesis, one could in principle construct an action of $SO(2)$ on $\mathbb{B}ord_{2}^{fr}$ and then chase it through the relevant zig-zag of equivalences described in this work. It is, unfortunately, rather likely that such an action would be difficult to compute with in practice. We believe one should be able to construct an action directly on the space $\mathbb{K}(\mathbb{M}^{fd})$. 

As explained in \cite{chris_dualizability_in_low_dim_category_theory}, such an action can be defined as a map of classifying spaces $BSO(2) \rightarrow BAut(\mathbb{K}(M^{fd}))$. Note that we have $BSO(2) \simeq \mathbb{CP}^{\infty}$, which admits a CW-structure with even-dimensional skeleta given by $\mathbb{CP}^{n}$. Observe that since $BAut(\mathbb{K}(M^{fd}))$ is a $3$-type, to give such a map is the same as to give an element in $S \in \pi_{2}(BAut(\mathbb{K}(M^{fd}))$ that extends to $\mathbb{CP}^{2}$, that is, is null-homotopic when composed with the Hopf fibration $S^{3} \rightarrow S^{2}$. 

By unwinding definitions, one sees that to specify $S$ is to give an invertible endohomomorphism of $\mathbb{K}(\mathbb{M}^{fd})$, this element in our case should be assembled from Serre autoequivalences. On the other hand, $S$ composed with the Hopf fibration will give a datum of an invertible modification of the identity of $\mathbb{K}(\mathbb{M}^{fd})$, the condition that this element is null-homotopic will then translate into an equation at the level of $2$-cells. The relevant equation is written down explicitly in \cite{chris_dualizability_in_low_dim_category_theory}.

We believe that details of the argument given above can be worked out directly in the setting of symmetric monoidal bicategories, similarly to the way we explicitly verified some other expected properties of Serre autoequivalence without invoking the Cobordism Hypothesis. The bicategory $\mathcal{C}ohFullyDualPair(\mathbb{M})$ of coherent fully dual pairs could be of possible help here, as on one hand it is equivalent to $\mathbb{K}(\mathbb{M}^{fd})$ in a very explicit way, on the other, objects of $\mathcal{C}ohFullyDualPair(\mathbb{M})$ come encoded with a choice of a Serre autoequivalence.
\end{so2_action_on_fully_dualizable_objects*}

\begin{cobordism_hypothesis_for_oriented_bordism*}
One of the great accomplishments of Jacob Lurie was to describe not only the space of framed topological field theories, which was conjectured by Baez-Dolan, but in fact classify all fully extended topological field theories defined on $G$-manifolds, where $G$ is some fixed type of tangential structure, see \cite{lurietqfts}. 

In the setting of symmetric monoidal bicategories, the expected answer would take the form of an equivalence

\begin{center}
$\textbf{SymMonBicat}(\mathbb{B}ord_{2}^{G}, \mathbb{M}) \simeq (\mathbb{K}(\mathbb{M}^{fd}))^{hG}$,
\end{center}
where $G$ is a topological group over $O(2)$, $\mathbb{B}ord_{2}^{G}$ is the bicategory of manifolds equipped with a $G$-structure and $(\mathbb{K}(\mathbb{M}^{fd}))^{hG}$ is the space of homotopy fixed points of $\mathbb{K}(\mathbb{M}^{fd})$ with respect to the action of $G$ derived from the $O(2)$-action.

It is difficult to say whether arguments similar to those in \cite{lurietqfts} could be used directly in our setting to derive all of these "$G$-Cobordism Hypotheses" at once. 

One accessible variant should be the oriented and unoriented cases, as for these types of bordism bicategories we have explicit presentations derived by Christopher Schommer-Pries in \cite{chrisphd}. Moreover, as the current work is an extension of these results, both presentations are rather similar to the one we have for the framed bordism bicategory. We believe it is possible to prove these cases of the two-dimensional Cobordism Hypothesis directly, using methods similar to ours. 

The oriented bordism bicategory $\mathbb{B}ord_{2}^{or}$ looks particularly promising, as one sees that its presentation given in \cite{chrisphd} is essentially the definition of a coherent fully dual pair with Serre autoequivalences $q, q^{-1}$ artificially forced to be identities.

On the other hand, the Cobordism Hypothesis predicts that the space of oriented two-dimensional topological field theories is equivalent to $SO(2)$-homotopy fixed points of the space of fully dualizable objects. As the $SO(2)$-action essentially encodes the Serre autoequivalence, taking homotopy fixed points should have a similar effect of forcing it to be trivial. 
\end{cobordism_hypothesis_for_oriented_bordism*}

\FloatBarrier
\appendix

\section{Freely generated bicategories}
\label{appendix_free_monoidal_bicategories}

In this chapter we present the theory of \emph{freely generated monoidal bicategories}, which we use extensively in the parts of the current work related to coherence of dualizable objects. Our approach is modeled on the language of \emph{computadic symmetric monoidal bicategories} developed by Christopher Schommer-Pries and is formally analogous. 

We also give a few technical results related to what we call \emph{bicategories of $G$-shapes}, which are used to model bicategories of homomorphisms out of freely generated bicategories. These are rather simple-minded in nature and are used in several places of the current work to simplify book-keeping. 

In a pinch, one can also use the ideas presented here to gain some understanding of the symmetric monoidal case, which we don't include as it is virtually identical. A definite reference for the subject is \cite{chrisphd}, but one should take note that we use a slightly different language, with a dictionary given in \textbf{Table \ref{tab:comparison_of_notions}} below. 

\begin{table}[!htbp]
\caption{Dictionary between languages related to freely generated bicategories}
\begin{tabular}{c || c}
	\textbf{Current work} & \cite{chrisphd} \\

	\hline
	\hline
	generating datum & computad \\
	\hline 
	$i$-truncated generating datum & $i$-computad \\
	\hline
	freely generated $\Box$ bicategory & computadic $\Box$ bicategory \\
	\hline
	$G$-shape & $G$-data
\end{tabular}
\label{tab:comparison_of_notions}
\end{table}
Note that $\Box$ could stand for any of the kinds of bicategories considered in the current work, which are \emph{monoidal}, \emph{symmetric monoidal} and \emph{unbiased semistrict symmetric monoidal}. 

\subsection{Freely generated monoidal bicategories}

There are different kinds of data one can use to generate a monoidal bicategory, like sets, category-enriched graphs and bicategories themselves. Many of these have been studied in the literature, see \cite{gurski_coherence_in_threedimensional_cat_theory}. 

These approaches are, however, not quite sufficient for our purposes. We would like to allow the generating 1-cells to have domains and codomains that are only \emph{consequences} of the generating data at the level of objects and similarly for 2-cells. This can make the presentations of bicategories defined using such data much smaller and more readable.  

For example, in the presentation of the unoriented bordism bicategory one has only one generating object, a point, even though the domain of the generating morphism given by a right elbow consists of a disjoint union of two such points. Using this approach, this disjoint union, which corresponds to the tensor product of a point with itself, does not have to be added separately as a generating object, avoiding redundancy. 

\begin{defin}
A \textbf{$0$-truncated generating datum for a monoidal bicategory} consists of a set $G_{0}$. Given such a generating datum, we inductively define the set $BW(G_{0})$ of \textbf{binary words in $G_{0}$} by declaring that 

\begin{itemize}
\item the symbol $I$ is a binary word,
\item the symbol $X$ is a binary word for all $X \in G_{0}$,
\item if $X, Y$ are binary words, then so is $(X \otimes Y)$.
\end{itemize}
\end{defin}

\begin{notation*}
In what follows we will talk and construct objects with \emph{source} and \emph{target} maps into some other kind of objects. We will then use the function notation to talk about the values of these maps, even when there's no category in sight, so that for example "$f: A \rightarrow B$" will just mean "$f$ has source $A$ and target $B$".
\end{notation*}

\begin{defin}
A \textbf{$1$-truncated generating datum for a monoidal bicategory} consists of a tuple of sets $(G_{0}, G_{1})$ together with \textbf{source} and \textbf{target} maps $s, t: G_{1} \rightarrow BW(G_{0})$. If $(G_{0}, G_{1})$ is a $1$-truncated generating datum, then we inductively define the set $BW(G_{1})$ of \textbf{binary words in $G_{1}$}, also with source and target in $BW(G_{0})$, by declaring that

\begin{itemize}
\item if $X$ is a binary word in $G_{0}$, then the symbol $id_{X}$ is a binary word in $G_{1}$ with source and target $X$,
\item If $X, Y, Z$ are binary words in $G_{0}$ and $x$ is a symbol taken from \textbf{Table \ref{tab:symbols_for_binary_words_in_1_morphisms}}, then $x$ is a binary word in $G_{1}$ with source and target given according to the table,
\item if $X, Y, Z$ are binary words in $G_{0}$ and $x$ is a symbol taken from table \textbf{Table \ref{tab:symbols_for_binary_words_in_1_morphisms}}, then $x^{\bullet}$ is a binary word in $G_{1}$ with source and target opposite to the one given in the table,
\item If $f \in G_{1}$, then the symbol $f$ is a binary word with source $s(f)$ and target $t(f)$. 
\end{itemize}
Moreover, we inductively define the set $BS(G_{1})$ of \textbf{binary sentences in $G_{1}$}, again with source and target in $BW(G_{0})$, by declaring that 

\begin{itemize}
\item If $f \in BW(G_{1})$ is a binary word, then $f$ is a binary sentence with the same source and target,
\item If $f, g \in BS(G_{1})$ are binary sentences such that the target of $f$ matches the source of $g$, then $(g) \circ (f)$ is a binary sentence with source $s(f)$ and target $t(g)$,
\item If $f, g \in BW(G_{1})$ are binary sentences, then $f \otimes g$ is a binary sentence with source $s(f) \otimes s(g)$ and target $t(f) \otimes t(g)$.
\end{itemize}
\end{defin}

\begin{table}[!htbp]
\caption{Symbols for binary words in $G_{1}$}
\begin{tabular}{l || c | c}
	\textbf{Symbol} & \textbf{Source} & \textbf{Target} \\

	\hline
	$l_{X}$ & $I \otimes X$ & $X$ \\
	\hline 
	$r_{X}$ & $X \otimes I$ & $X$ \\
	\hline
	$a_{X, Y, Z}$ & $(X \otimes Y) \otimes Z$ & $X \otimes (Y \otimes Z)$
\end{tabular}
\label{tab:symbols_for_binary_words_in_1_morphisms}
\end{table}

\begin{defin}
A \textbf{free generating datum for a monoidal bicategory} consists of a tuple $(G_{0}, G_{1}, G_{2})$, where $(G_{0}, G_{1})$ is a $1$-truncated generating datum, together with \textbf{source} and \textbf{target} maps $s, t: G_{2} \rightarrow BS(G_{1})$ satisfying the globularity condition $s(s(\zeta)) = s(t(\zeta)$ and $t(s(\zeta)) = t(t(\zeta))$ for all $\zeta \in G_{2}$. 

If $(G_{0}, G_{1}, G_{2})$ is a free generating datum, then we inductively the set of \textbf{binary words in $G_{2}$} with source and target in $BS(G_{1})$, by declaring that

\begin{itemize}
\item if $X, Y, Z$ are binary words in $G_{0}$ and $x$ is a symbol taken from \textbf{Table \ref{tab:symbols_for_binary_words_in_1_morphisms}}, then $\mu_{x}$ is a binary word in $G_{2}$ with source $id _{s(x)}$ and target $x^{\bullet} \circ x$ and $\epsilon_{x}$ is a binary word in $G_{2}$ with source $x \circ x^{\bullet}$ and target $id_{t(x)}$,
\item if $f: A \rightarrow B$, $f^\prime: B \rightarrow C$, $f^{\prime \prime}: C \rightarrow D$ are binary sentenes in $G_{1}$ and $x$ is a bicategorical constraint symbol from \textbf{Table \ref{tab:bicategorical_constraint_symbols_for_binary_words_in_2_morphisms}}, then $x$ is a binary word in $G_{2}$ with source and target as given in the table and $x^{-1}$ is a binary word with source and target opposite to the given ones,
\item if $f: B \rightarrow C$, $f^\prime: A \rightarrow B$, $g: Y \rightarrow Z$, $g^\prime: X \rightarrow Y$ and $h: P \rightarrow Q$ are binary sentences in $G_{1}$ and $x$ is a monoidal morphism constraint symbol from \textbf{Table \ref{tab:monoidal_morphism_constraint_symbols_for_binary_words_in_2_morphisms}}, then $x$ is also a binary word in $G_{2}$ with source and target as given in the table and $x^{-1}$ is a binary word with source and target opposite to the given ones,
\item if $A, B, C, D$ are binary words in $G_{0}$ and $x$ is a monoidal object constraint symbol from \textbf{Table \ref{tab:object_constraint_symbols_for_binary_words_in_2_morphisms}}, then $x$ is also a binary word in $G_{2}$ with source and target as given in the table and $x^{-1}$ is a binary word with source and target opposite to the given ones,
\item if $\zeta$ in $G_{2}$, then the symbol $\zeta$ is binary word.
\end{itemize}
Consequently, we inductively define the set $BS(G_{2})$ of \textbf{binary sentences in $G_{2}$} with source and target maps in $BS(G_{1})$ by declaring that

\begin{itemize}
\item if $\alpha$ is a binary word in $G_{2}$, then it's also a binary sentence with the same source and target,
\item if $\alpha, \beta$ are binary sentences in $G_{2}$ such that $t(t(\alpha)) = s(s(\beta))$, then $\beta * \alpha$ is a binary sentence in $G_{2}$ with source $s(\beta) \circ s(\alpha)$ and target $t(\beta) \circ t(\alpha)$,
\item if $u, v$ are binary sentences in $G_{2}$, then $u \otimes v$ is binary sentence with source $s(u) \otimes s(v)$ and target $t(u) \otimes t(v)$, 
\item if $\alpha _{0}, \alpha_{1}, \ldots, \alpha _{k}$ is a composable sequence of sentences, that is, we have $s(\alpha_{i}) = t(\alpha_{i-1})$ for all $1 \leq i \leq k$, then $p_{k} p_{k-1} \ldots p_{0} = p_{k} \circ p_{k-1} \circ \ldots \circ p_{0}$ is a binary sentence in $G_{2}$ with target $t(p_{k})$ and source $s(p_{0})$.
\end{itemize}
\end{defin}

\begin{table}[!htbp]
\caption{Bicategorical constraint symbols for binary words in $G_{2}$}
\begin{tabular}{l || c | c}
	\textbf{Symbol} & \textbf{Source} & \textbf{Target} \\

	\hline
	$id_{f}$ & f & f \\
	\hline
	$a^{c} _{f, f^\prime, f^{\prime \prime}}$ & $(f \circ f^\prime) \circ f^{\prime \prime}$ & $f \circ (f^\prime \circ f^{\prime \prime})$ \\
	\hline
	$ r^{c}_{f}$ & $f \circ id_{A}$ & $f$ \\
	\hline
	$ l^{c}_{f}$ & $id_{B} \circ f$  & $f$
\end{tabular}
\label{tab:bicategorical_constraint_symbols_for_binary_words_in_2_morphisms}
\end{table}

\begin{table}[!htbp]
\caption{Morphism constraint symbols for binary words in $G_{2}$}
\begin{tabular}{l || c | c}
	\textbf{Symbol} & \textbf{Source} & \textbf{Target} \\
	
	$ \phi ^{\otimes} _{(f, g), (f^\prime, g^\prime)}$ & $(f \otimes g) \circ (f^\prime \circ g^\prime)$ & $(f \circ f^\prime) \otimes (g \circ g^\prime)$ \\
	\hline
	$ \alpha _{f, g, h} $ & $a_{C, Z, Q} \circ (f \otimes g) \otimes h$ & $f \otimes (g \otimes h) \circ a_{B, Y, P}$ \\
	\hline
	$ l_{f} $ & $l_{C} \circ (id_{I} \otimes f)$ & $f \circ l_{B}$ \\
	\hline
	$ r_{f} $ & $r_{C} \circ(f \otimes id_{I})$ & $f \circ l_{B}$
\end{tabular}
\label{tab:monoidal_morphism_constraint_symbols_for_binary_words_in_2_morphisms}
\end{table}

\begin{table}[!htbp]
\caption{Object constraint symbols for binary words in $G_{2}$}
\begin{tabular}{l || c | c}
	\textbf{Symbol} & \textbf{Source} & \textbf{Target} \\

	$\phi ^{\otimes} _{(A, B)}$ & $id_{A \otimes B}$ & $id_{A} \otimes id_{B}$ \\
	\hline
	$\pi _{A, B, C D}$ & $((id_{A} \otimes a_{B, C, D}) \circ a_{A, B \otimes C, D}) \circ (a_{A, B, C} \otimes id)$ & $a_{A, B, C \otimes D} \circ a_{A \otimes B, C ,D}$ \\
	\hline
	$\mu _{A, B}$ & $((id_{A} \otimes l_{B}) \circ a_{A, I, B}) \circ (r_{A} \otimes id_{B})$ & $id_{A \otimes B}$ \\
	\hline 
	$\lambda _{A, B}$ & $l_{A} \otimes id_{B}$ & $l_{A \otimes B}\circ_{1, A, B}$ \\
	\hline
	$ \rho _{A, B}$ & $id_{A} \otimes r_{B}$ & $a _{A, B, 1} \circ_{A \otimes B}$ 
\end{tabular}
\label{tab:object_constraint_symbols_for_binary_words_in_2_morphisms}
\end{table}

We are now ready to define a free monoidal bicategory generated by a free generating datum.

\begin{defin}
Let $G = (G_{0}, G_{1}, G_{2})$ be a free generating datum for a monoidal bicategory. We define $\mathbb{F}(G)$, the \textbf{monoidal bicategory generated by $G$} as follows.

\begin{itemize}
\item The objects of $\mathbb{F}(G)$ are precisely the binary words in $G_{0}$. 
\item The morphisms of $\mathbb{F}(G)$ are precisely the binary sentences in $G_{1}$, with source and target as defined.
\item The 2-cells of $\mathbb{F}(G)$ are \emph{equivalence classes} of sentences in $G_{2}$, with source and target as defined.
\end{itemize}
The equivalence relation $\sim$ on binary sentences is the smallest equivalence relation such that

\begin{itemize}
\item if $x$ is a symbol from \textbf{Tables \ref{tab:bicategorical_constraint_symbols_for_binary_words_in_2_morphisms}}, \textbf{\ref{tab:monoidal_morphism_constraint_symbols_for_binary_words_in_2_morphisms}} or \textbf{\ref{tab:object_constraint_symbols_for_binary_words_in_2_morphisms}}, then $xx^{-1} \sim id_{t(x)}$ and $x^{-1} x \sim id_{s(x)}$,
\item the $2$-cells $a^{c}, r^{c}, l^{c}, \phi^{\otimes}_{(f, g), (f^\prime, g^\prime)}, \phi^{\otimes} _{(X, Y)}, a_{f, g, h}, l_{f}, r_{f}$ are components of a natural transformation, that is, in the relevant naturality pasting diagrams the two different composites are equivalent
\item the axioms of a monoidal bicategory hold,
\item the equivalence relation is closed under $\otimes$, horizontal composition $*$ and vertical composition $\circ$. 
\end{itemize}
The structure of a monoidal bicategory is defined formally. More specifically, composite of two binary words $f, f^\prime$ is $f \circ f^\prime$, horizontal composite of equivalence classes $\alpha, \beta$ is $\beta * \alpha$, their vertical composite is $\beta \circ \alpha$, similarly one defines the monoidal product of objects, morphisms and $2$-cells using the symbol $\otimes$. All the required axioms that this structure has to satisfy are enforced by the equivalence relation on binary sentences in $G_{2}$. 
\end{defin}
We will now discuss relations on monoidal bicategories, the only kind we will consider is of relations between paralell $2$-cells. 

\begin{defin}
A \textbf{class of relations} $\mathcal{R}$ on a monoidal bicategory $\mathbb{M}$ consists of a relation on the set of its $2$-cells such that if $\alpha \sim _{\mathcal{R}} \beta$, then the source and target of $\alpha, \beta$ coincide.  If $\mathcal{R}$ is a class of relations, then we define its \textbf{closure} $c(\mathcal{R})$ to be the smallest equivalence class of relations that contains $\mathcal{R}$ and also 

\begin{itemize}
\item is closed under vertical composition, that is, if $\alpha \sim _{c(\mathcal{R})} \alpha ^\prime$ and $\beta \sim _{c(\mathcal{R})} \beta ^\prime$ and $\alpha, \beta$ are vertically composable, then $\beta \alpha \sim _{c(\mathcal{R})} \beta^\prime \alpha^\prime$,
\item is closed under whiskering, that is, if $\alpha \sim _{c(\mathcal{R})} \alpha^\prime$, then also $f * \alpha \sim _{c(\mathcal{R})} f * \alpha^\prime$ and $\alpha * g \sim _{c(\mathcal{R})} \alpha^\prime * f$ whenever this makes sense, ie. when $f$ and $s(f)$ or $s(f)$ and $g$ are composable,
\item is closed under tensor product, that is, if $\alpha \sim _{c(\mathcal{R})} \alpha^\prime$ then also $f \otimes \alpha \sim _{c(\mathcal{R})} f \otimes \alpha^\prime$ and $\alpha \otimes f \sim _{c(\mathcal{R})} \alpha^\prime \otimes f$.
\end{itemize}
\end{defin}

\begin{defin}
Suppose $\mathbb{M}$ is a monoidal bicategory and $\mathcal{R}$ is a class of relations on it. Then we define the \textbf{quotient monoidal bicategory} $\mathbb{M} / \mathcal{R}$ to have the same objects and morphisms as $\mathbb{M}$ and with $2$-cells given by equivalence classes of $2$-cells of $\mathbb{M}$ under the closure $c(\mathcal{R})$. 
\end{defin}
Observe that there is a canonical strict quotient homomorphism $\pi_{\mathcal{R}}: \mathbb{M} \rightarrow \mathbb{M} / \mathcal{R}$ given by identity on objects and morphisms and by passing to equivalence classes on $2$-cells. We will now establish its universal property.

\begin{prop}
For any monoidal bicategory $\mathbb{N}$, the induced by precomposition strict homomorphism 
\[ \pi^{*}: \textbf{MonBicat}(\mathbb{M} / \mathcal{R}, \mathbb{N}) \rightarrow \textbf{MonBicat}(\mathbb{M}, \mathbb{N}) \]
between bicategories of monoidal homomorphisms identifies $\textbf{MonBicat}(\mathbb{M} / \mathcal{R}, \mathbb{N})$ as the full subbicategory of $\textbf{MonBicat}(\mathbb{M}, \mathbb{N})$ on those homomorphism $\phi: M \rightarrow N$ satisfying $\phi(\alpha) = \phi(\beta)$ for all $\alpha, \beta$ such that $\alpha \sim _{\mathcal{R}} \beta$. 
\end{prop}

\begin{proof}
Clearly any homomorphism coming from precomposition with $\pi$ will have this property, as we have $\pi(\alpha) = \pi(\beta)$ if $\alpha \sim _{\mathcal{R}} \beta$. Conversely, if $\phi: \mathbb{M} \rightarrow \mathbb{N}$ is a monoidal homomorphism such that $\phi(\alpha) = \phi(\beta)$ if $\alpha \sim _{\mathcal{R}} \beta$, then it must also necessarily satisfy $\phi(\alpha) = \phi(\beta)$ if $\alpha \sim _{c(\mathcal{R})} \beta$, as is reasonably easy to verify. This shows that $\pi^{*}$ is an inclusion onto the right subbicategory.

The homomorphism $\pi^{*}$ is fully faithful by direct inspection of a definiton of a monoidal transformation and modification, which have no components related to $2$-cells at all.
\end{proof}
We can now define a generating datum for a monoidal bicategory, which will consists of a free generating datum together with a class of relations. The freely generated monoidal bicategory will be the obvious quotient. 

\begin{defin}
A \textbf{generating datum for a monoidal bicategory} or a \textbf{presentation} consists of a tuple $(G, \mathcal{R})$, where $G$ is free generating datum $G = (G_{0}, G_{1}, G_{2})$ and $\mathcal{R}$ is a class of relations on $\mathbb{F}(G)$. The \textbf{free symmetric monoidal bicategory} $\mathbb{F}(G, \mathcal{R})$ generated by it is the quotient $\mathbb{F}(G) / \mathcal{R}$. 
\end{defin}

\begin{rem}
Observe that the usage of the word \emph{free} is different than in say, the theory of groups. Here we say that a monoidal bicategory is freely generated if it admits a presentation at all, the presentation doesn't have to necessarily be \emph{relation-free}. This is motivated by the fact that such monoidal bicategories are still rather special, unlike "groups admitting a presentation", which is simply all of them.
\end{rem}

\subsection{Bicategories of shapes}

The purpose of this section is to present the theory of \emph{shapes} associated to a generating datum. In more detail, we will see that to any generating datum $(G, \mathcal{R})$ and any ambient monoidal bicategory $\mathbb{M}$ one can associate a bicategory $\mathbb{M}(G, \mathcal{R})$ of \textbf{$(G, \mathcal{R})$-shapes with values in $\mathbb{M}$}. 

This bicategory serves as a compact model for $\textbf{MonBicat}(\mathbb{F}(G, \mathcal{R}), \mathbb{M})$, the bicategory of monoidal homomorphisms $\mathbb{F}(G, \mathcal{R}) \rightarrow \mathbb{M}$ out of the freely generated monoidal bicategory on the given datum. 

More precisely, $\mathbb{M}(G, \mathcal{R})$ is by construction \emph{isomorphic} to a suitably defined bicategory $\textbf{MonBicat}_{str.}(\mathbb{F}(G, \mathcal{R}), \mathbb{M})$ of strict homomorphisms. The two homomorphism bicategories are equivalent by the so called \emph{cofibrancy theorem}.

It follows that intuitively a $(G, \mathcal{R})$-shape encodes a strict homomorphism from the free monoidal bicategory into some other monoidal bicategory $\mathbb{M}$. A strict homomorphism is given by its values on the objects, morphisms and $2$-cells, since all the coherence data is by definition trivial. However, in the case of mapping from a freely generated monoidal bicategory, it is enough to give the values of the strict homomorphism on the generators, this list of values is what constitutes a $(G, \mathcal{R})$-shape.

Thus, a $(G, \mathcal{R})$-shape $P$ will be given by a triple of maps $P_{i}: G_{i} \rightarrow \mathbb{M}_{i}$, where $\mathbb{M}_{i}$ is the set of $i$-cells of $\mathbb{M}$, the triple is required to be globular, that is, respect source and targets. 

To make sense of this, one needs to observe that given such a triple there is a canonical extension of $P_{0}$ to all of $BW(G_{0})$, the set of binary words in generating objects of $G$, given by inductively declaring that $P_{0}(I) = I_{\mathbb{M}}$, the unit of $\mathbb{M}$ and that $P_{0}(X \otimes Y) = P_{0}(X) \otimes _{\mathbb{M}} P_{0}(Y)$. Similarly, there are extensions of $P_{1}$ to the set $BS(G_{1})$ of binary sentences in $G_{1}$ and of $P_{2}$ to the set $BS(G_{2})$ of binary sentences in $G_{2}$ given by evaluation of expressions using the structure of the monoidal bicategory $\mathbb{M}$. 

Additionally, if $\mathcal{R}$ is non-empty, there are additional conditions that are needed to ensure that the homomorphism will in fact factor through $\mathbb{F}(G, \mathcal{R})$, which as the reader recalls was defined as a quotient.

\begin{defin}
Let $G$ be a free generating datum for a monoidal bicategory. A triple $P$ of maps $P_{i}: G_{i} \rightarrow \mathbb{M}_{i}$, where $\mathbb{M}_{i}$ is the set of $i$-cells of $\mathbb{M}$ is a \textbf{$G$-shape with values $\mathbb{M}$} if the globularity conditions

\begin{itemize}
\item $s(P_{1}(f)) = P_{0}(s(f))$, $t(P_{1}(f)) = P_{0}(t(f))$ for all $f \in G_{1}$
\item $s(P_{2}(\alpha)) = P_{1}(s(\alpha))$, $t(P_{2}(\alpha)) = P_{1}(t(\alpha))$ for all $\alpha \in G_{2}$
\end{itemize}
are satisfied after the canonical extension of $P_{0}$ to the set of binary words in $G_{0}$, of $P_{1}$ to the set of binary sentences in $G_{1}$ and of $P_{2}$ to the set of binary sentences in $G_{2}$. If $P$ is a $G$-shape with values in $\mathbb{M}$, then one defines its \textbf{associated strict homomorphism} $\mathbb{F}(G) \rightarrow \mathbb{M}$, by abuse of language also denoted by $P$, to be the unique strict homomorphism given on $i$-cells by $P_{i}$.
\end{defin}

\begin{defin}
If $(G, \mathcal{R})$ is a generating datum for a monoidal bicategory, then a $G$-shape is called a \textbf{$(G, \mathcal{R})$-shape} if and only if the associated strict homomorphism factors through $\mathcal{F}(G) / \mathcal{R}$, that is, exactly when we have $P(\alpha) = P(\beta)$ for $2$-cells $\alpha, \beta$ such that $\alpha \sim _{\mathcal{R}} \beta$. 
\end{defin}
Observe that \emph{being a $(G, \mathcal{R})$-shape} is a \emph{property}, the bicategory of $(G, \mathcal{R})$-shapes will be defined as the full subbicategory of the bicategory of all $G$-shapes. Thus, when defining maps and $2$-cells, we will restrict to the special case of $G$-shapes, which corresponds to $\mathcal{R} = \emptyset$, as this is sufficient.  

The intuition is that a map of $G$-shapes encodes a natural transformation between their associated strict homomorphisms. Againm, since the source is a free monoidal bicategory, it is enough to give the values on the generators.

\begin{defin}
A \textbf{homomorphism} $w: P \rightarrow P^\prime$ of $G$-shapes with values in $\mathbb{M}$ consists of a tuple of maps $w_{0}: G_{0} \rightarrow \mathbb{M}_{1}$ and $w_{1}: G_{1} \rightarrow \mathbb{M}_{2}$. We require that $w_{0}(X): P(X) \rightarrow P^{\prime}  (X)$ for all $X \in G_{0}$ and that $w_{1}(f)$ is an isomorphism fitting the diagram
\begin{center}

\end{center}
for all $f \in G_{1}$, $f: A \rightarrow B$, where $w_{0}$ is again implicitly extended to the set of binary words in $G_{0}$. Additionally,we require naturality with respect to all $\alpha \in G_{2}$. In other words, that for all such $\alpha: f_{1} \rightarrow f_{2}$ we have

\begin{center}
	%
.
\end{center}
Here $w_{1}$ was uniquely extended to all binary paragrahs in $G_{1}$ in such a way that this family of $2$-cells plays a role of constraint isomorphisms for a natural transformation between $\widetilde{P}$ and $\widetilde{P}^{\prime}$. By abuse of language, we denote the constructed transformation also by $w$ and call it the \textbf{associated semi-strict transformation}. 
\end{defin}

We call the associated transformation \emph{semi-strict}, because it is not strict as a natural transformation between homomorphisms of bicategories, but the additional data making it into a monoidal transformation \emph{is} trivial. We could in principle model only strict transformations by requiring additionally that all $w_{1}(f)$ from the definition above are identities, but the bicategory we would obtain in this way would be too small. 

The theme of considering only strict homomorphisms together with non-strict transformations is a recurrent one in category theory. For example, the natural enrichment of the category $2-Cat$ of strict two-categories in itself using the notion of the internal $2$-category of strict functors and all transformations leads to the \textbf{Gray} tensor product, which is already "weak enough" to model all tricategories.

\begin{defin}
If $w, v: P \rightarrow P^\prime$ are homomorphisms of $G$-shapes, then a \textbf{transformation} $\zeta: w \rightarrow v$ consists of a map $\zeta_{0}: G_{0} \rightarrow \mathbb{M}_{2}$ such that for all $X \in G_{0}$ we have $\zeta_{0}(X): w(X) \rightarrow v(X)$. Additionally, we require naturality with respect to all $f \in G_{1}$. In other words, that for all such $f: A \rightarrow B$ we have

\begin{center}
	\begin{tikzpicture}
		\node (PA) at (0,0) {$ P(A) $};
		\node (PpA) at (0, 2) {$ P^\prime(A) $};
		\node (PB) at (2,0) {$ P(B) $};
		\node (PpB) at (2,2) {$ P^\prime(B) $};
		
		\draw [->] (PB) to node[right] {$ w(B) $} (PpB);
		\draw [->] (PA) to node[auto] {$ P(f) $} (PB);
		\draw [->] (PpA) to node[auto] {$ P^\prime(f) $} (PpB);
		\draw [->] (PA) to [out=135, in=-135] node[auto] {$ v(A) $} (PpA);
		\draw [->] (PA) to [out=45, in=-45] node[right] {$ w(A) $} (PpA);
		
		\node at (0, 1.3) {$ \Rightarrow $};
		\node at (0,0.9) {$ \zeta(A) $};
		\node at (1.3, 1.5) {$ \simeq w(f) $};
		
		\node (rPA) at (4.5,0) {$ P(A) $};
		\node (rPpA) at (4.5, 2) {$ P^\prime(A) $};
		\node (rPB) at (6.5,0) {$ P(B) $};
		\node (rPpB) at (6.5,2) {$ P^\prime(B) $};
		
		\draw [->] (rPB) to [out=45, in=-45] node[right] {$ w(B) $} (rPpB);
		\draw [->] (rPB) to [out=135, in=-135] node[left] {$ v(B) $} (rPpB);
		\draw [->] (rPA) to node[auto] {$ P(f) $} (rPB);
		\draw [->] (rPpA) to node[auto] {$ P^\prime(f) $} (rPpB);
		\draw [->] (rPA) to node[left] {$ v(A) $} (rPpA);
		
		\node at (5.2, 1.5) {$ \simeq v(f) $};
		\node at (6.5, 1.3) {$ \Rightarrow $};
		\node at (6.5,0.9) {$ \zeta(B) $};
		
		\node at (3.25, 1) {$ = $};
	\end{tikzpicture},
\end{center}
where the map $\zeta$ was trivially extended to the set of all binary words in $G_{0}$ so that its components give the \textbf{associated modification} between natural transformations presented by $w, v$. 
\end{defin}

\begin{defin}
If $(G, \mathcal{R})$ is a generating datum and $\mathbb{M}$ is a monoidal bicategory, then we define the \textbf{bicategory $\mathbb{M}(G, \mathcal{R})$ of $(G, \mathcal{R})$-shapes with values in $\mathbb{M}$} to have $(G, \mathcal{R})$-shapes as objects, their homomorphisms as $1$-cells and transformations between homomorphisms as $2$-cells. 
\end{defin}
Observe that by construction we have a strict inclusion functor $i: \mathbb{M}(G, \mathcal{R}) \hookrightarrow \textbf{MonBicat}(\mathbb{F}(G, \mathcal{R}), \mathbb{M})$ given by passing to the associated homomorphisms, transformations and modifications. The main property of that homomorphism is that it is an equivalence, we call this result the \emph{cofibrancy theorem} and give it without proof. 

\begin{thm}[Cofibrancy theorem]
\label{thm:cofibrancy_theorem}
The inclusion of bicategories
\[ i: \mathbb{M}(G, \mathcal{R}) \hookrightarrow \textbf{MonBicat}(\mathbb{F}(G, \mathcal{R}), \mathbb{M}) \]
from the bicategory of $(G, \mathcal{R})$-shapes into the bicategory of monoidal homomorphisms is an equivalence of bicategories whose image is exactly the full subbicategory on strict homomorphisms and semi-strict natural transformations, that is, those transformations whose monoidal constraint cells are given by identities.
\end{thm}

\begin{rem}
A detailed proof of the cofibrancy theorem in the symmetric monoidal setting is given in \cite{chrisphd}, an analogous statement in the setting of ordinary bicategories was proven by Stephen Lack in \cite{lack_quillen_model_structure_for_bicategories}.
\end{rem}

The following corollary establishes a particularly nice property of freely generated monoidal bicategories.
\begin{cor}
Any monoidal homomorphism from a freely generated monoidal bicategory is equivalent to a strict one.
\end{cor}

\subsection{Promotion between bicategories of shapes}

In this section we will present a few results concerning the "truncation" of $G$-shapes. We will be mainly interested in singling out some favourable conditions under which some truncated set of data can be "promoted" to more complete one. The results are simple-minded in nature, but are used throughout the current work to simplify book-keeping.

Everything we say here applies \emph{both} to the theory of $G$-shapes in monoidal bicategories we presented above and to the theory of shapes in symmetric monoidal bicategories from \cite{chrisphd}. This is not surprising, as the definitions of the relevant bicategories of shapes are essentially the same.

Let $G$ be a free generating datum for a monoidal category and let $G^\prime$ be its 1-truncation. More precisely, let $G^\prime$ be the datum defined by $G^{\prime}_{0} = G_{0}, G^{\prime} _{1} = G_{1}$ and $G^\prime _{2} = \emptyset$. The inclusions $G^\prime _{i} \subseteq G_{i}$ yield a strict homomorphism $i: \mathbb{F}(G^\prime) \rightarrow \mathbb{F}(G)$ which is, as one sees immediately from the construction, an inclusion of monoidal bicategories that is bijective on objects and morphisms. 

Dually, if $\mathbb{M}$ is any monoidal category, then any $G$-shaped diagram in $\mathbb{M}$ yields a $G^\prime$-shaped simply be neglectance of data, similarly one can "truncate" homomorphisms and transformations. This assembles a strict forgetful homomrophism $\pi_{1}: G(\mathbb{M}) \rightarrow G^\prime(\mathbb{M})$ which in our description of the bicategories of diagrams as bicategories of semi-strict monoidal homomorphisms can be understood as beind induced by precomposition with $i$. 

\begin{prop}
\label{prop:comparison_of_categories_of_shapes_with_first_truncation}
The forgetful homomorphism $\pi_{1}: \mathbb{M}(G) \rightarrow M(G^\prime)$ is locally on Hom-categories an inclusion of a subcategory closed under isomorphisms.
\end{prop}

\begin{proof}
One sees directly from the definition of the shape bicategories that if $P_{1}, P_{2} \in \mathbb{M}(G)$, then to give a map $P_{1} \rightarrow P_{2}$ of $G$-shapes is \emph{exactly the same data} as to give a map of $G^\prime$-shapes, as it concerns only objects and 1-cells of the generating datum, only the conditions are different. In different words, to be a map of $G$-shapes - compared to being a map of $G^\prime$-shapes - is a \emph{property} and not a structure. This gives local injectivity on morphisms of $\pi_{1}$. 

Similarly, directly from the definitions one sees that a 2-cell in the bicategory of shapes is given by data concerning only the objects of the generating datum and axioms only concerning the 1-cells. These both coincide for $G$ and $G^\prime$, this yields that $\pi_{1}$ is locally bijective on $2$-cells.

We now only have to verify closure under isomorphisms. Let $s: P_{1} \rightarrow P_{2}$ be a map of $G$-shapes, $s: P_{1} \rightarrow P_{2}$ be a map of $G^\prime$-shapes and suppose we have an isomorphism $\omega: s \simeq t$. We have to show that then $s: P_{1} \rightarrow P_{2}$ is also a map in $\mathbb{M}(G)$, which amounts to verifying that it satisfies naturality with respect to $2$-cells in $G_{2}$. 

Let $\alpha: w_{1} \rightarrow w_{2}$ be an element of $G_{2}$, where $w_{i}$ are sentences in $G_{1}$ with $w_{i}: A \rightarrow B$.  We have to verify that 

\begin{center}
.
\end{center}
However, since $\omega$ is an isomorphism of maps of $G^\prime$-shapes, we have for each sentence $w$ in $G_{1}$

\begin{center}
	%
.
\end{center}
Pasting these decompositions of $t_{w_{i}}$ in terms of $s_{w_{i}}$ into the naturality equation above we see that it reduces to naturality for $s$, which we assumed.
 
\end{proof}

\begin{lem}[Promoting equivalences of $G^{\prime}$-shapes]
\label{lem:promoting_equivalences_of_g_shapes}
Let $P \in \mathbb{M}(G)$, $S \in \mathbb{M}(G^\prime)$ and suppose we are given an equivalence $s: \pi_{1}(P) \rightarrow S$. Then there is a unique $G$-shape $\widetilde{S}$ such that there exists an equivalence $\widetilde{s}: P \rightarrow \widetilde{S}$ of $G$-shapes with $\pi_{1}(\widetilde{s}) = s$. 
\end{lem}
We can think of the lemma as allowing as to promote equivalence of $G^\prime$-shapes to equivalence of $G$-shapes by transport of structure. 

\begin{cor}
If $P, S$ are $G$-shapes and $s: \pi_{1}(P) \rightarrow \pi_{1}(S)$ is an equivalence, then it is also an equivalence of $G$-shapes if and only if we have $S = \widetilde{\pi(S)}$. This follows immediately from the uniqueness part of the lemma.
\end{cor}

\begin{proof}
To promote $S$ to a $G$-shape we have to define it on the 2-cells in $G_{2}$. To say that $s: P \rightarrow \widetilde{S}$ is an equivalence of $G$-shapes is to say that for all $\alpha \in G_{2}$ with $\alpha: w_{1} \rightarrow w_{2}$ and $w_{i}: A \rightarrow B$ the naturality equation

\begin{center}
	\begin{tikzpicture}
		\node (PA) at (0,0) {$ P(A) $};
		\node (PB) at (2.5,0) {$ P(B) $};
		\node (SA) at (0,2) {$ S(A) $};
		\node (SB) at (2.5,2) {$ S(B) $};
		
		\draw [->] (PA) to node[left] {$ s_{A} $} (SA);
		\draw [->] (PB) to node[right] {$ s_{B} $} (SB);
		\draw [->] (PA) to [out=35, in=145] node[above] {$ P(w_{1}) $} (PB);
		\draw [->] (PA) to [out=-35, in=-145] node[below] {$ P(w_{2}) $} (PB);
		\draw [->] (SA) to node[auto] {$ S(w_{1}) $} (SB);
		
		\node at (1.25,0) {$ \Downarrow P(\alpha) $};
		\node at (1.2, 1.5) {$ \simeq s_{w_{1}} $};
		
		\node (rPA) at (5.5,-0.3) {$ P(A) $};
		\node (rPB) at (8,-0.3) {$ P(B) $};
		\node (rSA) at (5.5,1.7) {$ S(A) $};
		\node (rSB) at (8,1.7) {$ S(B) $};
		
		\draw [->] (rPA) to node[left] {$ t_{A} $} (rSA);
		\draw [->] (rPB) to node[right] {$ t_{B} $} (rSB);
		\draw [->] (rSA) to [out=35, in=145] node[above] {$ S(w_{1}) $} (rSB);
		\draw [->] (rSA) to [out=-35, in=-145] node[below] {$ S(w_{2}) $} (rSB);
		\draw [->] (rPA) to node[below] {$ P(w_{2}) $} (rPB);
		
		\node at (6.7, 0.2) {$ \simeq s_{w_{2}} $};
		\node at (6.75, 1.7) {$ \Downarrow S(\alpha) $};
		\node at (4.1, 1) {$ = $};
	\end{tikzpicture}.
\end{center}
holds. However, since $s$ is an equivalence, the vertical maps are all equivalences and we see that this equation uniquely defines $\widetilde{S}(\alpha)$. 
\end{proof}
The above result has the following corollary which can be pleasantly phrased in the language of the model structure on bicategories introduced by Stephen Lack, for details see \cite{lack_quillen_model_structure_for_bicategories}.
\begin{cor}
The restriction homomorphism $\pi_{1}$ is a fibration of underlying bicategories in the Lack model structure on bicategories and strict homomorphisms.
\end{cor}

\begin{proof}
Notice that $\pi_{1}$ has an invertible 2-cells lifting property as it is locally fully faithful and it has an equivalence lifting property by the lemma.
\end{proof}

The following little result establishes a bit more control over the transported structure of a $G$-shape.
\begin{lem}
\label{lem:promotion_respects_classes_of_relations}
Let $P \in \mathbb{M}(G), S \in \mathbb{M}(G^\prime)$, let $s: \pi(P) \rightarrow S$ be an equivalence and $\widetilde{s}: P \rightarrow \widetilde{S}$ be the promotion constructed using the \textbf{Lemma \ref{lem:promoting_equivalences_of_g_shapes}} above. Let $\mathcal{R}$ be a set of relations on $2$-cells on the free bicategory $\mathbb{F}(G)$ and suppose that $P$ is in fact a $(G, \mathcal{R})$-shapes. Then so is $\widetilde{S}$.
\end{lem}

\begin{proof}
If $\alpha = \beta$ is such relation, for $\alpha, \beta$ some parallel $2$-cells in $\mathbb{F}(G)$, then one sees from the construction of the transported structure that $\widetilde{S}(\alpha), \widetilde{S}(\beta)$ are uniquely defined by the same diagram, since we have $P(\alpha) = P(\beta)$.
\end{proof}

\begin{rem}
This also proves, and is equivalent to, the fact that given two equivalent $G$-shapes, one of them is a $(G, \mathcal{R})$-shape if and only if the other is.
\end{rem}

Promotion of equivalences is very useful in constructing equivalences between $G$-shapes, as it sometimes allows one to cut down on the number of equations that need to be satisfied.

We will now need a "2-cell analogue" of the promotion result. It will be useful to cut down on the number of equations that need to be satisfied when one tries construct a 2-cell in the bicategory of shapes. Again, we start by comparing the category of shapes to its truncation.

By $G^{\prime \prime}$ we will denote the $0$-truncation of a free generating datum $G$. More precisely, $G^{\prime \prime}$ has the same generating objects as $G$ and has no generating $1$-cells or $2$-cells.

\begin{prop}
The restriction homomorphism $\pi_{0}: \mathbb{M}(G^\prime) \rightarrow \mathbb{M}(G^{\prime \prime})$ is locally faithful on $2$-cells.
\end{prop}

\begin{proof}
Observe that to give a 2-cell in the category of shapes is to give data concerning only the objects of the generating datum. The two agree for $G^{\prime}$ and $G^{\prime \prime}$ and the conclusion follows. 
\end{proof}

\begin{lem}[Promotion of invertible $2$-cells]
\label{lem:promotion_of_isomorphisms_between_maps_of_gpp_shapes}
Let $P_{1}, P_{2} \in \mathbb{M}(G^\prime)$ and suppose we are given a homomorphism $s: P_{1} \rightarrow P_{2}$ of $G^{\prime}$-shapes, a homomorphism $t: \pi_{0}(P_{1}) \rightarrow \pi_{0}(P_{2})$ of $G^{\prime \prime}$-shapes and an isomorphism $\omega: \pi_{0}(s) \rightarrow t$. Then, there is a unique map of $G^{\prime}$-shapes $\widetilde{t}$ with $\pi_{0}(\widetilde{t}) = t$ such that there exists an isomorphism $\widetilde{\omega}: s \rightarrow \widetilde{t}$ with $\pi_{0}(\widetilde{\omega}) = \omega$. 
\end{lem}

\begin{cor}
If $s, t: P_{1} \rightarrow P_{2}$ are both maps of $G^{\prime}$-shapes and $\omega: \pi_{0}(s) \rightarrow \pi_{0}(t)$ is an isomorphism of maps of $G^{\prime \prime}$-shapes, then it is also an isomorphism in $\mathbb{M}(G^\prime)$ if and only if we have $t = \widetilde{\pi_{0}(t)}$. This follows immediately from the uniqueness part.
\end{cor}

\begin{proof}
To promote $t$ to a map of $G^{\prime}$-shapes we have to define the constraint isomorphisms for all $1$-cells in $G^{\prime} _{1} = G_{1}$. Let $f: A \rightarrow B \in G_{1}$ be one, where $A, B$ are some words in $G_{0}$. If we want $\omega$ to become an isomorphism 2-cell in $\mathbb{M}(G^{\prime \prime})$, then for each such $f$ the naturality equation 

\begin{center}
	\begin{tikzpicture}
		\node (P1A) at (0,0) {$ P_{1}(A) $};
		\node (P2A) at (0, 2) {$ P_{2}(A) $};
		\node (P1B) at (2,0) {$ P_{1}(B) $};
		\node (P2B) at (2,2) {$ P_{2}(B) $};
		
		\draw [->] (P1B) to node[right] {$ t_{B} $} (P2B);
		\draw [->] (P1A) to node[auto] {$ P_{1}(f) $} (P1B);
		\draw [->] (P2A) to node[auto] {$ P_{2}(f) $} (P2B);
		\draw [->] (P1A) to [out=135, in=-135] node[auto] {$ s_{A} $} (P2A);
		\draw [->] (P1A) to [out=45, in=-45] node[right] {$ t_{A} $} (P2A);
		
		\node at (0,1) {$ \Rightarrow \omega _{A} $};
		\node at (1.4, 1.4) {$ \simeq t_{f} $};
		
		\node (rP1A) at (4.5,0) {$ P_{1}(A) $};
		\node (rP2A) at (4.5, 2) {$ P_{2}(A) $};
		\node (rP1B) at (6.5,0) {$ P_{1}(B) $};
		\node (rP2B) at (6.5,2) {$ P_{2}(B) $};
		
		\draw [->] (rP1B) to [out=45, in=-45] node[right] {$ t_{B} $} (rP2B);
		\draw [->] (rP1B) to [out=135, in=-135] node[left] {$ s_{B} $} (rP2B);
		\draw [->] (rP1A) to node[auto] {$ P_{1}(f) $} (rP1B);
		\draw [->] (rP2A) to node[auto] {$ P_{2}(f) $} (rP2B);
		\draw [->] (rP1A) to node[left] {$ s_{A} $} (rP2A);
		
		\node at (5.1, 1.4) {$ \simeq s_{f} $};
		\node at (6.5,1) {$ \Rightarrow \omega _{B} $};
		
		\node at (3.25, 1) {$ = $};
	\end{tikzpicture}
\end{center}
must hold. However, since the components of $\omega$ are isomorphisms, this equation clearly has a unique solution $t_{f}$. This allows one to define the needed constraint $2$-cells and shows that they are unique.
\end{proof}

\begin{rem}
By the closure under isomorphisms part of \textbf{Proposition \ref{prop:comparison_of_categories_of_shapes_with_first_truncation}}, in the lemma above if $P_{1}, P_{2}$ are in fact $G$-shapes and $s$ is a map of $G$-shapes, the promotion $\widetilde{t}$ will also be.
\end{rem}

We end the section with the following observation, which is almost trivial to prove, but we use it repeatedly in the current work and so give it as a separate statement.
\begin{lem}[Reflection of equations]
\label{lem:reflection_of_equations}
Let $(G, \mathcal{R})$ be a generating datum for a monoidal bicategory and suppose that $\phi: \mathbb{M} \rightarrow \mathbb{M}^{\prime}$ is a strict monoidal homomorphism which is locally injective on 2-cells. Then the induced homomorphism
\[ \phi_{*}: \mathbb{M}(G) \rightarrow \mathbb{M}^\prime(G) \]
has the property that $P$ is a $(G, \mathcal{R})$ shape if and only if $\phi_{*}(P)$ is. 
\end{lem}

\begin{proof}
If we identify $P \in \mathbb{M}(G)$ with a strict homomorphism $\mathbb{F}(G) \rightarrow \mathbb{M}$, then $\phi_{*}(P) = \phi \circ P$. In this context the assertion follows immediately from the fact that $\phi$ reflects equations on 2-cells.
\end{proof}

\bibliographystyle{amsalpha}
\bibliography{on_dualizable_objects_in_monoidal_bicategories_bibliography}

\end{document}